\documentclass[11pt]{article}

\usepackage{latexsym}
\usepackage[all]{xy}
\usepackage{amsmath}
\usepackage{amssymb}
\usepackage{amsfonts}
\usepackage{color}
\usepackage{theorem}

\usepackage{color}
\usepackage{ifthen}

\newcommand {\Map} {\mathbb{R}\mathbf{Map}}

\newcommand {\OO} {\mathcal{O}}

\newcommand {\DR}{\mathbf{DR}}
\newcommand {\Pol}{\mathbf{Pol}}
\newcommand {\W}{\mathbf{W}}
\newcommand {\A} {\mathcal{A}}
\newcommand {\E} {\mathcal{E}}
\newcommand {\LL} {\mathcal{L}}
\newcommand {\B} {\mathcal{B}}
\newcommand {\V} {\mathbb{V}}
\newcommand {\C} {\mathcal{M}}
\newcommand {\T} {\mathcal{T}}
\newcommand {\Spec} {\mathbf{Spec}}
\newcommand  {\dg}     {\mathbf{dg}}
\newcommand  {\dglie}     {\mathbf{dgLie}}
\newcommand  {\cdg}     {\mathbf{codg}}
\newcommand  {\mdg}     {\mathbf{dg}^{gr}}
\newcommand  {\edg}     {\epsilon-\mathbf{dg}}
\newcommand  {\medg}     {\epsilon-\mathbf{dg}^{gr}}
\newcommand {\scat} {\infty-\mathbf{Cat}}
\newcommand  {\mecdga}     {\epsilon-\mathbf{cdga}^{gr}}

\newcommand  {\cdga}     {\mathbf{cdga}}

\newcommand  {\dAff}     {\mathbf{dAff}}
\newcommand  {\dFSt}     {\mathbf{dFSt}}

\newcommand  {\dSt}   {\mathbf{dSt}}

\newcommand{\s}{\infty}
\newcommand {\D} {\mathbb{D}}

\newtheorem{thm}{Theorem}[subsection]
\newtheorem{prop}[thm]{Proposition}
\newtheorem{lem}[thm]{Lemma}
\newtheorem{sublem}[thm]{Sub-Lemma}
\newtheorem{df}[thm]{Definition}
\newtheorem{cor}[thm]{Corollary}
\newtheorem{conj}[thm]{Conjecture}

{\theorembodyfont{\rmfamily} \newtheorem{rmk}[thm]{Remark}}
\newtheorem{ex}[thm]{Example}
\newtheorem{Q}[thm]{Question}

\newtheorem{thmi}{Main results}

\newcommand{\Appendix}[1]{%
  \refstepcounter{section}%
  \addtocontents{toc}{\protect\setcounter{tocdepth}{1}}
  \addcontentsline{toc}{section}%
    {\bfseries\appendixname~\thesection\ #1}%
    {\medskip\noindent \Large\bfseries\appendixname\ \thesection\ #1}%
\sectionmark{#1}\smallskip\noindent
\renewcommand{\theequation}{{\bf 
{{\thesection}}.{\arabic{equation}}}}
}
\begin{document}

\title{\textbf{
Shifted Poisson Structures and \\ deformation quantization}}  

\author{D. Calaque, T. Pantev, B. To\"en, M. Vaqui\'e, G. Vezzosi}  

\date{}

\maketitle

\begin{abstract}
  This paper is a sequel to \cite{ptvv}. We develop a general and
  flexible context for differential calculus in derived geometry,
  including the de Rham algebra and the mixed algebra of polyvector
  fields. We then introduce the formalism of formal derived stacks and
  prove formal localization and gluing results. These allow us to
  define shifted Poisson structures on general derived Artin stacks,
  and to prove that the non-degenerate Poisson structures correspond
  exactly to shifted symplectic forms. Shifted deformation
  quantization for a derived Artin stack endowed with a shifted
  Poisson structure is discussed in the last section. This paves the
  road for shifted deformation quantization of many interesting derived
  moduli spaces, like those studied in \cite{ptvv} and  many
  others.
\end{abstract}

\tableofcontents

\section*{Introduction} 
\addcontentsline{toc}{section}{Introduction}

This work is a sequel of \cite{ptvv}. We introduce the notion of a
\emph{shifted Poisson structure} on a general derived Artin stack,
study its relation to the shifted symplectic structures from
\cite{ptvv}, and construct a deformation quantization of it.  As a
consequence, we construct a deformation quantization of any derived
Artin stack endowed with an $n$-shifted symplectic structure, as soon
as $n\neq 0$. In particular we quantize many derived moduli spaces
studied in \cite{ptvv}. In a nutshell the results of this work
are summarized as follows.

\begin{thmi}\label{ti1}
\begin{enumerate}
\item There exists a meaningful notion of \emph{$n$-shifted Poisson
    structures} on derived Artin stacks locally of finite
  presentation, which recovers the usual notion of Poisson structures on
  smooth schemes when $n=0$.
\item For a given derived Artin stack $X$, the space of $n$-shifted
  symplectic structures on $X$ is naturally equivalent to the space of
  non-degenerate $n$-shifted Poisson structures on $X$.
\item Let $X$ be any derived Artin stack locally of finite
  presentation endowed with an $n$-shifted Poisson structure $\pi$. For
  $n\neq 0$ there exists a canonical deformation quantization of $X$
  along $\pi$, realized as an $E_{|n|}$-monoidal $\s$-category
  $\mathsf{Perf}(X,\pi)$, which is a deformation of the symmetric monoidal
  $\s$-category $\mathsf{Perf}(X)$ of perfect complexes on $X$.
\end{enumerate}
\end{thmi}

\

\noindent
As a corollary of these, we obtain the existence of deformation
quantization of most derived moduli stacks studied in \cite{ptvv},
e.g.~of the derived moduli of $G$-bundles on smooth and proper
Calabi-Yau schemes, or the derived moduli of $G$-local systems on 
compact oriented topological manifolds. The
existence of these deformation quantizations is a completely new
result, which is a far reaching generalization of the construction of
deformation quantization of character varieties heavily studied in
topology, and provides a new world of quantized moduli spaces to
explore in the future.

The above items are not easy to achieve. Some ideas of what
\emph{$n$-shifted Poisson structures} should be have been available in
the literature for a while (see \cite{mel,toenems,toenicm}), but up
until now no general definition of $n$-shifted Poisson structures on
derived Artin stacks existed outside of the rather restrictive case of
Deligne-Mumford stacks. The fact that Artin stacks have affine covers
only up to smooth submersions is an important technical obstacle which
we have to deal with already when we define shifted Poisson structures
in this general setting. Indeed, in contrast to differential forms,
polyvectors do not pull-back along smooth morphisms, so the well
understood definition in the affine setting (see \cite{mel,toenems})
can not be transplanted to an Artin stack without additional effort,
and such a transplant requires a new idea. A different complication
lies in the fact that the comparison between non-degenerate shifted
Poisson structures and their symplectic counterparts requires keeping
track of non-trivial homotopy coherences even in the case of an affine
derived scheme.  One reason for this is that non-degeneracy is only
defined up to quasi-isomorphism, and so converting a symplectic
structure into a Poisson structure by dualization can not be performed
easily. Finally, the existence of deformation quantization requires
the construction of a deformation of the globally defined
$\s$-category of perfect complexes on a derived Artin stack. These
$\s$-categories are of global nature, and their deformations are not
easily understood in terms of local data.

In order to overcome the above mentioned technical challenges we
introduce a new point of view on derived Artin stacks by developing
tools and ideas from formal geometry in the derived setting. This new approach is one of the technical hearts of the paper, and
we believe it will be an important general tool in derived geometry, even outside the applications to shifted Poisson and symplectic structures discussed in this work. The key
new idea here is to understand a given derived Artin stack $X$ by
means of its various formal completions $\widehat{X}_{x}$, at all of
its points $x$ in a coherent fashion. For a smooth algebraic variety,
this idea has been already used with great success in the setting of
deformation quantization (see for instance \cite{fe,ko,bezka}), but
the extension we propose here in the setting of derived Artin stacks
is new. By \cite{lu}, the geometry of a given formal completion
$\widehat{X}_{x}$ is controlled by a dg-Lie algebra, and our approach,
in a way, rephrases many problems concerning derived Artin stacks in
terms of dg-Lie algebras. In this work we explain how shifted
symplectic and Poisson structures, as well as $\s$-categories of
perfect complexes, can be expressed in those terms.  Having this
formalism at our disposal is what makes our Main statement \ref{ti1}
accessible. The formalism essentially allows us to reduce the problem
to statements concerning dg-Lie algebras over general base rings and
their Chevalley complexes. The general formal geometry results we
prove on the way are of independent interest and will be useful for
many other questions related to derived Artin stacks. 
For a slightly different approach to formal derived geometry, we recommend \cite[Part IV]{garo}.\\

Let us now discuss the mathematical content of the paper in more
detail. To start with, let us explain the general strategy and the
general philosophy developed all along this manuscript. For a given
derived Artin stack $X$, locally of finite presentation over a base
commutative ring $k$ of characteristic $0$, we consider the
corresponding de Rham stack $X_{DR}$ of \cite{si1,si2}. As an
$\s$-functor on commutative dg-algebras, $X_{DR}$ sends $A$ to
$X(A_{red})$, the $A_{red}$-points of $X$ (where $A_{red}$ is defined
to be the reduced ordinary commutative ring $\pi_{0}(A)_{red}$). The
natural projection $\pi : X \longrightarrow X_{DR}$ realizes $X$ as a
family of formal stacks over $X_{DR}$: the fiber of $\pi$ at a given
point $x\in X_{DR}$, is the formal completion $\widehat{X}_{x}$ of $X$
at $x$. By \cite{lu} this formal completion is determined by a dg-Lie
algebra $l_x$. However, the dg-Lie algebra $l_x$ itself does not
exist globally as a sheaf of dg-Lie algebras over $X_{DR}$, simply
because its underlying complex is $\mathbb{T}_{X}[-1]$, the shifted
tangent complex of $X$, which in general does not have a flat
connection and thus does not descend to $X_{DR}$. However, the
Chevalley complex of $l_x$, viewed as a graded mixed commutative
dg-algebra (\emph{cdga} for short) can be constructed as a global
object $\B_{X}$ over $X_{DR}$. To be more precise we construct
$\B_{X}$ as the derived de Rham complex of the natural inclusion
$X_{red} \longrightarrow X$, suitably sheafified over $X_{DR}$. One of
the key insights of this work is the following result, expressing
global geometric objects on $X$ as sheafified notions on $X_{DR}$
related to $\B_X$.

\begin{thmi}\label{ti2}
With the notation above:

\begin{enumerate}
\item The $\s$-category $\mathsf{Perf}(X)$, of perfect complexes on
  $X$, is naturally equivalent, as a symmetric monoidal
  $\s$-category, to the $\s$-category of perfect sheaves of graded
  mixed $\B_X$-dg-modules on $X_{DR}$:
$$
\mathsf{Perf}(X) \simeq \B_X - Mod_{\medg}^{\mathsf{Perf}}.
$$
\item There is an equivalence between the space of $n$-shifted
  symplectic structures on $X$, and the space of closed and
  non-degenerate $2$-forms on the sheaf of graded mixed cdgas $\B_X$.
\end{enumerate}
\end{thmi}

\

The results above state that the geometry of $X$ is largely recovered
from $X_{DR}$ together with the sheaf of graded mixed cdgas $\B_X$, and
that the assignment $X \mapsto (X_{DR},\B_X)$ behaves in a faithful
manner from the perspective of derived algebraic geometry. In the last part of the paper, we 
take advantage of this in order to define the deformation quantization problem for objects belonging to general categories over $k$. In particular, we study and quantize shifted
Poisson structures on $X$, by considering compatible brackets on the
sheaf $\B_X$. Finally, we give details for three relevant quantizations and compare them to the existing literature. The procedure of replacing $X$ with $(X_{DR},\B_X)$ is crucial for derived Artins stacks because it essentially reduces statements and notions to the
case of a sheaf of graded mixed cdgas. As graded mixed cdgas can also
be understood as cdgas endowed with an action of a derived group
stack, this further reduces statements to the case of (possibly
unbounded) cdgas, and thus to an affine situation.

\

\smallskip

\noindent \textbf{Description of the paper.}

In the first section, we start with a very general and flexible
context for (relative) differential calculus. We introduce the
\emph{internal cotangent complex} $\mathbb{L}^{int}_A$ and
\emph{internal de Rham complex} $\mathbf{DR}^{int}(A)$ associated with
a commutative algebra $A$ in a good enough symmetric monoidal stable
$k$-linear $\infty$-category $\C$ (see Section~\ref{1.1} and
Section~\ref{infty} for the exact assumptions we put on $\C$).  The
internal de Rham complex $\mathbf{DR}^{int}(A)$ is defined as a graded
mixed commutative algebra in $\C$. Next we recall from \cite{ptvv} and
extend to our general context the spaces $\mathcal A^{p,cl}(A,n)$ of
(closed) $p$-forms of degree $n$ on $A$, as well as of the space
$\mathsf{Symp}(A,n)$ of $n$-shifted symplectic forms on $A$. We
finally introduce (see also \cite{ptvv,mel,toenems,toenicm}) the
\emph{object} $\mathbf{Pol}^{int} (A,n)$ of \emph{internal $n$-shifted
  polyvectors on $A$}, which is a graded $n$-shifted Poisson algebra
in $\C$. In particular, $\mathbf{Pol}^{int}(A,n)[n]$ is a graded Lie
algebra object in $\C$. We recall from \cite{mel} that the space
$\mathsf{Pois}(A,n)$ of graded $n$-shifted Poisson structures on $A$
is equivalent to the mapping space from $\mathbf{1}(2)[-1]$ to
$\mathbf{Pol}^{int}(A,n+1)[n+1]$ in the $\infty$-category of graded
Lie algebras in $\C$, and we thus obtain a reasonable definition of
non-degeneracy for graded $n$-shifted Poisson structures. Here
$\mathbf{1}(2)[-1]$ denotes the looping of the monoidal unit of $\C$
sitting in pure weight $2$ (for the grading). We finally show that

\

\noindent {\bfseries Corollary~\ref{cpoissnd3}} \emph{ If
  $\mathbb{L}_A^{int}$ is a dualizable $A$-module in $\mathcal M$,
  then there is natural morphism
$$
\mathsf{Pois}^{nd}(A,n) \longrightarrow \mathsf{Symp}(A,n)
$$
from the space $\mathsf{Pois}^{nd}(A,n)$ of non-degenerate $n$-shifted
Poisson structures on $A$ to the space $\mathsf{Symp}(A,n)$ of
n-shifted symplectic structures on $A$. }

\

We end the first part of the paper by a discussion of what happens
when $\C$ is chosen to be the $\infty$-category $\epsilon-(k-
mod)^{gr}$ of graded mixed complexes, which will be our main case of
study in order to deal with the sheaf $\B_X$ on $X_{DR}$ mentioned
above. We then describe two lax symmetric monoidal functors
$|-|,|-|^t:\epsilon-\C^{gr}\to \C$, called standard realization and
Tate realization.  We can apply the Tate realization to all of the
previous internal constructions and get in particular the notions of
Tate $n$-shifted symplectic form and non-degenerate Tate $n$-shifted
Poisson structure. We prove that, as before, these are equivalent as
soon as $\mathbb{L}_A^{int}$ is a dualizable $A$-module in $\mathcal
M$.

\medskip

One of the main difficulties in dealing with $n$-shifted polyvectors
(and thus with $n$-shifted Poisson structures) is that they do not
have sufficiently good functoriality properties.  Therefore, in
contrast with the situation with forms and closed forms, there is no
straightforward and easy global definition of $n$-shifted
polyvectors and $n$-shifted Poisson structures.  Our strategy is to
use ideas from formal geometry and define an $n$-shifted Poisson
structure on a derived stack $X$ as a flat family of $n$-shifted
Poisson structures on the family of all formal neighborhoods of points
in $X$.  The main goal of the second part of the paper is to make
sense of the previous sentence for general enough derived stacks,
i.e. for locally almost finitely presented derived Artin stacks over $k$.
In order to achieve this, we develop a very general theory of \emph{derived formal localization} that will be certainly very useful in other applications of derived geometry, as well.
\medskip

We therefore start the second part by introducing various notions of
formal derived stacks: formal derived stack, affine formal derived
stack, good formal derived stack over $A$, and perfect formal derived
stack over $A$.  It is important to note that if $X$ is a derived
Artin stack, then
\begin{itemize}
\item the formal completion $\widehat{X}_f:X\times_{X_{DR}}F_{DR}$ along
  any map $f:F\to X$ is a formal  
derived stack. 
\item the formal completion $\widehat{X}_x$ along a point 
$x:\mathbf{Spec}(A)\to X$ is an affine formal 
derived stack. 
\item each fiber $X\times_{X_{DR}}\mathbf{Spec}(A)$ of $X\to X_{DR}$
  is a good formal derived stack over $A$, which is moreover perfect
  if $X$ is locally of finite presentation.
\end{itemize} 
Our main result here is the following 

\

\noindent {\bfseries Theorem \ref{t1}} \emph{There exists an
  $\infty$-functor $\mathbb{D}$ from affine formal derived stacks to
  commutative algebras in $\C =\epsilon-(k-mod)^{gr}$, together
  with a conservative $\infty$-functor
$$
\phi_X:\mathsf{QCoh}(X)\to\mathbb{D}(X)-mod_{\C},
$$
which becomes fully faithful when restricted to perfect modules. }

\

Therefore, $\mathsf{Perf}(X)$ is identified with a full
sub-$\infty$-category $\mathbb{D}(X)- mod_{\mathcal M}^{perf}$ of
$\mathbb{D}(X)-mod_{\mathcal M}$ that we explicitly determine.  We then prove that the de Rham
theories of $X$ and of $\mathbb{D}(X)$ are equivalent for a perfect
algebraisable formal derived stack over $A$. Namely we show that:
$$
\mathbf{DR}\big(\mathbb{D}(X)/\mathbb{D}(\mathbf{Spec}\,A)\big)
\simeq\mathbf{DR}^t\big(\mathbb{D}
(X)/\mathbb{D}(\mathbf{Spec}\,A)\big)
\simeq\mathbf{DR}(X/\mathbf{Spec}\,A)
$$
as commutative algebras in $\epsilon-(A-mod)^{gr}$.  We finally extend
the above to the case of families $X\to Y$ of algebraisable perfect
formal derived stacks, i.e. families for which every fiber
$X_A:=X\times_Y\mathbf{Spec}\,A\to\mathbf{Spec}\,A$ is an
algebraisable perfect formal derived stack. We get an equivalence
of symmetric monoidal $\infty$-categories
$\phi_X:\mathsf{Perf}(X)\simeq\mathbb{D}_{X/Y}-mod_{\mathcal
  M}^{perf}$ as well as equivalences:
$$
\Gamma\Big(Y,\mathbf{DR}\big(\mathbb{D}_{X/Y}/
\mathbb{D}(Y)\big)\Big)\simeq\Gamma\Big(Y,
\mathbf{DR}^t\big(\mathbb{D}_{X/Y}/
\mathbb{D}_Y)\big)\simeq\mathbf{DR}(X/Y)
$$
of commutative algebras in
$\C$

In particular, whenever $Y=X_{DR}$, we get a description of the de Rham
(graded mixed) algebra $\mathbf{DR}(X)\simeq\mathbf{DR}(X/X_{DR})$ by
means of the global sections of the relative Tate de Rham (graded
mixed) algebra  $\mathcal B_X:=\mathbb{D}_{X/X_{DR}}$ over
$\mathbb{D}_{X_{DR}}$. Informally speaking, we prove that a (closed)
form on $X$ is a float family of (closed) forms on the family of all
formal completions of $X$ at various points.

\medskip

The above justifies the definitions of shifted polyvector fields and
shifted Poisson structures that we introduce in the third part of the
paper.  Namely, the $n$-shifted Poisson algebra $\mathbf{Pol}(X/Y,n)$
of $n$-shifted polyvector fields on a family of algebraisable perfect
formal derived stacks $X\to Y$ is defined to be
$$
\Gamma\big(Y,\mathbf{Pol}^t(\mathbb{D}_{X/Y}/\mathbb{D}_X,n)\big)
$$
The space of $n$-shifted Poisson structures $\mathsf{Pois}(X/Y,n)$ is
then defined as the mapping space from $k(2)[-1]$ to
$\mathbf{Pol}(X/Y,n+1)[n+1]$ in the $\infty$-category of graded
Lie algebras in $\mathcal M$.  Following the affine case treated in
the first part (see also \cite{mel}), we again prove that this is
equivalent to the space of $\mathbb{D}_{Y}$-linear $n$-shifted Poisson
algebra structures on $\mathbb{D}_{X/Y}$.  We then
prove\footnote{Recently J. Pridham proved this comparison theorem for
  derived Deligne-Mumford stacks by a different approach
  \cite{pridham-compare}. In a later version of \cite{pridham-compare}, which appeared after our paper was put on the arXiv, the author modified his approach in order to treat also the case of derived Artin stacks.} the following

\

\noindent
{\bfseries Theorem~\ref{t4}} \emph{The subspace of non-degenerate
  elements in $\mathsf{Pois}(X,n):=\mathsf{Pois}(X/X_{DR},n)$ is
  equivalent to $\mathsf{Symp}(X,n)$ for any derived Artin stack that
  is locally of finite presentation.  }

\

We conclude the third Section by defining the deformation quantization problem for 
$n$-shifted Poisson structures, whenever $n\geq-1$. 
For every such $n$, we consider a $\mathbb{G}_m$-equivariant $\mathbb{A}^1_k$-family of $k$-dg-operads 
$\mathbb{BD}_{n+1}$ such that its $0$-fiber is the Poisson operad $\mathbb{P}_{n+1}$ and its generic 
fiber is the $k$-dg-operad $\mathbb{E}_{n+1}$ of chains of the little $(n+1)$-disk topological operad. 
The general deformation quantization problem can then be stated as follows: 

\

\noindent
\textbf{Deformation Quantization Problem. } \emph{Given a $\mathbb{P}_{n+1}$-algebra stucture on an object,  does it exist a family of 
$\mathbb{BD}_{n+1}$-algebra structures such that its $0$-fiber is the original $\mathbb{P}_{n+1}$-algebra structure?}

\

Let $X$ be a derived Artin stack locally of finite presentation over $k$, and equipped with an $n$-shifted Poisson structure. 
Using the formality of $\mathbb{E}_{n+1}$ for $n\geq1$, we can solve the deformation quantization problem 
for the $\mathbb{D}_{X_{DR}}$-linear $\mathbb{P}_{n+1}$-algebra structure on $\mathcal B_X$. 
This gives us, in particular, a $\mathbb{G}_m$-equivariant $1$-parameter family of $\mathbb{D}_{X_{DR}}$-linear 
$\mathbb{E}_{n+1}$-algebra structures on $\mathcal B_X$. Passing to perfect modules we get a $1$-parameter deformation 
of $\mathsf{Perf}(X)$ as an $\mathbb{E}_n$-monoidal $\infty$-category, which we call the $n$\emph{-quantization of} $X$. 

We work out three important examples in some details: 
\begin{itemize}
\item the case of an $n$-shifted symplectic structure on the formal neighborhood of a $k$-point in $X$: we recover Markarian's Weyl $n$-algebra from \cite{markarian}.  
\item the case of those $1$-shifted Poisson structure on $BG$ that are given by elements in $\wedge^3(\mathfrak g)^{\mathfrak g}$: 
we obtain a deformation, as a monoidal $k$-linear category, of the category $\mathsf{Rep}^{fd}(\mathfrak g)$ of finite dimensional representations of $\mathfrak g$. 
\item the case of $2$-shifted Poisson structures on $BG$, given by elements in $Sym^2(\mathfrak g)^{\mathfrak g}$: 
we obtain a deformation of $\mathsf{Rep}^{fd}(\mathfrak g)$ as a braided monoidal category. 
\end{itemize}

Finally, Appendices A and B contains some technical results used in Sections 1 and 3, respectively.

\

\smallskip

\noindent \textbf{Further directions and future works.} In order to 
finish this introduction, let us mention that the present
work does not treat some important questions related to quantization, which we hope to address in
the future. For instance, we introduce a general notion of coisotropic
structures for maps with an $n$-shifted Poisson target, analogous to
the notion of Lagrangian structures from \cite{ptvv}. However, the
definition itself requires a certain additivity theorem, whose proof
has been announced recently by N. Rozenblyum but is not available yet.
Also, we did not address the question of comparing Lagrangian
structure and co-isotropic structures that would be a relative version
of our comparison between shifted symplectic and non-degenerate
Poisson structures. Neither did we address the question of
quantization of coisotropic structures.  In a different direction, our
deformation quantizations are only constructed under the restriction
$n\neq 0$. The case $n=0$ is presently being investigated, but at the
moment is still open. In the same spirit, when $n=-1$ and $n=-2$,
deformation quantization admits an interpretation different from our
construction (see for example \cite[Section~6.2]{toenicm}). We believe
that our present formal geometry approach can also be applied to these
two specific cases.

\

\smallskip

\noindent \textbf{Acknowledgments.}  We are thankful to D. Kaledin for
suggesting to us several years ago that formal geometry should give a
flexible enough setting for dealing with shifted polyvectors and
Poisson structures.  We would also like to thank K. Costello and
O. Gwilliam for their explanations about the Darboux lemma
\cite[Lemma~11.2.0.1]{factorization} in the setting of minimal
$\LL_{\s}$-algebras, which can be found in a disguised form in the
proof of our Lemma \ref{l12}. We are grateful to M. Harris, D. Joyce,
M. Kontsevich, V. Melani, M. Porta, P. Safronov, D. Tamarkin, and N.
Rozenblyum for illuminating conversations on the subject of this
paper. V. Melani and M. Porta's questions were very useful in order to clarify some obscure or implicit points in
a previous version of the paper. 
It is a pleasure to thank once again C. Simpson, for having
brought to us all along these years the importance of the de Rham
stack $X_{DR}$, which is a central object of the present work.

Damien Calaque acknoweldges the support of Institut Universitaire 
de France and ANR SAT. 
During the preparation of this work Tony Pantev was partially
supported by NSF research grant DMS-1302242. 
Bertrand To\"en and Michel Vaqui\'e have been partially supported by 
ANR-11-LABX-0040-CIMI within the program ANR-11-IDEX-0002-02. 
Gabriele Vezzosi is a member of the GNSAGA-INDAM group (Italy) and of PRIN-Geometria delle
variet\'{a} algebriche (Italy). In addition Gabriele Vezzosi would
like to thank the IAS (Princeton) and the IHES (Bures-sur-Yvette),
where part of this work was carried out, for providing excellent
working conditions.  

\

\bigskip

\noindent \textbf{Notation.} 

\begin{itemize}
\item Throughout this paper $k$ will denote a Noetherian commutative
  $\mathbb{Q}$-algebra. 
\item We will use $(\infty, 1)$-categories (\cite{lutop}) as our model
  for $\infty$-categories. They will be simply called
  $\infty$-categories. 
\item For a model category $N$, we will denote by $L(N)$ the
  $\infty$-category defined as the homotopy coherent nerve of the
  Dwyer-Kan localization of $N$ along its weak equivalences.  
\item The $\infty$-category $\T:= L(\textbf{sSets})$ will be called
  the $\infty$-category of \emph{spaces}. 
\item All symmetric monoidal categories we use will be symmetric monoidal
  (bi)closed categories. 
\item $\cdga_{k}$ will denote the $\s$-category of non-positively
  graded differential graded $k$-algebras (with differential
  increasing the degree). For $A\in \cdga_{k}$, we will write
  $\pi_{i}\, A = H^{-i}(A)$ for any $i \geq 0$. 
\item For $A \in \cdga_{k}$, we will write either $\mathrm{L}(A)$
  or $\mathrm{L}_{\mathrm{QCoh}}(A)$ for the $\s$-category of
  $A$-dg-modules 
\item For $A \in \cdga_{k}$, we will denote by
  $\mathrm{L}_{\mathsf{Perf}}(A)$ the full sub-$\s$-category of
  $\mathrm{L}(A)$ consisting of perfect $A$-dg-modules. 
\item If $X$ is a derived geometric stack, we will denote by either
  $\mathsf{QCoh}(X)$ or $\mathrm{L}_{\mathsf{QCoh}}(X)$ the $k$-linear
  symmetric monoidal dg-category of quasi-coherent complexes on $X$.
\item If $X$ is a derived geometric stack, we will denote by either
  $\mathsf{Perf}(X)$ or $\mathrm{L}_{\mathsf{Perf}}(X)$ the symmetric
  monoidal sub-dg-category of $\mathsf{QCoh}(X)$  consisting of
  dualizable objects. 
\item If $X$ is a derived geometric stack, we will denote by either
  $\mathsf{Coh}(X)$ or $\mathrm{L}_{\textrm{Coh}}(X)$ the full sub-dg
  category of $\mathsf{QCoh}(X)$ consisting of complexes whose
  cohomology sheaves are coherent over the truncation $\mathsf{t}_{0}
  X$.
  \item For a derived stack $X$, $\Gamma(X, -)$ will always denote the \emph{derived} functor of global sections on $X$ (i.e. hypercohomology).
\end{itemize}

\section{Relative differential calculus}\label{RDCsection}

In this section we describe the basics of differential calculus
\emph{inside} any reasonable $k$-linear symmetric monoidal
$\infty$-category. In particular, we introduce cotangent complexes, De
Rham mixed dg-algebras, shifted (closed) forms and 
polyvectors, and two different \emph{realizations} (standard and Tate)
of such objects over $k$.

\subsection{Model categories setting}\label{1.1} 

Let $k$ be a Noetherian commutative $\mathbb{Q}$-algebra, and let
$C(k)=\dg_{k}$ be the category of (unbounded, cochain)
$k$-dg-modules. We endow $C(k)$ with its standard model category
structure whose equivalences are quasi-isomorphisms and whose
fibrations are epimorphisms (\cite[Theorem~2.3.11]{hov}).  The natural
tensor product $-\otimes_{k}-$ of dg-modules endows $C(k)$ with the
structure of a symmetric monoidal model category
(\cite[Proposition~4.2.13]{hov}). As a monoidal model category $C(k)$
satisfies the monoid axiom of \cite[Definition~3.3]{ss}, and moreover,
since $k$ is a $\mathbb{Q}$-algebra, $C(k)$ is freely-powered in the
sense of \cite[Definition~4.5.4.2]{lualg}.

Suppose next that $M$ is a symmetric monoidal model category that is
combinatorial as a model category (\cite[Definition~A.2.6.1]{lutop}).
Assume furthermore that $M$ admits a $C(k)$-enrichment (with tensor
and cotensor) compatible with both the model and the monoidal
structures, i.e. $M$ is a symmetric monoidal $C(k)$-model algebra as
in \cite[Definition~4.2.20]{hov}. As a consequence (see our
Proposition~\ref{Mstable}) such an $M$ is a \emph{stable} model category,
i.e. it is pointed and the suspension functor is a self equivalence of
its homotopy category.

\noindent All along this first section, and as a reference for the
rest of the paper, we make the following  \textsf{standing
  assumptions on} $M$

\begin{enumerate}
\item The unit $\mathbf{1}$ is a cofibrant object in $M$.
\item For any cofibration $j : X \rightarrow Y$ in $M$, any object $A
  \in M$, and for any morphism $u: A\otimes X \rightarrow C$ in $M$
  the push-out square in $M$
$$\xymatrix{
C \ar[r] & D \\
A\otimes X \ar[u]^-{u} \ar[r]_-{id\otimes j} & A \otimes Y, \ar[u]}$$
is a homotopy push-out square.  
\item For a cofibrant object $X \in M$, the functor $X \otimes - : M
  \longrightarrow M$ 
preserves equivalences (i.e. \emph{cofibrant objects in $M$ are
  $\otimes$-flat}). 
\item $M$ is a \emph{tractable} model category, i.e. there are
  generating sets of cofibrations $I$, and trivial cofibrations $J$ in
  $M$ with cofibrant domains.
\item Equivalences are stable under filtered colimits and finite
  products in $M$. 
 \end{enumerate}

 We note that conditions $(2)-(5)$ together imply that $M$ satisfies
 the monoid axiom of \cite[Definition~3.3]{ss} In particular
 (\cite[Theorem~4.1 (2)]{ss}), for any commutative monoid $A\in
 Comm(M)$, the category of $A$-modules in $M$, denoted by $A-Mod_{M}$,
 is endowed with the structure of a symmetric monoidal combinatorial
 model category, for which the equivalences and fibrations are defined
 in $M$, and it again satisfies the monoid axiom. Moreover,
 $A-Mod_{M}$ comes with an induced compatible $C(k)$-enrichment (with
 tensor and cotensor). Moreover, as shown in Proposition
 \ref{modinvariance}, the conditions $(2)-(5)$ on $M$ imply that if $A
 \longrightarrow A'$ is an equivalence in $Comm(M)$, then the induced
 restriction-extension of scalars Quillen adjunction
$$A-Mod_{M} \longleftrightarrow A'-Mod_{M}$$
is a Quillen equivalence.

As $k$ is a $\mathbb{Q}$-algebra, $M$ is itself a $\mathbb{Q}$-linear
category. This implies that $M$ is freely-powered in the sense of
\cite[Definition~4.5.4.2]{lualg}, since quotients by finite group
actions are split epimorphisms in characteristic $0$. As a
consequence, the category $Comm(M)$ of commutative and unital monoids
in $M$, is again a combinatorial model category for which the
equivalences and fibrations are defined via the forgetful functor to
$M$, and whose generating (trivial) cofibrations are given by
$\mathrm{Sym}(I)$ (respectively, $\mathrm{Sym}(J)$), where $I$
(respectively $J$) are generating (trivial) cofibrations in $M$
(\cite[Proposition~4.5.4.6]{lualg}). 

\

Let $B$ be a $k$-linear commutative Hopf dg-algebra. We let $B-\cdg_{M}$ be the category of $B$-comodules in
$M$, i.e. the category whose
\begin{itemize} 
\item objects are objects $P$ in $M$ equipped with a morphism
  $\rho_{P}: P \to P\otimes_{k} B$ in $M$ ($\otimes_{k}: M \times C(k)
  \to M$ being the tensor product given by the
  $C(k)$-enrichment\footnote{Note that this slightly abusive notation
    for the tensor enrichment $\otimes_{k}:= \otimes_{C(k)}$ is
    justified by the fact that the properties of the enrichment give a
    canonical isomorphism $P \otimes_{C(k)} (B\otimes_{k} B) \simeq
    (P\otimes_{C(k)} B) \otimes_{C(k)} B$.}) satisfying the usual
  identities
$$
\begin{aligned}
(\rho_{P} \otimes_{k} \mathrm{id}_{B})\circ \rho_{P} & =
(\mathrm{id}_{P} \otimes_{k} \Delta_{B})\circ \rho_{P} \\
(\mathrm{id}_{P}\otimes_{k} \varepsilon_{B}) \circ \rho_{P} & =
  \mathrm{id}_{P}
\end{aligned}
$$ 
where $\Delta_{B}$ (respectively
  $\varepsilon_{B}$) denotes the comultiplication (respectively the
  counit) of $B$, and we have implicitly identified $P$ with $P\otimes
  k$ via the obvious $M$-isomorphism $P\otimes_{k} k \to P$;
\item morphisms are given by $M$-morphisms commuting with the
  structure maps $\rho$. 
\end{itemize}
The category \mbox{$B-\cdg_{M}$} comes equipped with a left adjoint
forgetful functor  $B-\cdg_{M} \longrightarrow M$, whose
right adjoint sends an object $X\in M$ to $X\otimes B$ endowed with
its natural $B$-comodule structure. The multiplication in $B$ endows
$B-\cdg_{M}$ with a natural symmetric monoidal structure for which the
forgetful functor $B-\cdg_{M} \longrightarrow M$ becomes a symmetric
monoidal functor.

%\begin{rmk}
%Whenever $M=C(k)$ there exists model structure on $B-\cdg_{M}$ (see for instance \cite[Theorem 2.5.17]{hov}, 
%where this is done for chain complexes an can be adapted to cochain complexes without problem) such that 
%the above adjunction is a Quillen adjunction (\cite[Proposition 2.5.20]{hov}). 
%It is probably still true for a general $M$ but we don't need it in this paper as in the example we consider 
%$B-\cdg_{M}$ can be identified with the category of graded mixed modules, for which an explicit model structure is given below. 
%\end{rmk}

We will be especially interested in the case where
$B=k[t,t^{-1}]\otimes_{k}k[\epsilon]$ defined as follows.  Here
$k[\epsilon]:=Sym_{k}(k[1])$ is the free commutative $k$-dg-algebra
generated by one generator $\epsilon$ in cohomological degree $-1$,
and $k[t,t^{-1}]$ is the usual commutative algebra of functions on
$\mathbb{G}_{m}$ (so that $t$ sits in degree $0$).  The
comultiplication on $B$ is defined by the dg-algebra map
$$ 
\xymatrix@R-1.5pc{ \Delta_{B} : \hspace{-1pc} &  B \ar[r] & B \otimes_{k}B \\ 
  & t\equiv t\otimes 1 \ar@{|->}[r] & (t\otimes 1) \otimes (t \otimes 1)\equiv
  t\otimes t \\ 
& \epsilon \equiv 1\otimes \epsilon \ar@{|->}[r] &
  (1\otimes \epsilon) \otimes (1\otimes 1) + (t\otimes 1) \otimes
  (1\otimes \epsilon)\equiv \epsilon \otimes 1 + t\otimes \epsilon }
$$
where $\equiv$ is used as a concise, hopefully clear notation for
canonical identifications.
%sending $\epsilon$ to $\epsilon \otimes 1 + 1 \otimes t.\epsilon$ and 
%$t$ to $t\otimes t$. 
Together with the counit dg-algebra map 
$$
\varepsilon_{B}: B
\longrightarrow k \,\,, \, t\mapsto 1\, ,\, \epsilon \mapsto 0 \, ,
$$
$B$ becomes a commutative $k$-linear Hopf dg-algebra.

\begin{rmk} Note that $B$ can be identified geometrically with the 
dg-algebra of functions on the affine group stack $\mathbb{G}_{m}
\ltimes \Omega_{0}\mathbb{G}_{a}$, semi-direct product of
$\mathbb{G}_{m}$ with $\Omega_{0}\mathbb{G}_{a} = K(\mathbb{G}_{a},
-1) = \mathbb{G}_{a}[-1]$ induced by the natural action of the
multiplicative group on the additive group. This is similar to \cite[Remark~1.1]{ptvv} where the algebra of functions on
$\mathbb{G}_{m} \ltimes B\mathbb{G}_{a} = \mathbb{G}_{m} \ltimes
\mathbb{G}_{a}[1]$ was used instead. In fact these two Hopf dg-algebras have
equivalent comodule theories and can be used interchangeably (see
Remark~\ref{r1}). This observation will not be used in any essential way in the rest of the paper.
\end{rmk}

Recall from above that $k[t,t^{-1}]$ co-acts on $k[\epsilon]$ and that $B$ is their semi-direct (co-)product. 
Therefore $k[\epsilon]$ is a coalgebra in $\mathbb{G}_m$-modules and $B-\cdg_{C(k)}$ is identified with the 
category of $k[\epsilon]$-comodules in $\mathbb{G}_m$-modules. 
Since $k[\epsilon]$ is dualisable, then we have a dual algebra $k[\epsilon]^\vee=k[e]$ in $\mathbb{G}_m$-modules: 
it is the free commutative algebra on a generator $e$ of both degree and $\mathbb{G}_m$-weight one. 
Therefore the category of $B$-comodules $B-\cdg_{C(k)}$ identifies
naturally with the category of \emph{graded mixed complexes of $k$-dg-modules}.

\begin{rmk}
Note that $k[\epsilon]^\vee$ is the universal enveloping algebra of the abelian (and thus nilpotent) dg-Lie algebra $k[-1]$ in $\mathbb{G}_m$-mod.
Its exponentiation is the group $\mathbb{G}_m$-stack $\mathbb{G}_a[-1]$. Therefore the category $B-\cdg_{C(k)}$ is a model for 
the $\s$-category of representations of $\mathbb{G}_{m}\ltimes \mathbb{G}_{a}[-1]$. This observation will not be used in any essential way in the rest of the paper.
\end{rmk}

More precisely, by a \emph{graded mixed complexe of
  $k$-dg-modules} we mean a family of $k$-dg-modules
$\{E(p)\}_{p\in \mathbb{Z}}$, together with families of morphisms
$$
\epsilon : E(p) \longrightarrow E(p+1)[1],
$$ 
such that
$\epsilon^2=0$. The identification between $B-\cdg_{C(k)}$ and the category of graded mixed complexes, actually an isomorphism of
categories, can be made directly by observing that the co-restriction functor 
$$
p_*:B-\cdg_{C(k)} \to k[t, t^{-1}]-\cdg_{C(k)} 
$$ 
along the
coalgebra map $p: B \to k[t, t^{-1}]$ (sending $\epsilon$ to $0$),
yields the usual $C(k)$-isomorphism $\oplus_{p\in \mathbb{Z}}E(p) \to
E$, where
$$E(p):= \rho_{p_{*}E}^{-1}(E \otimes_{k} k\cdot t^{p})\, ,
\, p\in \mathbb{Z}
$$ 
or, equivalently, 
$$
E(p):= \rho_{E}^{-1}(E
\otimes_{k} (k\cdot t^{p} \oplus k[t,t^{-1}]\epsilon)\, , \, p\in
\mathbb{Z}.
$$ 
Note that the morphism $\epsilon : E(p) \longrightarrow
E(p+1)[1]$ is then defined by sending $x_i \in E(p)^{i}$ to the image
of $x_i$ under the composite map 
$$
\xymatrix{E \ar[r]^-{\rho} &
  E\otimes_{k}B \ar[r]^-{\textrm{pr}} & E\otimes_{k} k\cdot
  t^{p+1}\epsilon}.
$$ 
Therefore, objects in $B-\cdg_{C(k)}$ will be
often simply denoted by $E=\oplus_{p}E(p)$, and the corresponding
mixed differential by $\epsilon$.

In order to avoid confusion, we will refer to the decomposition
$E=\oplus_{p}E(p)$ as the \emph{weight decomposition}, and refer to
$p$ as the \emph{weight degree} in order to distinguish it from the
\emph{cohomological} or \emph{internal degree}.

\begin{rmk}\label{r1}
Note that here we have adopted a convention opposite to the one in 
\cite[1.1]{ptvv}: the category $B-\cdg_{C(k)}$ of graded mixed
complexes introduced above, is naturally equivalent to the category of
graded mixed complexes used in \cite[1.1]{ptvv} where the mixed
structures decrease the cohomological degrees by one.  An explicit
equivalence is given by sending an object $\oplus_{p}E(p)$ in
$B-\cdg_{C(k)}$ to $\oplus_{p}(E(p)[2p])$ together with its natural
induced mixed structure (which now decreases the cohomological degree
by $1$).
\end{rmk}

\

More generally, the category of \emph{graded mixed objects in} $M$ is
defined to be $B-\cdg_{M}$, the category of $B$-comodules in $M$, with
$B=k[t,t^{-1}]\otimes_{k}k[\epsilon]$, and will be denoted by
$\epsilon-M^{gr}$. Its objects consist of

\begin{itemize}
\item $\mathbb{Z}$-families  $\{E(p)\}_{p\in \mathbb{Z}}$ of objects of $M$,
\item together with morphisms in $M$ 
$$
\epsilon\equiv  \left\{\epsilon_{p} : E(p) \longrightarrow E(p+1)[1]
\right\}_{p\in \mathbb{Z}},
$$ 
where for $P \in M$ and $n\in \mathbb{Z}$ we define $P[n] :=
P^{k[-n]}$ using the (cotensored) $C(k)$-enrichment, and we require
that 
$\epsilon^2=0$, i.e. that the composition
$$
\xymatrix{E(p) \ar[r]^-{\epsilon_{p}} & E(p+1)[1]
  \ar[rr]^-{\epsilon_{p+1}[1]} & & E(p+2)[2] }
$$ 
is zero for any $p\in \mathbb{Z}$.
\end{itemize}
Note that, by adjunction, $\epsilon_{p}$ can also be specified by
giving a map $E(p)\otimes_{k} k[-1] \to E(p+1)$ in $M$ or,
equivalently, a map $k[1] \to \underline{\textrm{Hom}}(E(p), E(p+1))$
in $C(k)$, (where $\underline{\textrm{Hom}}$ denotes the
$C(k)$-enriched hom in $M$). The morphisms $\epsilon$ will 
sometimes be called \emph{mixed maps} or \emph{mixed differentials},
following the analogy with the case $M=C(k)$.

We will consider on $\epsilon-M^{gr}=B-\cdg_{M}$ the following model structure.
First of all, the category $M^{gr}:= \prod_{p\in \mathbb{Z}} M$ is naturally a
symmetric monoidal model category with weak equivalences (respectively
cofibrations, respectively fibrations) defined component-wise, and a
monoidal structure defined by
$$
(E\otimes 
E')(p):=\bigoplus_{i+j=p}E(i)\otimes E'(j)
$$ 
where $\oplus$ denotes
the coproduct in $M$, and the symmetry constraint \emph{does not}
involve signs, and simply consists in exchanging the two factors in
$E(i)\otimes E'(j)$. It is easy to check, 
using our standing assumptions $(1)-(5)$ on $M$, that
$\epsilon-M^{gr}$ comes then equipped with a combinatorial symmetric
monoidal model category structure for which the equivalences and
cofibrations are defined through the forgetful functor
$$
\epsilon-M^{gr} \longrightarrow M^{gr}.
$$
Again the symmetric monoidal structure
on $\epsilon-M^{gr}$ can be described on the level of graded objects by
the formula
$(E\otimes E')(p):=\oplus_{i+j=p}E(i)\otimes E'(j)$ where $\oplus$
denotes the coproduct in $M$, and 
again the symmetry constraint \emph{does not}
involve signs, and simply consist of the exchange of the two factors
in $E(i)\otimes E'(j)$.  
The mixed differentials on $E\otimes E'$ are then defined by the usual
formula, taking the sums (i.e. coproducts) of  
all maps 

$$
\xymatrix@R-1.1pc{
\epsilon \otimes 1 + 1 \otimes \epsilon' : \hspace{-1.5pc} &
E(i) \otimes E'(j) \ar[r] & (E(i+1)[1]\otimes E'(j)) \bigoplus
(E(i)\otimes E'(j+1)[1]) \ar@{=}[d] \\
& & \left( (E(i+1)\otimes E'(j)) \bigoplus
(E(i)\otimes E'(j+1)) \right)[1]
}
$$
As a symmetric monoidal model category $\epsilon-M^{gr}$ again
satisfies all of our standing assumptions $(1)-(5)$, 
and the forgetful functor $\epsilon-M^{gr} \longrightarrow M^{gr}$
comes equipped with a natural  
symmetric monoidal structure. 

Note that  $\epsilon-M^{gr}$ is also an  $\epsilon-C(k)^{gr}$-enriched
symmetric monoidal model category.  Let us just briefly define the
graded mixed complex $\underline{\mathrm{Hom}}^{gr}_{\epsilon}(E,F)$,
for $E, F \in \epsilon-M^{gr}$, leaving the other details and
properties of this enrichment to the reader. We define  
\begin{itemize}
\item $\underline{\mathrm{Hom}}^{gr}_{\epsilon}(E,F)(p) := \prod_{q
  \in \mathbb{Z}}  \underline{\mathrm{Hom}}_{k}(E(q), F(q+p))$, for
  any $p \in \mathbb{Z}$ 
\item the mixed differential $\epsilon_{p}:
  \underline{\mathrm{Hom}}^{gr}_{\epsilon}(E,F)(p) \to
  \underline{\mathrm{Hom}}^{gr}_{\epsilon}(E,F)(p+1)[1]$ as the map
  whose $q$-component 
$$
\prod_{q' \in \mathbb{Z}}
  \underline{\mathrm{Hom}}_{k}(E(q'), F(q'+p))\longrightarrow
  \underline{\mathrm{Hom}}_{k}(E(q), F(p+q+1))[1]\simeq
  \underline{\mathrm{Hom}}_{k}(E(q), F(p+q+1)^{k[-1]} 
$$ 
is given by
  the sum $\alpha + \beta$  
where 
$$
\xymatrix@R-0.6pc{ 
\prod_{q' \in \mathbb{Z}}
  \underline{\mathrm{Hom}}_{k}(E(q'), F(q' + p)) \ar[dr]_{\alpha}
  \ar[r]^-{\textrm{pr}} 
  &  \underline{\mathrm{Hom}}_{k}(E(q), F(q + p)) \ar[d]^-{\alpha'}
  \\
&
  \underline{\mathrm{Hom}}_{k}(E(q), F(p+q+1))[1]
}
$$ 
$\alpha'$ being adjoint to the composite 
$$
\xymatrix@1{\underline{\mathrm{Hom}}_{k}(E(q), F(q + p))
  \otimes E(q) \ar[r]^-{\textrm{can}} & F(q+p) \ar[r]^-{\epsilon_{F}}
  & F(q+p+1)^{k[-1]},
}
$$ 
and 
$$
\xymatrix@R-0.6pc{
\prod_{q' \in
    \mathbb{Z}}  \underline{\mathrm{Hom}}_{k}(E(q'), F(q' +
  p)) 
\ar[dr]_-{\beta}\ar[r]^-{\textrm{pr}} &  
\underline{\mathrm{Hom}}_{k}(E(q+1), F(q
  + 1+ p)) \ar[d]^-{\beta'} \\
& \underline{\mathrm{Hom}}_{k}(E(q),
  F(p+q+1))[1]
}
$$ 
$\beta'$ being adjoint to the composite 
$$
\xymatrix@R-0.6pc{
  \underline{\mathrm{Hom}}_{k}(E(q+1), F(q + 1+ p)) \otimes_{k}
  (E(q)\otimes_{k} k[-1]) \ar[d]_-{id \otimes \epsilon_{E}} & \\
  \underline{\mathrm{Hom}}_{k}(E(q+1), F(q+p+1)) \otimes_{k} E(q+1) 
\ar[r]^-{\textrm{can}} & F(q+p+1). 
} 
$$ 
\end{itemize}
Therefore, as already observed for $M$, the category
$Comm(\epsilon-M^{gr})$, of commutative and unital monoids
in graded mixed objects in $M$, is a combinatorial model category
whose weak equivalences and fibrations are defined through the
forgetful functor $Comm(\epsilon-M^{gr}) \longrightarrow
\epsilon-M^{gr}$ (\cite[Proposition~4.5.4.6]{lualg}).

\subsection{$\s$-Categories setting}\label{infty}
We will denote by $\C:=L(M)$ the $\s$-category obtained from $M$ by
inverting the equivalences (see \cite[\S 2.1]{toenems}).  Since $M$ is
a stable model category (Proposition~\ref{Mstable}), $\C$ is
automatically a stable $\s$-category. Moreover, as explained in
\cite[\S 2.1]{chern}, $\C$ possesses a natural induced symmetric
monoidal structure. An explicit model for $\C$ is the simplicial
category of fibrant and cofibrant objects in $M$, where the simplicial
sets of morphisms are defined by applying the Dold-Kan construction to
the truncation in non-negative degrees of the complexes of morphisms
coming from the $C(k)$-enrichment (see \cite{tab}).  The symmetric
monoidal structure on $\C$ is harder to describe explicitly, and we
will not discuss it here since it will not be used in an essential way
in what follows. Parallel results hold for $\C^{gr}:=L(M^{gr})$. We
refer to \cite[\S 2.1]{chern} for more about localization of symmetric
monoidal model categories.

We recall from Section~\ref{1.1} that $Comm(M)$ is the model category of
commutative monoids in $M$, and we  
let 
$$
\cdga_{\C}:=L(Comm(M)),
$$ 
to be the $\s$-category obtained by localizing $Comm(M)$ along weak
equivalences. Note that  
our notation suggests that  $\cdga_{\C}$  is the $\s$-category of
\emph{commutative dg-algebras 
internal to $\C$} in the sense of \cite{lualg}. This is justified by
the existence of a natural  
equivalence of $\s$-categories 
$$
L(Comm(M)) \simeq Comm(LM).
$$
This equivalence is a consequence of \cite[Theorem~4.5.4.7]{lualg},
since by Proposition~\ref{cofibforgettocofib} the forgetful functor
$Comm(M) \to M$ preserves fibrant-cofibrant objects.

The Quillen adjunction $\epsilon-M^{gr} \longleftrightarrow M^{gr}$
(see Section~\ref{1.1}) 
induces an adjunction of $\s$-categories 
$\epsilon-\C^{gr}= L(\epsilon-M^{gr}) \longleftrightarrow \C^{gr}:=L(M^{gr})$.

\begin{df}\label{d3}
  The symmetric monoidal $\s$-category $\epsilon-\C^{gr}$ of
  \emph{graded mixed objects in $\C$} is defined as $\epsilon-\C^{gr}
  := L(\epsilon-M^{gr})$.  The $\s$-category
  $\epsilon-\cdga^{gr}_{\C}$ of \emph{graded mixed commutative
    dg-algebras in $\C$} is defined as $\epsilon-\cdga^{gr}_{\C} :=
  L(Comm(\epsilon-M^{gr}))$.
\end{df}

\bigskip 

\noindent
Note that, again,  \cite[Theorem~4.5.4.7]{lualg} and
Proposition~\ref{cofibforgettocofib} imply that we have a natural
equivalence of $\s$-categories 
$$
Comm(\epsilon-\C^{gr}) \simeq L(Comm(\epsilon-M^{gr})),
$$
and so $\epsilon-\cdga_{\C}^{gr}$ can also be considered as the
$\s$-category  
of commutative monoid objects in the symmetric monoidal $\s$-category 
$\epsilon-\C^{gr}$.
We have an adjunction of $\s$-categories
$$\epsilon-{\C}^{gr} \longleftrightarrow \epsilon-\cdga_{\C}^{gr},$$
where the right adjoint forgets the algebra structure.

At a more concrete level, 
objects in $\epsilon-\cdga^{gr}_{\C}$ can be described as
commutative monoids in $\epsilon-M^{gr}$, i.e. 
as the following collections of data
\begin{enumerate}
\item a family of objects $\{A(p) \in M\}_{p\in \mathbb{Z}}$. 
\item a family of morphisms $\epsilon\equiv  \{\epsilon_{p} : A(p)
  \longrightarrow A(p+1)[1] \}_{p\in \mathbb{Z}}$, satisfying  
$\epsilon^2=0$.
\item a family of multiplication maps
$$
\{ A(p) \otimes A(q) \longrightarrow A(p+q)\}_{(p,q)\in \mathbb{Z}
  \times \mathbb{Z}},
$$
which are associative, unital, graded commutative, and compatible with
the maps  
$\epsilon$ above.  
\end{enumerate}

\

\begin{rmk} 
Since $\C$ is stable, we have equivalences in $\C$ 
$$
\Sigma X
\simeq X\otimes_{k} k[1] \simeq X[1] = X^{k[-1]} \simeq \Omega^{-1}
X
$$
 where the the tensor and cotensor products are to be understood in
the $\s$-categorical sense (i.e. in the derived sense when looking at
$M$). These equivalences are natural in $X \in \C$. In particular
there is no ambiguity about what $X[n]$ means in $\C$, for any $n \in
\mathbb{Z}$: $X[n] \simeq X\otimes_{k} k[n] \simeq X^{k[-n]}$. Beware
that these formulas are not correct, on the nose, in $M$, unless $X$
is fibrant and cofibrant.
\end{rmk}

\subsection{De Rham theory in a relative setting}

Let $M$ be a symmetric monoidal model category satisfying the
conditions from Section~\ref{1.1}.  We denote the corresponding
$\s$-category by $\C$. As above we have the category $\epsilon-M^{gr}$
of graded mixed objects in $M$ and the corresponding $\s$-category
$\epsilon-\C^{gr}$ of graded mixed objects in $\C$.

Since $\mathbf{1}_M$ is cofibrant in $M$, there is a natural Quillen
adjunction 
$$
-\otimes \mathbf{1}_M: C(k) \longleftrightarrow M :
\underline{Hom}(\mathbf{1}_M ,-),
$$
where the left adjoint sends an object $x \in C(k)$ to 
$x\otimes \mathbf{1} \in M$ (tensor enrichment of $M$ over $C(k)$), 
while the right adjoint is given by the $C(k)$-hom enrichment.
The induced adjunction
on the corresponding $\s$-categories will be denoted by 
$$
-\otimes \mathbf{1}_M: \dg_k=\mathrm{L}(k) \longleftrightarrow \C :
|-|:=\mathbb{R}\underline{Hom}(\mathbf{1}_M,-).
$$
Since $\mathbf{1}_M$ is a comonoid object in $M$, the right Quillen
functor $\underline{Hom}(\mathbf{1}_M ,-)$ is lax symmetric
monoidal. Therefore, 
 we get similar adjunctions at the commutative monoids and graded
 mixed level (simply denoted through the corresponding right adjoints) 
$$
\xymatrix@R-1.5pc{
\cdga_{k} \ar@{<->}[r] &  \cdga_{\C} : |-| \\
\epsilon-\cdga_{k} \ar@{<->}[r] &  \epsilon-\cdga_{\C} : |-| \\
\epsilon-\dg^{gr}_{k} \ar@{<->}[r] &   \epsilon-{\C}^{gr} : |-| \\
\epsilon-\cdga^{gr}_{k}\ar@{<->}[r] &  \epsilon-\cdga_{\C}^{gr} :
|-|
}
$$

\

\begin{df}\label{dreal}
The right adjoint $\s$-functors $|-|$ defined above will be called the
\emph{realization $\s$-functors}.
\end{df}

\

\begin{rmk}\label{enhanced} Note that 
if $A\in \cdga_{\C}$ and $P \in A-Mod_{\C}$, then $|P| \in
|A|-\dg_{k}$, and we get a refined realization functor 
$$
|-|: A-Mod_{\C} \longrightarrow |A|-\dg_{k}.
$$
\end{rmk}

\subsubsection{Cotangent complexes.}
We start with the notion of a cotangent complex for a commutative
dg-algebra 
inside $\C$. For $A \in \cdga_{\C}$ we have an $\s$-category
$A-Mod_{\C}$ of $A$-modules in $\C$. If the object $A$ corresponds to
$A\in Comm(M)$, the $\s$-category $A-Mod_{\C}$ can be defined as the
localization of the category $A-Mod_{M}$, of $A$-modules in $M$, along
the equivalences.  The model category $A-Mod_{M}$ is a stable model
category and thus $A-Mod_{\C}$ is itself a presentable stable
$\s$-category. As $A$ is commutative, $A-Mod_{M}$ is  a
symmetric monoidal category in a natural way, for the tensor product
$-\otimes_{A}-$ of $A$-modules. This makes $A-Mod_{M}$ a symmetric
monoidal model category which satisfies again the conditions $(1)-(5)$
(see Proposition \ref{1-5formod}). The corresponding $\s$-category
$A-Mod_{\C}$ is thus itself a symmetric monoidal presentable and
stable $\s$-category.

For an $A$-module $N \in A-Mod_{M}$, we endow $A \oplus N$ with the
trivial square zero structure, as in \cite[1.2.1]{hagII}. We 
denoted the coproduct in $M$ by $\oplus$; note however that since
$A-Mod_{M}$ is stable, any finite coproduct is identified with the
corresponding finite product. The projection $A\oplus N \rightarrow A$
defines an object $A\oplus N \in Comm(M)/A$, as well as an object in
the comma $\s$-category $A\oplus N \in \cdga_{\C}/A$ of commutative
monoids in $\C$ augmented to $A$.

\begin{df}\label{d1}
In the notations above, the \emph{space of derivations from $A$ to
  $N$} is defined by 
$$
Der(A,N):=\mathsf{Map}_{\cdga_{\C}/A}(A,A\oplus N) \in \T.
$$
\end{df}

\

For a fixed $A \in \cdga_{\C}$, the construction $N \mapsto Der(A,N)$
can be naturally promoted to  
an $\s$-functor 
$$ 
Der(A,-) : A-Mod_{\C} \longrightarrow \T.
$$

\begin{lem}\label{l1}
For any $A\in A-Mod_{\C}$, the $\s$-functor $Der(A,-)$ is corepresentable 
by an object $\mathbb{L}^{int}_{A} \in A-Mod_{\C}$.
\end{lem}
\textbf{Proof:}  This is a direct application of
\cite[Proposition~5.5.2.7]{lutop}, since $A-Mod_{\C}$ and $\T$ are both
presentable $\infty$-categories, 
and the $\s$-functor $Der(A,-)$ is accessible and commutes with small limits. 
\hfill $\Box$ 

\

\begin{df}\label{d2}
Let $A \in \cdga_{\C}$. 
\begin{enumerate} 
\item The object $\mathbb{L}^{int}_{A} \in A-Mod_{\C}$ is called the
\emph{cotangent complex of $A$, internal to $\C$}. 
\item The \emph{absolute cotangent complex} (or simply
\emph{the cotangent complex} of $A$) is 
$$
\mathbb{L}_{A}:=|\mathbb{L}_{A}^{int}| \in \dg_k,
$$
where $|-| : \C \longleftrightarrow \dg_k$ 
is the realization $\s$-functor of definition \ref{dreal}.
\end{enumerate}
\end{df}

\

\begin{rmk} Both $A-Mod_{\C}$ and $\cdga_{\C}/A$ are presentable
  $\infty$-categories, and the $\s$-functor $N \mapsto A\oplus N$ is
  accessible and preserves limits, therefore  
(\cite[Corollary~5.5.2.9]{lutop}) it admits a 
left adjoint  $\mathbf{L}^{int}: \cdga_{\C}/A \to A-Mod_{\C}$, and we
have $\mathbb{L}^{int}_{A} = \mathbf{L}^{int}(A)$. 
\end{rmk}

The construction $A \mapsto \mathbb{L}^{int}_{A}$ possesses all
standard and expected properties. For a morphism $A \longrightarrow B$
in $\cdga_{\C}$, we have an adjunction of $\s$-categories
$$
B\otimes_{A} - : A-Mod_{\C} \longleftrightarrow B-Mod_{\C}  :
\textsf{forg} 
$$
 where $\textsf{forg}$ is the forgetful $\s$-functor, and we have a
natural morphism $B\otimes_{A} \mathbb{L}^{int}_{A} \longrightarrow
\mathbb{L}^{int}_{B}$ in $B-Mod_{\C}$.  The cofiber of this morphism,
in the $\s$-category $B-Mod_{\C}$, is denoted by
$\mathbb{L}^{int}_{B/A}$, and is called the \emph{relative cotangent
  complex of $A\rightarrow B$ internal to $\C$}. We have, by
definition, a fibration-cofibration sequence of $B$-modules
$$
\xymatrix@1{B\otimes_{A} \mathbb{L}^{int}_{A} \ar[r] & \mathbb{L}^{int}_{B} 
\ar[r] & \mathbb{L}^{int}_{B/A}.}
$$ 
Moreover, the internal cotangent complex is compatible with
push-outs in $\cdga_{\C}$, in the following sense. For a cocartesian
square of objects in $\cdga_{\C}$
$$
\xymatrix{
A \ar[r] \ar[d] & B \ar[d] \\
C \ar[r] & D,}
$$
the induced square of objects in $B-Mod_{\C}$
$$
\xymatrix{
D\otimes_{A} \mathbb{L}^{int}_{A} \ar[r] \ar[d] & D\otimes_{B} \mathbb{L}^{int}_{B} \ar[d] \\
D\otimes_{C} \mathbb{L}^{int}_{C} \ar[r] & \mathbb{L}^{int}_{D}}
$$
is again cocartesian. 
%[\texttt{We should give a sketch of proof here or loosely refer to hagII}]% \\

\begin{rmk}\label{r2}
The above definition of an internal cotangent complex gives the usual
cotangent complex of commutative dg-algebras $A$ over $k$ when one
takes $M=C(k)$. More precisely, for $M=C(k)$, the $\s$-functor $|-|$
is isomorphic to the forgetful functor $\mathsf{forg}: A-Mod \to
C(k)$, and we have $\mathsf{forg}(\mathbb{L}_{A}^{int}) \simeq
\mathbb{L}_{A}$ in $C(k)$.
\end{rmk}

\subsubsection{De Rham complexes.}\label{1.3.2}

We have defined, for any object $A\in \cdga_{\C}$ a cotangent complex
$\mathbb{L}^{int}_{A} \in A-Mod_{\C}$. We will now show how to
associate to any $A\in \cdga_{\C}$ its de Rham complex. As for
cotangent complexes we will have two versions, an internal de Rham
complex $\DR^{int}(A)$, and an absolute one $\DR(A)$, respectively
related to $\mathbb{L}_{A}^{int}$ and $\mathbb{L}_{A}$. The first
version, $\DR^{int}(A)$ will be a graded mixed cdga in $\C$, whereas
$\DR(A)$ will be a graded mixed cgda in $\dg_k$. These of
course will be related by the formula
$$
\DR(A)=|\DR^{int}(A)|
$$
where $|-| : \C \longrightarrow \dg_k$ (or equivalently, $|-|:
\epsilon-\cdga_{\C}^{gr} \longrightarrow \epsilon-\cdga^{gr}_{k} $) is
the realization $\s$-functor of Definition~\ref{dreal}.

We recall from Section~\ref{infty}
that a mixed graded commutative dg-algebra $A$ in $\C$ can 
be described as the following data

\begin{enumerate}
\item a family of objects $\{A(p) \in M\}_{p\in \mathbb{Z}}$. 
\item a family of morphisms $\epsilon\ = \{\equiv \epsilon_{p} : A(p)
  \longrightarrow A(p+1)[1] \}_{p\in \mathbb{Z}}$, satisfying
  $\epsilon_{p+1}[1] \circ \epsilon_{p}=0$.

\item a family of multiplication maps
$$
\{ A(p) \otimes A(q) 
\longrightarrow A(p+q)\}_{(p,q)\in \mathbb{Z} \times \mathbb{Z}},
$$ 
which are associative, unital, graded commutative, and compatible
with the maps $\epsilon$.
\end{enumerate}

The (formal) decomposition $A=\oplus A(p)$ will be called the
\emph{weight decomposition}, and $A(p)$ the \emph{weight $p$ part} of
$A$.

By point 3. above, for $A \in \epsilon-\cdga_{\C}^{gr}$, the weight
$0$ object $A(0)\in \C$ comes equipped with an induced commutative
monoid structure and thus defines an object $A(0)\in \cdga_{\C}$. This
defines an $\s$-functor
$$
(-)(0): \epsilon-\cdga^{gr}_{\C} \longrightarrow \cdga_{\C}
$$
which picks out the part of  weight degree $0$ only.
The compatibility of the multiplication with the mixed  structure
$\epsilon$ expresses in particular
that the property that the morphism $A(0) \longrightarrow A(1)[1]$ is
a derivation of the commutative monoid $A(0)$ 
with values in $A(1)[1]$.
%[\texttt{give a proof sketch of this statement}]% 
We thus have a natural induced morphism in the stable $\s$-category 
of $A(0)$-modules
$$\varphi_{\epsilon} : \mathbb{L}^{int}_{A(0)}[-1] \longrightarrow A(1).$$

\begin{prop}\label{l2}
The $\s$-functor 
$$ 
(-)(0): \epsilon-\cdga^{gr}_{\C} \longrightarrow \cdga_{\C},
$$
has a left adjoint
$$
\mathbf{DR}^{int} :   \cdga_{\C} \longrightarrow
\epsilon-\cdga^{gr}_{\C}.
$$
\end{prop}
\textbf{Proof.} This is an application of the adjoint functor theorem
(\cite[Corollary~5.5.2.9]{lutop}). We just need to show that the
$\s$-functor $A \mapsto A(0)$ is accessible and preserves limits. For
this, we use the commutative diagram of $\s$-categories
$$
\xymatrix{
\epsilon-\cdga^{gr}_{\C} \ar[r] \ar[d] & \cdga_{\C} \ar[d] \\
\epsilon-\C^{gr} \ar[r] & \C,
}
$$
where the vertical $\s$-functors forget the commutative monoid
structures and the horizontal $\s$-functors select the parts of
weight $0$. These vertical $\s$-functors are conservative and commute
with all limits. We are thus reduced to checking that the bottom
horizontal $\s$-functor $\epsilon-\C^{gr} \longrightarrow \C$
preserves limits. This last $\s$-functor has in fact an explicit left
adjoint, obtained by sending an object $X \in \C$, to the graded mixed
object $E$ defined by
$$
E(0)=X \qquad E(1)=X[-1] \qquad E(i)=0 \,\,\,\,  \forall i\neq 0,1,
$$
and with $\epsilon : E(0)\rightarrow E(1)[1]$ being the canonical
isomorphism $X[-1][1]\simeq   X$.  
\hfill $\Box$ 

\

\begin{df}\label{d4}
Let $A \in \cdga_{\C}$ be a commutative dg-algebra in $\C$. 
\begin{enumerate}
\item 
The \emph{internal de Rham object of $A$}
is the graded mixed commutative dg-algebra over $\C$ defined by
$$
\mathbf{DR}^{int}(A) \in \epsilon-\cdga^{gr}_{\C}.
$$
\item The
\emph{absolute de Rham object of $A$} (or simply the \emph{de Rham object})
is the graded mixed commutative dg-algebra over $k$ defined by
$$
\mathbf{DR}(A):=|\DR^{int}(A)| \in \epsilon-\cdga^{gr}_{k}
$$ 
where $|-|: \epsilon-\cdga_{\C}^{gr} \longrightarrow
\epsilon-\cdga^{gr}_{k} $ is the realization $\s$-functor of
Definition~\ref{dreal}. 
\end{enumerate}
\end{df}

\

\begin{rmk}\label{r3}
Abusing the language we will often refer to the de Rham objects
$\DR^{int}(A)$ and $\DR(A)$ as the (internal or absolute) \emph{de
  Rham complexes of $A$}, even though they are not just complexes but
a rather objects of $\epsilon-\cdga^{gr}_{\C}$ or of
$\epsilon-\cdga^{gr}_{k}$.
\end{rmk}

\

We will also need the following

\begin{df}\label{grmonmixalg}
Let $Comm(M)^{gr}$ be the category with objects $\mathbb{Z}$-indexed
families $\{A(n)\}_{n\in \mathbb{Z}}$ of objects in $Comm(M)$, and
morphisms $\mathbb{Z}$-indexed families $\{A(n) \to B(n)\}_{n\in
  \mathbb{Z}}$ of morphisms in $Comm(M)$. 

$Comm(M)^{gr}$ has a model structure with fibrations, weak
equivalences (and cofibrations) defined levelwise. Its localization
$L(Comm(M)^{gr})$ along weak equivalences will be denoted by
$\cdga_{\C}^{gr}$ and called the \emph{$\s$-category of graded
  (non-mixed) commutative dg-algebras in} $\C$.
\end{df}

\

By definition, the de Rham object $\DR^{int}(A)$ comes equipped with
an adjunction morphism $A \longrightarrow \DR^{int}(A)(0)$ in
$\cdga_{\C}$.  Moreover, the structure of a mixed graded cdga on
$\DR^{int}(A)$ defines a derivation $\DR^{int}(A)(0) \longrightarrow
\DR^{int}(A)(1)[1]$, and thus a canonical morphism in the
$\s$-category of $\DR^{int}(A)(0)$-modules
$$\mathbb{L}^{int}_{A}\otimes_{A}\DR^{int}(A)(0) 
\longrightarrow \mathbb{L}^{int}_{\DR^{int}(A)(0)} \longrightarrow
\DR^{int}(A)(1)[1].
$$
Note that this is the same as a morphism 
$$
\mathbb{L}^{int}_{A}[-1] \longrightarrow \DR^{int}(A)(1)
$$
in the stable $\s$-category of $A$-modules.

This extends to a morphism in $\cdga_{\C}^{gr}$
$$
\phi_{A} : Sym_{A}(\mathbb{L}^{int}_{A}[-1]) \longrightarrow
\DR^{int}(A),
$$ 
where the grading on the left hand side is defined by letting
$\mathbb{L}^{int}_{A}[-1]$ be pure of weight $1$. Note that, by
construction, the morphism $\phi_{A} $ is natural in $A$.

\begin{prop}\label{p1}
For all $A \in \cdga_{\C}$ the above morphism
$$
\phi_{A} : Sym_{A}(\mathbb{L}^{int}_{A}[-1]) \longrightarrow
\DR^{int}(A)
$$
is an equivalence in $\cdga_{\C}^{gr}$. 
\end{prop}
\textbf{Proof.} The morphism $\phi_{A}$ is functorial in $A$, and
moreover, any commutative dg-algebra in $\C$ is a colimit of free
commutative dg-algebras (see, e.g. \cite[3.2.3]{lualg}).  It is
therefore enough to prove the following two assertions

\begin{enumerate}
\item The morphism $\phi_{A} : Sym_{A}(\mathbb{L}^{int}_{A}[-1])
  \longrightarrow \DR^{int}(A)$ is an equivalence when $A=Sym(X)$ is
  the free commutative dg-algebra over an object $X \in \C$.
\item The two $\s$-functors $A \mapsto
  Sym_{A}(\mathbb{L}^{int}_{A}[-1])$ and $A \mapsto \DR^{int}(A)$,
  from commutative dg-algebras in $\C$ to graded commutative algebras
  in $\C$, commute with all colimits.
\end{enumerate}

\

Proof of $1.$ \ Let $A=Sym(X) \in \cdga_{\C}$ be a free object. Explicitly
its de Rham object $\DR^{int}(A)$ can be described as follows. Let us
denote by $Y \in \epsilon-\C^{gr}$ the free graded mixed object over
$X$, the free graded mixed object functor being left adjoint to the
forgetful functor $\epsilon-\C^{gr} \longrightarrow \C$.  As already
observed, we have $Y(0)=X$, $Y(1)=X[-1]$, $Y(i)=0$ if $i\neq 0,1$, and
with the canonical mixed structure $X \simeq X[-1][1]$. The de Rham
object $\DR^{int}(A)$, is then the free commutative monoid object in
$\epsilon-\cdga^{gr}_{\C}$ over $Y$. We simply denote by $X\oplus
X[-1]$ the graded object in $\C$ obtained by forgetting the mixed
differential in $Y$. As forgetting the mixed structure is a symmetric
monoidal left adjoint, the graded commutative algebra underlying
$\DR^{int}(A)$ is thus given by
$$
\begin{aligned}
\DR^{int}(Sym(X)) & \simeq Sym(X\oplus X[-1]) \simeq 
Sym(X) \otimes Sym(X[-1]) \simeq 
Sym_{Sym(X)}(A\otimes X[-1]) \\ 
& \simeq Sym_{A}(\mathbb{L}^{int}_{A}[-1]),
\end{aligned}
$$ 
where $Sym : \C^{gr} \to \cdga^{gr}_{\C}$ denotes the left adjoint to
the forgetful functor. Note that, for $Y \in \C^{gr}$ sitting entirely
in weight $0$, $Sym(Y)$ sits entirely in weight $0$. On the other hand
if $Z \in \C$, and we write $Z(1) \in \C^{gr} $ for $Z$ sitting in
degree $1$, then $Sym (Z(1))$ coincides with $Sym (Z)$ with its
``usual'' full $\mathbb{N}$-weight-grading (with $Z$ sitting in weight
$1$). This proves $1.$.

\

Proof of $2.$ This follows because both $\s$-functors are obtained by
composition of various left adjoint $\s$-functors. Indeed, for the
case of $A \mapsto \DR^{int}(A)$ this is the composition of the
$\s$-functor $\DR^{int}$ from lemma \ref{l2} with the forgetful
$\s$-functor from $\epsilon-\cdga_{\C}^{gr} \longrightarrow
\cdga_{\C}^{gr}$ which are both left adjoints. For the second
$\s$-functor, we have, for any $B \in \cdga_{\C}^{gr}$, a natural
morphism of spaces
$$
\mathsf{Map}_{\cdga_{\C}^{gr}}(Sym_{A}(\mathbb{L}^{int}_{A}[-1]),B) 
\longrightarrow \mathsf{Map}_{\cdga_{\C}}(A,B(0)).
$$
The fiber of this map at a given morphism $A \rightarrow B(0)$, is
naturally equivalent to 
$\mathsf{Map}_{A-Mod_{\C}}(\mathbb{L}^{int}_{A}[-1],B(1))$.  By the
definition of the cotangent complex this fiber is also naturally
equivalent to $\mathsf{Map}_{\cdga_{\C}/B(0)}(A,B(0)\oplus B(1)[1])$.
This implies that, for a fixed $B\in \cdga_{\C}^{gr}$, the
$\s$-functor $A \mapsto
\mathsf{Map}_{\cdga_{\C}^{gr}}(Sym_{A}(\mathbb{L}^{int}_{A}[-1]),B)$
transforms colimits into limits, and thus that $A \mapsto
Sym_{A}(\mathbb{L}^{int}_{A}[-1])$, as an $\s$-functor $\cdga_{\C}
\rightarrow \cdga_{\C}^{gr}$ preserves colimits.  \
 \mbox{\ } \hfill 
$\Box$ 

\

\begin{rmk} Observe that $\phi_{A} : Sym_{A}(\mathbb{L}^{int}_{A}[-1])
  \longrightarrow \DR^{int}(A)$ is actually an equivalence in the
  under-category $A/ \cdga_{\C}^{gr}$ (where $A$ sits in pure weight
  $0$), simply by inducing the map $A\to \DR^{int}(A)$ using $\phi(A)$
  and the canonical map $A \to Sym_{A}(\mathbb{L}^{int}_{A}[-1])$.
\end{rmk}

\

An important corollary of the previous proposition is the existence of
a \emph{de Rham differential}, for any object $A \in \cdga_{\C}$.

\begin{cor}\label{c1}
For any object $A\in \cdga_{\C}$, the graded commutative dg-algebra
$Sym_{A}(\mathbb{L}^{int}_{A}[-1])$
possesses a canonical mixed structure making it into a mixed graded
commutative dg-algebra in $\C$. The corresponding mixed differential
is called the \emph{de Rham  
differential}.
\end{cor}

\

\begin{rmk} Note that, from the point of view of $\s$-categories 
(which is the point of view adopted in its statement), Corollary~\ref{c1}
  is almost tautological.  In fact, from this point of view, for a
  graded cdga $B$ in $M$, a mixed structure on $B$ means a \emph{weak}
  mixed structure, i.e.  a pair $(B',u)$, where $B'$ is a graded mixed
  cdga in $M$ and $u : B' \simeq B$ is an equivalence of graded
  cdga. This is the exact content of Cororllary~\ref{c1}.
\end{rmk}

\

\noindent \textbf{Relative $\DR^{int}$.} We conclude this subsection with the \emph{relative} version of
$\DR^{int}$. Let $A \in \cdga_{\C}$, and consider the $\s$-functor 
$$
(-)(0): A/\mecdga_{\C} \longrightarrow A/\cdga_{\C} \,\,\,\, C
\longmapsto C(0)
$$ 
where in $A/\mecdga_{\C}$, $A$ is considered as concentrated in pure
weight $0$ (hence with trivial mixed differential). 

\begin{prop}\label{DREL} For any $A \in \cdga_{\C}$, the $\s$-functor 
$$
(-)(0) : A/\mecdga_{\C} \longrightarrow A/\cdga_{\C}
$$ 
has a left adjoint, denoted as 
$$
\DR^{int}(-/A): A/\cdga_{\C} \longrightarrow A/\mecdga_{\C} \,\,\,\,\,
(A\to B) \longmapsto \DR^{int}(B/A)$$ 
\end{prop}
\noindent \textbf{Proof.} This is an application of the adjoint
functor theorem (\cite[Corollary~5.5.2.9]{lutop}), completely
analogous to the proof of Proposition~\ref{l2}. We leave the details
to the reader \hfill $\Box$

\

\noindent
Proceeding as in Proposition~\ref{p1}, we also get

\begin{prop}\label{p1rel}
For all $A \in \cdga_{\C}$ there is a canonical morphism
$$
\phi_{/A} : Sym_{B}(\mathbb{L}_{B/A}^{int}[-1]) \longrightarrow
\DR^{int}(B/A)
$$
is an equivalence in $A/\cdga_{\C}^{gr}$.  
\end{prop}

\

\noindent
Consider the $\s$-functor 
$$
\DR^{int} : \mathsf{Mor}(\cdga_{\C}) \longrightarrow \mecdga_{\C},
$$
sending a morphism $A \rightarrow B$ to $\DR^{int}(B/A)$.  This
$\s$-functor can be explicitly constructed as the localization along
equivalences of the functor
$$
DR^{str} : \mathsf{Cof}(cdga_{\C}) \longrightarrow
\epsilon-cdga_{\C}^{gr},
$$
from the category of cofibrations between cofibrant cdga to 
the category of graded mixed cdga, sending  
a cofibration $A \rightarrow B$ to 
$DR^{str}(B/A)=Sym_{B}(\Omega_{B/A}^{1}[-1])$, 
with mixed structure given by the de Rham differential. 
The following result gives a useful description of $\DR^{int}(B/A)$.

\begin{lem}\label{relativeDR} 
For the $\s$-functor 
$$
\DR^{int} : \mathsf{Mor}(\cdga_{\mathcal{M}})
\longrightarrow \mecdga_{\mathcal{M}} \,\,\,\,\,\,\, (A\to B)
\longmapsto \DR^{int}(B/A),
$$ 
we have an equivalence in
$A/\mecdga_{\mathcal{M}}$
$$
\DR^{int}(B/A) \simeq \DR^{int}(B)\otimes_{\DR^{int}(A)} A
$$ 
where $A$ is concentrated in weight $0$ (hence, with trivial mixed
differential), and the rhs denotes the obvious pushout in the category
$\mecdga_{\mathcal{M}}$.
\end{lem}
\textbf{Proof.} We have to prove that the $\s$-functor 
$$
A/\cdga_{\mathcal{M}} \longrightarrow A/\mecdga_{\mathcal{M}}
\,\,\,\,\,\,\, (A\to B) \longmapsto \DR^{int}(B)\otimes_{\DR^{int}(A)}
A
$$
 is left adjoint to the functor sending $C$ to $C(0)$.
Now,
 $$\mathsf{Map}_{A/\mecdga_{\mathcal{M}}}(X, C) \simeq
\mathsf{Map}_{\mecdga_{\mathcal{M}}}(X, C)
\times_{\mathsf{Map}_{\mecdga_{\mathcal{M}}}(A,C)} \{*\}
$$ 
where the map $\{*\} \to \mathsf{Map}_{\mecdga_{\mathcal{M}}}(A,C)$ is
induced by the structure map $\rho: A \to C$, defining $C$ as an
object in $A/\mecdga_{\mathcal{M}}$. Taking $X=
\DR^{int}(B)\otimes_{\DR^{int}(A)} A$, and using the shortcut notation
$\mathsf{Map}:=
\mathsf{Map}_{\cdga_{\C}}$, we thus get 
$$
\begin{aligned}
\mathsf{Map}_{A/\mecdga_{\mathcal{M}}}(& 
\DR^{int}(B)\otimes_{\DR^{int}(A)} A, C) \\
& \simeq 
(\mathsf{Map}(B,C(0)) \times_{\mathsf{Map}(A,C(0))}
\mathsf{Map}_{\mecdga_{\mathcal{M}}}(A,C))
\times_{\mathsf{Map}_{\mecdga_{\mathcal{M}}}(A, C)} \{*\} \\
& \simeq \mathsf{Map}(B,C(0)) \times_{\mathsf{Map}(A,C(0))}\{*\} 
\end{aligned}
$$ 
where the map $\{*\} \to \mathsf{Map}(A, C(0))$ is induced by the
weight $0$ component $\rho(0)$ of
$\rho$. Therefore 
$$
\begin{aligned}
\mathsf{Map}_{A/\mecdga_{\mathcal{M}}}(
\DR^{int}(B)\otimes_{\DR^{int}(A)} A, C) & \simeq \mathsf{Map}(B,C(0))
\times_{\mathsf{Map}(A,C(0))}\{*\} \\ 
\simeq
\mathsf{Map}_{A/\mecdga_{\mathcal{M}}}(B,C(0))
\end{aligned}
$$ 
as we wanted. \hfill $\Box$

\subsubsection{Strict models.} \label{1.3.3}

For future reference we give here \emph{strict models} for both the
cotangent complex $\mathbb{L}^{int}_{A}$ and the de Rham object
$\DR^{int}(A)$. For $A \in \cdga_{\C}$, corresponding to an object $A
\in Comm(M)$, we can consider the functor
$$
Der^{str}(A,-) : A-Mod_{M} \longrightarrow Set,
$$
sending an $A$-module $M$ to the set $Hom_{Comm(M)/A}(A,A\oplus
M)$. This functor commutes with limits and thus is corepresentable by
an $A$-module $\Omega^{1}_{A} \in A-Mod_{M}$.

Let $Q(A) \longrightarrow A$ be a cofibrant replacement inside
$Comm(M)$. As this is an  
equivalence it induces an equivalence of homotopy categories
$$
Ho(A-Mod_{\C}) \simeq Ho(A-Mod_{M})\simeq Ho(Q(A)-Mod).
$$
Through these identifications, we have a natural 
isomorphism in $Ho(A-Mod_{\C})$
$$
\Omega^{1}_{Q(A)} \simeq \mathbb{L}^{int}_{A}.
$$
In particular, when $A$ is cofibrant $\Omega^{1}_{A}$ is a model for
the cotangent complex of $A$. 

De Rham complexes also possess similarly defined strict models. We
have the functor 
$$
Comm(\epsilon-M^{gr}) \longrightarrow Comm(M),
$$
sending a graded mixed commutative monoid 
$A$ to its part of weight zero $A(0)$. 

This functor commutes with limits and thus possesses a left adjoint
$$
DR^{str} : Comm(M) \longrightarrow Comm(\epsilon-M^{gr}).
$$
For the same formal reasons, the analogue of the Lemma \ref{p1}
remains correct, and for any $A \in Comm(M)$, we have a functorial
isomorphism of graded commutative monoids in $M$
$$
Sym_{A}(\Omega_{A}^{1}[-1]) \simeq DR^{str}(A).
$$
In particular, $Sym_{A}(\Omega_{A}^{1}[-1])$ has a uniquely defined
mixed structure compatible 
with its natural grading and multiplicative structure. 
This mixed structure is given by a map in $M$
$$
\epsilon : \Omega^{1}_{A} \longrightarrow \wedge^{2}\Omega_{A}^{1}
$$
which is called the \emph{strict de Rham differential}. 

If $Q(A)$ is a cofibrant model for $A$ in $Comm(M)$, we have a natural
equivalence of mixed graded commutative dg-algebras in $\C$
$$
DR^{str}(Q(A)) \simeq \DR^{int}(A).
$$
Therefore, the explicit graded mixed commutative monoid
$Sym_{Q(A)}(\Omega_{Q(A)}^{1}[-1])$ 
is a model for $\DR^{int}(A)$.

\begin{rmk}\label{r4}
When $M=C(k)$, and $A$ is a commutative dg-algebra over $k$,
$\DR^{int}(A)$ coincides with the de Rham object $\DR(A/k)$
constructed in \cite{tove}.
\end{rmk}

\subsection{Differential forms and polyvectors}

Next we describe the notions of \emph{differential forms},
\emph{closed differential forms} and \emph{symplectic structure}, as
well as the notion of \emph{$\mathbb{P}_n$-structure} on commutative
dg-algebras over a fixed base $\s$-category $\C$.  We explain a first
relation between Poisson and symplectic structures, by constructing
the symplectic structure associated to a \emph{non-degenerate} Poisson
structure.

\subsubsection{Forms and closed forms.}\label{closedforms} 

Let $A\in \cdga_{\C}$ be a commutative dg-algebra over $\C$.  As
explained in Section~\ref{1.3.2} we have the associated de Rham object
$\DR^{int}(A) \in \epsilon-\cdga_{\C}^{gr}$. We let $\mathbf{1}$ be
the unit object in $\C$, considered as an object in $\epsilon-\C^{gr}$
in a trivial manner (pure of weight zero and with zero mixed
structure). We let similarly $\mathbf{1}(p)$ be its twist by $p\in
\mathbb{Z}$: it is now pure of weight $p$ again with the zero mixed
structure. Finally, we have shifted versions $\mathbf{1}[n](p)\equiv
\mathbf{1}(p)[n] \in \epsilon-\C^{gr}$ for any $n\in \mathbb{Z}$.

For $q\in \mathbb{Z}$, we will denote the weight-degree shift by $q$
functor as 
$$
(-)((q)): \epsilon-\C^{gr} \longrightarrow
\epsilon-\C^{gr} \, E \longmapsto E((q)) \,\, ;
$$ 
it sends $E= \{
E(p), \epsilon \}_{p\in \mathbb{Z}}$ to the graded mixed object in
$\C$ having $E(p+q)$ in weight $p$, and with the obvious induced mixed
structure (with no signs involved). Note that $(-)((q))$ is an
equivalence for any $q \in \mathbb{Z}$, it commutes with the
cohomological-degree shift, and that, in our previous notation, we
have $\mathbf{1}(p)= \mathbf{1}((-p))$.

We will also write $\mathrm{Free}_{\epsilon , 0}^{gr}: \C \to
\epsilon-\C^{gr}$ for the left adjoint to the weight-zero functor
$\epsilon-\C^{gr} \to \C$ sending $E= \{ E(p), \epsilon \}_{p\in
  \mathbb{Z}}$ to its weight-zero part $E(0)$. Note that, then, the
functor $\epsilon-\C^{gr} \to \C$ sending $E= \{ E(p), \epsilon
\}_{p\in \mathbb{Z}}$ to its weight-$q$ part $E(q)$ is right adjoint
to the functor $X \mapsto (\mathrm{Free}_{\epsilon , 0}^{gr}
(X))((-q))$.

Below we will not distinguish notationally between $\DR^{int}(A)$ and
its image under the forgetful functor $\epsilon-\cdga_{\C}^{gr} \to
\epsilon-\C^{gr}$, for $A\in \cdga_{\C}$. The same for $\DR(A)$ and
its image under the forgetful functor $\epsilon-\cdga^{gr}_{k} \to
\epsilon-\dg_{k}^{gr}$, and for $\wedge_{A}^{p}\mathbb{L}^{int}_{A}$
and its image under the forgetful functor $A-Mod_{\C} \to \C$.

\begin{df}\label{d5}
For any $A\in \cdga_{\C}$, and any integers $p\geq 0$ and $n\in
\mathbb{Z}$, we define 
the \emph{space of closed $p$-forms of degree $n$ on $A$} by
$$
\mathcal{A}^{p,cl}(A,n):=
\mathsf{Map}_{\epsilon-\C^{gr}}(\mathbf{1}(p)[-p-n],\DR^{int}(A))
\in \T.
$$
The \emph{space of $p$-forms of degree $n$ on $A$} is defined by
$$
\mathcal{A}^{p}(A,n)
:=\mathsf{Map}_{\C}(\mathbf{1}[-n],\wedge_{A}^{p}\mathbb{L}^{int}_{A}) 
\in \T.
$$
\end{df}

\

\begin{rmk} \label{NB} Note that by definition of realization functors
  (Definition~\ref{dreal}), we have natural identifications 
$$
\begin{aligned}
\mathcal{A}^{p,cl}(A,n) & =\mathsf{Map}_{\medg_k}(k(p)[-p-n],\DR(A)) \\
\mathcal{A}^{p}(A,n) & =
\mathsf{Map}_{\dg_k}(k[-n],\wedge_{|A|}^{p}\mathbb{L}_{A})
\end{aligned}
$$
where $|A| \in \cdga_{k}$. Note also that 
$|\wedge_{A}^{p}\mathbb{L}^{int}_{A}| \simeq \wedge_{|A|}^{p}\mathbb{L}_{A}$.
\end{rmk}

By Proposition~\ref{p1}, we
have 
$$
\mathcal{A}^{p}(A,n)=\mathsf{Map}_{\C}(\mathbf{1}[-n], 
\wedge_{A}^{p}\mathbb{L}^{int}_{A}) \simeq \mathsf{Map}_{\epsilon-\C^{gr}}
(\mathrm{Free}_{\epsilon ,
  0}^{gr}(\mathbf{1})((-p)),\DR^{int}(A)[p+n])
$$ 
and the identity map $\mathbf{1} \to \mathbf{1}$ induces a map
$\mathrm{Free}_{\epsilon , 0}^{gr}(\mathbf{1})((-p)) \rightarrow
\mathbf{1}((-p))$ in $\epsilon-\C^{gr}$ (where, in the target we abuse
notation and write $\mathbf{1}$ for the object $\mathbf{1}$ in pure
weight zero). In particular, we get an induced canonical map
$$
\mathcal{A}^{p,cl}
(A,n) \longrightarrow \mathcal{A}^{p}(A,n)
$$ 
which should be thought of as the map assigning to a closed $p$-form
its \textit{underlying $p$-form}.

In order to gain a better understanding of the spaces 
$\mathcal{A}^{p, cl}(A,n)$, we observe that the object 
$\mathbf{1} \in \epsilon-\C^{gr}$ possesses a natural cell decomposition 
consisting of a sequence of push-outs in $\epsilon-\C^{gr}$
$$
\xymatrix{
X_{m} \ar[r]  & X_{m+1} \\
L_{m+1}[-1] \ar[u] \ar[r] & \ar[u] 0,}
$$
with the following properties

\begin{enumerate}
\item $X_{-1}\simeq 0$.
\item $L_{m} \in \epsilon-\C^{gr}$ is the free graded mixed object in
  $\C$ generated by $\mathbf{1} \in \C$, and weight-shifted by $(-m)$,
  i.e. $L_{m}:= (\mathrm{Free}_{\epsilon , 0}^{gr}(\mathbf{1}))
  ((-m))$.  Note that $L_{m}$ is \textbf{not} concentrated in one
  single weight.
\item There is a natural equivalence $colim_{m} \, X_m \simeq
  \mathbf{1}$.
\end{enumerate}

\

We can give a completely explicit description of this cell
decomposition, by first studying the case of the enriching category
$M=C(k)$. In $\epsilon-C(k)^{gr}$ there is a natural cell model for
$k=k(0)$, considered as a trivial graded mixed complex pure of weight
zero.  The underlying $k$-module is generated by a countable number of
variables $\{x_n,y_n\}_{n\geq 0}$, where $x_n$ is of cohomological
degree $0$ and $y_n$ of cohomological degree $1$, and the
cohomological differential is defined by $d(x_n)=y_{n-1}$ (with the
convention $y_{-1}=0$).  The weight-grading is defined by declaring
$x_n$ to be pure of weight $n$ and $y_n$ pure of weight
$(n+1)$. Finally, the mixed structure is defined by $\epsilon
(x_n)=y_n$.  This graded mixed complex will be denoted by
$\widetilde{k}$ and is easily seen to be equivalent to $k$ via the
natural augmentation $\widetilde{k} \rightarrow k$ sending $x_0$ to
$1$ and all other generators to zero. Note that while $k$ is cofibrant
in the injective model structure on $\epsilon-C(k)^{gr}$ (where
cofibrations and weak equivalences are detected through the forgetful
functor $\mathrm{U}_{\epsilon}: \epsilon-C(k)^{gr} \to C(k)^{gr}$), it
is not cofibrant in the \emph{projective} model structure on
$\epsilon-C(k)^{gr}$ (where fibrations and weak equivalences are
detected through the same forgetful functor
$\mathrm{U}_{\epsilon}$). In fact the map $\widetilde{k} \rightarrow
k$ is a cofibrant replacement of $k$ in the projective model structure
on $\epsilon-C(k)^{gr}$.  Moreover, the graded mixed complex
$\widetilde{k}$ comes naturally endowed with a filtration by
sub-objects $\widetilde{k} = \cup_{m\geq -1} Z_{m}$, where $Z_m$ is
the sub-object spanned by the $x_n$'s and $y_n$'s, for all $n\leq
m$. 

For a general symmetric monoidal model category $M$, enriched over
$C(k)$ as in Section~\ref{1.1}, we can consider $\widetilde{k}\otimes
\mathbf{1}$ as a graded mixed object in $M$. Since $(-)\otimes_{k}
\mathbf{1} $ is left Quillen, the cell decomposition of
$\widetilde{k}$ defined above, induces the required cell decomposition
in $\epsilon-\C^{gr}$
$$
colim_{m} \, X_m \simeq \mathbf{1},
$$
where $X_m := Z_m \otimes \mathbf{1}$.

In particular, we have, for all $m \geq -1$ ($X_{-1}:=0$), a
cofibration sequence in $\epsilon-\C^{gr}$ 
$$
\xymatrix{X_{m} \ar[r]  & X_{m+1} \ar[r] & L_{m+1}.}
$$
Passing to mapping spaces, we obtain, for all graded mixed object $E
\in \epsilon-\C^{gr}$, a tower decomposition
$$
\mathsf{Map}_{\epsilon-\C^{gr}}(\mathbf{1}, E) \simeq lim_{m} \,
\mathsf{Map}_{\epsilon-\C^{gr}}(X_m, E),
$$
together with fibration sequences 
$$
\xymatrix{\mathsf{Map}_{\epsilon-\C^{gr}}(L_{m+1},E) 
\simeq \mathsf{Map}_{\C}(\mathbf{1},E(m+1)) \ar[r] &
\mathsf{Map}_{\epsilon-\C^{gr}}(X_{m+1},E) \ar[r] &
\mathsf{Map}_{\epsilon-\C^{gr}}(X_{m}, E).}
$$
Note that , for any $(n,q)\in \mathbb{Z}^2$, the degree-shift and
weight-shift functors 
$$
[n]\, , \, ((q)): \epsilon-\C^{gr} \to
\epsilon-\C^{gr}
$$ 
are equivalences, hence commute with colimits.
Therefore by taking $E$ to be the graded mixed object
$\DR^{int}(A)[n+p]((p))$, we have the following decomposition of the
space of closed $p$-forms of degree $n$
$$
\mathcal{A}^{p,cl}(A,n) 
\simeq  lim_{m}\, \mathcal{A}^{p,cl}(A,n)(\leq m),$$
where 
$$
\mathcal{A}^{p,cl}(A,n)(\leq m) := \mathsf{Map}_{\epsilon-\C^{gr}}(X_{m},
\DR^{int}(A)[n+p]((p))).
$$ 
These data are all packaged in fibration sequences
$$
\xymatrix{
\mathsf{Map}_{\C}(\mathbf{1},(\wedge_{A}^{p+m+1}\mathbb{L}^{int}_{A})[n-m-1]) 
\ar[r] & \mathcal{A}^{p,cl}(A,n)(\leq m+1) \ar[r] & 
\mathcal{A}^{p,cl}(A,n)(\leq m)}
$$
where we have used Proposition~\ref{p1} to identify 
$$
\DR^{int}(A)[n+p](m+1+p) \simeq Sym_{A}^{m+p+1}(\mathbb{L}^{int}_A
   [-1]) [n+p] \simeq (\wedge_{A}^{p+m+1}\mathbb{L}^{int}_{A})[n-m-1].
$$
These successive fibration sequences embody the \emph{Hodge
  filtration} on the de Rham complex of $A$. Note that $L_0 \simeq
X_0$ so that $\mathcal{A}^{p,cl}(A,n)(\leq -1) \simeq
\mathcal{A}^{p}(A,n)$. In particular, the canonical map
$\mathcal{A}^{p,cl}(A,n) \longrightarrow \mathcal{A}^{p}(A,n) $ from
closed $p$-forms to $p-$forms, defined above, can be re-obtained as
the canonical map 
$$ 
lim_{m}\, \mathcal{A}^{p,cl}(A,n)(\leq m)
\longrightarrow \mathcal{A}^{p,cl}(A,n)(\leq -1)
$$ 
from the limit to
the level $(\leq -1)$ of the tower.

We are now ready to define the notion of a shifted symplectic
structure on a commutative dg-algebra in $\C$. Let $A \in \cdga_{\C}$
and $A-Mod_{\C}$ be the symmetric monoidal $\s$-category of
$A$-modules in $\C$. The symmetric monoidal $\s$-category $A-Mod_{\C}$
is closed, so any object $M$ possesses a dual
$$
M^{\vee}:=\underline{Hom}_{\C}(M,A) \in A-Mod_{\C}.
$$
For an object $M \in A-Mod_{\C}$, and a morphism
$w : A \longrightarrow M\wedge_{A} M[n]$, we have an adjoint morphism
$$
\Theta_w : M^{\vee} \longrightarrow M[n]
$$
where $M^{\vee}$ is the dual object of $M$. 

\begin{df}\label{dtang}
For $A\in \cdga_{\C}$ the 
\emph{internal tangent complex of $A$ is 
defined by}
$$
\mathbb{T}_{A}^{int}:=(\mathbb{L}_{A}^{int})^{\vee} \in A-Mod_{\C}.
$$
\end{df}
Note that the space of (non-closed) $p$-forms of degree $n$ on $A$ can
be canonically identified as the mapping space
$$
\A^{p}(A,n)\simeq
\mathsf{Map}_{A-Mod_{\C}}(A,\wedge^{p}\mathbb{L}^{int}_{A}[n]).
$$
In particular, when $p=2$ and when $\mathbb{L}^{int}_{A}$ is a
dualizable $A$-module, any $2$-form $\omega_0$ of degree $n$ induces a
morphism of $A$-modules
$$
\Theta_{\omega_0} : \mathbb{T}_{A}^{int} \longrightarrow
\mathbb{L}_{A}^{int}[n].
$$

\begin{df}\label{dsymp}
Let $A \in \cdga_{\C}$. We assume that $\mathbb{L}^{int}_{A}$ is
a dualizable object in the symmetric monoidal $\s$-category of $A$-modules
in $\C$. 
\begin{enumerate}
\item A closed $2$-form $\omega \in \pi_{0}(\A^{2,cl}(A,n))$ of degree
  $n$ on $A$ is \emph{non-degenerate} if the underlying $2$-form $\omega_0 \in 
\pi_{0}(\A^{2}(A,n))$ induces an equivalence of $A$-modules
$$
\Theta_{\omega_0} : \mathbb{T}^{int}_{A} \simeq
\mathbb{L}^{int}_{A}[n].
$$
\item The \emph{space} $\mathsf{Symp}(A;n)$ 
\emph{of $n$-shifted symplectic structures on $A$} is
the subspace of $\A^{2,cl}(A,n)$ consisting of the union of connected
components corresponding to non-degenerate elements.
\end{enumerate}
\end{df}

\

De Rham objects have strict models, as explained in our previous
subsection, so the same is true for the space of forms and closed
forms.  Let $A \in \cdga_{\C}$ be a commutative dg-algebra in $\C$,
and choose a cofibrant model $A' \in Comm(M)$ for $A$. Then, the space
of closed $p$-forms on $A$ can be described as follows. We consider
the unit $\mathbf{1} \in M$, and set
$$
|-| : M \longrightarrow C(k)
$$
the functor defined by sending $x \in M$ to
$\underline{Hom}_{k}(\mathbf{1},R(x)) \in C(k)$, where $R(x)$ is a
(functorial) fibrant replacement of $x$ in $M$ and
$\underline{Hom}_{k}$ is the enriched hom of $M$ with values in
$C(k)$. The graded mixed object $\DR^{int}(A)$ can be represented by
$DR^{str}(Q(A))$, and $\DR(A)$ by $|DR^{str}(Q(A))|$. We have by
construction
$$
\mathcal{A}^{p,cl}(A,n) \simeq
\mathsf{Map}_{\epsilon-C(k)^{gr}}(k(p)[-p-n], |DR^{str}(Q(A))|).
$$
In order to compute this mapping space we observe that the injective
model structure on $\epsilon-C(k)^{gr}$ (where cofibrations and weak
equivalences are detected through the forgetful functor
$\mathrm{U}_{\epsilon}: \epsilon-C(k)^{gr} \to C(k)^{gr}$) is Quillen
equivalent to the projective model structure on $\epsilon-C(k)^{gr}$
(where fibrations and weak equivalences are detected through the same
forgetful functor $\mathrm{U}_{\epsilon}$), therefore the
corresponding mapping spaces are equivalent objects in $\T$. It is
then convenient to compute
$\mathsf{Map}_{\epsilon-C(k)^{gr}}(k(p)[-p-n], |DR^{str}(Q(A))|)$ in the
projective model structure, since any object is fibrant here, and we
have already constructed an explicit (projective) cofibrant
resolution $\widetilde{k}$ of $k$.  This
way, we get the following explicit strict model for the space of
closed forms on $A$
$$
\begin{aligned}
\mathcal{A}^{p,cl}(A,n) & \simeq
\mathsf{Map}_{C(k)}(k[-n],\prod_{j\geq
  p}|\wedge_{A'}^{j}\Omega^{1}_{A'}|[-j]) \\ 
& =\mathsf{Map}_{C(k)}(k[-n],\prod_{j\geq p}
\DR(A)(j)).
\end{aligned}
$$
Here $\prod_{j\geq p}|\wedge_{A'}^{j}\Omega^{1}_{A'}|[-j]$ is the
complex with the total differential, which is sum of the cohomological
differential and mixed structure as in \cite[\S 5]{toenems}.
%\texttt{[ give strict models also for $p$-forms and the underlying-form map ]}%

\subsubsection{Shifted polyvectors.} 

We will now introduce the dual notion to differential forms, namely
\emph{polyvector fields}. Here we start with strict models, as the
$\s$-categorical aspects are not totally straightforward and will be
dealt with more conveniently in a second step.

\

\noindent \textbf{Graded dg shifted Poisson algebras in $\C$}.  Let us
start with the case $M=C(k)$, $n \in \mathbb{Z}$, and consider the
\emph{graded n-shifted Poisson operad} $\mathbb{P}^{gr}_{n} \in
Op(C(k)^{gr})$ defined as follows. As an operad in $C(k)$ (i.e. as an
ungraded dg-operad), it is freely generated by two operations
$\cdot,[-,-]$, of arity $2$ and respective cohomological degrees $0$
and $(1-n)$
$$
\cdot \in \mathbb{P}^{gr}_{n}(2)^{0} \qquad [-,-] \in
\mathbb{P}^{gr}_{n}(2)^{1-n},
$$
with the standard relations expressing the conditions that $\cdot$ is
a graded commutative product, and that $[-,-]$ is a biderivation of
cohomological degree $1-n$ with respect to the product $\cdot$.

A $\mathbb{P}^{gr}_{n}$-algebra in $C(k)$ is  just a
commutative dg-algebra $A$ endowed with a compatible Poisson bracket
of degree $(1-n)$
$$
[-,-] : A \otimes_{k} A \longrightarrow A[1-n].
$$
The weight-grading on $\mathbb{P}^{gr}_{n}$ is then defined by letting
$\cdot$ be of weight $0$ and $[-,-]$ be of weight $-1$. When $n>1$,
the operad $\mathbb{P}_{n}$ is also the operad $H_{\bullet}(E_n)$ of
homology of the topological little $n$-disks or $E_n$-operad, endowed
with its natural weight-grading for which $H_0$ is of weight $0$ and
$H_{n-1}$ of weight $-1$ (see \cite{cohen} or \cite{sinha} for a very
detailed account).

We consider $M^{gr}$, the category of $\mathbb{Z}$-graded objects in
$M$, endowed with its natural symmetric monoidal structure. With
fibrations and equivalences defined levelwise, $M^{gr}$ is a symmetric
monoidal model category satisfying our standing assumptions $(1)-(5)$
of \ref{1.1}.  We can then consider $\mathsf{Op}(M^{gr})$ the category
of (symmetric) operads in $M^{gr}$.  As already observed, the category
$M^{gr}$ is naturally enriched over $C(k)^{gr}$, via a symmetric
monoidal functor $C(k)^{gr} \to M^{gr}$.  This induces a functor
${\mathsf{Op}}((C(k))^{gr}) \to \mathsf{Op}(M^{gr})$, and we will denote
by $\mathbb{P}^{gr}_{M, n} \in \mathsf{Op}(M^{gr})$ the image of
$\mathbb{P}^{gr}_{n}$ under this functor.  The category of
$\mathbb{P}^{gr}_{M, n}$-algebras will be denoted by
$\mathbb{P}_{n}-\cdga_{M}^{gr}$, and its objects will be called
\emph{graded $n$-Poisson commutative dg-algebras in $M$}. Such and algebra
consists of the following data.

\begin{enumerate}
\item A family of objects $A(p) \in M$, for $p\in \mathbb{Z}$. 
\item A family of multiplication maps
$$
A(p) \otimes A(q) \longrightarrow A(p+q),
$$
which are associative, unital, and graded commutative.
\item A family of morphisms 
$$
[-,-] : A(p)\otimes A(q) \longrightarrow A(p+q-1)[1-n].
$$
\end{enumerate}
These data are furthermore required to satisfy the obvious
compatibility conditions for a Poisson algebra (see \cite[\S
  1.3]{gejo} for the ungraded 
dg-case).  
%\texttt{[perhaps we could give explicit formulas here].}%
We just recall that, in particular, $A(0)$ should be  a commutative monoid in
$M$, and that the morphism 
$$
[-,-] : A(1)\otimes A(1) \longrightarrow A(1)[1-n]
$$
has to make $A(1)$ into a $n$-Lie algebra object in $M$, or equivalently, 
$A(1)[n-1]$ has to be a Lie algebra object in $M$ when endowed with
the induced pairing 
$$
A(1)[n-1]\otimes A(1)[n-1] \simeq (A(1) \otimes A(1))[2n-2] \longrightarrow
A(1)[n-1].
$$
Since the bracket is a derivation with respect to the product, this
Lie algebra object acts naturally on $A(0)$ by derivations, making the
pair 
$(A(0),A(1)[n-1])$ into a Lie algebroid object in $M$ (see 
\cite{vez2}). Moreover, $A[n-1]$ is a Lie algebra object in $M^{gr}$.

\begin{rmk} Note that for any dg-operad $\mathcal{O}$ over $k$, and for any symmetric monoidal $E: C(k)$-model category $M$ as in \S \ref{1.1},
the symmetric monoidal functor $C(k) \to M$ has a natural extension to a symmetric monoidal functor $E^{gr }: C(k)^{gr} \to M^{gr}$. Moreover $\mathcal{O}$ has a \emph{naive} extension to an operad $\mathcal{O}^{gr, \textrm{naive}}$ in $C(k)^{gr}$, and via $E^{gr}$ there is an induced operad $\mathcal{O}_{M}^{gr, \textrm{naive}}$ on $M^{gr}$. However, for $\mathcal{O}= \mathbb{P}_n$ the $n$-Poisson operad, and for $\mathcal{O}= \mathsf{Lie}$ the Lie operad, the graded versions $\mathbb{P}^{gr}_{n}$ and $\mathsf{Lie}^{gr}$ we are considering here are \emph{not} the naive versions, due to the non-zero weight of the bracket operation. The same is true for our operads $\mathbb{P}^{gr}_{M, n}$ and $\mathsf{Lie}_{M}^{gr}$. 
\end{rmk}

\begin{df}\label{d7}
The \emph{$\s$-category of graded $n$-Poisson commutative dg-algebras
  in the $\s$-category $\C$} is defined to be
$$\mathbb{P}_{n}-\cdga_{\C}^{gr}:=L(\mathbb{P}_{n}-\cdga_{M}^{gr}).$$ 
\end{df}

\bigskip

\

\noindent \textbf{Shifted polyvectors.}  Let $A \in \cdga_{M}$ be a
commutative monoid in $M$. We define a graded $\mathbb{P}_{n}$-algebra
of $n$-shifted polyvectors on $A$ as follows.  As in the case of
forms, we will have an internal and external version of shifted
polyvectors on $A$.  We consider the $A$-module $\Omega_{A}^1$
corepresenting derivations (see \ref{1.3.3}), and we write
$$
T(A,n):=\underline{Hom}_{A}(\Omega_{A}^{1},A[n]) \in A-Mod_{M}
$$ 
for the $A$-module object of derivations from $A$ to to the $A$-module
$A[n]$ (note that $T(A,n)$ is a model for $\mathbb{T}_{A}^{int}[n]$ of
Definition \ref{dtang} only when $A$ is cofibrant and fibrant object
in $\cdga_M$).

Note that $T(A,n)$ can also be identified as follows.  Consider the
canonical  map $$\alpha: \underline{Hom}_{M}(A,A[n]) \otimes A
\longrightarrow A[n]$$ in $M$, adjoint to the identity of
$\underline{Hom}_{M}(A,A[n])$, and the multiplication map $$m:
A\otimes A \longrightarrow A.$$ Then, we consider the following three
maps

\begin{itemize}
\item $\mu_{1}'$ defined 
as the composition
$$
\xymatrix@1{\underline{Hom}_{M}(A,A[n]) \otimes A \otimes A
  \ar[r]^-{\textrm{id}_{A}\otimes m} & \underline{Hom}_{M}(A,A[n])
  \otimes A \ar[r]^-{\alpha} & A[n] }
$$
\item $u_{1}'$ defined as the composition 
$$
\xymatrix{\underline{Hom}_{M}(A,A[n]) \otimes A \otimes A
  \ar[r]^-{\alpha \otimes \textrm{id}_{A}} & A[n] \otimes A
  \ar[r]^-{r} & A[n] }
$$ 
where $r$ is the right $A$-module structure on $A[n]$;
\item $v_{1}'$ defined as the composition 
$$
\xymatrix{\underline{Hom}_{M}(A,A[n]) \otimes A \otimes A
  \ar[r]^-{\sigma \otimes \textrm{id}_{A}} & A \otimes
  \underline{Hom}_{M}(A,A[n]) \otimes A
  \ar[r]^-{\textrm{id}_{A}\otimes \alpha} & A\otimes A[n] \ar[r]^-{l}
  & A[n] }
$$ 
where $l$ is the left $A$-module structure on $A[n]$, and $\sigma$ is
the symmetry for $\underline{Hom}_{M}(A,A[n]) \otimes A$; 
\end{itemize}

If we denote by $\mu_{1}, u_{1}, v_{1} : \underline{Hom}_{M}(A,A[n])
\longrightarrow \underline{Hom}_{M}(A\otimes A,A[n])$ the adjoint maps
to $\mu_{1}', u_{1}', v_{1}'$, then the object $T(A,n)$ is the kernel
of the morphism
$$
\mu_{1} -u_{1}-v_{1} : \underline{Hom}_{M}(A,A[n]) \longrightarrow
\underline{Hom}_{M}(A^{\otimes 2},A[n]).
$$

More generally, for any $p\geq 0$, we define $T^{(p)}(A,n)$ the
$A$-module of \emph{$p$-multiderivations} from $A^{\otimes p}$ to
$A[np]$. This is the $A$-module of morphisms $A^{\otimes p}
\longrightarrow A[np]$ which are derivations in each variable
separately. More precisely, let us consider the canonical
map $$\alpha_p: \underline{Hom}_{M}(A^{\otimes\, p},A[np]) \otimes
A^{\otimes\, p} \longrightarrow A[np]$$ in $M$, adjoint to the
identity of $\underline{Hom}_{M}(A^{\otimes\, p} ,A[n])$, the
multiplication map $m: A\otimes A \longrightarrow A$, and, for any
pair $(P,Q)$ of $A$-modules, let us denote by $\sigma(P,Q)$ the
symmetry map  $P\otimes Q \to Q\otimes P$. Then, for any $1
\leq i\leq p$, we can define the following three morphisms

\begin{itemize}
\item $\mu_{i}'$ defined as the composition 
$$
\xymatrix{\underline{Hom}_{M}(A^{\otimes\, p},A[np]) \otimes
  A^{\otimes\, p+1} \ar[rr]^-{\textrm{id}\otimes m \otimes
    \textrm{id}} & & \underline{Hom}_{M}(A^{\otimes\, p},A[np])
  \otimes A^{\otimes\, p} \ar[r]^-{\alpha_{p}} & A[n] }
$$ 
where $m$ is the multiplication map 
$A_{(i)}\otimes A_{(i+1)} \to A$ on the $(i,i+1)$ factors of
$A^{\otimes\, p+1}$; 
\item $u_{i}'$ defined as the composition 
$$
\xymatrix@C+1.5pc{\underline{Hom}_{M}(A^{\otimes\, p}, A[np]) \otimes
  A^{\otimes\, p+1} \ar[r]^-{\textrm{id} \otimes \sigma_{(i+1)}} &
  \underline{Hom}_{M}(A^{\otimes\, p},A[np]) \otimes A^{\otimes\, p+1}
  \ar[r]^-{\alpha \otimes \textrm{id}_{A}} & A[n] \otimes A 
  \ar[d]^-{r} \\ 
&& A[n] 
}
$$ 
where  $\sigma_{(i+1)}:= \sigma (A_{(i+1)}, A^{\otimes\, p-i})$, and
$r$ is the right $A$-module structure on $A[n]$; 
\item $v_{i}'$ defined as the composition 
$$
\xymatrix@C+1.5pc{\underline{Hom}_{M}(A^{\otimes\, p},A[np]) \otimes
  A^{\otimes\, p+1} \ar[r]^-{\tau_{(i)} \otimes \textrm{id}} & A
  \otimes \underline{Hom}_{M}(A^{\otimes\, p},A[np]) \otimes
  A^{\otimes\, p} \ar[r]^-{\textrm{id}_{A}\otimes \alpha} & A\otimes
  A[n] \ar[d]^-{l} \\
& & A[n] }
$$ 
where $\tau_{(i)}:= \sigma(\underline{Hom}_{M}(A^{\otimes\, p},A[np])
\otimes A^{\otimes\, i-1}, A_{(i)})$, and $l$ is the left $A$-module
structure on $A[n]$.
\end{itemize}

\

We denote by $\mu_{i}, u_{i}, v_{i} : \underline{Hom}_{M}(A^{\otimes\,
  p},A[np]) \longrightarrow \underline{Hom}_{M}(A^{\otimes\,
  p+1},A[np])$ the adjoint maps to $\mu_{i}', u_{i}', v_{i}'$.

We have, for each $1\leq i\leq p$ a sub-object in $M$
$$
Ker(\mu_{i} - u_i -v_i) \subset \underline{Hom}_{M}(A^{\otimes
  p},A[np]).
$$
The intersection of all these sub-objects defines 
$$
T^{(p)}(A,n):=\cap Ker(\mu_{i} - u_i -v_i) \subset
\underline{Hom}_{M}(A^{\otimes p},A[np]). 
$$
The symmetric group $\Sigma_{p}$ acts on
$\underline{Hom}_{M}(A^{\otimes p},A[np])$,  
by its standard action on $A^{\otimes p}$, and by $(-1)^{n} \cdot
\textrm{Sign}$ on $A[np]$ 
which is the natural action when $A[np]$ is identified with
$A[n]^{\otimes_{A} p}$. This action  
stabilizes the sub-object $T^{(p)}(A,n)$ and thus induces a
$\Sigma_{p}$-action on  
$T^{(p)}(A,n)$. We set\footnote{Since we work in characteristic $0$,
  we could have used coinvariants instead of invariants.}
$$
\mathsf{Pol}^{int}(A,n):=\bigoplus_{p\geq
  0}(T^{(p)}(A,-n))^{\Sigma_{p}} \in M, $$
and call it the \emph{object of internal $n$-shifted polyvectors on
  $A$}. 

\

The object $\mathsf{Pol}^{int}(A,n)$ is naturally endowed with a
structure of a graded $(n+1)$-Poisson commutative dg-algebra in $M$ as
follows.

\begin{itemize}
\item The weight $\mathbb{Z}$-grading is the usual one, with
  $(T^{(p)}(A,-n))^{\Sigma_{p}}$ being of weight $p$ by
  definition. The multiplication morphisms
$$
(T^{(p)}(A,-n))^{\Sigma_{p}} \otimes (T^{(q)}(A,-n))^{\Sigma_{q}}
\longrightarrow (T^{(p+q)}(A,-n))^{\Sigma_{p+q}}
$$
are induced by composing the natural morphisms
$$
\underline{Hom}_{M}(A^{\otimes p},A[-np]) \otimes 
\underline{Hom}_{M}(A^{\otimes q},A[-nq]) \longrightarrow 
\underline{Hom}_{M}(A^{\otimes p+q},A[-np]\otimes A[-nq]),
$$
with the multiplication in the monoid $A$:
$$
A[-np]\otimes A[-nq] \simeq (A\otimes A) [-n(p+q]] \longrightarrow
A[-n(p+q)],
$$
and then applying the symmetrization with respect to $\Sigma_{p+q}$. 
This endows the object $\mathsf{Pol}^{int}(A,n)$ with the structure of
a graded commutative monoid object in $M$. 
\item The Lie structure, shifted by $-n$, on $\mathsf{Pol}^{int}(A,n)$ 
is itself a version of the Schouten-Nijenhuis bracket on
polyvector fields. One way to define it categorically is to consider
the graded object  
$\mathsf{Pol}^{int}(A,n)[n]$ as a sub-object of 
$$
Conv(A,n):=\bigoplus_{p\geq 0}\underline{Hom}_{M}(A^{\otimes
  p},A[-np])^{\Sigma_{p}}[n].
$$
The graded object $Conv(A,n)$ is a graded Lie algebra in $M$, where
the Lie bracket is given by natural explicit formulas given by generalized
commutators (the notation $Conv$  here refers to the convolution
Lie algebra of the operad $Comm$ with the endomorphism operad of $A$,
see \cite{lv}). We refer to \cite[10.1.7]{lv} and \cite[\S 2]{mel} for
more details.  This Lie bracket restricts to a graded Lie algebra
structure on $\mathsf{Pol}^{int}(A,n)[n]$.
\end{itemize}

The Lie bracket $\mathsf{Pol}^{int}(A,n)$ is easily seen to be
compatible with the graded algebra structure, i.e.
$\mathsf{Pol}^{int}(A,n)$ is a graded $\mathbb{P}_{n+1}$-algebra
object in $M$.

\begin{df}\label{d8}
  Let $A \in\cdga_{M}$ be a commutative monoid in $M$. The
  \emph{graded $\mathbb{P}_{n+1}$-algebra of $n$-shifted polyvectors
    on $A$} is defined to be
$$
\mathsf{Pol}^{int}(A,n) \in \mathbb{P}_{n+1}-\cdga_{M}^{gr}
$$
described above.
\end{df}

\

For a commutative monoid $A \in Comm(M^{gr})$, the graded
$\mathbb{P}_{n+1}$-algebra $\mathsf{Pol}^{int}(A,n)$ is related to the
set of (non graded) $\mathbb{P}_{n}$-structures on $A$ in the
following way. The commutative monoid structure on $A$ is given by a
morphism of (symmetric) operads in $C(k)$
$$
\phi_{A} : Comm \longrightarrow \underline{Hom}_{k}(A^{\otimes
  \bullet},A),
$$
where the right hand side is the usual endomorphism operad of $A \in
M$ (which is an operad in $C(k)$). We have a natural morphism of
operads $Comm \longrightarrow \mathbb{P}_{n}$ , inducing the forgetful
functor from $\mathbb{P}_{n}$-algebras to commutative monoids, by
forgetting the Lie bracket. The set of $\mathbb{P}_{n}$-algebra
structures on $A$ is by definition the set of lifts of $\phi_{A}$ to a
morphism $\mathbb{P}_{n} \longrightarrow
\underline{Hom}_{C(k)}(A^{\otimes \bullet},A)$
$$
\mathbb{P}^{str}_{n}(A):=Hom_{Comm/\mathsf{Op}}( 
\mathbb{P}_{n},\underline{Hom}_{k}(A^{\otimes  
\bullet},A)).
$$
The superscript \emph{str} stands for \emph{strict}, and is used to
distinguish this operad from its $\s$-categorical version that will be
introduced below.  Recall that $\mathsf{Pol}^{int}(A,n)[n]$ is a Lie
algebra object in $M^{gr}$, and consider another Lie algebra object
$\mathbf{1}(2)[-1]$ in $M^{gr}$ given by $\mathbf{1}[-1] \in M$ with
zero bracket and pure weight grading equal to $2$.

\begin{prop}\label{valerio} There is a natural bijection
$$
\mathbb{P}_{n}^{str}(A) \simeq Hom_{Lie_{M}^{gr}}(\mathbf{1}(2)[-1],
\mathsf{Pol}^{int}(A,n)[n])
$$
where the right hand side is the set of morphisms of Lie algebra
objects in $M^{gr}$.
\end{prop}
\noindent \textbf{Proof.} Recall that $M^{gr}$ is
$C(k)^{gr}$-enriched, and let us consider the corresponding symmetric
lax monoidal functor $R:=\underline{Hom}^{gr}_{k}(\mathbf{1},-) :
M^{gr} \longrightarrow C(k)^{gr}$, where $\mathbf{1}$ sits in pure
weight $0$.  From a morphism $f: \mathbf{1}(2)[-1] \longrightarrow
\mathsf{Pol}^{int}(A,n)[n]$ of graded Lie algebras in $M$, we get a
morphism of graded Lie algebras in $C(k)$
$$
R(f): k(2)[-1] \longrightarrow R(\mathsf{Pol}^{int}(A,n)[n]).
$$
Now, the image under $R(f)$ of the degree $1$-cycle $1\in k$ is then a
morphism
$$
\varphi:= R(f)(1): \mathbf{1} \longrightarrow
T^{(2)}(A,-n)[n+1]^{\Sigma_{2}}
$$ 
in $M$. By definition of $T^{(2)}(A,-n)[n+1]^{\Sigma_{2}}$, the shift
$\varphi[2(n-1)]$ defines a morphism in $M$
$$
[-,-] : A[n-1]\otimes A[n-1] \longrightarrow A[n-1],
$$
which is a derivation in each variable and is $\Sigma_2$-invariant. The
fact that the Lie bracket is zero on $k[-1]$ implies that this bracket
yields a Lie structure on $A$.  This defines a
$\mathbb{P}_{n}$-structure on $A$ and we leave to the reader to verify
that this is a bijection (see also \cite[Proof of
Theorem~3.1]{mel}). \newline \mbox{\quad} \hfill $\Box$

\

\smallskip

Later on we will need the $\s$-categorical version of the previous
proposition, which is a much harder statement. For future reference we
formulate this $\s$-categorical version below but we refer the reader
to \cite{mel} for the details of the proof.  Let $A \in \cdga_{\C}$ be
a commutative dg-algebra in $\C$. We consider the forgetful
$\s$-functor
$$
\mathrm{U}_{\mathbb{P}_{n}} : \mathbb{P}_{n}-\cdga_{\C}
\longrightarrow \cdga_{\C}
$$
sending a $\mathbb{P}_{n}$-algebra in $\C$ to its underlying
commutative monoid in $\C$. The fiber at $A \in \cdga_{\C}$ of this
$\s$-functor is an $\s$-groupoid and thus corresponds to a space
$$
\mathbb{P}_{n}(A):=\mathrm{U}_{\mathbb{P}_{n}}^{-1}(\{A\}) \in \T.
$$
\

\begin{thm} \emph{\cite[Thm. 3.2]{mel}}\label{valerio2} Suppose that $A$ is 
  fibrant and cofibrant in $\cdga_{M}$. There is a natural equivalence
  of spaces
$$
\mathbb{P}_{n}(A) \simeq
Map_{Lie_{\C}^{gr}}(\mathbf{1}_{M}(2)[-1],\mathsf{Pol}^{int}(A,n)[n])
$$
where the right hand side is the mapping space of morphisms of inside
the $\s$-category of Lie algebra objects in $\C^{gr}$.
\end{thm}

\

\begin{rmk} Theorem 3.2 in \cite{mel} is stated for $M$ the model category of non-positively graded dg-modules over $k$, but the same proof extends immediately to our general $M$. The original statement seems moreover to require a restriction to those cdga's having a dualizable cotangent complex. This is due to the fact that the author uses the tangent complex (i.e. the dual of the cotangent complex) in order to identify derivations. However, the actual proof produces an equivalence between (weak, shifted) Lie brackets
and (weak) biderivations. Therefore if one identifies derivations using the linear dual of the symmetric algebra of the cotangent complex, the need to pass to the tangent complex disappears, and the result holds with the same proof and without the assumption of the cotangent complex being dualizable. This is the main reason we adopted Def. \ref{d8} as our definition of internal polyvectors.
\end{rmk}

Now we give a slight enhancement of Theorem \ref{valerio2} and, as a corollary, we will get a strictification result (Cor. \ref{strett}) that will be used in \S \ref{dim}.\\

Let $\mathsf{Poiss}^{\textrm{eq}}_{M, n}$ be the category whose objects are pairs $(A, \pi)$ where $A$ is a fibrant-cofibrant object in $\cdga_{M}$, and $\pi$ is a map $\mathbf{1}_M[-1](2) \to \mathsf{Pol}^{int}(A,n+1)[n+1] $ in the homotopy category of  $Lie_{M}^{gr}$, and whose morphisms $(A,\pi) \to (A', \pi')$ are weak equivalences $u: A \to A'$ in $\cdga_{M}$ such that the diagram $$\xymatrix{ & \mathsf{Pol}^{int}(A,n+1)[n+1] \ar[dd]^-{\mathsf{Pol}^{int}(u, n)[n]} \\ \mathbf{1}_\mathsf{M}[-1](2) \ar[ru]^-{\pi} \ar[rd]_-{\pi'} & \\  & \mathsf{Pol}^{int}(A',n+1)[n+1] } $$ is commutative in the homotopy category of  $Lie_{M}^{gr}$. We denote the nerve of $\mathsf{Poiss}^{\mathrm{eq}}_{M, n}$ by $\mathbf{Poiss}^{\mathrm{eq}}_{M, n}$.  \\

There is an obvious (strict) functor $w$ from the category cofibrant-fibrant objects in $\mathbb{P}_{n+1}-\cdga_{M}$ and weak equivalences, to $\mathsf{Poiss}^{\mathrm{eq}}_{M, n}$, sending a strict $\mathbb{P}_{n+1}$-algebra $B$ in $\mathsf{M}$ to the pair $(B, \pi)$, where $\pi$ is induced, in the standard way, by the (strict) Lie bracket on $B$ (since the bracket is strict, it is a strict biderivation on $B$, and the classical construction carries over). Restriction to weak equivalences (between cofibrant-fibrant objects) in $\mathbb{P}_{n}-\cdga_{M}$ ensures this is a functor, and note that objects in the image of $w$ are, by definition, \emph{strict} pairs, i.e. maps $\pi: \mathbf{1}_{M}[-1](2) \to \mathsf{Pol}^{int}(A,n+1)[n+1] $ are actual morphism in $Lie_{M}^{gr}$ (rather than just maps in the homotopy category). The functor $w$ is compatible with the forgetful functors $p: \mathbb{P}_{n}-\cdga_{M} \to \cdga_{M}$, and $q: \mathsf{Poiss}^{\mathrm{eq}}_{M, n} \to \cdga_{M}$, and by passing to the nerves, we thus obtain a commutative diagram in $\T$ (where we have kept the same name for the maps) $$\xymatrix{  \mathcal{I}(\mathbb{P}_{n+1}-\cdga_{\C}) \ar[dr]_-{p} \ar[rr]^-{w} & &\mathbf{Poiss}^{\mathrm{eq}}_{M, n} \ar[dl]^-{q} \\ & \mathcal{I}(\cdga_{\C}) & }$$ where $\mathcal{I}(\mathcal{C})$ denotes the maximal $\infty$-subgroupoid of an $\s$-category $\mathcal{C}$, i.e. the classifying space of $\mathcal{C}$. Note that $p$, and $q$ are both surjective, since they both have a section given by choosing the trivial bracket or the trivial strict map $\pi$.

\begin{thm}\label{melext} The map of spaces $w: \mathcal{I}(\mathbb{P}_{n+1}-\cdga_{\C}) \to \mathbf{Poiss}^{\mathrm{eq}}_n $ is an equivalence.
\end{thm}

\noindent \textbf{Proof}. It is enough to prove that for any cofibrant $A \in \cdga_{M}$, the map induced by $w$ between $q$ and $p$ fibers over $A$ is an equivalence. But this is exactly Theorem \ref{valerio2}. \hfill $\Box$ \\

As an immediate consequence, we get the following useful strictification result. An arbitrary object $(A, \pi)$ in $\mathsf{Poiss}^{\textrm{eq}}_{M, n}$ will be called a \emph{weak pair}, and we will call it a \emph{strict pair} if $\pi$ is strict, i.e. is an actual morphism $\pi: \mathbf{1}_M[-1](2) \to \mathsf{Pol}^{int}(A,n+1)[n+1] $ in $Lie_{M}^{gr}$.

\begin{cor}\label{strett} Any weak pair is equivalent, inside $\mathsf{Poiss}^{\textrm{eq}}_{M, n}$,  to a strict pair.
\end{cor}

\noindent \textbf{Proof}. By Theorem \ref{melext},  an object $(A, \pi) \in \mathsf{Poiss}^{\textrm{eq}}_{M, n}$ (i.e. an a priori weak pair), is equivalent to a pair of the form $w(B)$, where $B \in  \mathbb{P}_{n+1}-\cdga_{M}$ (i.e. is a strict $\mathbb{P}_{n+1}$-algebra in $M$), whose underlying commutative algebra is weakly equivalent  to $A$ in  $\cdga_{M}$. We conclude by observing that objects in the image of $w$ are always strict pairs. \hfill $\Box$ \\

\smallskip

\noindent \textbf{Functoriality}. The assignment $A \mapsto
\mathsf{Pol}^{int}(A,n)$ is not quite functorial in $A$, and it is
therefore not totally obvious how to define its derived version. We
will show however that it can be derived to an $\s$-functor from a
certain sub-$\s$-category of formally \'etale morphisms
$$
\mathsf{Pol}^{int}(-,n) : \cdga_{\C}^{fet} \longrightarrow
\mathbb{P}_{n+1}-\cdga_{\C}^{gr}.
$$

We start with a (small) category $I$ and consider the model category
$M^{I}$ of diagrams of shape $I$ in $M$. It is endowed with the model
category structure for which the cofibrations and equivalences are
defined levelwise. As such, it is a symmetric monoidal model category
which satisfies again our conditions $(1)-(5)$ of \ref{1.1}. For $$(\underline{A}: i \ni I \longmapsto A_i \in Comm(M))
\in Comm(M^{I})\simeq Comm(M)^{I}$$ an $I$-diagram of commutative
monoids in $M$, we have its graded $\mathbb{P}_{n+1}$-algebra of
polyvectors $\mathsf{Pol}^{int}(\underline{A},n) \in
\mathbb{P}_{n+1}-\cdga_{M^{I}}^{gr} \simeq
(\mathbb{P}_{n+1}-\cdga_{M}^{gr})^I$.

\begin{lem}\label{l3}
  With the above notation, assume that $\underline{A}$ satisfies the following
  conditions
\begin{itemize}
\item $\underline{A}$ is a fibrant and cofibrant object in $Comm(M)^{I}$. 
\item For every morphism $i \rightarrow j$ in $I$, the  morphism
$A_i \rightarrow A_j$ induces an equivalence in $Ho(M)$
$$\mathbb{L}_{A_{i}} \otimes^{\mathbb{L}}_{A_{i}}A_{j} \simeq 
\mathbb{L}_{A_j}.$$
\end{itemize}
Then, we have: 
\begin{enumerate}
\item for every object $i \in I$ there is a natural equivalence of
  graded $\mathbb{P}_{n+1}$-algebras
$$
\xymatrix{\mathsf{Pol}^{int}(\underline{A},n)_i  \ar[r]^-{\sim} &
  \mathsf{Pol}^{int}(A_i,n),}
$$
\item for every morphism $i \rightarrow j$ the induced morphism 
$$
\mathsf{Pol}^{int}(\underline{A},n)_i \longrightarrow \mathsf{Pol}^{int}(\underline{A},n)_j
$$
is an equivalence of graded $\mathbb{P}_{n+1}$-algebras.
\end{enumerate}
\end{lem}
\noindent \textbf{Proof.} Since $\underline{A}$ is fibrant and cofibrant as an
object of $Comm(M)^I$, we have that for all $i\in I$ the object $A_i$
is again fibrant and cofibrant in $Comm(M)$. As a consequence, for all
$i \in I$, the $A_i$-module $\mathbb{L}_{A_{i}}$ can be represented by
the strict model $\Omega_{A_{i}}^{1}$. Moreover, the second assumption
implies that for all $i\rightarrow j$ in $I$ the induced morphism
$$
\Omega_{A_{i}}^{1} \otimes_{A_{i}}A_j \longrightarrow
\Omega_{A_{j}}^{1}
$$
is an equivalence in $M$. 

As $\underline{A}$ is cofibrant, so is the $\underline{A}$-module $\Omega_{\underline{A}}^{1} \in
\underline{A}-Mod_{M^{I}}$. This implies that $(\Omega_{\underline{A}}^{1}) ^{\otimes_{\underline{A}} p}$
is again a cofibrant object in $\underline{A}-Mod_{M^{I}}$. The graded object
$\mathsf{Pol}^{int}(\underline{A},n)$ in $M^{I}$ of $n$-shifted polyvectors on $A$
is thus given by
$$
\bigoplus_{p\geq 0} \underline{Hom}_{\underline{A}-Mod_{M^{I}}}((\Omega_{\underline{A}}^{1})
^{\otimes_{\underline{A}} p},\underline{A}[-np])^{\Sigma_p}.
$$
For all $i\in I$, and all $p\geq 0$, we have a natural
evaluation-at-$i$ morphism
$$ 
\underline{Hom}_{\underline{A}-Mod_{M^{I}}}((\Omega_{\underline{A}}^{1}) ^{\otimes_{\underline{A}}
  p}, \underline{A}[-np])^{\Sigma_p}\longrightarrow
\underline{Hom}_{A_i-Mod_{M}}((\Omega_{A_{i}}^{1}) ^{\otimes_{A_{i}}
  p},A_{i}[-np])^{\Sigma_p}.
$$
We now use the following sublemma

\begin{sublem}\label{1.6.6}
Let $\underline{A}$ be a commutative monoid in $M^{I}$. Let $E$ and $F$ be two 
$\underline{A}$-module objects, with $E$ cofibrant and $F$ fibrant. We assume that 
for all $i\rightarrow j$ in $I$ the induced morphisms
$$
E_i \longrightarrow E_j \qquad F_i \longrightarrow F_j
$$
are equivalences in $M$. Then, for all $i\in I$, the evaluation
morphism
$$
\underline{Hom}_{\underline{A}-Mod_{M^{I}}}(E,F)_{i} \longrightarrow 
\underline{Hom}_{A_{i}-Mod_{M}}(E_{i},F_{i}) 
$$
is an equivalence in $M$.
\end{sublem}
\noindent \textit{Proof of sub-lemma \ref{1.6.6}}. 
%\texttt{[ I'll fix this proof. ]}% 
For $i\in I$, we have a natural isomorphism
$$
\underline{Hom}_{\underline{A}-Mod_{M^{I}}}(E,F)_{i} \simeq
\underline{Hom}_{M}(E_{|i},F_{|i}),
$$
where $(-)_{|i} : M^{I} \longrightarrow M^{i/I}$ denotes the
restriction functor, and $\underline{Hom}_{M}$ now denotes the natural
enriched $Hom$ of $M^{i/I}$ with values in $M$.  This restriction
functor preserves fibrant and cofibrant objects, so $E_{|i}$ and
$F_{|i}$ are cofibrant and fibrant $A_{|i}$-modules. By assumption, if
we denote by $E_{i} \otimes A_{|i}$ the $A_{|i}$-module sending
$i\rightarrow j$ to $E_{i}\otimes_{A_{i}}A_{j} \in A_{j}-Mod_{M}$, the
natural adjunction morphism
$$
E_{i} \otimes A_{|i} \longrightarrow E_{|i}
$$
is an equivalence of cofibrant $A_{|i}$-modules. This implies that the
induced morphism
$$
\underline{Hom}_{M}(E_{|i},F_{|i}) \longrightarrow
\underline{Hom}_{M}(E_{i}\otimes A_{|i},F_{|i})\simeq
\underline{Hom}_{M}(E_i,F_i)
$$
is an equivalence in $M$. \hfill $\Box$ 

\

\noindent Sublemma \ref{1.6.6} implies that the evaluation morphism
$\mathsf{Pol}^{int}(\underline{A},n)_{i} \longrightarrow
\mathsf{Pol}^{int}(A_i,n)$ is an equivalence. As this morphism is a
morphism of graded $\mathbb{P}_{n+1}$-algebras, this proves assertion
$(1)$ of the lemma. Assertion $(2)$ is proven in the same manner.
\hfill $\Box$

\

While it is not true that an arbitrary morphism $A \longrightarrow B$
in $Comm(M)$ induces a morphism $\mathsf{Pol}^{int}(A,n)
\longrightarrow \mathsf{Pol}(B,n)$ (i.e.  polyvectors are not
functorial for arbitrary morphisms), Lemma \ref{l3} provides a way to
understand a restricted functoriality of the construction $A\mapsto
\mathsf{Pol}^{int}(A,n)$.  In fact, let $I$ be the sub-category of
morphisms in $\cdga_{M}$ consisting of all morphisms $A \rightarrow B$
which are \emph{formally \'etale} i.e. morphisms for which the induced map
$$
\mathbb{L}^{int}_A \otimes^{\mathbb{L}}_{A}B \longrightarrow
\mathbb{L}^{int}_{B}
$$
is an isomorphism in $Ho(M)$, or equivalently in $B-Mod_M$. The category $I$ is not small but things
can be arranged by fixing universes, or bounding the cardinality of
objects.  We have a natural inclusion functor $I \longrightarrow
\cdga_{M}$, and we choose a fibrant and cofibrant model for this
functor, denoted as
$$
\mathcal{A} : I \longrightarrow \cdga_{M}.
$$
This functor satisfies the conditions of Lemma~\ref{l3} above, and
thus induces an $\s$-functor after inverting equivalences
$$\mathsf{Pol}^{int}(\mathcal{A},n) : L(I) \longrightarrow 
L(\mathbb{P}_{n+1}-\cdga_{M}^{gr})=\mathbb{P}_{n+1}-\cdga_{\C}^{gr}.$$
The $\s$-category $L(I)$ is naturally equivalent to the 
(non-full) sub-$\s$-category of $L(\cdga_{M})=\cdga_{\C}$
consisting of formally \'etale morphisms. We denote this 
$\s$-category by $\cdga_{\C}^{fet} \subset \cdga_{\C}$.
We thus have constructed 
an $\s$-functor
$$
\Pol^{int}(-,n):=\mathsf{Pol}^{int}(\mathcal{A},n) : \cdga_{\C}^{fet}
\longrightarrow \mathbb{P}_{n+1}-\cdga_{\C}^{gr}.
$$
\

\begin{df}\label{d9}
The $\s$-functor 
$$
\Pol^{int}(-,n) : \cdga_{\C}^{fet} \longrightarrow
\mathbb{P}_{n+1}-\cdga_{\C}^{gr}
$$
is called the functor of \emph{graded $\mathbb{P}_{n+1}$-algebras of
  internal $n$-shifted polyvectors} in $\C$.
\begin{enumerate}
\item
If $A \in\cdga_{\C}$ is a commutative dg-algebra in $\C$, the 
\emph{graded $\mathbb{P}_{n+1}$-algebra of internal $n$-shifted
polyvectors on $A$} is its value 
$\Pol^{int}(A,n) \in \mathbb{P}_{n}-\cdga_{\C}^{gr}$
at $A$.
\item If $A \in\cdga_{\C}$ is a commutative dg-algebra in $\C$, the 
\emph{graded $\mathbb{P}_{n+1}$-algebra of $n$-shifted
polyvectors on $A$} is 
$
\Pol(A,n):=|\Pol^{int}(A,n)| \in \mathbb{P}_{n}-\cdga_{k}^{gr}.
$
\end{enumerate}
\end{df}

\

\begin{rmk}\label{notable}Note that, by lemma \ref{l3}, we know that the values of the
$\s$-functor $\Pol^{int}$ at $A \in \cdga_{M}$ is naturally equivalent, inside $\mathbb{P}_{n}-\cdga_{\C}^{gr}$,
to the graded $\mathbb{P}_{n+1}$-algebra $\mathsf{Pol}^{int}(QR(A),n)$,
where $QR(A)$ is a fibrant and cofibrant model for $A$ in $\cdga_{M}$.
\end{rmk}
\subsubsection{$\mathbb{P}_{n}$-structures and symplectic forms.} \label{comparazione}

In this section we explain how the standard relation between Poisson
structures and differential forms manifests itself in our setting.\\

\noindent \textbf{Construction $\phi_{\pi}$.} Let $A' \in \cdga_{M}$ be a commutative dg-algebra over $M$. We fix an
integer $n \in \mathbb{Z}$, and we consider on one side
$\mathsf{Pol}^{int}(A',n)$, the $n$-shifted polyvectors on $A'$, and on
the other side, $DR^{str}(A')$, the strict de Rham complex of $A'$. By
Proposition~\ref{valerio}, a (strict) $\mathbb{P}_{n}$-structure on $A'$ is
nothing else than a morphism of graded dg-Lie algebras in $M$
$$
\pi : \mathbf{1}(2)[-1] \longrightarrow \mathsf{Pol}^{int}(A',n)[n]. 
$$ 
Assume that one such $\mathbb{P}_{n}$-structure $\pi$ is
fixed on $A'$.  We can use $\pi$ in order to define a structure of a
graded mixed object on $\mathsf{Pol}^{int}(A',n)$, as follows.  Recall
that the weight $q$ part of $\mathsf{Pol}(A',n)$ is the object
$T^{(q)}(A',-n)^{\Sigma_{q}}$ of $\Sigma_{q}$-invariant
multiderivations $A'^{\otimes q} \longrightarrow A'[-nq]$.  Consider the
symmetric lax monoidal functor
$R:=\underline{Hom}^{gr}_{k}(\mathbf{1},-) : M^{gr} \longrightarrow
C(k)^{gr}$ (where $\mathbf{1}$ sits in weight $0$).  Then $R(\pi):
k(2)[-1] \longrightarrow R(\mathsf{Pol}^{int}(A',n))[n]$ is a morphism
of graded Lie algebras in $C(k)$.  The image under $R(\pi)$ of
the degree $1$ cycle $1\in k$ is then a morphism
$$
\underline{\pi}:= R(\pi)(1): \mathbf{1} \longrightarrow
T^{(2)}(A',-n)[n+1]^{\Sigma_{2}}
$$ 
in $M$.  The composite map 
$$
\xymatrix{ \epsilon_{\pi}: \mathbf{1}
  \otimes T^{(q)}(A', -n)^{\Sigma_{q}} \ar[r]^-{\underline{\pi} \otimes
    \textrm{id}} & T^{(2)}(A', -n)[n+1]^{\Sigma_{2}} \otimes
  T^{(q)}(A',-n)^{\Sigma_{q}} \ar[r]^-{[-,-]} & T^{(q+1)}(A',
  -n)[1]^{\Sigma_{q+1}} }
$$ 
(where $[-,-]$ denotes the Lie bracket part of the graded
$\mathbb{P}_{n+1}$-structure on $\mathsf{Pol}^{int}(A',n)$) defines
then a mixed structure on the graded object $\mathsf{Pol}^{int}(A',n)$,
making it into a graded mixed object in $M$. This graded mixed
structure is also compatible with the multiplication and endows
$\mathsf{Pol}^{int}(A',n)$ with a \emph{graded mixed commutative
  dg-algebra} structure in $M$.

Since in weight $0$ we have 
$\mathsf{Pol}^{int}(A',n)(0)=A'$, the identity map $A' \to A'$ induces, by
Section~\ref{1.3.3}, a morphism 
$$
\phi_{\pi, A'} : DR^{str}(A') \longrightarrow \mathsf{Pol}^{int}(A',n)
$$
of graded mixed commutative algebras
in $M$.\\

\begin{rmk}
Here is an equivalent way of constructing $\phi_{\pi, A'} : DR^{str}(A') \to \mathsf{Pol}^{int}(A',n)$.  The morphism $\underline{\pi}$
defines a morphism of $A'$-modules $\wedge^{2}_{A'}\Omega_{A'}^{1}
\longrightarrow A'[1-n]$, and, by duality, a morphism of $A$-modules
$$ 
\Omega_{A'}^{1}[-1] \longrightarrow
\underline{Hom}_{A'-Mod}(\Omega_{A'}^{1}, A'[-n])\simeq T^{(1)}(A',-n)
$$
Since $\mathsf{Pol}^{int}(A',n)\in \cdga^{gr}_{M}$, by composing it
with the map $T^{(1)}(A',-n) \to \mathsf{Pol}^{int}(A',n)$, and using
adjunction, we get and induced map
$$
Sym_{A'}(\Omega^{1}_{A'}[-1]) \longrightarrow \mathsf{Pol}^{int}(A',n)
$$ 
of graded commutative algebras in $M$. Now it is enough to invoke
the isomorphism $DR^{str}(A') \simeq Sym_{A}(\Omega^{1}_{A'}[-1])$ (see
Section~\ref{1.3.3}), to obtain a map of graded commutative
algebras
$$
\phi_{\pi, A'} : DR^{str}(A') \longrightarrow \mathsf{Pol}^{int}(A',n)
$$ 
that can be verified to strictly preserve with the mixed differentials
on both sides. Thus $\phi_{\pi, A'} $ is a map of graded mixed commutative algebras in $M$.
\end{rmk}

\

Let now $A \in \mathbb{P}_{n}-\cdga_{\C}$. Since $\mathbb{P}_{n}-\cdga_{\C} = L(\mathbb{P}_{n}-\cdga_{M})$, we may choose $A'$ fibrant-cofibrant in $\mathbb{P}_{n}-\cdga_{M}$ (a strict $\mathbb{P}_{n}$-algebra in $M$) which is equivalent to $A$ inside $\mathbb{P}_{n}-\cdga_{\C}$.  Since we have equivalences  $$\DR^{int}(A) \simeq \DR^{int}(A') \simeq DR^{str}(A')$$ in the $\s$-category of graded mixed commutative algebras in $\C$, and equivalences (see Remark \ref{notable}) $$\Pol^{int}(A,n)\simeq \Pol^{int}(A',n) \simeq \mathsf{Pol}^{int}(A',n)$$ in $\mathbb{P}_{n+1}-\cdga_{\C}^{gr}$, we may run the above Construction $\phi_{\pi}$ on $A'$, and use Definition \ref{d9} in order to :
\begin{itemize}
\item turn $\Pol^{int}(-,n)$ into an $\s$-functor $$
 \Pol^{int}(-,n) : 
(\mathbb{P}_{n}-\cdga_{\C})^{eq} \longrightarrow
\epsilon-\cdga^{gr}_{M} \, \, ; 
$$
\item consider the functor $$\DR^{int} : 
(\mathbb{P}_{n}-\cdga_{\C})^{eq} \longrightarrow
\epsilon-\cdga^{gr}_{M}$$ as the composition of the restriction $\DR^{int}: (\cdga_{\C})^{eq} \to \epsilon-\cdga^{gr}_{M}$ of the usual $\DR^{int}$ functor, with the forgetful functor $(\mathbb{P}_{n}-\cdga_{\C})^{eq} \to (\cdga_{\C})^{eq}$ ;
\item promote  the collection of all $\phi_{\pi, A}$'s to a morphism 
$$\phi_{\pi}: \DR^{int} \longrightarrow \Pol^{int}(-,n)$$
that is well defined in the $\s$-category of $\s$-functors from 
$(\mathbb{P}_{n}-\cdga_{\C})^{eq}$ 
to $\epsilon-\cdga^{gr}_{M}$.
\end{itemize}

\

%Since we have equivalences  $$\DR^{int}(A) \simeq \DR^{int}(A') \simeq DR^{str}(A') \,, \,\,\,\,\,\,\,\,\,  \Pol^{int}(A,n)\simeq \Pol^{int}(A',n) \simeq \mathsf{Pol}^{int}(A',n)$$ inside the $\s$-category of graded mixed commutative algebras in $\C$, we get an induced morphism
%$$
%\phi_{\pi, A} : \DR^{int}(A) \longrightarrow \Pol^{int}(A,n)
%$$
%of graded mixed cdga's in $\C$, that is functorial with respect to equivalences in $A$. More precisely, we have two 
%$\s$-functors
%$$
%\DR^{int}, \Pol^{int}(-,n) : 
%(\mathbb{P}_{n}-\cdga_{\C})^{eq} \longrightarrow
%(\epsilon-\cdga^{gr}_{M})^{eq}, 
%$$
%and
%\

\begin{df}\label{dpoissnd}
  We say that $A \in \mathbb{P}_{n}-\cdga_{\C}$ is \emph{non-degenerate} if the
  previous morphism
$$
\phi_{\pi, A} : \DR^{int}(A) \longrightarrow \Pol^{int}(A,n)
$$
is an equivalence.
\end{df}

\

\begin{rmk} This definition is obviously independent of the choice of a strict model $A'$, and hence of the corresponding map $\pi$ (which is then uniquely identified by $A'$, by Thm. \ref{valerio}).
\end{rmk}

\noindent
For an $n$-Poisson commutative cdga $A \in
\mathbb{P}_{n}-\cdga_{\C}$, we consider, as above, a strict cofibrant-fibrant model $A'\in \mathbb{P}_{n}-\cdga_{M}$, together with the corresponding strict map
$\pi : \mathbf{1}(2)[-1] \longrightarrow \mathsf{Pol}^{int}(A',n)[n]$ of graded dg-Lie algebras in $M$.
Such a $\pi$
defines a morphism of graded mixed objects in $M$:
$$
\omega_{\pi, A'} : \mathbf{1}(2) \longrightarrow
\mathsf{Pol}^{int}(A',n)[n+1]. 
$$
Since $\Pol^{int}(A,n) \simeq \mathsf{Pol}^{int}(A',n)$, 
we thus obtain 
a diagram of graded mixed objects in $\C$: 
$$
\xymatrix{
\DR^{int}(A)[n+1] \ar[rr]^-{\phi_{\pi, A}[n+1]} & & \Pol^{int}(A,n)[n+1] &
\ar[l]_-{\omega_{\pi, A}} 
\mathbf{1}(2),}
$$
for each $A \in \mathbb{P}_{n}-\cdga_{\C}$, which, upon realization, produces a diagram in graded mixed $k$-dg modules
$$ 
\xymatrix{ \DR(A)[n+1] \ar[rr]^-{\phi_{\pi, A}[n+1]} & & \Pol(A,n)[n+1] &
  \ar[l]_-{\omega_{\pi, A}} k(2).}
$$

We use $\phi_{\pi, A}[n+1]$, and $\omega_{\pi , A}$ to identify $\DR(A)[n+1]$ and $k(2)$ as objects in the $\infty$-over-category $\medg_{k}/\Pol(A,n+1)[n+1]$, and give the following

\begin{df}\label{dformpoiss}
  Let $A \in \mathbb{P}_{n}-\cdga_{\C}$.The space of \emph{closed
    $2$-forms compatible with the $\mathbb{P}_{n}$-structure on $A$}
  is the space
$$\mathsf{Map}_{\medg_{k}/\Pol(A,n+1)[n+1]}(k(2),\DR(A)[n+1]) \in \T .$$
\end{df}

\

\noindent
In other words, the space of closed $2$-forms
compatible with the $\mathbb{P}_{n}$-structure on $A$
consists of lifts $k(2) \longrightarrow \DR(A)[n+1]$
of the morphism $\omega_{\pi}$. There is a natural forgetful morphism
$$
\mathsf{Map}_{\medg_{k}/\Pol(A,n)[n+1]}(k(2),\DR(A)[n+1]) 
\longrightarrow \mathsf{Map}_{\medg_{k}}(k(2),\DR(A)[n+1]) \simeq
\mathcal{A}^{2,cl}(A,n-1),
$$ 
to the space of closed $2$-forms on $A$ of
degree $(n-1)$.

Note that, by definition, if a $\mathbb{P}_{n}$-algebra $A$ in $\C$ is
non-degenerate, then the space of closed $2$-forms compatible with the
$\mathbb{P}_{n}$-structure on $A$ is contractible. In particular, we
obtain in this case a well defined  (in $\pi_0(\mathcal{A}^{2,cl}(A,n-1))$) and canonical closed $2$-form $\omega$
of degree $(n-1)$ on $A$. Moreover, since $\pi$ is assumed to be non-degenerate,
then so is the corresponding underlying $2$-form. For reference, we
record this observation in the following

\begin{cor}\label{cpoissnd} \begin{itemize}
\item  Let $A\in \mathbb{P}_{n}-\cdga_{\C}$ be non-degenerate. Then there is a unique, up to a
  contractible space of choices, closed and non-degenerate $2$-form of
  degree $(n-1)$ compatible with the $\mathbb{P}_{n}$-structure on
  $A$.
  \item As a consequence, for any $A\in \cdga_{\C}$, there is a well-defined morphism of spaces $$
\W_{A} : \mathbb{P}_{n}(A)^{nd} \longrightarrow \mathsf{Symp}(A,n-1),$$ from the space of non-degenerate $\mathbb{P}_{n}$-structures on $A$ to
the space of $(n-1)$-shifted symplectic structures on $A$.
  \end{itemize} 
\end{cor}

\

\begin{rmk} Moreover, exactly as we did for Definition~\ref{d9}, we get that $\W_{A}$ is functorial in $A$, with respect to formally \'etale
maps in $\cdga_{\mathcal{M}}$.
\end{rmk}

\
%%%%%%%%%%%%%%%%%%%%%%%%% DA QUI 
\noindent
%
%\begin{df} \begin{itemize}
%\item Let $A \in \cdga_{\C}$ such that $\mathbb{L}_{A}^{int}$ is a
%  dualizable $A$-module in $\C$. A morphism of graded dg-Lie
%algebras
%$$
%k(2)[-1] \longrightarrow \mathsf{Pol}(A,n)[n]= |
%$$
%is \emph{non-degenerate} if the underlying morphism in $\C$
%$$ 
%\mathbf{1} \longrightarrow Sym_{A}^{2}(\mathbb{T}_{A}^{int}[-n])[n+1]
%$$
%induces an equivalence of $A$-modules
%$$
%\mathbb{L}_{A}^{int} \simeq \mathbb{T}_{A}^{int}[1-n].
%$$
%\item
%\end{itemize}
%\end{df}
% 
% Theorem~\ref{valerio2}
% 
%\begin{cor}\label{cpoissnd2}
%  Let $A \in \cdga_{\C}$ such that $\mathbb{L}_{A}^{int}$ is a
%  dualizable $A$-module in $\C$. Then, there is a natural morphism of
%  spaces, functorial in $A$ with respect to formally \'etale
%  morphisms:
%$$
%\mathsf{Map}_{\dglie_k^{gr}}^{nd}(k(2)[-1],\mathsf{Pol}(A,n)[n]) \longrightarrow
%\mathsf{Symp}(A,n-1),
%$$
%where
%$\mathsf{Map}_{\dglie_k^{gr}}^{nd}(k(2)[-1],\mathsf{Pol}(A,n)[n])$ is
%the subspace of
%$\mathsf{Map}_{\dglie_k^{gr}}(k(2)[-1],\mathsf{Pol}(A,n)[n])$ of
%connected components of non-degenerate elements.
%\end{cor}

We finish this paragraph with an important corollary that will be used in \S \ref{ndpoisson}. \\ In order to prepare for the next definition, first of all, observe that if $A \in \cdga_{\C}$ has the property that $\mathbb{L}_{A}^{int}$ is dualizable in $A-Mod_{\C}$ , then we have an equivalence $$\mathbf{Pol}^{int}(A,n) \simeq \bigoplus_{p\geq 0} Sym^{p}_{A}(\mathbb{T}_{A}^{int}[-n])[n]$$ in $\C^{gr}$. Then, recall that for any $A$-module $P$ in $\C$, maps $\mathbf{1}_{\C} \to P$ in $\C$ are in bijection with maps $A\simeq \mathbf{1}_{\C} \otimes_{\C} A \to P$ in  $A-Mod_{\C}$.

\begin{df} Let $A \in \cdga_{\C}$ such that $\mathbb{L}_{A}^{int}$ is a
  dualizable $A$-module in $\C$.  \begin{itemize}
\item A morphism in $Lie_{\C}^{gr}$ $$
k(2)[-1] \longrightarrow \mathbf{Pol}(A,n)[n]= |\mathbf{Pol}^{int}(A,n)[n] |
$$
is \emph{non-degenerate} if the map in $\C$ (induced by the adjunction data for $|-|$) 
$$ 
\mathbf{1}_{\C} \longrightarrow Sym_{A}^{2}(\mathbb{T}_{A}^{int}[-n])[n+1]
$$
yields, by adjunction, an equivalence of $A$-modules
$$
\mathbb{L}_{A}^{int} \simeq \mathbb{T}_{A}^{int}[1-n].
$$
\item We denote by $\mathsf{Map}_{\dglie_k^{gr}}^{nd}(k(2)[-1],\mathbf{Pol}(A,n)[n])$ is
the subspace of
$\mathsf{Map}_{\dglie_k^{gr}}(k(2)[-1],\mathbf{Pol}(A,n)[n])$ of
connected components of non-degenerate morphisms.
\end{itemize}
\end{df}
 
By Theorem~\ref{valerio2} and Corollary \ref{cpoissnd}, we get
 
\begin{cor}\label{cpoissnd2} Let $A \in \cdga_{\C}$ such that $\mathbb{L}_{A}^{int}$ is a
  dualizable $A$-module in $\C$
\begin{enumerate}
\item The map in Thm.~\ref{valerio2} restrict to an equivalence $$\mathbb{P}_{n}(A)^{nd} \simeq \mathsf{Map}_{\dglie_k^{gr}}^{nd}(k(2)[-1],\mathbf{Pol}(A,n)[n]).$$
\item  There is a natural morphism (induced by $\W_{A} $ of Cor. \ref{cpoissnd}) in $\T$ 
$$
\mathsf{Map}_{\dglie_k^{gr}}^{nd}(k(2)[-1],\mathbf{Pol}(A,n)[n]) \longrightarrow
\mathsf{Symp}(A,n-1),
$$
functorial in $A$ with respect to formally \'etale
  morphisms.
\end{enumerate}
\end{cor}

%%%%%%%%%%%%%%%%%%%%%%% A QUI

\subsection{Mixed graded modules: Tate realization}

One of the most important situations in which we will use the above
formalism of de Rham objects and shifted polyvectors is when $\C$ is
itself the $\s$-category of graded mixed $k$-dg-modules, or more
generally diagrams of such. The situation gets complicated because
several different graded mixed structures interact in this
setting. The language of relative differential calculus developed in
the previous section comes handy here and allows us to avoid
confusion.

Throughout this subsection, $M=\epsilon-dg^{gr}_{k}$.  $M$ is a
symmetric monoidal category. Recall that, unless otherwise stated, it
will be endowed with the injective model structure, for which
cofibrations and weak equivalences are defined on the underlying
graded complexes of $k$-modules; as such is a symmetric monoidal model
category satisfying our standing assumptions (see
Section~\ref{1.1}). We let $\C=\medg_k$ be the corresponding
$\s$-category.  Recall that for $M=\epsilon-dg^{gr}_{k}$, and $E, F
\in M$, the $\dg_{k}$-enriched hom object is explicitly given by
$$
\underline{Hom}_{k} (E,F)\equiv \underline{Hom} (E,F) :=
\mathrm{Z}_{\epsilon}(\underline{\mathrm{Hom}}^{gr}_{\epsilon}(E,F)
(0)) \in \dg_{k}
$$
where $\underline{\mathrm{Hom}}^{gr}_{\epsilon}$ denotes the internal
hom object in $M$ (see Section~\ref{1.1}), and, for $X \in M$, we denoted
by $\mathrm{Z}_{\epsilon}(X(0)) \in \dg_{k}$ the kernel of the map of
dg-modules $\epsilon : X(0) \to X(1)[1]$.  The corresponding
$\dg_{k}$-tensor structure is given by $$V\otimes E := V(0)
\otimes_{M} E$$ where $V(0)$ is the mixed graded dg-module
concentrated in weight $0$ and with trivial mixed differentials, and
$\otimes_{M}$ is the monoidal structure in $M$ (Section~\ref{1.1}). Note
that the functor  $\dg_{k} \to M$ sending $V$ to $V(0)$ (in the
notation just introduced) is exactly the symmetric monoidal left
Quillen functor defining the $\dg_{k}$-algebra model structure on $M$.

The category of commutative monoids in $M$ is simply the category
$\epsilon-cdga^{gr}_{k}$ of graded mixed cdgas, whose corresponding
$\s$-category is then $\cdga_{\C}=\mecdga_k$. As already observed
earlier in this section, we have a forgetful $\s$-functor
$$
\mathrm{U}_{\epsilon} : \medg_k \longrightarrow \dg_{k}^{gr}
$$
forgetting the mixed structure.
This $\s$-functor is induced by a left Quillen symmetric monoidal
functor and thus induces a functor 
$$  
\mathrm{U}_{\epsilon} : \mecdga_k \longrightarrow \cdga_{k}^{gr}
$$
It is easy to see that this $\s$-functor preserves de Rham objects, in
the sense that, for any $A \in \cdga_{\C}=\mecdga_k$, the natural
morphism\footnote{About the target, recall that $\mathbb{L}^{int}_{B}
  \simeq \mathbb{L}_{B}$ in $\dg_k$ (respectively in $\dg^{gr}_{k}$) for any
  $B \in \cdga_{k}$ (respectively $B \in \cdga^{gr}_{k}$). }
$$
\mathrm{U}_{\epsilon}(\mathbb{L}^{int}_{A}) \longrightarrow
\mathbb{L}^{int}_{\mathrm{U}_{\epsilon}(A)}
$$
induces an equivalence 
$$
\mathrm{U}_{\epsilon}(\DR^{int}(A)) \simeq
\DR^{int}(\mathrm{U}_{\epsilon}(A)),
$$
of graded mixed cdga inside the $\s$-category $\dg_k^{gr}$ of graded
dg-modules (note that on the left hand side the functor
$\mathrm{U}_{\epsilon}$ sends $\epsilon-\cdga^{gr}_{\C}$ to
$\epsilon-\cdga^{gr}_{\dg^{gr}_{k}}$). At the level of strict models
this is even simpler, as for $A$ a graded mixed cdga, the graded mixed
$A$-module $\Omega^1_{A}$ is canonically isomorphic, as a graded
$A$-module, to $\Omega^1_{\mathrm{U}_{\epsilon}(A)}$. In other words,
in order to compute $\Omega^1_A$ as a graded mixed $A$-module we
simply compute it as a graded $A$-module, and then endow it with the
natural mixed structure coming from the one on $A$.

Recall (Definition~\ref{dreal} with $M=\epsilon-dg^{gr}_{k}$) that we have
defined a realization functor $$|-|: \mathcal{M}=\medg_k
\longrightarrow \dg_k$$ as the $\s$-functor
$\mathbb{R}\underline{Hom}(1_M ,-)$ associated to the right derived
functor of the Quillen right adjoint to the functor $- \otimes 1_{M}:
dg_{k} \to M$ (here $1_M = k(0)$ is $k$ sitting in weight $0$, degree
$0$, with trivial differential and trivial mixed differential). As
above, $M$ is endowed here with the injective model
structure, for which the monoidal unit $1_{M}$ is cofibrant. However,
$M$ can also be given the \emph{projective} model structure
$M^{\textrm{proj}}$ where fibrations and weak equivalences are defined
on the underlying graded complexes of $k$-modules. In
$M^{\textrm{proj}}$ the monoidal unit $1_{M}$ is no longer cofibrant,
and we have already constructed in \ref{closedforms} an explicit
cofibrant replacement $\widetilde{k} \to 1_{M}$ in
$M^{\textrm{proj}}$. Moreover, $\widetilde{k}$ is a counital comonoid
object in $M$, therefore we have a Quillen pair 
$$
-\otimes
\widetilde{k} : dg_{k} \longleftrightarrow M :
\underline{Hom}(\widetilde{k},-)
$$ 
where the right adjoint is lax
symmetric Quillen monoidal. The identity functor on $M$ induces an
identification (equivalence) on the associated $\s$-categories, and
the realization functor $|-|$ is equivalent, under this
identification, to the $\s$-functor induced by the right derived
Quillen functor $\mathbb{R}\underline{Hom}(\widetilde{k} ,-)$,
i.e. the derived
functor with respect to the projective model structure on $M$. Since
in $M^{\textrm{proj}}$, unlike in the injective model structure on
$M$, every object is fibrant, we have
$\mathbb{R}\underline{Hom}(\widetilde{k} ,-)\simeq
\underline{Hom}(\widetilde{k} ,-)$. Thus we conclude that as
$\s$-functors we have an equivalence
$$
\mathbb{R}\underline{Hom}(k(0),-) := |-| \simeq
\underline{Hom}(\widetilde{k} ,-) \ : \ \medg_k \longrightarrow \dg_k
$$

\

\begin{prop}\label{realnotenuff}
For any $E \in M$, there is a canonical isomorphism of $k$-dg modules
$$ 
\prod_{p\geq 0}E(p) \simeq \underline{Hom}_{k}(\widetilde{k}, E)
$$
where the source is endowed with the total differential, i.e. the sum of the
cohomological and the mixed differentials.
\end{prop}

\noindent \textbf{Proof.} The complex 
$\underline{\mathrm{Hom}}^{gr}_{\epsilon}(\widetilde{k},E) (0) \in
C(k)$ is given in degree $n$ by $$E(0)^{n} \times \prod_{p>0}
(E(p)^{n} \times E(p)^{n+1}).$$ The map $f: \prod_{p\geq 0}E(p) \to
\underline{\mathrm{Hom}}^{gr}_{\epsilon}(\widetilde{k},E) (0)$ defined
(with obvious notations) in degree $n$ by 
$$
f^n : \{x_{0},
(x_{p})_{p>0} \} \longmapsto \{x_{0}, (x_{p},
-\epsilon_{E}(x_{p-1}))_{p>0} \} 
$$ 
is a map of complexes, and the
composite 
$$
\xymatrix{ \prod_{p\geq 0}E(p) \ar[r]^-{f} &
  \underline{\mathrm{Hom}}^{gr}_{\epsilon}(\widetilde{k},E) (0)
  \ar[rr]^{\epsilon_{\underline{Hom}_{M}}} & &
  \underline{\mathrm{Hom}}^{gr}_{\epsilon}(\widetilde{k},E) (1)[1]}
$$
is zero. A computation now shows that the induced
map $\overline{f}: \prod_{p\geq 0}E(p) \simeq
\underline{Hom}_{k}(\widetilde{k}, E)$ is an isomorphism of
$k$-dg-modules. \hfill $\Box$

By Proposition~\ref{realnotenuff}, we get that the $\s$-functor
$$
|-| : \medg_k \longrightarrow \dg_k
$$ 
has a canonical strict model given by 
$$
E \longmapsto \prod_{p\geq 0}E(p),
$$
where the right hand side is endowed with the total differential $=$ sum of
the cohomological differential and the mixed structure.

Since for any $i \in \mathbb{Z}$ the $(-i)$-weight shift
$\widetilde{k}((-i))$ is a cofibrant resolution of $k(i)$ (i.e. of
$k[0]$ concentrated in weight $i$) in $M^{\textrm{proj}}$, the above
computation yields the following equivalences in $\dg_k$
$$
\mathbb{R}\underline{Hom}_{k}(k(i),k(i+1)) \simeq k.
$$
We thus have a canonical morphism 
$u_i : k(i) \longrightarrow k(i+1)$ in $\edg_k$ for
all $i\in \mathbb{Z}$, corresponding to $1\in k$ in the above
formula. In particular,  
we get a pro-object in $\edg^{gr}_k$
$$
k(-\s):=\{ \cdots \to k(-i) \rightarrow k(-i+1) \rightarrow \cdots \to
k(-1)
\rightarrow k(0)\}.
$$

\begin{df}\label{dtate}
  The \emph{Tate} or \emph{stabilized realization $\s$-functor} is
  defined to be
$$|-|^t := \mathbb{R}\underline{Hom}_{k}(k(-\s),-) : \medg_k
\longrightarrow \dg_k,$$ 
sending $E \in \edg^{gr}_k$ to 
$$|E|^t = \textrm{colim}_{i\geq 0} \mathbb{R}\underline{Hom}_{k}(k(-i),E) \simeq 
\textrm{colim}_{i\geq 0} \prod_{p\geq -i} E(p).$$
\end{df}

%Explicitely, we have
%$$|E|^t \simeq \textrm{colim}_{i\geq 0} \prod_{p\geq -i} E(p).$$
The natural map $k(-\s) \longrightarrow k(0)$ of pro-objects in
$\edg_k$ (where $k(0)$ is considered as a constant pro-object)
provides a natural transformation
$$|-| \longrightarrow |-|^t$$ 
from the standard realization to the Tate realization. By definition,
we see that this natural transformation induces an equivalence
$|E| \simeq |E|^t$ in $\dg_{k}$, as soon as $E(p)=0$ for all $p<0$.

The $\s$-functor $|-|$ is lax symmetric monoidal, and this endows
$|-|^t$ with a canonical structure of a lax symmetric monoidal
$\s$-functor. This follows, for instance, from the fact that the
pro-object $k(-\s)$ defined above is a cocommutative and counital
coalgebra object, which is the dual of the commutative and unital
algebra $\textrm{colim}_{i\geq 0}k(i)$.  Therefore the Tate
realization induces an $\s$-functor on commutative algebra objects in
$\mathcal{M}=\medg_k$, and more generally on all kind of algebra-like
structures in $\mathcal{M}$. In particular, we have Tate realization
functors, denoted with the same symbol, for graded mixed cdgas over
$\medg_k$, as well as for graded $\mathbb{P}_{n+1}$-cdgas
$$|-|^t : \mecdga_{\medg_k} \longrightarrow \mecdga_{k}$$
$$|-|^t : \mathbb{P}_{n+1}-\cdga^{gr}_{\medg_k} \longrightarrow 
\mathbb{P}_{n+1}-\cdga^{gr}_{k}.$$
This way we get Tate versions of the de Rham and shifted polyvectors objects 
introduced in Def. \ref{d4} and \ref{d9}. 

\begin{df}\label{drpoltate}
Let $A \in \cdga_{\medg_{k}}$ be commutative cdga in the $\s$-category
of graded mixed complexes (i.e. a graded mixed cdga over $k$). 
\begin{enumerate}

\item The \emph{Tate de Rham complex of $A$} is defined by
$$\DR^{t}(A):=|\DR^{int}(A)|^t \in \mecdga_k.$$

\item The \emph{Tate $n$-shifted polyvectors cof $A$} is defined by
$$\Pol^{t}(A,n):=|\Pol^{int}(A,n)|^t \in \mathbb{P}_{n+1}-\cdga_{k}^{gr}.$$
\end{enumerate}
\end{df}

Note that we have natural induced morphisms
$$\DR(A) \longrightarrow \DR^t(A) \qquad \Pol(A,n) \longrightarrow
\Pol^{t}(A,n)$$ 
which are not always equivalences. More precisely, if $A(p)=0$ for all
$p< 0$, then $\mathbb{L}_{A}^{int}$ is itself only positively
weighted, and we get $\DR(A) \simeq \DR^t(A)$ by the natural
morphism. On the other hand, $\mathsf{Pol}(A,n)$ has in general both
positive and non-positive weights, as the weights of
$\mathbb{T}^{int}_{A}$ are dual to that of $A$.  So, except in some
very degenerate cases, $\mathsf{Pol}(A,n) \longrightarrow
\mathsf{Pol}^t(A,n)$
will typically \emph{not} be an equivalence. 

\

To finish this section we mention the Tate analogue of the
morphism constructed in Corollary~\ref{cpoissnd2} from the space of
non-degenerate $n$-shifted Poisson structures to the space of
$n$-shifted symplectic structures.

The notion of Tate realization functor, can be interpreted as a
standard realization functor for a slight modification of the base
$\s$-category $\C=\medg$. The same is true for the objects $\DR^t(A)$
and $\Pol^t(A)$ at least under some mild finiteness conditions on
$A$. In order to see this, we let $\C':=\mathsf{Ind}(\C)$ be the
$\s$-category of Ind-objects in $\C$. The $\s$-category $\C'$ is again
symmetric monoidal and possesses as a model the model category
$\mathsf{Ind}(M)$ of Ind-objects in $M$ (see \cite[Thm. 1.5]{ind}):
$$\mathsf{Ind}(\C) \simeq L(\mathsf{Ind}(M)).$$

We consider the following Ind-object in $\C$
%$$k(-\s) := \xymatrix{k(0) \ar[r] & k(-1) \ar[r] & \dots k(-i) \ar[r]
%& k(-i-1) \dots}$$ 
$$k(\s) := \{\xymatrix{k(0) \ar[r] & k(1) \ar[r] & \dots k(i) \ar[r] &
  k(i+1) \ar[r] & \cdots }\}$$ 
which is objectwise dual to the pro-object $k(-\s)$ we have considered above. 
Now, the standard realization $\s$-functor $|-|: \C' \to \dg_{k}$ for $\C'$ 
recovers the Tate realization on $\C$, since we have a naturally
commutative diagram of $\s$-functors
$$\xymatrix{
\C \ar[rr]^-{-\otimes k(\s)} \ar[rd]_-{|-|^t} & &  \C' \ar[ld]^-{|-|} \\
 & \dg_k & }$$

Moreover, the natural equivalences $k(i)\otimes k(j) \simeq k(i+j)$
makes $k(\s)$ into a commutative cdga in $\C'=\mathsf{Ind}(\C)$. For
any $A\in \cdga_{\C}$, viewed as a constant commutative cdga in $\C'$
via the natural functor $\mathcal{M} \to
\mathsf{Ind}(\mathcal{M})=\C'$, we thus have a natural object obtained
by base change
$$A(\s):=A\otimes k(\s) \in \cdga_{\C'}.$$ 
Note that, as an Ind-object in $\C$, we have
$$A(\s) = \{ \xymatrix{A\otimes k(0) \ar[r] & A\otimes  k(1) \ar[r] & \dots 
A\otimes k(i) \ar[r] & A\otimes  k(i+1) \ar[r] & \cdots} \}$$
The cdga $A(\s)$ will be considered as 
a $k(\s)$-algebra object in $\C'$
$$A(\s) \in k(\s)-\cdga_{\C'}= k(\s)/\cdga_{\C'} .$$
It therefore has the corresponding relative de Rham 
and polyvector objects
$$\DR^{int}(A(\s)/k(\s)) \in \mecdga_{\mathcal{M'}} \qquad
\Pol^{int}(A(\s)/k(\s),n) \in
\mathbb{P}_{n+1}-\cdga_{\mathcal{M'}}^{gr} ,$$ 
and, as usual, we will denote by  
$$\DR(A(\s)/k(\s)) \in \mecdga_{k} \qquad \Pol(A(\s)/k(\s),n) \in
\mathbb{P}_{n+1}-\cdga_{k}^{gr}$$ 
the corresponding images under the standard realization $|-|: \C' \to
\dg_{k}$, that, recall, is lax symmetric monoidal so it sends $\mecdga_{\mathcal{M'}}$ to $\mecdga_{k}$, and  $\mathbb{P}_{n+1}-\cdga_{\mathcal{M'}}^{gr}$ to $\mathbb{P}_{n+1}-\cdga_{k}^{gr}$.

The following lemma compares de Rham and polyvectors objects of $A\in
\cdga_{\C}$, and of $A(\s)$ relative to $k(\s)$, under suitable
finiteness hypotheses on $A$.

\begin{lem}\label{l(-s)}
  If $A \in \cdga_{\C}$ is such that $\mathbb{L}_{A}^{int}$ is a
  perfect (i.e. dualizable) $A$-module, then there are natural
  equivalences of graded mixed cdgas over $k$ and, respectively, of
  graded $\mathbb{P}_{n+1}$-algebras over $k$
$$\DR^{t}(A) \simeq \DR(A(\s)/k(\s))$$
$$\Pol^t(A,n) \simeq \Pol(A(\s)/k(\s),n).$$
\end{lem}

\noindent \textbf{Proof.} Without any assumptions on $A$, we have
$$\DR^{int}(A)\otimes k(\s) \simeq \DR^{int}(A(\s)/k(\s)).$$
Since, as already observed, $|-\otimes k(\s)| \simeq |-|^{t}$, this
shows that $\DR^{t}(A) \simeq \DR(A(\s)/k(\s))$.

For polyvectors, the dualizability condition on $\mathbb{L}_{A}^{int}$
implies that the natural morphism
$$\Pol^{int}(A,n)\otimes k(\s) \longrightarrow \Pol^{int}(A(\s)/k(\s),n)$$
is an equivalence. So, again, we have 
$$\Pol^t(A,n) \simeq \Pol(A(\s)/k(\s),n).$$
\hfill $\Box$

We can therefore state a Tate version of Corollary \ref{cpoissnd2}, by
working in $\C'$, for $A \in \cdga_{\C}$ with dualizable
$\mathbb{L}_{A}^{int}$. In the corollary below the non-degeneracy
conditions is required in $\C'$, that is after tensoring with
$k(\s)$. This modifies the notion of shifted symplectic structures as
follows. If $\A^{2,cl}(A,n)$ is the space of closed $2$-forms of
degree $n$ on $A$, we say that an element $\omega \in \pi_{0}
\,\A^{2,cl}(A,n)$ is \emph{Tate non-degenerate} if the underlying
adjoint morphism in $\C$
$$\Theta_{\omega_0} : \mathbb{T}_{A}^{int} \longrightarrow
\mathbb{L}_{A}^{int}[n]$$ 
induces an equivalence in $\C'$ 
$$\Theta_{\omega_0}(\s) : \mathbb{T}_{A}^{int}(\s) \longrightarrow 
\mathbb{L}_{A}^{int}(\s)[n]$$
i.e. after tensoring with $k(\s)$. The space $\mathsf{Symp}^t(A,n)$ of
\emph{$n$-shifted Tate symplectic structures} 
on $A$ is then the  subspace of $\A^{2,cl}(A,n)$ consisting of
connected components  
of Tate non-degenerate elements. Note that by Lemma \ref{l(-s)} we have
$$\mathsf{Symp}^t(A,n) \simeq \mathsf{Symp}(A(\s)/k(\s),n),$$
where the right hand side is the space of n-shifted 
symplectic structures on $A(\s)$ relative to $k(\s)$, computed
in $\C'=\mathsf{Ind}(\C)$.

\begin{cor}\label{cpoissnd3}
Let $A \in \cdga_{\C}$ such that $\mathbb{L}_{A}^{int}$ is a
dualizable $A$-module in $\C$. Then, there is a natural morphism of
spaces, functorial in $A$ with respect to formally \'etale morphisms
$$ 
\mathsf{Map}_{\dglie_k^{gr}}^{nd}(k(2)[-1],\mathsf{Pol}^t(A,n)[n])
\longrightarrow \mathsf{Symp}^t(A,n-1),
$$
where
$\mathsf{Map}_{\dglie_k^{gr}}^{nd}(k(2)[-1],\mathsf{Pol}^t(A,n)[n])$
is the subspace of
$\mathsf{Map}_{\dglie_k^{gr}}(k(2)[-1],\mathsf{Pol}^t(A,n)[n])$
consisting of connected components of non-degenerate elements.
\end{cor}

\section{Formal localization}\label{FLsection}

A commutative dg-algebra (in non-positive degrees) $A$ over $k$ is
\emph{almost finitely presented} if $H^0(A)$ is a $k$-algebra of
finite type, and each $H^{i}(A)$ is a finitely presented
$H^{0}(A)$-module. Notice that, in particular, such an $A$ is
\emph{Noetherian} i.e. $H^{0}(A)$ is a Noetherian $k$-algebra (since
our base $\mathbb{Q}$-algebra $k$ is assumed to be Noetherian), and
each $H^{i}(A)$ is a finitely presented $H^{0}(A)$-module.

We let $\dAff_{k}$ be the opposite $\s$-category of almost finitely
presented commutative dg-algebras over $k$ concentrated in
non-positive degrees.  We will simply refer to its objects as
\emph{derived affine schemes} without mentioning the base $k$ or the
finite presentation condition. When writing $\Spec\, A$, we implicitly
assume that $\Spec\, A$ is an object of $\dAff_{k}$, i.e.  that $A$ is
almost finitely presented commutative $k$-algebra concentrated in
non-positive degrees.  The $\s$-category $\dAff_{k}$ is equipped with
its usual \'etale topology of \cite[Def. 2.2.2.3]{hagII}, and the
corresponding $\s$-topos of stacks will be denoted by $\dSt_{k}$. Its
objects will simply be called \emph{derived stacks} (even though they
should be, strictly speaking, called \emph{locally almost finitely
  presented derived stacks over $k$}).

With these conventions, an algebraic derived $n$-stack will have a
smooth atlas by objects in $\dAff_{k}$, i.e. by objects of the form
$\Spec\, A$ where $A$ is almost finitely presented over $k$.
Equivalently, all our algebraic derived $n$-stacks will be derived
$n$-stacks according to \cite[\S 2]{hagII}, that is such stacks are
defined on the category of all commutative dg-algebra concentrated in
non-positive degrees. Being locally almost of finite presentation
these stacks $X$ have cotangent complexes which are in $\mathsf{Coh}(X)$
and bounded on the right.

\subsection{Derived formal stacks}

We start by a zoology of derived stacks with certain infinitesimal
properties.

\begin{df}\label{d10-}
A \emph{formal derived stack} is
an object $F \in \dSt_{k}$
satisfying the following conditions.

\begin{enumerate}

\item The derived stack $F$ is \emph{nilcomplete} i.e. for all $\Spec\, B \in \dAff_{k}$, the canonical map 
$$F(B)\longrightarrow \lim_{k}F(B_{\leq k}),$$
where $B_{\leq k}$ denotes the $k$-th Postnikov truncation of $B$, is
an equivalence in $\T$.

\item The derived stack $F$ is \emph{infinitesimally cohesive}
  i.e. for all cartesian squares of almost finitely presented
  $k$-cdgas in non-positive degrees
$$\xymatrix{
  B \ar[r] \ar[d] & B_{1} \ar[d] \\
  B_{2} \ar[r] & B_{0},}$$ such that each $\pi_{0}(B_i)
\longrightarrow \pi_{0}(B_0)$ is surjective with nilpotent kernel,
then the induced square
$$\xymatrix{
F(B) \ar[r] \ar[d] & F(B_{1}) \ar[d] \\
F(B_{2}) \ar[r] & F(B_{0}),}$$
is cartesian in $\T$.
\end{enumerate}
\end{df}

\begin{rmk} Note that if one assumes that a derived stack $F$ has a
  cotangent complex (\cite[\S 1.4]{hagII}), then $F$ is a formal
  derived stack if and only if it is nilcomplete and satisfies the
  infinitesimally cohesive axiom where \emph{at least one} of the two
  $B_i \to B_0$ is required to have $\pi_{0}(B_i) \longrightarrow
  \pi_{0}(B_0)$ surjective with nilpotent kernel
  (\cite[Proposition~2.1.13]{luXIV}). We also observe that, even if we
  omit the nilpotency condition on the kernels but keep the
  surjectivity, we have that the diagram obtained by applying $\Spec$
  to the square of cdgas in \ref{d10-} $(2)$ is a homotopy push-out in
  the $\s$-category of derived schemes, hence in the $\s$-category of
  derived \emph{algebraic} stacks (say for the \'etale topology). This
  is a derived analog of the fact that pullbacks along surjective maps
  of rings induce pushout of schemes. In particular, any derived
  algebraic stack $F$ sends any diagram as in \ref{d10-} $(2)$, with
  the nilpotency condition possibly omitted, to pullbacks in $\T$,
  i.e. is actually \emph{cohesive} (\cite[DAG IX, Corollary~6.5]{lu2}
  and \cite[Lemma 2.1.7]{luXIV}).
\end{rmk}

There are various sources of examples of formal derived stacks. 

\begin{itemize}

\item Any algebraic derived $n$-stack $F$, in the sense of \cite[\S
  2.2]{hagII}, is a formal derived stack. Nilcompleteness of $F$ is
  (the easy implication of) \cite[Theorem~c.9 (c)]{hagII}, while the
  infinitesimally cohesive property follows from nilcompleteness, the
  existence of a cotangent complex for $F$, and the general fact that
  any $B_i\to B_0$ with $\pi_{0}(B_i) \to \pi_{0}(B_0)$ surjective
  with nilpotent kernel can be written as the limit in
  $\cdga_{k}/B_{0}$ of a tower $\cdots \to C_n \to \cdots C_1 \to
  C_0:=B_0$ where each $C_{n}$ is a square-zero extension of $C_{n-1}$
  by some $C_{n-1}$-module $P_{n}[k_n]$, where $k_{n} \to +\infty$ for
  $n\to +\infty$ (see \cite[Lemma 2.1.14]{luXIV}, or
  \cite[Proposition~2.1.13]{luXIV} for a full proof of the infinitesimal
  cohesive property for a stack that is nilcomplete and has a
  cotangent complex). Alternatively, one can observe (\cite[Lemma
  2.1.7]{luXIV}) that any derived algebraic stack is actually cohesive
  (hence infinitesimally cohesive).

\item For all $\Spec\, A \in \dAff_{k}$ we let $QCoh^{-}(A)$ be the
full sub-$\s$-groupoid of $\mathrm{L}(A)$ consisting of $A$-dg-modules
$M$ with $H^{i}(M)=0$ for $i>>0$. The $\s$-functor 
$A \mapsto QCoh^{-}(A)$ defines a derived stack which can be 
checked to be a formal derived stack. 

\item Any (small) limit, in $\dSt_{k}$, of formal derived stacks is
  again a formal derived stack. This follows from the fact that (by
  Yoneda), for any $A \in \cdga_{k}$, the functor $\dSt_{k} \to \T$
  given by evaluation at $A$ commutes with (small) limits, and that
  both convergence and infinitesimal cohesiveness are expressed by
  conditions on objectwise limits.
\end{itemize}

Let us consider the inclusion functor $i:
\mathbf{alg}^{\textrm{red}}_{k} \longrightarrow \cdga_{k}$ of the full
reflective sub $\s$-category of \emph{reduced discrete} objects
(i.e. $R\in \cdga_{k}$ such that $R$ is discrete and $R \simeq H^0(R)$
is a usual reduced $k$-algebra).  The functor $i$ has a left
adjoint 
$$(-)^{red}: \cdga_{k} \longrightarrow
\mathbf{alg}^{\textrm{red}}_{k} \, , \, \, A \longmapsto
A^{red}:=H^0(A)/\textrm{Nilp}(H^0 (A).$$ 
Moreover, it is easy to
verify that $i$ is both continuous and cocontinuous for the \'etale
topologies on $\cdga_{k}^{op}$, and
$(\mathbf{alg}^{\textrm{red}}_{k})^{op}$. If we denote by
$\mathbf{St}_{\textrm{red}, k}$ the $\s$-category of stacks on
$(\mathbf{alg}^{\textrm{red}}_{k})^{op}$ for the \'etale topology, we
thus get an induced $\s$-functor 
$$i^*: \dSt_{k} \longrightarrow
\mathbf{St}_{\textrm{red}, k}$$ 
that has both a right adjoint $i_*$,
and a left adjoint $i_!$, obeying the following
properties:\begin{itemize}
\item $i_* \simeq ((-)^{red})^*$ (thus $i^*\Spec \, A \simeq
  \mathrm{Spec} (A^{red})$).
\item $i_!$ and $i_*$ are fully faithful (equivalently, the adjunction
  maps $\mathrm{Id} \to i^*i_!$ and $i^*i_* \to \mathrm{Id}$ are
  objectwise equivalences).
\item $i^*\simeq ((-)^{red})_{!}$.
\item $i_!i^*$ is left adjoint to $i_*i^*$.
\end{itemize}

%For any $F\in \dSt_{k}$ we can consider the left (respectively,
%right) Kan extension $\mathsf{Lan}(F\circ i, i)$
%(resp. $\mathsf{Ran}(F\circ i, i)$)  of $F\circ i$ along $i$. It is
%easy to verify that $$\mathsf{Ran}(F\circ i, i)(A) \simeq
%F(A^{red})$$ while $$\mathsf{Lan}(\Spec A \circ i, i) \simeq \Spec \,
%(A^{red}).$$ 

\begin{df}\label{red&dR} \begin{enumerate}
\item The functor 
$$(-)_{DR} := i_*i^* : \dSt_{k} \longrightarrow
  \dSt_{k}$$ 
is called the \emph{de Rham stack functor}. By
  adjunction,  for any $F\in \dSt_{k}$, we have a canonical natural
  map $\lambda_F: F \mapsto F_{DR}$. 
\item The functor 
$$(-)_{\textrm{red}} := i_! i^* : \dSt_{k} \longrightarrow \dSt_{k}$$  
is called the \emph{reduced stack functor}. By adjunction,  for any
$F\in \dSt_{k}$, we have a canonical natural map $\iota_{F} :
F_{\textrm{red}} \mapsto F$.  
\item Let $f : F \longrightarrow G$ be a morphism
in $\dSt_{k}$. We define the \emph{formal completion} $\widehat{G}_{f}$ 
\emph{of} $G$ \emph{along the morphism} $f$ as the fibered product in
$\dSt_{k}$: 
$$\xymatrix{\widehat{G}_{f} \ar[r]^-{\beta_{f}} \ar[d] &
  F_{DR} \ar[d]^-{f_{DR}} \\ G \ar[r]_-{\lambda_{G}} & G_{DR}.}$$  
\end{enumerate}
\end{df}

\

\medskip

\noindent Since the left adjoint to $i$ is $(-)^{red}$, then it is
easy to see that 
$$F_{DR}(A) \simeq F(A^{red}) \,\,\,\,\,\,\,\,
\textrm{and} \,\,\,\,\,\,\,\,(\Spec \, A)_{\textrm{red}} \simeq \Spec
\, (A^{red}),$$ 
for any $A\in \cdga_{k}$. Therefore
$\widehat{G}_{f}(A)=G(A) \times_{G(A^{red})}F(A^{red})$, for $f: F \to
G$ in $\dSt_{k}$. We already observed that $(-)_{DR}$ is right adjoint
to $(-)_{\textrm{red}}$, as functors $\dSt_{k}\to \dSt_{k}$.

%\begin{df}\label{red&dR} \begin{enumerate}
%\item For $F\in \dSt_{k}$, we put $$F_{DR}:=\mathsf{Ran}(F\circ i, i)
%\,\, \in \dSt_{k},
%\,\,\,\,\,\,\,\,\,\,\,\,\,\,\,\,\,\,\,\,\,\,\,\,\,\,\,\,\,\,\,\,\,\,\,\,\,
%F_{\textrm{red}}:= \mathsf{Lan}(F\circ i, i)\,\, \in \dSt_{k}.$$
%$F_{DR}$ will be called the \emph{de Rham stack} of $F$, and
%$F_{\textrm{red}}$ the \emph{reduced stack} of $F$. Note that,  by
%definition of left and right Kan extensions, we have functorial maps
%in $\dSt_{k}$ $$F \longrightarrow F_{DR} \,\,\,\,\,\,\,\,\,\,\,\,\,
%\,\,\,\,\,\,\,\,\,\,\,\,\,  F_{\textrm{red}} \longrightarrow F $$ 
%\item Let $f : F \longrightarrow G$ be a morphism
%in $\dSt_{k}$. We define the \emph{formal completion} $\widehat{G}_{f}$ 
%\emph{of} $G$ \emph{along the morphism} $f$ as the fibered product in $\dSt_{k}$ $$\widehat{G}_{f}:= G\times_{G_{DR}} F_{DR}$$ i.e.  as the $\s$-functor 
%sending $X=\Spec\, A\in \dAff_{k}$ to
%$$\widehat{G}_{f}(A):=G(A) \times_{G(A^{red})}F(A^{red}),$$
%where $A^{red}:=H^{0}(A)/\textrm{Nilp}(H^0 (A))$. 
%\end{enumerate}
%\end{df} 

Since taking the reduced algebra is a projector, we have that the
canonical map  $F_{DR} \to (F_{DR})_{DR}$ is an equivalence; the same
holds for $(F_{\textrm{red}})_{\textrm{red}} \to
F_{\textrm{red}}$. Moreover, for any $F \in \dSt_{k}$, we have $F_{DR}
\simeq \widehat{(\Spec \, k)}_{f}$, where $f: F \to \Spec \, k$ is the
structure morphism (note that the canonical map $\Spec \, k \to
(\Spec \, k)_{DR}$ is an equivalence). We list below a few elementary
properties of de Rham stacks and reduced stacks.

\begin{prop}\label{easy}
\begin {enumerate}
\item $F_{DR}$ is a formal derived stack for any $F \in \dSt_{k}$.
\item If $G$ is a formal derived stack, the formal completion
  $\widehat{G}_{f}$, along any map $f:F\to G$ in $\dSt_{k}$, is a
  formal derived stack. 
\item For any $F \in \dSt_{k}$, the canonical map $\lambda_{F}: F\to
  F_{DR}$ induces an equivalence  $F_{\textrm{red}} \to
  (F_{DR})_{\textrm{red}}$. 
\item For any map $f:F\to G$ in $\dSt_{k}$, the canonical map
  $\alpha_{f}: F \to \widehat{G}_{f} $ induces an equivalence
  $F_{\textrm{red}} \to (\widehat{G}_{f})_{\textrm{red}} $.
\item For any $F \in \dSt_{k}$, the canonical map $F_{\textrm{red}}
  \to F $ induces an equivalence $(F_{\textrm{red}})_{DR} \to F_{DR}$.
\item For any $F \in \dSt_{k}$, if $j: \mathsf{t}_0 F \to F$ denotes
  the canonical map in $\dSt_{k}$ from the truncation of $F$ to $F$,
  then the canonical map $\widehat{F}_{j} \to F$ is an equivalence.
\item If $f:F\to G$ is a map in $\dSt_{k}$ such that
  $f_{\textrm{red}}$ is an equivalence, then the canonical map
  $\widehat{G}_{f} \to G$ is an equivalence.
\end{enumerate}
\end{prop} 
\textbf{Proof.} Since $F_{DR} (A) = F(A^{red})$, and $(-)^{red}$ sends
cartesian squares as in Definition~\ref{d10-} $(2)$ to cartesian
squares of isomorphisms, $(1)$ follows. $(2)$ follows from $(1)$ and
the fact that formal derived stacks are closed under small
limits. $(3)$ follows from the fact that $F$ and $F_{DR}$ agree when
restricted to $\mathbf{alg}^{\textrm{red}}_{k}$ i.e. $i^*(\lambda_F) :
i^*F \to i^* (F_{DR})$ is an equivalence and hence $i_!i^*(\lambda_F)$
is an equivalence as well.  Let us prove $(4)$. Both $G$ and $G_{DR}$
agree on $\mathbf{alg}^{\textrm{red}}_{k}$, $i^*(\lambda_G) : i^*G \to
i^* (G_{DR})$ is an equivalence, and $i^*$ is right (and left)
adjoint, and also $i^*(\beta_{f}): i^*(\widehat{G}_{f}) \to
i^*(F_{DR})$ is an equivalence. Furthermore $i_!$ is fully faithful,
so $\beta_{f,\textrm{red}}: (\widehat{G}_{f})_{\textrm{red}} \to
(F_{DR})_{\textrm{red}}$ is an equivalence. Now, the composite
$\alpha_{f} \circ \beta_{f}$ is equal to $\lambda_{F}$, so we get
$(4)$ from $(3)$. In order to prove $(5)$ it is enough to observe that
the adjunction map $i^*i_! \to \mathrm{Id}$ is an objectwise
equivalence. The assertion $(6)$ follows immediately by observing that
$\mathsf{t}_0 F(S) = F(S)$ for any discrete commutative $k$-algebra
$S$, therefore $j_{DR}$ is an equivalence. Finally, by $(5)$ we get
that if $f_{\textrm{red}}$ is an equivalence, then so is $f_{DR}$, and
so $(7)$ follows.  \hfill $\Box$

\

\begin{df}\label{d10}
\begin{enumerate}
\item A formal derived stack $F$ according to Definition \ref{d10-} is
  called \emph{almost affine} if $F_{red} \in \dSt_{k}$ is an affine
  derived scheme.
\item An almost affine formal derived stack $F$ in the sense above is
  \emph{affine} if $F$ has a cotangent complex in the sense of
  \cite[\S 1.4]{hagII}, and if, for all $\Spec \, B \in \dAff_{k}$ and
  all morphism $u : \Spec\, B \longrightarrow F$, the $B$-dg-module
  $\mathbb{L}_{F,u} \in \mathsf{L}(B)$ (\cite[Definition~1.4.1.5]{hagII}) is
  coherent and cohomologically bounded above.
\end{enumerate}
\end{df}

Recall our convention throughout this section, that all derived affine
schemes are automatically assumed to be almost of finite
presentation. Therefore, any derived affine scheme is an affine formal
derived stack according to Definition~\ref{d10}.

Note that when $F$ is any affine formal derived stack, there is a
globally defined quasi-coherent complex $\mathbb{L}_{F} \in
\mathsf{L}_{\mathsf{Qcoh}}(F)$ such that for all $u : \Spec\, B
\longrightarrow F$, we have a natural equivalence of $B$-dg-modules
$$u^{*}(\mathbb{L}_{F}) \simeq \mathbb{L}_{F,u}.$$
The quasi-coherent complex $\mathbb{L}_{F}$ is then itself coherent,
with cohomology bounded above.

Since $(\Spec\, A)_{\textrm{red}} \simeq \Spec (A^{red})$, we get by
Proposition~\ref{easy} (4), that for any algebraic derived stack $F$,
and any morphism in $\dSt_{k}$
$$f : \Spec\, A \longrightarrow F\, ,$$ 
the formal completion $\widehat{F}_{f}$ of $F$ along $f$ is an affine
formal derived stack in the sense of Definition \ref{d10}
above. Moreover, the natural morphism $v : \widehat{F}_{f}
\longrightarrow F$ is formally \'etale, i.e.  the natural morphism
$$v^{*}(\mathbb{L}_{F}) \longrightarrow \mathbb{L}_{\widehat{F}_{f}}$$
is an equivalence in $\mathsf{L}_{\mathsf{Qcoh}}(\widehat{F}_{f})$.  

\

\noindent
This formal completion construction along a map from an affine will be
our main source of examples of affine formal derived stacks.

We will ultimately be concerned with affine formal derived stacks
\emph{over affine bases}, which we proceed to discuss.

\begin{df}\label{d10+}
  Let $X:=\Spec\, A \in \dAff_{k}$. A \emph{good formal derived stack
    over $X$} is an object $F \in \dSt_{k}/X$ satisfying the following
  two conditions.
\begin{enumerate}
\item The derived stack $F$ is an affine formal derived stack.
\item The induced morphism $F_{red} \longrightarrow (\Spec\,
  A)_{\textrm{red}} = \Spec\, A^{red}$ is an equivalence.
\end{enumerate}

The full sub-$\s$-category of $\dSt_{k}/X$ consisting of good formal
derived stacks over $X =\Spec\, A$ will be denoted as $\dFSt_{X}^{g}$,
or equivalently as $\dFSt_{A}^{g}$.

Finally, a \emph{perfect formal derived stack $F$ over $\Spec\, A$} is
a good formal derived stack over $\Spec\, A$ such that moreover its
cotangent complex $\mathbb{L}_{F/\Spec\, A} \in
\mathsf{L}_{\mathsf{Qcoh}}(F)$ is a perfect complex.
\end{df}

\

\begin{rmk} Since $i_{!}$ is fully faithful, it is easy to see that if
  $F \to \Spec \, A$ is a good (respectively, perfect) formal derived
  stack, then for any $\Spec \, B \to \Spec \, A$, the base change
  $F_B \to \Spec \, B$ is again a good (respectively, perfect) formal
  derived stack. In this sense, good (respectively, perfect) formal
  derived stacks are stable under arbitrary derived affine base
  change.
\end{rmk}

The fundamental example of a good formal derived stack is given by an
incarnation of the so-called \emph{Grothendieck connection} (also
called \emph{Gel'fand connection} in the literature). It consists, for
an algebraic derived stack $F \in \dSt_{k}$ which is locally almost of
finite presentation, of the family of all formal completions of $F$ at
various points. This family is equipped with a natural flat
connection, or in other words, is a crystal of formal derived stacks.

Concretely, for $F \in \dSt_{k}$ we consider the canonical map $F \to
F_{DR}$ whose fibers can be described as follows.

\begin{prop}\label{fibersDR}
  Let $F \in \dSt_{k}$, $\Spec\, A \in \dAff_{k}$, and
  $\overline{u}:\Spec\, A \longrightarrow F_{DR}$, corresponding (by
  Yoneda and the definition of $F_{DR}$) to a morphism $u : \Spec\,
  A^{red} \longrightarrow F$. Then the derived stack
  $F\times_{F_{DR}}\Spec\, A$ is equivalent to the formal completion
  $\widehat{(\Spec\, A \times F)}_{(i,u)}$ of the graph morphism
$$(i,u) : \Spec\, A^{red} \longrightarrow \Spec\, A \times F,$$
where $i : \Spec\, A^{red} \longrightarrow \Spec\, A$ is the natural 
closed embedding.
\end{prop} 
\textbf{Proof.}  Let $X:= \Spec \, A$. By Proposition~\ref{easy}, we
have $(X_{red})_{DR} \simeq X_{DR}$. Therefore the formal completion
$\widehat{(X \times F)}_{(i,u)}$ is in fact the pullback of the
following diagram $$\xymatrix{ \widehat{(X \times F)}_{(i,u)} \ar[rr]
  \ar[d] & & X \times F \ar[d]^-{\lambda_X \times \lambda_F} \\ X_{DR}
  \ar[rr]_-{(\textrm{id}, u_{DR})} & & X_{DR} \times F_{DR}}.$$ But
$(i,u)$ is a graph, so the following diagram is
cartesian $$\xymatrix{\widehat{(X \times F)}_{(i,u)} \ar[r] \ar[d] &
  X_{DR} \ar[rr]^-{u_{DR}} \ar[d] & &F_{DR} \ar[d]^-{\Delta} \\ F
  \times X \ar[r]_-{\lambda_{F} \times \lambda_{X}} & F_{DR \times
    X_{DR}} \ar[rr]_-{id \times u_{DR}} & & F_{DR} \times F_{DR}. }$$
Now recall that, in any $\s$-category with products, a
diagram $$\xymatrix{D \ar[r] \ar[d] & A \ar[d]^-{g} \\ C\ar[r]_-{h} &
  B}$$ is cartesian iff the diagram $$\xymatrix{ D \ar[d] \ar[r] & B
  \ar[d]^-{\Delta} \\ A\times C \ar[r]_{g\times h} & B \times B}$$ is
cartesian, thus, in our case we conclude that $$\xymatrix{\widehat{(X
    \times F)}_{(i,u)} \ar[d] \ar[r] & F \ar[d]^-{\lambda_{F}} \\ X
  \ar[r]_-{u_{DR} \circ \lambda_{X}} & F_{DR} }$$ is cartesian.

% so it is equivalent to $X \times_{X_{DR}} G$, where $G$ is the
% pullback $$\xymatrix{G \ar[r] \ar[d] & F \ar[d]^-{\lambda_{F}} \\
% X_{DR} \ar[r]_-{u_{DR}} & F_{DR}}.$$ Now, since $u_{DR} \circ
% \lambda_X = \overline{u}$, the fiber $X\times_{F_{DR}}F$ of
% $\lambda_{F}$ at $\overline{u} : X \to F_{DR}$ fits into the
% following sequence of cartesian
% squares $$\xymatrix{X\times_{F_{DR}}F \ar[r] \ar[d] & G \ar[r]
% \ar[d] & F \ar[d]^-{\lambda_F} \\ X \ar[r]_-{\lambda_X } & X_{DR}
% \ar[r]_-{u_{DR}} & F_{DR} }$$ and we conclude.

\ \hfill $\Box$

By Proposition~\ref{easy} (4), we get the following corollary of
Proposition~\ref{fibersDR}

\begin{cor}\label{ffdrisgood} If $F$ is algebraic, then each fiber
  $F\times_{F_{DR}}\Spec\, A$ of $F \to F_{DR}$ is a  
  good formal derived stack over $A$, according to Definition~\ref{d10+},
  which is moreover perfect when $F$ is locally of finite
  presentation.
\end{cor} 

Let us remark that in most of our applications $F$ will indeed be
locally of finite presentation (so that its cotangent complex will be
perfect). 

By Proposition~\ref{fibersDR}, the fiber $F\times_{F_{DR}}\Spec\, A$
of $F \to F_{DR}$ when $A=K$ is a field, is simply the formal
completion $\hat{F}_{x}$ of $F$ at the point $x: \Spec\, K
\longrightarrow F$, and corresponds to a dg-Lie algebra over $K$ by
\cite[Theorem~5.3]{lu} or \cite{luX}.  This description tells us that $F
\longrightarrow F_{DR}$ is a family of good formal derived stacks over
$F_{DR}$, and is thus classified by a morphism of derived stacks
$$F_{DR} \longrightarrow \dFSt^{g}_{-},$$
where the right hand side is the $\s$-functor $A \mapsto
\dFSt_{A}^{g}$. We will come back to this point of view in
Section~\ref{GLOB}.

We conclude this section with the following easy but important
observation

\begin{lem}\label{crucial} Let $X$ be a derived Artin stack, and $q: X
  \to X_{DR}$ the associated map. Then $\mathbb{L}_{X}$ and
  $\mathbb{L}_{X/X_{DR}}$ both exist in
  $\mathrm{L}_{\mathrm{QCoh}}(X)$, and we have
$$\mathbb{L}_{X} \simeq \mathbb{L}_{X/X_{DR}}.$$
\end{lem}

\noindent \textbf{Proof.} The cotangent complex $\mathbb{L}_{X}$
exists because $X$ is Artin. The cotangent complex
$\mathbb{L}_{Y_{DR}}$ exists (in the sense of \cite[1.4.1]{hagII}),
for \emph{any} derived stack $Y$, and is indeed trivial. In fact, if
$A$ is a cdga over $k$, and $M$ a dg-module, then $$Y_{DR}(A\oplus M)
\simeq Y((A\oplus M)^{\mathrm{red}}) = Y(A^{\mathrm{red}}) \simeq
Y_{DR}(A).$$ Hence, we may conclude by the transitivity
sequence $$0\simeq q^*\mathbb{L}_{X_{DR}} \to \mathbb{L}_{X} \to
\mathbb{L}_{X/X_{DR}}.$$
\hfill $\Box$ 

\subsection{Perfect complexes on affine formal derived stacks}

For any formal derived stack $F$, we have its $\s$-category of
quasi-coherent complexes $\mathsf{L}_{\mathsf{Qcoh}}(F)$. Recall that
it can described as the following limit (inside the $\s$-category of
$\s$-categories)
$$\mathsf{L}_{\mathsf{Qcoh}}(F):=\lim_{\Spec\, B \longrightarrow
  F}\mathsf{L}(B) \in \scat.$$ 
We can define various full
$\s$-categories of $\mathsf{L}_{\mathsf{Qcoh}}(F)$ by imposing
appropriate finiteness conditions. We will be interested in two of
them, $\mathsf{L}_{\mathsf{Perf}}(F)$ and
$\mathsf{L}_{\mathsf{Qcoh}}^{-}(F)$, respectively of perfect and
cohomologically bounded on the right objects. They are simply defined
as
$$
\mathsf{L}_{\mathsf{Perf}}(F):=\lim_{\Spec\, B \longrightarrow
  F}\mathsf{L}_{\mathsf{Perf}}(B) \qquad
\mathsf{L}_{\mathsf{Qcoh}}^{-}(F):=\lim_{\Spec\, B \longrightarrow
  F}\mathsf{L}_{\mathsf{Qcoh}}^{-}(B).
$$

\begin{df}\label{d10++}
\begin{itemize}
\item Let $\dFSt_{k}^{\mathsf{aff}}$ be the full sub-$\s$-category of
$\dSt_{k}$ consisting of all affine formal derived stacks in the sense
of Definition~\ref{d10}.
\item An affine formal derived stack $F\in \dFSt_{k}^{\mathsf{aff}}$
  is \emph{algebraisable} if there exists $n\in \mathbb{N}$, an
  algebraic derived $n$-stack $F'$, and a morphism $f : F_{red}
  \longrightarrow F'$ such that $F$ is equivalent to the formal
  completion $\hat{F'}_{f}$.
\item A good formal derived stack over $X:=\Spec \, A$
  (Definition~\ref{d10+}) is \emph{algebraisable over $X$} if there
  exists $n\in \mathbb{N}$, an algebraic derived $n$-stack $G
  \longrightarrow \Spec\, A$, locally of finite presentation over
  $\Spec\, A$, together with a morphism $f : \Spec\, A_{red}
  \longrightarrow G$ over $\Spec\, A$, such that $F$ is equivalent, as
  a derived stack over $\Spec\, A$, to the formal completion
  $\widehat{G}_f$.

\end{itemize}
\end{df}

\

\bigskip

In the statement of the next theorem, for $F \in
\dFSt_{k}^{\mathsf{aff}}$, we will denote by $A_{F}$ any $k$-cdga such
that $F_{red} \simeq \Spec \, A_F$: such an $A_F$ exists for any
almost affine derived formal stack $F$, and is unique up to
equivalence.

\noindent The rest of this subsection will be devoted to prove the
following main result

\begin{thm}\label{t1}
There exists an $\s$-functor
$$\mathbb{D} :\dFSt_{k}^{\mathsf{aff}} \longrightarrow (\mecdga_{k})^{op}$$
satisfying the following properties
\begin{enumerate}
\item If $F\in \dFSt_{k}^{\mathsf{aff}}$ is algebraisable, then we
  have an equivalence of (non-mixed) graded cdga
$$\mathbb{D}(F) \simeq Sym_{A_{F}}(\mathbb{L}_{F_{red}/F}[-1]),$$
\item For all $F \in \dFSt_{k}^{\mathsf{aff}}$, there exists an
  $\s$-functor
$$\phi_{F}: \mathsf{L}_{\mathsf{Qcoh}}(F) \longrightarrow
\mathbb{D}(F)-Mod^{\mathsf{gr}}_{\edg},$$ 
natural in $F$, which is conservative, and 
induces an equivalence of $\s$-categories
$$\mathsf{L}_{\mathsf{Perf}}(F) \longrightarrow
\mathbb{D}(F)-Mod_{\edg}^{\mathsf{gr}, \mathsf{perf}},$$ 
where the right hand side is the full sub-$\s$-category 
of  $\mathbb{D}(F)-Mod^{\mathsf{gr}}_{\edg}$ consisting of graded
mixed $\mathbb{D}(F)$-modules 
$E$ which are equivalent, as graded $\mathbb{D}(F)$-modules, to 
$\mathbb{D}(F)\otimes_{A_{F}}E_{0}$ for some $E_{0} \in
\mathsf{L}_{\mathsf{Perf}}(A_{F})$. 
\end{enumerate}
\end{thm}

We will first prove Thm \ref{t1} for $F$ a derived affine scheme, and
then proceed to the general case. \\

\noindent \textbf{Proof of Theorem \ref{t1}: the derived affine case.}
We start with the special case of the theorem for the
sub-$\s$-category $\dAff_{k} \subset \dFSt_{k}^{\mathsf{aff}}$ of
derived affine schemes (recall our convention that all derived affine
schemes are locally finitely presented), and construct the
$\s$-functor
$$\mathbb{D} : \dAff_{k}^{op} \longrightarrow \mecdga_k$$ 
as follows.
We start by sending an object $\Spec\, A \in \dAff_{k}$ to the
morphism $A \longrightarrow A^{red}$. This defines an $\s$-functor
$\dAff_{k} \longrightarrow \mathsf{Mor}(\dAff_{k})$, from derived
affine schemes to morphisms between derived affine schemes. We then
compose this with the $\s$-functor (see end of \S \ref{1.3.2}, with
$\C=\dg_{k}$)
$$\DR : \mathsf{Mor}(\cdga_{k}) \longrightarrow \mecdga_{k},$$
sending a morphism $A \rightarrow B$ to $\DR(B/A)$. Recall that this
second $\s$-functor can be explicitly constructed as the localization
along equivalences of the functor
$$DR^{str} : \mathsf{Cof}(cdga_{k}) \longrightarrow \epsilon-cdga_{k}^{gr},$$
from the category of cofibrations between cofibrant cdgas to the
category of graded mixed cdgas, sending a cofibration $A \rightarrow
B$ to $DR^{str}(B/A)=Sym_{B}(\Omega_{B/A}^{1}[-1])$, with mixed
structure given by the de Rham differential.

\begin{prop}\label{p2}
  The $\s$-functor defined above
$$\mathbb{D} : \dAff_{k}^{op} \longrightarrow \mecdga_{k}\, :
\,\,\,\,\, A \longmapsto \DR(A^{red}/A)$$ 
is fully faithful. Its essential image is contained inside the full
sub-$\s$-category of  
graded mixed cdgas $B$ satisfying the following three conditions.
\begin{enumerate}
\item The cdga $B(0)$ is concentrated in cohomological degree $0$, and
  is a reduced $k$-algebra of finite type.
\item The $B(0)$-dg-module $B(1)$ is almost finitely presented and has
  amplitude contained in $(-\infty,0]$.
\item The natural morphism
$$Sym_{B(0)}(B(1)) \longrightarrow B$$
is an equivalence of graded cdgas.
\end{enumerate}
\end{prop}
\noindent \textbf{Proof.} For $\Spec\, A \in \dAff_{k}$, we have
$$\mathbb{D}(A)=\DR(A_{red}/A)\simeq
Sym_{A_{red}}(\mathbb{L}_{A_{red}/A}[-1]),$$ 
showing that conditions $1$, $2$, and $3$ above are indeed satisfied
for $\mathbb{D}(A)$ (for $2$, recall that $A \rightarrow A_{red}$
being an epimorphism, we have $\pi_{0}(\mathbb{L}_{A_{red}/A})=0$).
The fact that $\mathbb{D}$ is fully faithful is essentially the
content of \cite{bh}, stating that the relative derived de Rham
cohomology of any closed immersion is the corresponding formal
completion. Indeed, here $X=\Spec\, A$ is the formal completion of
$X_{red}=(\Spec\, A)_{red}$ inside $X$. For the sake of completeness,
we will provide here a new proof of this fact, for the specific closed
immersion $X_{red} \longrightarrow X$.

Let $\Spec\, A$ and $\Spec\, B$ be two derived affine schemes, and
consider the induced morphism of mapping spaces
$$\mathsf{Map}_{\dSt_{k}}(\Spec\, A, \Spec\, B) \simeq 
\mathsf{Map}_{\cdga_{k}}(B,A) \longrightarrow 
\mathsf{Map}_{\mecdga_{k}}(\mathbb{D}(B),\mathbb{D}(A)).$$
By Lemma \ref{relativeDR}, we have 
$$\mathsf{Map}_{\mecdga_{k}}(\mathbb{D}(B),\mathbb{D}(A)) \simeq
\mathsf{Map}_{\cdga_{k}}(B_{red}, A_{red}) \times
_{\mathsf{Map}_{\cdga_{k}}(B, A_{red})} \mathsf{Map}_{\mecdga_{k}}(B,
\mathbb{D}(A))$$ 
where $B$ is considered as a graded mixed 
cdga in a trivial manner (pure of weight $0$ and with zero mixed
structure). But the canonical map $\mathsf{Map}_{\cdga_{k}}(B_{red},
A_{red}) \to\mathsf{Map}_{\cdga_{k}}(B, A_{red})$ is an equivalence,
hence 
$$\mathsf{Map}_{\mecdga_{k}}(\mathbb{D}(B),\mathbb{D}(A)) \simeq
\mathsf{Map}_{\mecdga_{k}}(B,\mathbb{D}(A)).$$ 
Finally, by adjunction we have
$$\mathsf{Map}_{\mecdga_{k}}(B,\mathbb{D}(A)) \simeq
\mathsf{Map}_{\mecdga_{k}}(k(0)\otimes_{k} B,\mathbb{D}(A)) \simeq
\mathsf{Map}_{\cdga_{k}}(B,|\mathbb{D}(A)|)$$ 
where $|-| : \mecdga_{k}
\longrightarrow \cdga_{k}$ is the realization $\s$-functor of
Definition~\ref{dreal} for commutative monoids in $\mathcal{M}=
\medg_k$. Note that the commutative $k-$dg-algebra $|\mathbb{D}(A)|$
is exactly the derived de Rham cohomology of $A_{red}$ over $A$. By
putting these remarks together, we conclude that, in order to prove
that $\mathbb{D}$ is fully faithful, it will be enough to show that,
for any $A \in \cdga_{k}$, the induced natural morphism $A
\longrightarrow |\mathbb{D}(A)|$ is an equivalence,
i.e. the statement is reduced to the following

\begin{lem}\label{l4}
  For any $\Spec\, B \in \dAff_{k}$ the natural morphism $B
  \longrightarrow \mathbb{D}(B)$ of graded mixed cdgas induces an
  equivalence in $\cdga_{k}$
$$B \longrightarrow |\mathbb{D}(B)|.$$
\end{lem}

%%%%%%%%%%%%%%%%%%%%%%%%%%   I L     Q U I    R  I G O R O S
%%%%%%%%%%%%%%%%%%%%%%%%%%   O    %%%%%%%%%%%%%%%%%%%%%%%%%%%%% 

\noindent \textit{Proof of Lemma.} We can assume that $B$ is a cell
non-positively graded commutative dg-algebra with finitely many cells
in each dimension. As a commutative graded algebra $B$ is a free
commutative graded algebra with a finite number of generators in each
degree. In particular $B^0$ is a polynomial $k$-algebra and $B^i$ is a
free $B^0$-module of finite rank for all $i$. In the same way, we
chose a cofibration $B \hookrightarrow C$ which is a model for $B
\longrightarrow B_{red}$.  We chose moreover $C$ to be a cell $B$-cdga
with finitely many cells in each dimension. As $B_{red}$ is quotient
of $\pi_0(B)$ we can also chose $C$ with no cells in degree $0$.

We let $L:=\Omega_{C/B}^1[-1]$, which is a cell $C$-dg-module with
finitely many cells in each degree, and no cells in positive
degrees. The commutative dg-algebra $|\mathbb{D}(B)|$ is by definition
the completed symmetric cdga $\widehat{Sym}_{C}(L)$, with its total
differential, sum of the cohomological and the de Rham differential.
Note that, because $L$ has no cells in positive degrees and only
finitely many cells in each degree, the cdga $|\mathbb{D}(B)|$ is
again non-positively graded. Note however that it is not clear a
priori that $|\mathbb{D}(B)|$ is almost of finite presentation and
thus not clear that $Spec\, |\mathbb{D}(B)| \in \dAff_{k}$.

We let $C^0$ be the commutative $k$-algebra of degree zero elements in
$C$, and $L^0$ of degree zero elements in $L$. We have a natural
commutative square of commutative dg-algebras, relating completed and
non-completed symmetric algebras
$$
\xymatrix{
Sym_{C}(L) \ar[r] & \widehat{Sym}_{C}(L) \\
Sym_{C^0}(L^0) \ar[u] \ar[r] & \widehat{Sym}_{C^0}(L^0). \ar[u]}
$$ 
In this diagram we consider $Sym_{C}(L)$ and $\widehat{Sym}_{C}(L)$
both equipped with the total differential, sum of the cohomological and
the de Rham differential (recall that $L=\Omega^1_{C/B}[-1]$).

By assumption $C^0$ is a polynomial $k$-algebra over a finite number
of variables, and $C^i$ is a free $C^0$-module of finite type. This
implies that the diagram above is a push-out of commutative
dg-algebras, and, as the lower horizontal arrow is a flat morphism of
commutative rings, this diagram is moreover a homotopy push-out of
cdgas. We thus have a corresponding push-out diagram of the
corresponding cotangent complexes, which base changed to $C$ provides
a homotopy push-out of $C$-dg-modules
$$\xymatrix{
\mathbb{L}_{Sym_{C}(L)}\otimes_{Sym_{C}(L)}C \ar[r] & 
\mathbb{L}_{\widehat{Sym}_{C}(L)}\otimes_{\widehat{Sym}_{C}(L)}C \\
\mathbb{L}_{Sym_{C^0}(L^0)}
\otimes_{Sym_{C^0}(L^0)}C \ar[u] \ar[r] & 
\mathbb{L}_{\widehat{Sym}_{C^0}(L^0)}\otimes_{\widehat{Sym}_{C^0}(L^0)}C. \ar[u]}
$$ 
As $C^0$ is a polynomial algebra over $k$, the lower horizontal
morphism is equivalent to  
$$\Omega^1_{Sym_{C^0}(L^0)} 
\otimes_{Sym_{C^0}(L^0)}C \to 
\Omega^1_{\widehat{Sym}_{C^0}(L^0)}\otimes_{\widehat{Sym}_{C^0}(L^0)}C,
$$
which is the base change along $C^0 \longrightarrow C$ of the morphism
$$\Omega^1_{Sym_{C^0}(L^0)}
\otimes_{Sym_{C^0}(L^0)}C^0 \to 
\Omega^1_{\widehat{Sym}_{C^0}(L^0)}\otimes_{\widehat{Sym}_{C^0}(L^0)}C^0.
$$
This last morphism is an isomorphism, and thus the induced morphism
$$
\mathbb{L}_{Sym_{C}(L)}\otimes_{Sym_{C}(L)}C \longrightarrow
\mathbb{L}_{\widehat{Sym}_{C}(L)}\otimes_{\widehat{Sym}_{C}(L)}C
$$
is an equivalence of $C$-dg-modules. To put things differently, 
the morphism of cdgas
$Sym_{C}(L) \longrightarrow \widehat{Sym}_{C}(L)$ is 
formally \'etale along the augmentation. 

We deduce from this the existence of a canonical identification of
$C$-dg-modules 
$$\mathbb{L}_{B_{red}/B} \simeq 
\mathbb{L}_{|\mathbb{D}(B)|}\otimes_{|\mathbb{D}(B)|}B_{red}.$$
This equivalence is moreover induced by the diagram of cdgas
$$\xymatrix{
  B \ar[rr] \ar[rd] & & |\mathbb{D}(B)| \ar[ld] \\
  & B_{red}. &}$$ Equivalently, the morphism $B \longrightarrow
|\mathbb{D}(B)|$ is formally \'etale at the augmentation over
$B_{red}$. By the infinitesimal lifting property, the morphism of
$B$-cdgas $|\mathbb{D}(B)| \longrightarrow B_{red}$ can be extended
uniquely to a morphism $|\mathbb{D}(B)| \longrightarrow
\pi_0(B)$. Similarly, using the Postnikov tower of $B$, this morphism
extends uniquely to a morphism of $B$-cdgas $|\mathbb{D}(B)|
\longrightarrow B$. In other words, the adjunction morphism $i : B
\longrightarrow |\mathbb{D}(B)|$ possesses a retraction up to homotopy
$r : |\mathbb{D}(B)| \longrightarrow B$.  We have $ri\simeq id$, and
$\phi:=ir$ is an endomorphism of $|\mathbb{D}(B)|$ as a $B$-cdga,
which preserves the augmentation $|\mathbb{D}(B)| \longrightarrow
B_{red}$ and is formally \'etale at $B_{red}$.

By construction, $|\mathbb{D}(B)|\simeq \lim_{n}|\mathbb{D}_{\leq
  n}(B)|$, where 
$$|\mathbb{D}_{\leq n}(B)|:=Sym^{\leq
  n}_{B_{red}}(\mathbb{L}_{B_{red}/B}[-1])$$ 
is the truncated de Rham
complex of $B_{red}$ over $B$. Each of the cdga $|\mathbb{D}_{\leq
  n}(B)|$ is such that $\pi_0(|\mathbb{D}_{\leq n}(B)|)$ is a finite
nilpotent thickening of $B_{red}$, and moreover
$\pi_i(|\mathbb{D}_{\leq n}(B)|)$ is a $\pi_0(|\mathbb{D}_{\leq
  n}(B)|)$-module of finite type. Again by the infinitesimal lifting
property we see that these imply that the endomorphism $\phi$ must
be homotopic to the identity.

This finishes the proof that the adjunction morphism $B
\longrightarrow |\mathbb{D}(B)|$ is an equivalence of cdgas, and thus
the proof Lemma \ref{l4}.
\hfill $\Box$ 

\

The lemma is proved, and thus Proposition \ref{p2} is proved as well.
\hfill $\Box$

\

\medskip

\noindent
One important consequence of Proposition \ref{p2} is the following corollary, 
showing that quasi-coherent complexes over $\Spec\, A \in \dAff_{k}$ can be
naturally identified with certain $\mathbb{D}(A)$-modules.

\begin{cor}\label{c2}
Let $\Spec\, A \in \dAff_{k}$ be an affine derived scheme, and
$\mathbb{D}(A):=\DR(A_{red}/A)$ be the corresponding 
graded mixed cdga. There exists a 
symmetric monoidal stable $\s$-functor
$$\phi_{A} : \mathrm{L}_{\mathrm{QCoh}}(A) \hookrightarrow
\mathbb{D}(A)-Mod_{\edg},$$ 
functorial in $A$, inducing an equivalence
of $\s$-categories
$$\mathrm{L}_{\mathsf{Perf}}(A) \simeq
\mathbb{D}(A)-Mod_{\edg}^{\mathsf{Perf}},$$ 
where $\mathbb{D}(A)-Mod_{\edg}^{\mathsf{Perf}}$ is the full
sub-$\s$-category consisting of   
mixed graded $\mathbb{D}(A)$-modules $M$ for which there exists $E \in
\mathrm{L}_{\mathsf{Perf}}(A_{red})$, and an equivalence 
of (non-mixed) graded modules 
$$M \simeq \mathbb{D}(A) \otimes_{A_{red}}E.$$ 
\end{cor}

\noindent \textbf{Proof.} The $\s$-functor $\phi_{A}$ is defined by sending
an $A$-dg-module $E\in \mathrm{L}(A)$ to  
$$\phi_{A}(E):=\mathbb{D}(A) \otimes_{A}E \in
\mathbb{D}(A)-Mod_{\epsilon-\mathrm{L}(k)^{gr}},
$$ 
using that $\mathbb{D}(A)=\DR(A_{red}/A)$ is, naturally, an $A$-linear
graded mixed cdga. This $\s$-functor sends $A$ to $\mathbb{D}(A)$
itself. In particular, we have 
$$
\begin{aligned}
\mathsf{Map}_{\mathbb{D}(A)-Mod_{\epsilon-\mathrm{L}(k)^{gr}}}(
\mathbb{D}(A),\mathbb{D}(A)) 
& \simeq 
\mathsf{Map}_{\epsilon-\mathrm{L}(k)^{gr}}(k(0),\mathbb{D}(A)) \\
& \simeq \mathsf{Map}_{\mathrm{L}(k)}(k,|\mathbb{D}(A)|) \\
& \simeq \mathsf{Map}_{A-Mod}(A,A).
\end{aligned}
$$
This shows that $\phi_{A}$ is fully faithful on the single object $A$, 
so, by stability, it is also fully faithful when restricted to
$\mathrm{L}_{\mathsf{Perf}}(A)$, the 
$\s$-category of perfect $A$-dg-modules. 
\hfill $\Box$ 

\begin{rmk}
We believe that the $\s$-functor $\phi_{A}=\mathbb{D}(A)\otimes_A-$ coincides with 
the $\s$-functor $\mathrm{QCoh}(A)\hookrightarrow\mathrm{IndCoh}(A)$. Let us provide 
a sketchy evidence in the case when $A_{red}=k$, that is to say when $A$ is an Artinian cdga. 
First of all let us observe that the graded mixed cdga $\mathbb{D}(A)$ is the Chevalley-Eilenberg 
graded mixed cdga of the Lie $A$-algebra $\mathrm{Der}_A(k,k)$. We thus get an equivalence of 
$\s$-categories between $\mathbb{D}(A)-Mod_{\edg}$ and $\mathrm{Der}_A(k,k)-Mod$. Under this equivalence 
$\phi_{A}$ simply becomes $k\otimes_A-$. 
Then recall (see \cite{lu} and \cite[\S2.2]{toenbourbaki}) that the Lie algebra $\mathrm{Der}_A(k,k)$ 
is equivalent to the Lie algebra $\ell_{\mathrm{Spf}(A)}$ associated with the formal moduli problem $\mathrm{Spf}(A)$. 
Finally, it is known after \cite{lu} that we have an equivalence of $\s$-categories 
$\ell_{\mathrm{Spf}(A)}-Mod\simeq \mathrm{IndCoh}(A)$. 
\end{rmk}

Proposition~\ref{p2} and Corollary~\ref{c2} together prove
Theorem \ref{t1} in the derived affine case. 
We now move to the general case.

\

\noindent \textbf{Proof of Theorem \ref{t1} : the general case.} We
will extend the above relation between derived affine schemes and
graded mixed cdgas to the case of affine formal derived stacks. In
order to do this, we start with the $\s$-functor
$$\dAff_{k}^{opp} \longrightarrow \mecdga_{k}$$
sending $A$ to $\mathbb{D}(A)=\DR(A_{red}/A)$. This $\s$-functor is 
a derived stack for the \'etale topology on $\dAff_{k}^{op}$, and thus
has a right Kan extension as an 
$\s$-functor defined on all
derived stacks
$$\mathbb{D} : \dSt_{k}^{op} \longrightarrow \mecdga_{k}, \qquad 
F \longmapsto \lim_{\Spec \, A \to F}(\mathbb{D}(A))
$$
(with the limit being taken in $\mecdga_{k}$), and sending colimits in
$\dSt_{k}$ to limits.  In general, there are no reasons to expect that
$\mathbb{D}(F)$ is free as a graded cdga, and it is a remarkable
property that this is the case when $F$ is an \emph{algebraisable}
affine formal derived stack (Definition~\ref{d10++}); we do not know if the
result still holds for a general affine formal derived stack. The
following Proposition establishes this, and thus point 1 of
Theorem~\ref{t1}.

\begin{prop}\label{p3}
  Let $F\in \dFSt_{k}^{\mathrm{aff}}$ be an algebraisable affine
  formal derived stack, and let $F_{red}\simeq \Spec\, A_{0}$. We have
  a natural equivalence of (non-mixed) graded cdga's
$$Sym_{A_{0}}(\mathbb{L}_{F_{red}/F}[-1]) \simeq \mathbb{D}(F).$$
\end{prop}
\noindent \textbf{Proof.} For all $\Spec\, A \longrightarrow F$, we have 
a commutative square 
$$
\xymatrix{
\Spec\, A_{red} \ar[r] \ar[d] & \Spec\, A \ar[d] \\
F_{red}=\Spec\, A_0 \ar[r] & F}
$$
and, therefore, an induced  
a natural morphism of $A_{0}$-dg-modules
$$
\mathbb{L}_{F_{red}/F} \longrightarrow \mathbb{L}_{A_{red}/A}.
$$
This yields a morphism of (non mixed) graded cdgas
$$
Sym_{A_{0}}(\mathbb{L}_{F_{red}/F}[-1]) \longrightarrow
\mathbb{D}(A).
$$
Taking the limit over $(\Spec\, A \to F) \in \dAff/F$, we obtain a
natural morphism of (non mixed) graded  
cdgas
$$ 
\phi_{F} : Sym_{A_{0}}(\mathbb{L}_{F_{red}/F}[-1]) \longrightarrow
\mathbb{D}(F) = \lim_{\Spec\, A \rightarrow F}\mathbb{D}(A).
$$
Since $F$ is algebraisable (Definition~\ref{d10++}), there exists an
algebraic derived $n$-stack (for some integer $n$) $G$, a morphism $f
: \Spec\, A \longrightarrow G$ and an equivalence $\widehat{G}_{f}
\simeq F$. We will prove that $\phi_{F}$ is an equivalence by
induction on $n$.

We first observe that the statement is local in the \'etale topology of 
$\Spec\, A_{0}$ in the following sense. Let $A_{0} \longrightarrow A_{0}'$ 
be an \'etale morphism and $X'=\Spec\, A_{0}' \longrightarrow X=\Spec\, A_0$
be the induced morphism. We let $F'$ be the formal completion of the
morphism $X' \longrightarrow F$ (or equivalently of $X' \longrightarrow G$) 
so that we have a commutative square
of derived stacks
$$
\xymatrix{
X' \ar[r] \ar[d] & F' \ar[d] \\
X \ar[r] & F.}
$$
By construction this square is moreover cartesian, and induces
a morphism of graded cdgas
$$
\mathbb{D}(F) \longrightarrow \mathbb{D}(F').
$$
Thus the assignment $X' \mapsto \mathbb{D}(F')$ defines a stack of
graded cdgas over the small \'etale site of $X$, and, in the same way,
$X' \mapsto Sym_{A'_{0}}(\mathbb{L}_{X'/F'}[-1])$ is a stack of graded
cdgas on the small \'etale site of $X$. The various morphism
$\phi_{F'}$ organize into a morphism of \'etale stacks on $X$. In
order to prove that $\phi_{F}$ is an equivalence it is enough to prove
that $\phi_{F'}$ is so after some \'etale covering $X' \longrightarrow
X$.

The above \'etale locality of the statement implies that we can assume
that there is an affine $Y \in \dAff$, a smooth morphism $Y
\longrightarrow G$, such that $X\longrightarrow G$ comes equipped with
a factorization through $Y$
$$\xymatrix{
X \ar[r] \ar[d] & G \\
Y \ar[ru] & }
$$
We let $Y_{*}$ be the nerve of the morphism $Y \longrightarrow G$,
which is a smooth Segal groupoid in derived stacks (see \cite[\S
S.3.4]{hagII}). Moreover, $Y_0=Y$ is affine and $Y_i$ is an algebraic
$(n-1)$-stack. We consider the chosen lifting $X \longrightarrow Y_0$
as a morphism of simplicial objects $X \longrightarrow Y_*$, where $X$
is considered as simplicially constant. We let $\widehat{Y}_*$ be the
formal completion of $Y_*$ along $X$, defined by
$$
\widehat{Y}_{i}:=\widehat{(Y_{i})}_{X \rightarrow Y_{i}}.
$$
The simplicial object $\widehat{Y}_{*}$ can be canonically identified
with the nerve of the induced morphism on formal completions
$\widehat{Y_0} \longrightarrow F=\widehat{G}$. Moreover, by
construction $\widehat{Y_0} \longrightarrow F$ is an epimorphism of
derived stacks, and we thus have a natural equivalence of derived
stacks
$$ 
|\widehat{Y}_*| = \mathrm{colim}_{i} \widehat{Y_i} \simeq F.
$$
As the $\s$-functor $\mathbb{D}$ sends colimits to limits we have
$$
\mathbb{D}(F) \simeq \lim_{i} \mathbb{D}(\widehat{Y_i}).
$$
Also, for each $i$ the morphism $\widehat{Y_i} \longrightarrow Y_i$
is formally \'etale, and thus we have
$$
\mathbb{L}_{X/\widehat{Y}_{i}} \simeq \mathbb{L}_{X/Y_{i}}.
$$
Smooth descent for differential forms on $G$ (see Appendix B) then
implies that we have 
equivalences of $A_0$-dg-modules
$$
\wedge^{p} \mathbb{L}_{X/F} \simeq \wedge^{p} \mathbb{L}_{X/G} \simeq 
\lim_{i} \wedge^{p} \mathbb{L}_{X/Y_{i}} \simeq
\lim_{i} \wedge^{p} \mathbb{L}_{X/\widehat{Y}_{i}}.
$$
Therefore 
$$
Sym_{A_{0}}(\mathbb{L}_{X/F}[-1]) \simeq 
\lim_{i} Sym_{A_{0}}(\mathbb{L}_{X/\widehat{Y}_{i}}[-1]).
$$
The upshot is that, in order to prove that $\phi_{F}$ is an
equivalence, it is enough to prove that all the $\phi_{Y_{i}}$'s are
equivalences. By descending induction on $n$ this allows us to reduce
to the case where $G$ is a derived algebraic stack, and by further
localization on $G$ to the case where $G$ is itself a derived affine
scheme. Moreover, by refining the smooth atlas $Y \rightarrow G$ in
the argument above, we may also assume that $X \longrightarrow G$ is a
closed immersion of derived affine schemes.

Suppose $G=Z \in \dAff_{k}$, and $X \longrightarrow Z$ be a closed
immersion; recall that this means that the induced morphism on
truncations $t_0(X)=X \longrightarrow t_{0}(Z)$ is a closed immersion
of affine schemes. We may present $X \longrightarrow Z$ by a cofibrant
morphism between cofibrant cdgas $B \longrightarrow A$, and moreover
we may assume that $A$ is a cell $B$-algebra with finitely many cells
in each degree, and that $B$ is a cell $k$-algebra with finitely many
cells in each degree. We let $B^0$ be the $k$-algebra of degree zero
elements in $B$ and $Z^0=\Spec\, B^0$. The formal completion
$\widehat{Z}=F$ of $X \longrightarrow Z$ sits in a cartesian square of
derived stacks
$$
\xymatrix{
\widehat{Z} \ar[r] \ar[d] & \widehat{Z^0} \ar[d] \\
Z \ar[r] & Z^{0},}
$$
where $Z \longrightarrow Z^{0}$ is the natural morphism induced by
$B^0 \subset B$, and $\widehat{Z^0}$ is the formal completion of $Z^0$
along the closed immersion corresponding to the quotient of algebras
$$
B^0 \longrightarrow \pi_{0}(B) \longrightarrow \pi_{0}(A)\simeq A_0.
$$
We let $I \subset B^0$ be the kernel of $B^0 \longrightarrow A_0$,
and we choose generators $f_1,\cdots,f_p$ for $I$. 
We set $B^0(j):=K(B^0,f_1^j,\dots,f_p^j)$ the Koszul cdga over $B^0$ 
attached to the sequence $(f_1,\cdots,f_p)$, 
$Z^0_j:=\Spec\, B^0(j)$ and $Z_j:=Z \times_{Z^0}Z^0_{j}$.
We have a natural equivalence of derived stacks
$$
F=\widehat{Z} \simeq \mathrm{colim}_{j} Z_{j}.
$$
By our Appendix B we moreover know that $\widehat{Z^0}$ is equivalent
to $\mathrm{colim}_{j} Z_{j}^{0}$ as derived prestacks, or in other
words, that the above colimit of prestacks is a derived stack.  By
pull-back, we see that the colimit $\mathrm{colim}_{j} Z_{j}$ can be
also computed in derived prestacks, and thus the equivalence
$\widehat{Z} \simeq \mathrm{colim}_{j} Z_{j}$ is an equivalence of
derived prestacks (i.e.  of $\s$-functors defined on $\dAff_{k}$).  As
$\mathbb{D}$ sends colimits to limits, we do have an equivalence of
graded mixed cdgas
$$
\mathbb{D}(F) \simeq \lim_{n} \mathbb{D}(Z_j).
$$
The proposition follows by observing that, for any $p\geq 0$, the
natural morphism 
$$
\wedge^{p}\mathbb{L}_{X/Z} \longrightarrow
\lim_{n} \wedge^{p}\mathbb{L}_{X/Z_j}
$$
is indeed an equivalence of dg-modules over $A_0$ (see Appendix
B). \hfill $\Box$  

\

\medskip

As a consequence of Proposition \ref{p3}, if $F$ is an algebraisable
affine formal derived stack, and if $\mathbb{L}_{F}$ is of amplitude
contained in $]-\infty,n]$ for some $n$, then the graded mixed cdga
$\mathbb{D}(F)$ satisfies the following conditions.

\begin{enumerate}
\item The cdga $A:=\mathbb{D}(F)(0)$ is concentrated in degree 
$0$ and is a reduced $k$-algebra of finite type.
\item The $A$-dg-module $\mathbb{D}(F)(1)$ is almost finitely
  presented and of amplitude contained in $]-\infty,n]$.
\item The natural morphism
$$Sym_{A}(\mathbb{D}(F)(1)) \longrightarrow \mathbb{D}(F)$$
is an equivalence of graded cdgas.
\end{enumerate}

We now move to the proof of point 2 in Theorem \ref{t1}, i.e. we
define the $\s$-functor 
$$
\phi_F : \mathrm{L}_{\mathrm{QCoh}}(F)
\longrightarrow \mathbb{D}(F)-Mod_{\edg}
$$ 
for a general $F \in
\dFSt_{k}^{\mathsf{aff}}$.  This was already defined when $F$ is an
affine derived stack in Corollary \ref{c2}, and for general $F$ the
$\s$-functor $\phi_{F}$ will be simply defined by left Kan
extension. More precisely, if $F \in \dSt_{k}$, we start with
$$
\lim_{\Spec\, A \rightarrow F}\phi_{A} : \lim_{\Spec\, A \rightarrow F} 
\mathrm{L}(A) \longrightarrow \lim_{\Spec\, A \rightarrow F}
\mathbb{D}(A)-Mod_{\edg},
$$
where for each fixed $A$ the $\s$-functor $\phi_{A} : \mathrm{L}(A)
\longrightarrow \mathbb{D}(A)-Mod_{\edg}$, is the one of our corollary
\ref{c2} and sends an $A$-dg-module $E$ to
$E\otimes_{A}\DR(A_{red}/A)$. Finally, as $\mathbb{D}(F)=\lim_{\Spec\,
  A \rightarrow F} \mathbb{D}(A)$ there is a natural limit
$\s$-functor
$$
lim : \lim_{\Spec\, A \rightarrow F} \mathbb{D}(A)-Mod_{\edg}
\longrightarrow \mathbb{D}(F)-Mod_{\edg}. 
$$
By composing these two functors, we obtain a natural $\s$-functor
$$
\phi_{F} : \mathrm{L}_{\mathrm{QCoh}}(F) \longrightarrow
\mathbb{D}(F)-Mod_{\edg},
$$ 
which is clearly functorial in $F \in \dSt_{k}$. Note that $\phi_{F}$
exists for any $F$, without any extra conditions. The fact that it
induces an equivalence on perfect modules only requires $F_{red}$ to
be an affine scheme, as shown in Proposition~\ref{p4} below. This
establishes, in particular, point 2 of Theorem \ref{t1}, and thus
concludes its proof.

If $B \in \mecdga$ is graded mixed cdga, a graded mixed $B$-dg-module
$M\in B-Mod_{\edg}$ is called \emph{perfect}, if, as a graded
$B$-dg-module, it is (equivalent to a graded $B$-dg-module) of the
form $B\otimes_{B(0)}E$ for $E \in
\mathrm{L}_{\mathsf{Perf}}(B(0))$. Note that $E$ is then automatically
equivalent to $M(0)$. In other words, $M$ is perfect if it is free
over its degree $0$ part, as a graded $B$-dg-module. We let
$B-Mod_{\edg}^{\mathsf{Perf}}$ be the full sub-$\s$-category of
$B-Mod_{\edg}$ consisting of perfect graded mixed $B$-dg-modules.

\begin{prop}\label{p4}
  Let $F \in \dSt_{k}$, and assume that $F_{red}=\Spec\, A_{0}$ is an
  affine reduced scheme of finite type over $k$. Then, the
  $\s$-functor 
$$
\phi_{F} : \mathrm{L}_{\mathsf{Perf}}(F) \longrightarrow
\mathbb{D}(F)-Mod_{\edg}^{\mathsf{Perf}}
$$
is an equivalence of $\s$-categories.
\end{prop}
\noindent \textbf{Proof.} By Corollary \ref{c2}, we have a natural
equivalence of $\s$-categories
$$
\mathrm{L}_{\mathsf{Perf}}(F) \simeq \lim_{\Spec\, A \rightarrow F}
\mathrm{L}_{\mathsf{Perf}}(A)  
\simeq \lim_{\Spec\, A \rightarrow
  F}\mathbb{D}(A)-Mod_{\edg}^{\mathsf{Perf}}.
$$
As $\mathbb{D}(F)=\lim_{\Spec\, A \rightarrow F}\mathbb{D}(A)$, 
we have a natural adjunction of $\s$-categories
$$\mathbb{D}(F)-Mod_{\edg} \longleftrightarrow 
\lim_{\Spec\, A \rightarrow F}\mathbb{D}(A)-Mod_{\edg},$$ where the
right adjoint is the limit $\s$-functor. Thus this adjunction induces
an equivalences on perfect objects.  \hfill $\Box$

\subsection{Differential forms and polyvectors on perfect formal
  derived stacks}\label{subssfpspfds} 

In the previous section, we have associated to any formal affine
derived stack $F$, a mixed graded cdga $\mathbb{D}(F)$ in such a way
that $\mathrm{L}_{\mathsf{Perf}}(F) \simeq
\mathbb{D}(F)-Mod_{\edg}^{\mathsf{Perf}}$.  We will now compare the de
Rham theories of $F$ (in the sense of \cite{ptvv}) and of
$\mathbb{D}(F)$ (in the sense of \S 1.3), and prove that they are
equivalent when appropriately understood.

\subsubsection{De Rham complex of perfect formal derived
  stacks} \label{drpform} We let $F \longrightarrow \Spec\, A$ be a
perfect formal derived stack (Definition~\ref{d10+}) and $\D(F)$ the
corresponding graded mixed cdga of Theorem~\ref{t1}. The projection $F
\longrightarrow \Spec\, A$ induces a morphism of graded mixed cdga
$\D(A) \longrightarrow \D(F)$ which allows us to view $\D(F)$ as a graded
mixed $\D(A)$-algebra. By taking $\C=\edg^{gr}_{k}$ in
Proposition~\ref{DREL}, we may consider in particular its relative de Rham
object $\DR^{int}(\D(F)/\D(A))$ which is a graded mixed cdga over the
$\s$Å-category of graded mixed $\D(A)$-dg-modules. There is an
equivalence
$$
\DR^{int}(\D(F)/\D(A)) \simeq
Sym_{\D(F)}(\mathbb{L}_{\D(F)/\D(A)}^{int}[-1])
$$ 
of (non mixed) graded cdgas over the $\s$Å-category of graded mixed
$\D(A)$-dg-modules.
We can consider its realization, as in Definition~\ref{dreal},
$$
\DR(\D(F)/\D(A)):=|\DR^{int}(\D(F)/\D(A))| 
$$
which is thus a graded mixed cdga over $|\D(A)|\simeq A$
(Remark~\ref{enhanced}, and Lemma \ref{l4}), i.e.  
an $A$-linear graded mixed cdga. Moreover, according to \S 1.5, we can
also consider its Tate realization 
$$
\DR^{t}(\D(F)/\D(A)):=|\DR^{int}(\D(F)/\D(A))|^{t},
$$
which is, again, a graded mixed $A$-linear cdga. 

On the other hand, it is natural to consider the following

\begin{df}\label{reldrstacks} The \emph{de Rham object of an arbitrary 
derived stack $F$ over an affine derived stack $\Spec\, B$} is 
$$\DR(F/B):=\lim_{\Spec\, C \rightarrow F} \DR(C/B) \in \mecdga_B$$
where the limit is taken in the category $\mecdga_B$ of graded mixed
$B$-linear cdgas, and over all morphisms $\Spec\, C \rightarrow F$ of
derived stacks over $\Spec \, B$.
\end{df}

\begin{prop}\label{trans1}
Let $f: F \to G$ be a map in $\dSt/ \Spec \, B$. There is an induced map $\DR(G/B) \to \DR(F/B)$ in  $\mecdga_B$.
\end{prop}

\noindent \textbf{Proof.}
Let $I_{F}$ (resp. $I_G$) the category on which $\DR(F/B)$ (resp. $\DR(G/B)$) is defined as a limit: $\DR(F/B)= lim_{I_X} \DR^{F}$, where $$\DR^F: I_X \to  \mecdga_B \, , \,\,\,\, I_F \ni (\Spec \, C \to F) \longmapsto \DR(C/B)$$ (resp. $\DR(G/B)= lim_{I_X} \DR^{G}$, where $$\DR^G: I_X \to  \mecdga_B \, , \,\,\,\, I_G \ni (\Spec \, C \to G) \longmapsto \DR(C/B) \,\, ).$$ There is an obvious functor $\alpha_f: I_F \to I_G$, induced by composition by $f$, hence an induced  morphism $\lim_{I_G} \DR^{G} \to \lim_{I_F} \alpha^*_f(\DR^{G})$. But $\alpha^*_f(\DR^{G}) \simeq \DR^{F}$, hence we get a morphism $\DR(G/B) \to \DR(F/B)$.
\hfill $\Box$

\

We now claim that the two de Rham complexes $\DR^{t}(\D(F)/\D(A))$ and
$\DR(F/A)$ are naturally equivalent, at least when $F$ is a perfect  formal derived stack over $\Spec \, A$ that is moreover 
\emph{algebraisable over $\Spec \, A$} as in
Definition~\ref{d10++}. More precisely, we have
 
\begin{thm}\label{t2}
  Let $F \longrightarrow \Spec\, A$ be a perfect formal derived
  stack. We assume that $F$ is moreover algebraisable over $\Spec \,
  A$ (Definition~\ref{d10++}). Then, there are natural morphisms
$$
\xymatrix{
\DR(\D(F)/\D(A)) \ar[r] & \DR^t(\D(F)/\D(A)) \ar[r] & 
\lim_{\Spec\, B \rightarrow F} \DR^t(\D(B)/\D(A)) & \ar[l] \DR(F/A)}
$$
that are all equivalences of graded mixed $A$-cdgas.
\end{thm}
\noindent \textbf{Proof.}  We start by defining the three natural
morphisms. The first morphism on the left is induced by the natural
transformation $|.| \rightarrow |.|^t$, from realization to Tate
realization (see \S 1.5). The second morphism on the left is induced
by functoriality.  It remains to describe the morphism on the right
$$
\DR(F/A) \longrightarrow  \lim_{\Spec\, B \rightarrow F}
\DR^t(\D(B)/\D(A)).
$$ 
By definition \ref{reldrstacks} 
$$
\DR(F/A) \simeq \lim_{\Spec\, B \rightarrow F} \DR(B/A),
$$
and we have a morphism of graded mixed cdgas
$B \longrightarrow \D(B)$, where $B$ is considered with its trivial
mixed structure of pure weight $0$. This morphism is the adjoint to the
equivalence $B \simeq |\D(B)|$ of Proposition \ref{p2}. By functoriality
it comes with a commutative square of graded mixed cdgas
$$\xymatrix{
B \ar[r] & \D(B) \\
A \ar[u] \ar[r] & \D(A), \ar[u]}
$$
and thus induces a morphism on de Rham objects
$$\DR(B/A) \longrightarrow \DR(\D(B)/\D(A)) \longrightarrow
\DR^t(\D(B)/\D(A)).$$ 
By taking the limit, we get the desired map 
$$
\DR(F/A) \longrightarrow  \lim_{\Spec\, B \rightarrow F}
\DR^t(\D(B)/\D(A)).
$$

To prove the statement of Theorem~\ref{t2}, we first observe that all
the graded mixed cdgas $\D(F)$ and $\D(B)$ are positively weighted, as
they are freely generated, as graded cdgas, by their weight $1$ part
(see Proposition \ref{p3}). The natural morphisms
$$
\DR(\D(B)/\D(A)) \longrightarrow \DR^t(\D(B)/\D(A)) \qquad
\DR(\D(F)/\D(A)) \longrightarrow \DR^t(\D(F)/\D(A))
$$
are then equivalence by trivial weight reasons. So, it will be enough
to check the following two statements

\begin{enumerate}
\item The descent morphism
$$
\DR(\D(F)/\D(A)) \longrightarrow 
\lim_{\Spec\, B \rightarrow F} \DR(\D(B)/\D(A))
$$
is an equivalence. 
\item For any $\Spec\, B \longrightarrow \Spec\, A$, the natural morphism
$$
\DR(B/A) \longrightarrow \DR(\D(B)/\D(A))
$$
is an equivalence. 
\end{enumerate}

\

\noindent
Statement $(1)$ is proved using the fact that $F$ is algebraisable completely
analogously to the proof of Proposition \ref{p3}.  We first note that
the assignment $\Spec\, B \mapsto \DR(\D(B)/\D(A))$ is a stack for the
etale topology, so the right hand side in (1) is simply the left Kan
extension of $\Spec\, B \mapsto \DR(\D(B)/\D(A))$ to all derived
stacks. In particular, it has descent over $F$.  We write
$F=\widehat{G}_{f}$, for a morphism $f : \Spec\, A_{red}
\longrightarrow G$, with $G$ an algebraic derived $n$-stack locally of
finite presentation over $A$.  By localizing with respect to the
\'etale topology on $\Spec\, A_{red}$, we can assume that there is an
affine derived scheme $U$ with a smooth map $U \longrightarrow G$,
such that $f$ factors through $U$. We let $\widehat{U_{*}}$ denote the formal
completion of the nerve of $U \rightarrow G$ along the morphism
$\Spec\, A_{red} \longrightarrow U_*$.  We now claim that the natural
morphism
$$
\DR(\D(F)/\D(A)) \longrightarrow 
\lim_{n\in \Delta} \DR(\D(\widehat{U_n})/\D(A))
$$
is an equivalence. We will actually prove the stronger statement that
the induced morphism 
$$
\wedge^{p}\mathbb{L}^{int}_{\D(F)/\D(A)} \longrightarrow 
\lim_{n\in \Delta}\wedge^{p}\mathbb{L}^{int}_{\D(\widehat{U_n})/\D(A)}
\qquad (*)
$$
is an equivalence of non-mixed graded complexes for all $p$. For this,
we use Proposition~\ref{p3}, which implies that we have equivalences
of graded modules
$$
\begin{aligned}
\wedge^{p}\mathbb{L}^{int}_{\D(F)/\D(A)} & \simeq
\D(F)\otimes_{A_{red}}\wedge^{p}f^{*} (\mathbb{L}_{G/A}) \\
\wedge^{p}\mathbb{L}^{int}_{\D(\widehat{U_n})/\D(A)} &
\simeq
\D(\widehat{U_n})\otimes_{A_{red}}\wedge^{p}f^{*}(\mathbb{L}_{U_{n}/A}).
\end{aligned}
$$
Since $\D(F)\simeq \lim_{n}\D(\widehat{U_n})$, and tensor product of
perfect modules preserves limits, we obtain $(*)$ as all
$f^{*}(\mathbb{L}_{U_{n}/A})$ and $f^{*} (\mathbb{L}_{G/A})$ are
perfect complexes of $A_{red}$-modules, and because differential forms
satisfy descent (see Appendix B), so that
$$
f^*(\mathbb{L}_{G/A})\simeq
\lim_{n}\wedge^{p}f^{*}(\mathbb{L}_{U_{n}/A}).
$$
By induction on the geometric level $n$ of $G$, we finally see that
statement $(1)$ can be reduced to the case where $G=\Spec\, B$ is
affine and $f : \Spec\, A_{red} \longrightarrow G$ is a closed
immersion.  In this case, we have already seen that $F$ can be written
as $\mathrm{colim}_{n} \Spec\, B_n$, for a system of closed immersions
$\Spec\, B_n \longrightarrow \Spec\, B_{n+1}$ such that $(B_n)_{red}
\simeq A_{red}$. This colimit can be taken in derived prestacks, so
Appendix~B \ref{lappend} applies.  This implies statement $(1)$, as we
have
$$
\begin{aligned}
\wedge^p\mathbb{L}^{int}_{\D(F)/\D(A)} & \simeq \D(F) \otimes_{B}
\wedge^p \mathbb{L}_{B/A} \\
\wedge^p\mathbb{L}^{int}_{\D(B_n)/\D(A)} & \simeq \D(B_n) \otimes_{B}
\wedge^p \mathbb{L}_{B_n/A}.
\end{aligned}
$$

\noindent
It remains to prove statement $(2)$. We need to show that 
the natural morphism $B \rightarrow \D(B)$ and 
$A \rightarrow \D(A)$ induces an equivalence
$$
\wedge^{p}\mathbb{L}_{B/A} \longrightarrow
|\wedge^{p}\mathbb{L}^{int}_{\D(B)/\D(A)}|.
$$  
This is the relative version of the following lemma, and
can be in fact deduced from it.

\begin{lem}\label{lt33}
If $F=\Spec\, A$ is an affine derived scheme then 
the natural morphism
$$\DR(A/k) \longrightarrow \DR(\D(A))$$
is an equivalence of graded cdgas.
\end{lem}
\textbf{Proof of lemma.} It is enough to show that the induced
morphism
$$\mathcal{A}^{p}(A)\simeq \wedge^{p}\mathbb{L}_{A} \longrightarrow
|\wedge^{p}\mathbb{L}_{\mathbb{D}(A)}|$$
is an equivalence of complexes, for any $p\geq 0$.  

The proof will now involve strict models. We choose a cell model for
$A$ with finitely many cells in each dimension, and a factorization
$$
\xymatrix{
A \ar[r] & A' \ar[r] & A_{red},}
$$
where $A' \longrightarrow A_{red}$ is an equivalence and 
$A'$ is a cell $A$-algebra with finitely many cells in each dimension. 
Moreover, as $\pi_{0}(A) \longrightarrow \pi_{0}(A_{red})$ is 
surjective, we can chose $A'$ having cells only in dimension $1$ and higher
(i.e. no $0$-dimensional cells). With such choices, the cotangent complex
$\mathbb{L}_{A_{red}/A}$ has a strict model 
$\Omega^{1}_{A'/A}$, and is itself a cell $A'$-module with 
finitely many cells in each dimension, and no $0$-dimensional cell. 
We let $L:=\Omega^{1}_{A'/A}$. 

The graded mixed cdga $\mathbb{D}(A)$ can then be represented (\S
\ref{1.3.3}) by the strict de Rham algebra
$\mathbb{D}^{str}(A):=Sym_{A'}(L[-1])$. We consider $B:=(A')^{0}=A^0$
the degree $0$ part of $A'$ (which is also the degree $0$ part of $A$
because $A'$ has no $0$-dimensional cell over $A$), and let
$V:=L^{-1}$ the degree $(-1)$ part of $L$. The $k$-algebra $B$ is just
a polynomial algebra over $k$, and $V$ is a free $B$-module whose rank
equals the number of $1$-dimensional cells of $A'$ over $A$.

For the sake of clarity, we introduce the following notations. 
For $E \in \medg$ a graded mixed $k$-dg-module, we
let 
$$
|E|:=\prod_{i\geq 0}E(i),
$$
the product of the non-negative weight parts of $E$, endowed with its natural
total differential sum of the cohomological differential and the mixed
structure.  
In the same way, we let 
$$
|E|^{\oplus}:=\oplus_{i\geq 0}E(i),
$$
to be the coproduct of the non-negative weight parts of $E$, with the
similar differential, so that $|E|^{\oplus}$ sits naturally inside
$|E|$ as a sub-dg-module. Note that $|E|$ is a model for
$\mathbb{R}\underline{Hom}_{\edg}(k(0),E)$, whereas $|E|^{\oplus}$ is
a rather silly functor which is not even invariant under
quasi-isomorphisms of graded mixed dg-modules.

As we have already seen in the proof of Lemma \ref{l4}, there exists a
strict push-out square of cdgas
$$
\xymatrix{
|Sym_{A'}(L[-1])|^{\oplus} \ar[r] & |Sym_{A'}(L[-1])| \\
Sym_{B}(V) \ar[r] \ar[u] & \widehat{Sym}_{B}(V) \ar[u],}
$$
where $\widehat{Sym}$ denotes the completed symmetric algebra, i.e the
infinite product of the various symmetric powers. This push-out is
also a homotopy push-out of cdgas because the bottom horizontal
morphism is a flat morphism of commutative rings.

We have the following version of the above push-out square for modules, too.  
Let $M \in \mathbb{D}^{str}(A)-Mod_{\medg}$ a graded mixed
$Sym_{A'}(L[-1])$-dg-module. We assume that, as a graded dg-module, 
$M$ is isomorphic to 
$$
M \simeq \mathbb{D}^{str}(A) \otimes_{A'}E,
$$
where $E$ is a graded $A'$-dg-module pure of some weight $i$, and
moreover, $E$ is a cell module with finitely many cells in each
non-negative dimension.  Under these finiteness conditions, it can be
checked that there is a natural isomorphism
$$
|M|^{\oplus} \otimes_{Sym_{B}(V)}\widehat{Sym}_{B}(V) \simeq 
|M|.
$$
The same is true for any graded mixed $\mathbb{D}^{str}(A)$-dg-module
$M$ which is (isomorphic to) a successive extension of graded mixed
modules as above. In particular, we can apply this to
$\Omega^{1}_{\mathbb{D}^{str}(A)}$ as well as to
$\Omega^{p}_{\mathbb{D}^{str}(A)}$, for any $p>0$. Indeed, there is a
short exact sequence of graded $Sym_{A'}(L[-1])$-modules
$$
\xymatrix@1{ 0 \ar[r] & 
\Omega^{1}_{A'} \otimes_{A'} Sym_{A'}(L[-1]) \ar[r] &
\Omega^{1}_{\mathbb{D}^{str}(A)} \ar[r] &  
L\otimes_{A'}Sym_{A'}(L[-1])[-1] \ar[r] & 0.}
$$
This shows that for all $p>0$, we have a canonical isomorphism
$$
|\Omega^{p}_{\mathbb{D}^{str}(A)}|^{\oplus}
\otimes_{Sym_{B}(V)}\widehat{Sym}_{B}(V) \simeq
|\Omega^{p}_{\mathbb{D}^{str}(A)}|.
$$
Now we notice that the natural morphism
$$
|\Omega^{p}_{\mathbb{D}^{str}(A)}|^{\oplus} \longrightarrow 
|\Omega^{p}_{\mathbb{D}^{str}(A)}|^{\oplus}
\otimes_{Sym_{B}(V)}\widehat{Sym}_{B}(V)
$$
is isomorphic to 
$$
|\Omega^{p}_{\mathbb{D}^{str}(A)}|^{\oplus} \longrightarrow 
|\Omega^{p}_{\mathbb{D}^{str}(A)}|^{\oplus} \otimes_{|\mathbb{D}^{str}(A)|^{\oplus}}
|\mathbb{D}^{str}(A)|.
$$
Let us show that

\begin{sublem}
For all $p\geq 0$ the above morphism
$$|\Omega^{p}_{\mathbb{D}^{str}(A)}|^{\oplus} \longrightarrow 
|\Omega^{p}_{\mathbb{D}^{str}(A)}|^{\oplus} \otimes_{|\mathbb{D}^{str}(A)|^{\oplus}}
|\mathbb{D}^{str}(A)|$$
is a quasi-isomorphism.
\end{sublem}
\noindent \textit{Proof of sub-lemma.} First of all, in the push-out
square of cdgas
$$
\xymatrix{
|\mathbb{D}^{str}(A)|^{\oplus} \ar[r] & |\mathbb{D}^{str}(A)| \\
Sym_{B}(V) \ar[r] \ar[u] & \widehat{Sym}_{B}(V) \ar[u],}
$$
the bottom horizontal arrow is flat. This implies that the tensor
product 
$$
|\Omega^{p}_{\mathbb{D}^{str}(A)}|^{\oplus} \otimes_{|\mathbb{D}^{str}(A)|^{\oplus}}
|\mathbb{D}^{str}(A)|
$$ 
is also a derived tensor product. The sub-lemma would then
follow from the fact that the inclusion
$$
|\mathbb{D}^{str}(A)|^{\oplus} \hookrightarrow |\mathbb{D}^{str}(A)|
$$
is a quasi-isomorphism. To see this, we consider the diagram of
structure morphism over $A$ 
$$
\xymatrix{
& A \ar[dl]_-{u} \ar[dr]^-{v} & \\
|\mathbb{D}^{str}(A)|^{\oplus}  \ar[rr] & & |\mathbb{D}^{str}(A)|.}
$$
The morphism $v$ is an equivalence by Proposition \ref{p2}
and lemma \ref{l4}. The morphism 
$u$ is the inclusion of $A$ into the non-completed derived de Rham 
complex of $A_{red}$ over $A$, and thus is also a quasi-isomorphism.
\hfill $\Box$ 

\

\medskip

Now we can prove that the above sub-lemma implies Lemma
\ref{lt33}. Indeed, the morphism
$$
\wedge^{p}\mathbb{L}_{A} \longrightarrow
|\wedge^{p}\mathbb{L}_{\mathbb{D}(A)}|
$$ 
can be represented by 
the composition of morphisms between strict models 
$$
\xymatrix{ \Omega^{p}_{A} \ar[r] &
  |\Omega^{p}_{\mathbb{D}^{str}(A)}|^{\oplus} \ar[r] &
  |\Omega^{p}_{\mathbb{D}^{str}(A)}|^{\oplus}
  \otimes_{|\mathbb{D}^{str}(A)|^{\oplus}} |\mathbb{D}^{str}(A)|
  \ar[r] & |\Omega^{p}_{\mathbb{D}^{str}(A)}|.}
$$
The two rightmost morphisms are quasi-isomorphisms by what we have
seen, while the leftmost one can simply be identified, up to a
canonical isomorphism, with the natural morphism
$$
\Omega^{p}_{A} \longrightarrow
\Omega^{p}_{|\mathbb{D}^{str}(A)|^{\oplus}}.
$$
This last morphism is again a quasi-isomorphism because it is induced
by the morphism
$$
A \longrightarrow |\mathbb{D}^{str}(A)|^{\oplus}
$$
which is a quasi-isomorphism of quasi-free, and thus cofibrant, cdgas.
\hfill $\Box$ 

\

\medskip

\noindent
Lemma \ref{lt33} is proven, and we have thus finished the proof of
Theorem \ref{t2}.  \hfill $\Box$

\

\bigskip

\noindent The following corollary is a consequence of the proof
Theorem \ref{t2}.  

\begin{cor}\label{ct2}
  Let $F \longrightarrow \Spec\, A$ be a perfect formal derived stack
  over $\Spec\, A$, and assume that $F$ is algebraisable. Let
$$
\phi_F : \mathrm{L}_{\mathsf{Perf}}(F) \longrightarrow
\D(F)-Mod_{\edg}^{\mathsf{Perf}}
$$
be the equivalence of Proposition~\ref{p4}. Then, there is a canonical 
equivalence of graded mixed $\D(F)$-modules
$$
\phi_F (\mathbb{L}_{F/A}) \simeq
\mathbb{L}^{int}_{\D(F)/\D(A)}\otimes_{k} k((1)).
$$
\end{cor}
\noindent \textbf{Proof.} First of all, as graded $\D(F)$-modules we have
(Proposition \ref{p4})
$$
\mathbb{L}^{int}_{\D(F)/\D(A)} \simeq \D(F)\otimes_{A_{red}}
f^{*}(\mathbb{L}_{F/A}),
$$
where $f : \Spec\, A \longrightarrow F$ is the natural morphism, 
and $f^{*}(\mathbb{L}_{F/A})$  sits in pure weight $1$, so that,
according to our conventions, we should rather write 
$$
\mathbb{L}^{int}_{\D(F)/\D(A)} \simeq \D(F)\otimes_{A_{red}}
f^{*}(\mathbb{L}_{F/A}) \otimes_k k((-1)).
$$
In particular, $\mathbb{L}^{int}_{\D(F)/\D(A)}\otimes_{k} k((1))$
belongs to $\D(F)-Mod_{\edg}^{\mathsf{Perf}}$, as it is now free over
its weight $0$  
part. 

Moreover, the same proof as in Theorem~\ref{t2} shows that for any
perfect complex $E \in \mathrm{L}_{\mathsf{Perf}}(F)$, we have a
natural equivalence, functorial in $E$
$$
\Gamma(F,E \otimes_{\OO_F}\mathbb{L}_{F/A})\simeq 
|\phi_F (E)\otimes_{\D(F)}\mathbb{L}^{int}_{\D(F)/\D(A)}|.
$$
We have a natural map $k=k((0)) \rightarrow k((-1))$ in the
$\s$-category of graded mixed complexes, represented by the map
$\tilde{k} \to k((-1))$ sending $x_1$ to $1$, in the notation of  \S
\ref{closedforms}. Its weight-shift by $1$ gives us a canonical map
$k((1)) \to k$ in the $\s$-category of graded mixed complexes,
inducing a morphism  
$$ 
\mathbb{L}^{int}_{\D(F)/\D(A)} \otimes_k k((1)) \longrightarrow
\mathbb{L}^{int}_{\D(F)/\D(A)}.
$$
Finally, this morphism
induces an equivalence
$$ 
|\phi_F (E)\otimes_{\D(F)}\mathbb{L}^{int}_{\D(F)/\D(A)} \otimes_k
k((1))| \simeq |\phi_F
(E)\otimes_{\D(F)}\mathbb{L}^{int}_{\D(F)/\D(A)}|. $$
We thus get an equivalence 
$$
\Gamma(F,E \otimes_{\OO_F}\mathbb{L}_{F/A})\simeq 
|\phi_F
(E)\otimes_{\D(F)}\mathbb{L}^{int}_{\D(F)/\D(A)}\otimes_{k}k((1))|, 
$$ 
functorial in $E$.  Observe now that
$\phi_F(E)\otimes_{\D(F)}\mathbb{L}^{int}_{\D(F)/\D(A)}\otimes_{k}k((1))$
is a perfect graded mixed $\D(F)$-module. Since $E$ is perfect, these
equivalence can also be re-written as
$$
\mathbb{R}\underline{Hom}(E^{\vee},\mathbb{L}_{F/A})
\simeq \mathbb{R}\underline{Hom}(\phi_{F}(E)^{\vee},
\mathbb{L}^{int}_{\D(F)/\D(A)}\otimes_{k}k((1))).
$$
Now, $\phi_F$ is an equivalence, and therefore Yoneda lemma implies
that $\phi_F(\mathbb{L}_{F/A})$ and
$\mathbb{L}^{int}_{\D(F)/\D(A)}\otimes_{k}k((1))$ are naturally
equivalent.
\hfill $\Box$

\subsubsection{Shifted polyvectors over perfect formal derived
  stacks} \label{polypfds}
We present here a version of Theorem \ref{t2} for shifted  polyvectors.

Let $F$ be a perfect formal derived stack over $\Spec\, A$.  We have
the corresponding graded mixed cdga $\D(F)$, which we consider as a
graded mixed $\D(A)$-algebra. By taking $\C= \medg_{\D(A)}$, we have
the corresponding the graded $\mathbb{P}_{n+1}$-dg-algebra of
$n$-shifted polyvectors $\Pol(\D(F),n)$ (Definition~\ref{d9} (2)), as
well as its Tate version $\Pol^{t}(\D(F),n)$ (Definition~\ref{drpoltate}
(2)). To emphasize the fact that such objects are defined relative to
$\D(A)$, we will more precisely denote them by $\Pol(\D(F)/\D (A),n)$,
and $\Pol^{t}(\D(F)/\D(A),n)$, respectively.

On the other hand, we can give the following general 

\begin{df}\label{polgeneral} Let $n\in \mathbb{Z}$, and $f: X
  \longrightarrow Y$ be a morphism of derived stacks, such that the
  relative cotangent complex $\mathbb{L}_{X/Y}$ is defined and is an
  object in $\mathrm{L}_{\mathsf{Perf}}(X)$. Then, we define
$$
\Pol(X/Y,n):=\bigoplus_{p}
(\underline{Hom}_{\mathrm{L}_{\mathrm{QCoh}}(X)}(
\otimes^{p}\mathbb{L}_{X/Y},\mathcal{O}_X[pn]))^{\Sigma_p} \in
\mdg_{k},
$$
where $\mathrm{L}_{\mathrm{QCoh}}(X)\simeq \lim_{\Spec\, A \rightarrow
  X} \mathrm{L}(A)$ is considered as a dg-category over $k$, and
$\underline{Hom}_{\mathrm{L}_{\mathrm{QCoh}}(X)}$ denotes its
$k$-dg-module of morphisms.
\end{df}

Note that, in particular, $\Pol(X/Y,n)$ is defined if $X$ and $Y$ are
derived Artin stacks locally of finite presentation over $k$, or if
$Y=\Spec \, A$ and $f: X \to Y$ is a perfect formal derived stack. 

\

\begin{thm}\label{t3}
If $F$ is a perfect formal derived stack over $\Spec\, A$, and
 $F$ is algebraisable, then 
there is a natural equivalence of graded $k$-dg-modules
$$\Pol^{t}(\D(F)/\D(A),n) \simeq \Pol(F/A,n).$$
\end{thm}
\noindent \textbf{Proof.} We have $\mathbb{L}_{F/A} \in
\mathrm{L}_{\mathsf{Perf}}(F)$, and we consider the equivalence of
Corollary \ref{c2}
$$
\phi_{F} : \mathrm{L}_{\mathsf{Perf}}(F) \longrightarrow
\D(F)-Mod_{\edg}^{\mathsf{Perf}}.
$$
By Corollary \ref{ct2}, there is a natural 
equivalence of graded mixed $\D(F)$-modules
$$
\phi_{F}(\mathbb{L}_{F/A}) \simeq
\mathbb{L}^{int}_{\D(F)/\D(A)}\otimes_{k}k((1)).
$$
As $\phi_{F}$ is a symmetric
monoidal equivalence, we get 
$$
\phi_{F}(Sym^{p}_{\OO_F}(\mathbb{T}_{F/A}[n]) \simeq
Sym^{p}(\mathbb{T}^{int}_{\D(F)/\D(A)}[n])\otimes_{k}k((-p)),
$$
for any $n$ and $p$. The result then follows from the fact that $\phi_F$
is an equivalence together with the fact that the Tate realization is
a stable realization, i.e. that, for any graded mixed $\D(F)$-module
$E$, there is a natural equivalence $|E|^t \simeq
|E\otimes_{k}k(1)|^t$.
\hfill $\Box$ 

\

\begin{rmk}
Note that corollary \ref{ct2} implies that  
$$
\mathbb{T}^{int}_{\D(F)/\D(A)}\simeq
\D(F)\otimes_{A_{red}}f^{*}(\mathbb{T}_{F/A})\otimes_k k((-1)),
$$ 
as a graded modules, where $f : \Spec\, A_{red} \longrightarrow F$ is
the natural morphism. The weight-shift on the right hand side gives no
chance for Theorem \ref{t3} to be true if the Tate realization
$|-|^{t}$ is replaced by the standard one $|-|$, while this is true in
the case of de Rham complexes.
\end{rmk}

\subsection{Global aspects and mixed principal parts}\label{GLOB} 

In this last part of Section $2$ we present the global aspects of what
we have seen so far, namely \emph{families} of perfect formal derived
stacks and their associated graded mixed cdgas.

\subsubsection{Families of perfect formal derived stacks}\label{fam}

\noindent We start by the notion of \emph{families} of perfect formal
derived stacks.

\begin{df}
  A morphism $X \longrightarrow Y$ of derived stacks is \emph{a family
    of perfect formal derived stacks over $Y$} if, for all $\Spec\, A
  \in \dAff_k$ and all morphism $\Spec\, A \longrightarrow Y$, the
  fiber
$$
X_{A}:=X \times_{Y}\Spec\, A \longrightarrow \Spec\, A
$$
is a perfect formal derived stack over $\Spec\, A$ in the sense of
Definition \ref{d10+}. 
\end{df}

\

Note that, in the above definition, all derived stacks $X_{A}$ have
perfect cotangent complexes, for all $\Spec\, A$ mapping to $Y$.  This
implies that the morphism $X \longrightarrow Y$ itself has a relative
cotangent complex $\mathbb{L}_{X/Y} \in \mathrm{L}_{\mathrm{QCoh}}(X)$
which is moreover perfect (see \cite[\S 1.4.1]{hagII}). In particular,
for any $n\in \mathbb{Z}$, the graded $k$-dg-module $\Pol(X/Y,n)$ is
well defined (Definition~\ref{polgeneral}).

%%%%%%%%%%%%%%%%%%%
\begin{df}\label{generalreldrstacks} 
Let $F \longrightarrow G$ be an arbitrary map of derived stacks.
The \emph{relative de Rham object of the 
derived stack $F$ over $G$} is 
$$\DR(F/G):=\lim_{\Spec\, A \rightarrow G} \DR(F_A/ A) \in \mecdga_k$$
where  $\DR(F_A/ A)$ is as in Def. \ref{reldrstacks}, and the limit is taken in the $\infty$-category $\mecdga_k$ over all morphisms $\Spec\, A \rightarrow G$.
\end{df}
%%%%%%%%%%%%%%%%%%%

\begin{prop}\label{trans2} Let $f: F \to G$ and $g: G \to H$ be maps of derived stacks. There are canonical induced maps $\DR(G/H) \to \DR(F/H) \to \DR (F/G)$ in $\mecdga_k$.
\end{prop}

\noindent \textbf{Proof.} The first map $\DR(G/H) \to \DR(F/H)$ follows easily from Prop. \ref{trans1}. In order to produce the second map, let $J_{F/H}$ be the category on which $\DR(F/H)$ is defined as a limit $\DR (F/H)= \lim _{J_{F/H}} \DR^{F/H}$ where $$\DR^{F/H}: J_{F/H} \to \mecdga_k \, , \,\,\,\,\, J_{F/H} \ni (\Spec \, A \to H) \longmapsto \DR(F_{gf, A} /A) \in \mecdga_k \, ,$$ where $F_{gf, A}$ denotes the base change of $g \circ f$ along $\Spec \, A \to H$.\\
Analogously, let $J_{F/G}$ be the category on which $\DR(F/G)$ is defined as a limit $\DR (F/G)= \lim _{J_{F/G}} \DR^{F/G}$ where $$\DR^{F/G}: J_{F/G} \to \mecdga_k \, , \,\,\,\,\, J_{F/H} \ni (\Spec \, B \to G) \longmapsto \DR(F_{f,B} /B) \in \mecdga_k \, ,$$ where $F_{f, B}$ denotes the base change of $f$ along $\Spec \, B \to G$. \\There is an obvious functor $\alpha: J_{F/G} \to J_{F/H}$ induced by composition with $g$, hence a morphism $c: \lim_{J_{F/H}}\DR^{F/H} \to \lim_{J_{F/G}} \alpha^{*} (\DR^{F/H})$ $\mecdga_k$. \\ Moreover, for any $(\Spec\, B \to G) \in J_{F/G}$, we have an induced canonical map $F_{f, B} \to F_{gf, B}$ over $\Spec \, B$, hence, by Prop. \ref{trans1}, a further induced map $\DR(F_{gf, B}/B) \to \DR(F_{f, B}/B)$. Thus, we get a morphism of functors $\varphi: \alpha^* (\DR^{F/H}) \to \DR^{F/G}$. The composition $$\xymatrix{\DR(F/H) \simeq \lim_{J_{F/H}}\DR^{F/H} \ar[r]^-{c}& \lim_{J_{F/G}} \alpha^{*} (\DR^{F/H})  \ar[rr]^-{\lim\, \varphi} &  & \lim_{J_{F/G}}  \DR^{F/G} \simeq \DR(F/G)}$$ gives us the second map.

\hfill $\Box$

\

\begin{rmk} Note that, when $\mathbb{L}_{F/H}$ and $\mathbb{L}_{G/H}$ (hence $\mathbb{L}_{F/G}$) exist, the sequence in Proposition \ref{trans2} becomes a fiber-cofiber sequence when considered inside $\dg^{gr}_k$. In fact, for a map $X \to Y$ between derived stacks  having cotangent complexes, Proposition \ref{p1} implies an equivalence $$\DR(X/Y) \simeq \bigoplus_{p \geq 0}\Gamma (X, Sym_{\OO_X}^p (\mathbb{L}_{X/Y} [-1]))$$ in $\cdga_{k}^{gr}$.
\end{rmk}

\

\begin{rmk}
Our main example and object of interest will be the following family of perfect formal derived stacks 
$$
q : X \longrightarrow X_{DR}
$$
for $X$ an Artin derived stack locally of finite presentation over
$k$. Corollary~\ref{ffdrisgood} shows that this is indeed a family of perfect formal derived stacks.
\end{rmk}

Let $X \longrightarrow Y$ be a perfect family of formal derived stacks as
above. The $\s$-category $\dAff_{k}/Y$ of derived affine schemes
over $Y$ comes equipped with a tautological prestack of cdgas
$$
\OO_{Y} : (\dAff_{k}/Y)^{op} \longrightarrow \cdga_{k} , \qquad
(\Spec\, A \rightarrow Y) \longmapsto A.
$$
For each $\Spec\, A \rightarrow Y$, we may associate to the good
formal derived stack $X_A$ its graded mixed cdga $\D(X_A) \in
A/\mecdga_k$ (Theorem~\ref{t1}).  Moreover, the morphism $X_A
\rightarrow \Spec\, A$ induces a natural $\D(A)$-linear structure on
$\D(X_A)$, and we will thus consider $\D(X_A)$ as on object in $\D (A)
/\mecdga_k$.

If $\Spec\, B \longrightarrow \Spec\, A$ is a morphism 
in $\dAff_k/Y$ we have an induced natural morphism of $\D(A)$-linear
graded mixed cdgas 
$$
\D(X_A) \longrightarrow \D(X_B).
$$
With a bit of care in the $\s$-categorical constructions (e.g. by using
strict models in model categories of diagrams), we obtain 
the following \emph{prestacks} of graded mixed cdgas on $\dAff_{k}/Y$:
$$
\begin{aligned}
\D_Y:=\D(\OO_Y) : &  (\dAff_{k}/Y)^{op} \longrightarrow  \mecdga_{k},
\qquad 
(\Spec\, A \rightarrow Y) \longmapsto \D(A), \\
\D_{X/Y} : & (\dAff_{k}/Y)^{op} \longrightarrow \mecdga_k, \qquad 
(\Spec\, A \rightarrow Y) \longmapsto  \D(X_A).
\end{aligned}
$$
The natural $\D(A)$-structure on $\D(X_A)$ gives a natural morphism of
prestacks of graded mixed cdgas
$$
\D_Y \longrightarrow \D_{X/Y},
$$
which we consider as the datum of a $\D_Y$-linear structure on
$\D_{X/Y}$.

\begin{rmk}
  The two prestacks $\D_Y$ and $\D_{X/Y}$ defined above, are not
  stacks for the induced \'etale topology on $\dAff_{k}/Y$. See
  however Remark \ref{notstbut} below.
\end{rmk}

By taking $\C$ as the $\s$-category of functors $(\dAff_{k}/Y)^{op}
\rightarrow \medg_{k}$, we may apply to the prestacks $\D_Y$ and
$\D_{X/Y}$ the constructions $\DR$, $\DR^t$ and $\Pol^t$ of \S
\ref{drpform} and \S \ref{polypfds}, and obtain the following
prestacks on $\dAff_{k}/Y$
$$
\DR(\D_{X/Y}/\D_Y) \qquad \DR^t(\D_{X/Y}/\D_Y) \qquad
\Pol^t(\D_{X/Y}/\D_Y,n).
$$
The first two are prestacks of graded mixed cdgas while the last one
is a prestack of graded $\mathbb{P}_{n+1}$-algebras.

The main results of Subsection \ref{subssfpspfds}, i.e. Theorem~\ref{t2},
Corollary~\ref{ct2}, and Theorem~\ref{t3}, imply the following result for
\emph{families} of perfect formal derived stacks

\begin{cor}\label{c3}
  Let $f : X \longrightarrow Y$ be a family of perfect formal derived
  stacks.  We assume that for each $\Spec\, A \longrightarrow Y$ the
  perfect formal derived stack $X_A$ is moreover algebraisable. Then
\begin{enumerate}
\item There is a natural equivalence of graded mixed cdga's over $k$
$$
\DR(X/Y) \simeq \Gamma(Y,\DR(\D_{X/Y}/\D_Y))\simeq
\Gamma(Y,\DR^t(\D_{X/Y}/ \D_Y)).
$$
\item  For each $n\in \mathbb{Z}$,  there is a natural equivalence of
  graded $k$-dg-modules 
$$
\Pol(X/Y,n) \simeq \Gamma(Y,\Pol^t(\D_{X/Y}/\D_Y,n)).
$$
\item There is a natural equivalence 
of $\s$-categories
$$
\mathrm{L}_{\mathsf{Perf}}(X) \simeq
\D_{X/Y}-Mod_{\medg}^{\mathsf{Perf}},
$$
where $\D_{X/Y}-Mod_{\medg}^{\mathsf{Perf}}$ consists of prestacks $E$
of graded mixed $\D_{X/Y}$-modules on $Y$ satisfying the following two
conditions:
\begin{enumerate}
\item For all $\Spec\, A \longrightarrow Y$, the graded mixed
  $\D_{X/Y}(A)$-module 
$E(A)$ is perfect in the sense of Theorem~\ref{t1} (2).
\item $E$ is quasi-coherent in the following sense: for all $\Spec\, B
  \longrightarrow \Spec\, A$ in $\dAff_{k}/Y$ the induced morphism
$$
E(A) \otimes_{\D_{X/Y}(A)}\D_{X/Y}(B) \longrightarrow E(B)
$$
is an equivalence.
\end{enumerate}
\end{enumerate}
\end{cor}

\

\noindent
Note that in the above corollary the $\s$-category
$\D_{X/Y}-Mod_{\medg}^{\mathsf{Perf}}$ can also be defined as the
limit of $\s$-categories
$$
\D_{X/Y}-Mod_{\medg}^{\mathsf{Perf}} := 
\lim_{\Spec\, A \rightarrow Y}
\D_{X/Y}(A)-Mod_{\medg}^{\mathsf{Perf}}.
$$

\begin{rmk}\label{directimage} Parts $(1)$ and $(2)$ of Corollary \ref{c3}  can be made a bit more precise. We have direct image  
  prestacks on $\dAff_{k}/Y$
$$
f_*(\DR(-/Y)) \quad \text{and} \quad  f_*(\Pol(-/Y,n)),
$$
defined by sending $\Spec\, A\longrightarrow Y$ to 
$$
\DR(X_A/A) \quad \text{and} \quad \Pol(X_A/A,n).
$$
These are prestacks of graded mixed cdgas and of graded
$\mathbb{P}_{n+1}$-algebras, respectively, and are indeed stacks for
the \'etale topology (being direct images of stacks). Corollary
\ref{c3} can be refined to the existence of equivalences of prestacks
over $\dAff_k /Y$
$$
f_*(\DR(-/Y)) \simeq \DR(\D_{X/Y}/\D_Y) \qquad
f_* (\Pol(-/Y,n)) \simeq \Pol^t(\D_{X/Y}/\D_Y,n)
$$ 
\emph{before} taking global sections (i.e. one recovers Corollary~\ref{c3}
$(1)$ and $(2)$ from these equivalences of prestacks by taking global
sections, i.e. by applying $\lim_{\Spec \, A \to Y}$).
\end{rmk}

\

\noindent
As a consequence of Remark~\ref{directimage}, we get the following
corollary

\begin{cor}\label{c3'}
The prestacks 
$\DR(\D_{X/Y}/\D_Y)$ and $ \Pol^t(\D_{X/Y}/\D_Y,n)$
are stacks over $\dAff_{k}/Y$. 
\end{cor}

\

\noindent
We have a similar refinement also for statement $(3)$ of Corollary
\ref{c3}.  The $\s$-category $\D_{X/Y}-Mod_{\medg}^{\mathsf{Perf}}$
can be localized to a prestack of $\s$-categories on $\dAff_{k}/Y$
$$
\D_{X/Y}-\underline{Mod}_{\medg}^{\mathsf{Perf}} : (\Spec\, A \rightarrow Y)
\mapsto \D_{X/Y}(A)-Mod_{\medg}^{\mathsf{Perf}}.
$$
And we have an equivalence of prestacks of $\s$-categories 
on $\dAff_{k}/Y$
$$
f_*(\mathrm{L}_{\mathsf{Perf}}(-)) \simeq
\D_{X/Y}-\underline{Mod}_{\medg}^{\mathsf{Perf}}.
$$

\begin{rmk}\label{notstbut}
  Even though the prestacks $\D_Y$ and $\D_{X/Y}$ are not stacks for
  the induced \'etale topology, the associated constructions we are
  interested in, namely their de Rham complex, shifted polyvectors and
  perfect modules, are in fact stacks.  In a sense, this shows that
  the defect of stackiness of $\D_Y$ and $\D_{X/Y}$ is somehow
  artificial, and irrelevant for our purposes.
\end{rmk}

\subsubsection{Mixed principal parts on a derived Artin stack.}

We will be mainly interested in applying the results of \S \ref{fam}
to the special family
$$
q : X \longrightarrow X_{DR},
$$
for $X$ an Artin derived stack locally of finite presentation over
$k$. As already observed, this is a family of perfect formal derived stacks by
Corollary~\ref{ffdrisgood}.

\begin{df}\label{d11}
  Let $X$ be a derived Artin stack locally of finite presentation over
  $k$, and $q: X \longrightarrow X_{DR}$ the natural projection.
\begin{enumerate}
\item The prestack $\D_{X_{DR}}$ of graded mixed cdgas on
  $\dAff_{k}/X_{DR}$ will be called the \emph{mixed crystalline
    structure sheaf of $X$}.
\item The prestack $\D_{X/X_{DR}}$ of graded mixed cdgas under $\D_{X_{DR}}$ 
  will be called the \emph{mixed principal parts of $X$}. It will be
  denoted by
$$\mathcal{B}_{X}:=\D_{X/X_{DR}}.$$
\end{enumerate}
\end{df}

\

\noindent
The prestack mixed crystalline structure sheaf $\D_{X_{DR}}$ (which
is not a stack) is a graded mixed model for the standard crystalline
structure sheaf $\OO_{X_{DR}}$ on $\dAff_{k}/X_{DR}$. Indeed, by
Corollary \ref{c3}, we have
$$
|\D_{X_{DR}}| \simeq \DR(\D_{X_{DR}}/\D_{X_{DR}}) \simeq 
\OO_{X_{DR}}.
$$
Analogously, $\D_{X/X_{DR}}$ is a graded mixed model for the
standard sheaf of principal parts. Indeed, we have 
$$
|\D_{X/X_{DR}}| \simeq q_*(\OO_X).
$$
The value of the sheaf $q_{*}(\OO_X)$ on $\dAff_{k}/X_{DR}$ on
$\Spec\, A \rightarrow X_{DR}$ is the ring of functions on $X_A$, and
recall (Proposition~\ref{fibersDR}) that $X_A$ can be identified with
the formal completion of $X \times \Spec\, A$ along the graph of the
morphism $\Spec\, A_{red} \rightarrow X$. When $X$ is a smooth scheme
over $\Spec\, k$, the sheaf $\pi_*(\OO_X)$ is the usual sheaf of
principal parts on $X$ (\cite[16.7]{egaiv}), endowed with its natural
crystalline structure (i.e. descent data with respect to the map $q:X
\to X_{DR}$).  We may view $\B_X$ as controlling the formal completion
of $X$ along the diagonal, together with its natural
Grothendieck-Gel'fand connection.

Also recall (Lemma \ref{crucial}) that for $q: X\to X_{DR}$, we
have 
$$
\mathbb{L}_{X} \simeq \mathbb{L}_{X/X_{DR}}.
$$
In the special case of the perfect family of formal derived stacks
$q:X \rightarrow X_{DR}$, Corollary \ref{c3} thus yields the following

\begin{cor}\label{c4}
  Let $X$ be an Artin derived stack locally of finite presentation
  over $k$.
\begin{enumerate}
\item There is a natural equivalence of graded mixed cdgas over $k$
$$
\DR(X/X_{DR})\simeq \DR(X/k) \simeq
\Gamma(X_{DR},\DR(\B_X/\D_{X_{DR}})) \simeq \Gamma(X_{DR},\DR^{t}(\B_X/\D_{X_{DR}})).
$$
\item For each $n\in \mathbb{Z}$, there is a natural equivalence of
  graded complexes over $k$
$$
\Pol(X/X_{DR},n)\simeq \Pol(X,n) \simeq
\Gamma(X_{DR},\Pol^t(\B_X/\D_{X_{DR}},n)).
$$
\item There is a natural equivalence 
of $\s$-categories
$$
\mathrm{L}_{\mathsf{Perf}}(X) \simeq \B_X-Mod_{\medg}^{\mathsf{Perf}}.
$$
\item The natural $\s$-functor 
$$
\B_X-Mod_{\medg}^{\mathsf{Perf}} \longrightarrow
\B_X(\s)-Mod_{k(\s)-Mod}^{\mathsf{Perf}}, 
$$ 
induced by the base change $(-)\otimes k(\s)$, 
is an equivalence.
\end{enumerate}
\end{cor}

\noindent \textbf{Proof.} The first equivalence in $(1)$ is just the statement that the natural map $u: \DR(X/k) \to \DR(X/X_{DR})$  (Prop. \ref{trans2}) is an equivalence of mixed graded cdga's over $k$. In fact, by Prop. \ref{p1} we have equivalences of graded cdga's over $k$ $$\DR(X/X_{DR}) \simeq \bigoplus_{p \geq 0}\Gamma (X, Sym_{\OO_X}^p (\mathbb{L}_{X/X_{DR}} [-1]))\, ,$$ $$\DR(X/k) \simeq \bigoplus_{p \geq 0}\Gamma (X, Sym_{\OO_X}^p (\mathbb{L}_{X} [-1]))$$
Hence the map $u$ becomes an equivalence in $\cdga_{k}^{gr}$, by  Prop. \ref{crucial}, and is therefore itself an equivalence. The other equivalences in $(1)$ follows immediately from Corollary \ref{c3}. The proof of $(2)$ is analogous to the proof of $(1)$. Point $(3)$  follows immediately from the corresponding result in Corollary \ref{c3}. 
Only point $(4)$ requires some further explanations,
and an explicit proof.  First of all $k(\s)$ is a cdga in the $\s$-category
$\mathsf{Ind}(\medg_k)$ of Ind-objects in graded mixed complexes over
$k$. The notation $\B_X(\s)$ stands for $\B_X \otimes_{k}k(\s)$, which
is a prestack on $X_{DR}$ with values in cdgas inside
$Ind(\medg_k)$. As usual $\B_X(\s)-Mod_{k(\s)-Mod}$ denotes the
$\s$-category of prestacks of $\B_X(\s)$-modules.  Finally,
$\B_X(\s)-Mod_{k(\s)-Mod}^{\mathsf{Perf}}$ is defined as for
$\B_X-Mod_{\medg}^{\mathsf{Perf}}$: it is the full sub-$\s$-category
of $\B_X(\s)$-modules $E$ satisfying the following two conditions
\begin{enumerate}
\item For all $\Spec\, A \longrightarrow X_{DR}$, the 
$\B_X(\s)$-module $E(A)$ is of the form
$$
E(A)\simeq E_A \otimes_{\B_X(A)}\B_X(\s)(A),
$$
for $E_A$ a perfect $\B_X(A)$-graded mixed module in the sense of
Theorem~\ref{t1}.
\item For all $\Spec\, B \longrightarrow \Spec\, A$
in $\dAff_{k}/Y$, the induced morphism 
$$
E(A) \otimes_{\B_{X}(\s)(A)}\B_{X}(\s)(B) \longrightarrow E(B)
$$
is an equivalence of Ind-objects in $\medg_k$
\end{enumerate}

From this description, the natural $\s$-functor of point $(4)$ is
obtained by a limit of $\s$-functors
$$
\lim_{\Spec\, A \rightarrow X_{DR}} (\B_X(A)-Mod_{\medg_k}^{\mathsf{Perf}}
\longrightarrow \B_X(\s)(A)-Mod_{k(\s)-Mod}^{\mathsf{Perf}}).
$$
We will now prove that, for each $A$, the $\s$-functor
$$
\B_X(A)-Mod_{\medg_k}^{\mathsf{Perf}} \longrightarrow
\B_X(\s)(A)-Mod_{k(\s)-Mod}^{\mathsf{Perf}}
$$
is an equivalence. It is clearly essentially surjective by definition.
As both the source and the target of this functor are rigid symmetric
monoidal $\s$-categories, and the $\s$-functor is symmetric monoidal,
fully faithfulness will follow from the fact that for any object $E
\in \B_X(A)-Mod_{\medg_k}^{\mathsf{Perf}}$ the induced morphism of
spaces
$$
\mathsf{Map}_{\B_X(A)-Mod_{\medg_k}^{\mathsf{Perf}}}(\mathbf{1},E)
\longrightarrow
\mathsf{Map}_{\B_X(\s)(A)-Mod_{k(\s)-Mod}^{\mathsf{Perf}}}(\mathbf{1},E(\s))
$$
is an equivalence. By definition, $E$ is perfect, so is freely
generated over $\B_X(A)$ by its weight $0$ part. By Proposition
\ref{p3} $\B_X(A)$ is free over its part of degree $1$, as a graded
cdga. Therefore, both $\B_X(A)$ and $E$ has no non-trivial negative
weight components. The natural morphism of Ind-objects
$$
E 
\longrightarrow E(\s)
$$
induces an equivalence on realizations $|E| \simeq |E(\s)|\simeq |E|^t$. This
achieves the proof of Corollary, as we have natural identifications
$$
\begin{aligned}
  \mathsf{Map}_{\B_X(A)-Mod_{\medg_k}^{\mathsf{Perf}}}(\mathbf{1},E) &
  \simeq \mathsf{Map}_{\dg_k}(\mathbf{1},|E|) \\
  \mathsf{Map}_{\B_X(\s)(A)-Mod_{k(\s)-Mod}^{\mathsf{Perf}}}(\mathbf{1},E(\s))
  & \simeq \mathsf{Map}_{\dg_k} (\mathbf{1},|E(\s)|).
\end{aligned}
$$
\hfill $\Box$

\begin{rmk} We describe what happens over a
  \emph{reduced} point $f : \Spec\, A_{red} = \Spec\, A
  \longrightarrow X$.  The graded mixed cdga $\D_{X_{DR}}(A)$ reduces
  here to $A$ (with trivial mixed structure and pure weight $0$).
  Therefore, $\B_X(A)$ is here an $A$-linear graded mixed cdga
  together with an augmentation $\B_X(A) \longrightarrow A$ (as a map
  of graded mixed cdgas). Moreover, as a graded cdga, we have
  (Proposition \ref{p3})
$$
\B_X(A)\simeq Sym_A(f^{*}\mathbb{L}_{X}).
$$
This implies that $f^{*}(\mathbb{T}_{X})[-1]$ is endowed with a
natural structure of a dg-Lie algebra over $A$. This is the tangent
Lie algebra of \cite{hen}. Moreover,
$\B_X-Mod_{\medg}^{\mathsf{Perf}}$ is here equivalent to the
$\s$-category of perfect Lie $f^{*}(\mathbb{T}_{X})[-1]$-dg-modules,
and we recover the equivalence
$$
\mathrm{L}_{\mathsf{Perf}}(X_A) \simeq
f^{*}(\mathbb{T}_{X})[-1]-Mod^{\mathsf{Perf}},
$$
between perfect complexes on the formal completion of $X\times \Spec\,
A$ along the graph  
$\Spec\, A \longrightarrow X \times \Spec\, A$, and
perfect $A$-dg-modules with an action of the dg-Lie algebra
$f^{*}(\mathbb{T}_{X})[-1]$ (see \cite{hen}).

The situation over \emph{non-reduced} points is more complicated. In
general, the graded mixed cdga $\B_X(A)$ has no augmentation to $A$,
as the morphism $X_A \longrightarrow \Spec\, A$ might have no section
(e.g. if the point $\Spec\, A \longrightarrow X_{DR}$ does not lift to
$X$ itself).  In particular $\B_X(A)$ cannot be the Chevalley complex
of an $A$-linear dg-Lie algebra anymore. It is, instead, more accurate
to think of $\B_X(A)$ as the Chevalley complex of a dg-Lie
\emph{algebroid} over $\Spec\, A_{red}$, precisely the one given by
the nerve groupoid of the morphism $\Spec\, A_{red} \longrightarrow
X_A$. However, the lack of perfection of the cotangent complexes
involved implies that this dg-Lie algebroid is not the kind of objects
studied in \cite{vez2}. Finally, the action of $\D(A)$ on $\B_{X}(A)$
for a non-reduced cdga $A$, encodes the action of the Grothendieck
connection on the formal derived stack $X_A$.
\end{rmk}

\noindent \textbf{Shifted symplectic structures on derived stacks.} In this paragraph we make a link between \cite{ptvv} and the setting of this paper.\\ Recall that for a derived Artin stack $X$, we have defined $\DR(X/k)$ (Def. \ref{reldrstacks}), a mixed graded cdga over $k$, and for a map $X \to Y$ of derived Artin stacks over $k$, we have $\DR(X/Y)$ (Def. \ref{generalreldrstacks}), again a mixed graded cdga over $k$. The definitions of \cite{ptvv}, can be rephrased as follows.

\begin{df} Let $X$ be a derived Artin stack over $k$, $p \in \mathbb{N}$, and $n \in \mathbb{Z}$.
\begin{itemize}
\item the \emph{space of closed $p$-forms of degree $n$ on $X$} is 
$$
\mathcal{A}^{p,cl}(X,n): =\mathsf{Map}_{\medg_k}(k(p)[-p-n],\DR(X/k)) 
\in \T.
$$
\item The \emph{space of $p$-forms of degree $n$ on $X$} is defined by
$$
\mathcal{A}^{p}(X,n) :=
\mathsf{Map}_{\dg_k}(k[-n],\Gamma (X,\wedge_{\OO_{X}}^{p} \mathbb{L}_{X})) 
\in \T.
$$
\item By Proposition \ref{p1}, for any $\Spec \, A \to X$, there is a natural map $$\mathcal{A}^{p, cl}(\Spec \, A,n)\simeq \mathsf{Map}_{\medg_k}(k(p)[-p-n],\DR(A)) \to \mathsf{Map}_{\dg^{gr}_k}(k(p)[-p-n],\DR(A)) \simeq $$ $$ \simeq \mathsf{Map}_{\dg^{gr}_k}(k(p)[-p-n], Sym_{A}(\mathbb{L}_A[-1]) \simeq \mathsf{Map}_{\dg_k}(k[-p-n], Sym^{p}_{A}(\mathbb{L}_A[-1]) \simeq     \mathcal{A}^{p}(\Spec \, A,n)$$ which induces, by passing to the limit, a map $$\mathcal{A}^{p,cl}(X,n) \to \mathcal{A}^{p}(X,n)$$ called the \emph{underlying-form map}.
\item If the cotangent complex $\mathbb{L}_X \in L_{Perf}(X)$, then the \emph{space of $n$-shifted symplectic structures} on $X$ is $\mathsf{Symp}(X,n)$ is the subspace of $\mathcal{A}^{2,cl}(X,n)$ consisting of non-degenerate forms, i.e. elements of $\pi_0 (\mathcal{A}^{2,cl}(X,n))$ whose underlying form induces an equivalence  $\mathbb{T}_X \to \mathbb{L}_X [n]$.
\end{itemize}  
\end{df}

\

\begin{prop} \label{yeah}
Let $X$ be a derived Artin stack over $k$, and $n \in \mathbb{Z}$. There are canonical equivalences in $\T$

$$\mathsf{Symp}(X,n) \simeq \mathsf{Symp}(\B_X/\D_{X_{DR}}, n) \simeq \mathsf{Symp}(\B_X(\s)/\D_{X_{DR}}(\s),n)$$
where $ \mathsf{Symp}(\B_X/\D_{X_{DR}}, n) $ and $\mathsf{Symp}(\B_X(\s)/\D_{X_{DR}}(\s),n)$ are defined  as in Def. \ref{dsymp}.

\end{prop}

\noindent \textbf{Proof.}  By Cor. \ref{c4} we have equivalences $$\DR(X/k) \simeq \Gamma (X_{DR},\DR(\B_X/\D_{X_{DR}})) \simeq \Gamma (X_{DR},\DR^{t}(\B_X/\D_{X_{DR}})),$$ and by Lemma \ref{l(-s)} we deduce the further equivalence $$\Gamma (X_{DR},\DR^{t}(\B_X/\D_{X_{DR}})) \simeq \Gamma (X_{DR},\DR(\B_X(\s)/\D_{X_{DR}}(\s))).$$
Moreover,  $$\mathsf{Map}_{\medg_\C}(1(p)[-p-n], \DR(\B_X/\D_{X_{DR}})) \simeq \mathsf{Map}_{\medg_k}(k(p)[-p-n],   |\DR(\B_X/\D_{X_{DR}})|), $$ and $$|\DR(\B_X/\D_{X_{DR}})| \simeq \Gamma (X_{DR},\DR(\B_X/\D_{X_{DR}})),$$ where $\C$ is the $\s$-category of functors from $\dAff/X_{DR}$ to  $\medg_k$. Analogously, we get $$\mathsf{Map}_{\medg_{\C'}}(1(p)[-p-n], \DR(\B_X/\D_{X_{DR}})) \simeq \mathsf{Map}_{\medg_k}(k(p)[-p-n],  \Gamma (X_{DR},\DR(\B_X(\s)/\D_{X_{DR}}(\s))))$$ where $\C'$ is the $\s$-category of functors from $\dAff/X_{DR}$ to  $Ind(\medg_k)$. Since the non-degeneracy conditions match, we conclude.

\hfill $\Box$

\section{Shifted Poisson structures and quantization}

\subsection{Shifted Poisson structures: definition and
  examples}\label{sec:3.1}

Let $X$ be a derived Artin stack locally of finite presentation. In
the previous section (see Definition~\ref{d11}) we constructed the
prestack $\D_{X_{DR}}$, the mixed crystalline structure sheaf on $X_{DR}$,
and the prestack $\B_X$ of mixed principal parts, which is a
prestack of graded mixed $\D_{X_{DR}}$-cdgas on $X_{DR}$.  This gives us a
prestack of $\OO_{X_{DR}}$-linear graded $\mathbb{P}_{n+1}$-algebras
$\Pol^t(\B_X/\D_{X_{DR}},n)$ defined in Remark~\ref{directimage} (see
also Corollary~\ref{c3'}):
$$\Pol^t(\B_X/\D_{X_{DR}},n): (\Spec \, A \to X_{DR}) \longmapsto \Pol^t(\B_X(A)/\D_{X_{DR}}(A),n). $$
 We will define $\Pol(X,n)$ as the graded
$\mathbb{P}_{n+1}$-algebra obtained by taking its global sections on
$X_{DR}$:
$$
\Pol(X,n):=\Gamma(X_{DR},\Pol^t(\B_X/\D_{X_{DR}},n)) = \lim_{ \Spec \, A \to X_{DR}} \Pol^t(\B_X/\D_{X_{DR}},n)(A) ,
$$
and call it the \emph{n-shifted polyvectors on} $X$. Note that, by
Theorem \ref{t3}, the underlying graded complex
is 
$$
\bigoplus_{p\geq 0}\Gamma(X,Sym^{p}_{\OO_X}(\mathbb{T}_{X}[-n])),
$$ 
so our
notation $\Pol(X,n)$ should be unambiguous. The reader should just
keep in mind that from now on, unless otherwise stated, we view
$\Pol(X,n)$ with its full structure of graded
$\mathbb{P}_{n+1}$-algebra over $k$. In particular, $\Pol(X,n+1)[n+1]$
is a graded dg-Lie algebra over $k$.

\begin{df}\label{dpoiss}
  In the notations above, the space of \emph{$n$-shifted Poisson
    structures on $X$} is
$$
\mathsf{Poiss}(X,n):=\mathsf{Map}_{\dglie^{gr}_k}(k(2)[-1],\Pol(X,n+1)[n+1]),
$$
where $\dglie_{k}^{gr}$ is the $\s$-category of graded $k$-linear
dg-Lie algebras.
\end{df}

As a direct consequence of this definition and of the main theorem of
\cite{mel}, we get the following important result (see \S 1.5 for the
relation between Tate realization and twists by $k(\s)$).  In the
theorem below, $\D_{X_{DR}}(\s)$ is a prestack of commutative monoids in the
$\s$-category of Ind-objects in graded mixed complexes, $\B_X(\s)$ is
a prestack of commutative monoids in the $\s$-category of Ind-objects
in graded mixed complexes, and we have a canonical morphism
$$
\D_{X_{DR}}(\s) \longrightarrow \B_X(\s).
$$

\begin{thm}\label{tmel}
There is a canonical equivalence of spaces 
$$
\mathsf{Poiss}(X,n) \simeq \mathbb{P}_{n+1}-(\B_X(\s)/\D_{X_{DR}}(\s)),
$$
where the right hand side is the space of
$\mathbb{P}_{n+1}$-structures on $\B_X(\s)$ compatible with its fixed
structure of commutative monoid in the $\s$-category of prestacks of
graded mixed $\D_{X_{DR}}(\s)$-dg-modules.
\end{thm}

\noindent \textbf{Proof.}
Let $\mathcal{M'}$ be the $\s$-category of prestacks on $\dAff/ X_{DR}$ with values in $\mathrm{Ind}(\medg_k)$, and $\mathcal{M}$ be the $\infty$-category of $\D_{X_{DR}}(\infty)$-modules inside $\mathcal{M'}$. Recall that  $\B_{X}(\infty)$ is a commutative monoid in $\mathcal{M}$.
Let us first consider the space $$\mathsf{Map}_{\mathbf{Lie}_{\mathcal{M}}^{gr}} (1_{\mathcal{M}}[-1](2), \Pol^{int}(\B_{X}(\infty)/\D_{X_{DR}}(\infty),n+1)[n+1]).$$ This space is, on one hand, equivalent to $\mathbb{P}_{n+1}-(\B_X(\s)/\D_{X_{DR}}(\s))$ (i.e. the rhs of Thm \ref{tmel}) by Thm \ref{valerio2}, and on the other hand, equivalent to $$\mathsf{Map}_{\mathbf{dgLie}_{k}^{gr}} (k[-1](2), | \Pol^{int}(\B_{X}(\infty)/\D_{X_{DR}}(\infty),n+1)[n+1] | )$$ by the definition of the realization functor $|-|: \mathbf{Lie}_{\mathcal{M}}^{gr} \to \mathbf{dgLie}_{k}^{gr}$ as a right adjoint. We want to show that  the previous space is equivalent to the space $$K:=\mathsf{Map}_{\dglie^{gr}_k}(k(2)[-1],\Gamma(X_{DR},\Pol^t(\B_X/\D_{X_{DR}},n+1))[n+1]).$$ Since the forgetful functor $\mathbb{P}_{n+2}-\cdga^{gr}_{k} \to \dglie^{gr}_{k}$ commutes with limits, we have $$K \simeq  \lim_{\Spec \, A \in \dAff/X_{DR}}\mathsf{Map}_{\dglie^{gr}_k}(k[-1](2), \Pol^t(\B_X(A)/\D_{X_{DR}}(A),n+1))[n+1]). $$ By Lemma \ref{l(-s)} $$\Pol^t(\B_X(A)/\D_{X_{DR}}(A),n+1))[n+1] \simeq | \Pol^{int}(\B_X(\infty)(A)/\D_{X_{DR}}(\s)(A),n+1))[n+1] | \, ,$$ so that $$K \simeq \mathsf{Map}_{\mathbf{dgLie}_{k}^{gr}} (k[-1](2),  \lim_{\Spec \, A \in \dAff/X_{DR}} | \Pol^{int}(\B_{X}(\infty)(A)/\D_{X_{DR}}(\infty)(A),n+1)[n+1] | ).$$ We are thus reduced to proving an equivalence $$ \lim_{\Spec \, A \in \dAff/X_{DR}} | \Pol^{int}(\B_{X}(\infty)(A)/\D_{X_{DR}}(\infty)(A),n+1)[n+1] | \simeq  | \Pol^{int}(\B_{X}(\infty)/\D_{X_{DR}}(\infty),n+1)[n+1] |, $$ and this follows immediately from the general fact that the enriched  Hom in a category of diagrams $I^{op} \to C$, in an enriched symmetric monoidal category $C$, satisfies $$\underline{Hom}_{C^{I^{op}}}(\mathbf{1}_{C^{I^{op}}}, F) \simeq \lim_{x \in I^{op}} \underline{Hom}_{C}(\mathbf{1}_C, F(x)), $$ since the monoidal unit $\mathbf{1}_{C^{I^{op}}}$ is given by the constant $I^{op}$-diagram at $\mathbf{1}_C$. It's enough to apply this to our case $I= \dAff/X_{DR}$, and $C =\mathrm{Ind}(\medg_k)$.

\hfill $\Box$

We describe below what shifted Poisson structures look like on smooth
schemes and classifying stacks of reductive groups. We will see more
advanced examples later on.\\

\noindent \textbf{Smooth schemes.}  Let $X$ be a smooth scheme over
$k$. The $(n+1)$-shifted polyvectors can be sheafified over $X_{Zar}$
in an obvious way, and yield a stack of graded dg-Lie algebras 
$\Pol(X,n+1)[n+1]$ on $X_{Zar}$. As a stack of graded
$\OO_{X}$-dg-modules, this is just 
$\oplus_{p}Sym_{\OO_{X}}(\mathbb{T}_{X}[-1-n])[n+1]$.  As the weight
grading is compatible with the cohomological grading, this stack of
graded dg-Lie algebras is formal, and coincides with the standard
sheaf of shifted polyvectors with its (shifted) Schouten bracket.  By
theorem \ref{tmel}, we know that the space of $n$-shifted Poisson
structures on $X$ as defined in definition \ref{dpoiss} is equivalent
to the space of $\mathbb{P}_{n+1}$-structures on the sheaf
$\OO_{X}$. When $n=0$, this recovers the standard notion of
algebraic Poisson structure on the smooth scheme $X$. \\

\noindent \textbf{Classifying stacks.} Let $G$ be a reductive group
over $k$ with Lie algebra $\mathfrak{g}$. Again, as a graded
$k$-dg-module $\Pol(BG,n+1)$ is
$$
\Pol(BG,n+1)[n+1] \simeq \bigoplus_{p}
Sym^{p}_{k}(\mathfrak{g}[-n])^{G}[n+1].
$$
Again because the weight grading is compatible with the cohomological
grading, $\Pol(BG,n+1)$ is formal as a graded dg-Lie algebra, and the
bracket is here trivial. Using the explicit formulas for the
description of $\mathsf{Map}_{\dglie^{gr}_{k}}(k(2)[-1],-)$, we get
$$
\begin{aligned}
\pi_{0}(\mathsf{Poiss}(BG,2)) & \simeq Sym_{k}^{2}(\mathfrak{g})^{G} \\
\pi_{0}(\mathsf{Poiss}(BG,1)) & \simeq \wedge^{3}_{k}(\mathfrak{g})^{G} \\
\pi_{0}(\mathsf{Poiss}(BG,n)) & \simeq * \qquad \mathrm{if} \; n\neq
1,2. 
\end{aligned}
$$

\subsection{Non-degenerate shifted Poisson structures} \label{ndpoisson}

Let $X$ be a derived Artin stack locally of finite presentation over
$k$, and $p \in \pi_{0}\mathsf{Poiss}(X,n)$ an $n$-shifted Poisson
structure on $X$ in the sense of Definition \ref{dpoiss}. So, $p$ is a
morphism
$$
p : k(2)[-1] \longrightarrow \Pol(X,n+1)[n+1],
$$ 
in the $\s$-category of graded dg-Lie algebras over $k$,
and, in particular, it induces a morphism in the $\s$-category of graded
$k$-dg-modules
$$
p_0 : k(2) \longrightarrow \Pol(X,n+1)[n+2].
$$
Since by Thm. \ref{t3}, $\Pol(X,n+1)[n+2] \simeq \oplus_{p}\Gamma(X,Sym^{p}_{\OO_X}(\mathbb{T}_{X}[-n-1])[n+2])$ in the $\s$-category of graded
$k$-dg-modules, $p_0$ defines an element in 
$$
p_0 \in H^{-n}(X,\Phi_n^{(2)}(\mathbb{T}_{X})),
$$
where  
$$ 
\Phi_{n}^{(2)} (\mathbb{T}_{X}) := \left \{
  \begin{array}{c}
    Sym^{2}_{\OO_{X}}\mathbb{T}_{X}, \,\, \mathrm{if \; n\; is \; odd} \\
    \wedge^{2}_{\OO_{X}}\mathbb{T}_{X},\,\,  \mathrm{if \; n\; is \; even.}
  \end{array}
\right. 
$$
Hence $p_0$ induces, by adjunction, a map $\Theta_{p_0} : \mathbb{L}_{X} \to \mathbb{T}_{X}[-n]$ of perfect complexes.

\begin{df}\label{d12}
With the notations above, the $n$-shifted Poisson structure $p$ is
called \emph{non-degenerate} if the induced map
$$\Theta_{p_0} : \mathbb{L}_{X} \to \mathbb{T}_{X}[-n]$$
is an equivalence of perfect complexes on $X$.
\end{df}

%%%%%%%%%%%%%%%%%%%% DA QUI ?
By Theorem \ref{tmel}, the datum of $p \in \pi_0\mathsf{Poiss}(X,n)$
is equivalent to the datum of a compatible
$\mathbb{P}_{n+1}$-structure on the prestack of Tate principal parts
$\B_X(\s)$ on $X_{DR}$, relative to $\D_{X_{DR}}(\s)$.  The bracket of this
induced $\mathbb{P}_{n+1}$-structure provides a bi-derivation,
relative to $\D_{X_{DR}}(\s)$ ,
$$
[\cdot,\cdot] : \B_X(\s) \otimes_{\D_{X_{DR}}(\s)}\B_X(\s)  \longrightarrow 
\B_X(\s),
$$
and thus a morphism of prestacks of $\B_X(\s)$-modules on $X_{DR}$
$$
\mathbb{T}^{int}_{\B_X(\s)/\D_{X_{DR}}(\s)} \otimes 
\mathbb{T}^{int}_{\B_X(\s)/\D_{X_{DR}}(\s)} \longrightarrow \B_X(\s).
$$
By Corollary~\ref{ct2} and \ref{c4}, we know that $\mathbb{T}^{int}_{\B_X(\s)/\D_{X_{DR}}(\s)}$
can be naturally identified with the image of 
$\mathbb{T}_{X}$ by the equivalence  
$$
\phi_{X}: \mathrm{L}_{\mathsf{Perf}}(X) \simeq
\B_X(\s)-Mod_{k(\s)-Mod}^{\mathsf{Perf}}
$$ 
of Corollary~\ref{c4}.
As a consequence, we obtain the following 

\begin{cor}\label{c5}
  Let $X$ be a derived Artin stack locally of finite presentation over $k$, and $n \in \mathbb{Z}$. An $n$-shifted Poisson structure $p \in \pi_{0}\mathsf{Poiss}(X,n)$
  is non-degenerate in the sense of Definition \ref{d12} if and only
  if the corresponding $\mathbb{P}_{n+1}$-structure on the
  $\D_{X_{DR}}(\s)$-cdga $\B_X(\s)$ is non-degenerate in the sense of
  Definition \ref{dpoissnd}.
\end{cor}
%%%%%%%%%%%%%%%%%%%%% A QUI ?
\

\begin{rmk} We note that a similar corollary applies to the symplectic
  case. More precisely, if $\omega \in \A^{2,cl}(X,n)$ is an
  $n$-shifted closed $2$-form on $X$, it defines a canonical
  $n$-shifted closed $2$-form $\omega'$ on $\B_X(\s)$ relative to
  $\D_{X_{DR}}(\s)$. Then, $\omega$ is non-degenerate if and only if
  $\omega'$ is non-degenerate.
\end{rmk}

\

We are now ready to state the main theorem of this section. Let
$\mathsf{Poiss}^{nd}(X,n)$ the subspace of $\mathsf{Poiss}(X,n)$ of
connected components of non-degenerate $n$-shifted Poisson structures
on $X$. By Corollary~\ref{c5}, we get that the equivalence of Theorem~\ref{tmel} induces an equivalence 
$$\mathsf{Poiss}^{nd}(X,n) \simeq \mathbb{P}_{n+1}^{nd}(\B_X(\s)/\D_{X_{DR}}(\s))$$ in $\T$, where  $\mathbb{P}_{n+1}^{nd}(\B_X(\s)/\D_{X_{DR}}(\s))$
is the space of non-degenerate $\mathbb{P}_{n+1}$-structures on $\B_X(\s)$ relative to
$\D_{X_{DR}}(\s)$. On the other hand, if we take $\mathcal{M}$ to be the category of prestacks over $\dAff_{k}/X_{DR}$ of Ind-objects in mixed graded $\D_{X_{DR}}(\s)$-modules, then Corollary   \ref{cpoissnd2} (2) applied to $A$ equal to the prestack
$\B_X(\s)$ of $\D_{X_{DR}}(\s)$-linear cdgas, provides a morphism of spaces
$$
\psi : \mathbb{P}_{n+1}^{nd}(\B_X(\s)/\D_{X_{DR}}(\s)) \longrightarrow \mathsf{Symp}(\B_X(\s)/\D_{X_{DR}}(\s)).
$$ 
By Proposition \ref{yeah}, $\psi$ then induces a well defined
morphism in $\T$
$$
\psi : \mathsf{Poiss}^{nd}(X,n) \longrightarrow \mathsf{Symp}(X,n).
$$

\begin{thm}\label{t4}
The morphism constructed above
$$
\psi: \mathsf{Poiss}^{nd}(X,n) \longrightarrow \mathsf{Symp}(X,n)
$$
is an equivalence in $\T$.  
\end{thm}

\

\noindent
{\bfseries Note:} A version of this theorem for Deligne-Mumford
derived stacks was recently proven by J. Pridham by a different method
\cite{pridham-compare}. In a later version of \cite{pridham-compare}, which appeared after our paper was put on the arXiv, the author modified his approach in order to treat also the case of derived Artin stacks.\\

This theorem will be a consequence of the following finer statement, 
which implies Theorem \ref{t4} by taking global sections. 

\begin{thm}\label{t4'}
Let $q : X \longrightarrow X_{DR}$ be the natural projection.
Then, the induced morphism
$$
\psi : q_*(\mathsf{Poiss}^{nd}(-,n)) \longrightarrow 
q_*(\mathsf{Symp}(-,n))
$$
is an equivalence of stacks on $\dAff_{k}/X_{DR}$.
\end{thm}

\

\noindent 
Note that the stacks $q_*(\mathsf{Poiss}^{nd}(-,n))$
 and $q_*(\mathsf{Symp}(-,n))$ have values, respectively, 
 $$\mathbb{P}_{n+1}^{nd}(\B_X(\s)(A)/\D_{X_{DR}}(\s)(A))\simeq \mathsf{Poiss}^{nd}(X_A,n),$$ and  $$\mathsf{Symp}(\B_X(\s)(A)\D_{X_{DR}}(\s)(A),n)\simeq \mathsf{Symp}(X_A,n)$$
 on $\Spec \, A \to X_{DR}$.
 
 \ 
 
\noindent
The proof of Theorem \ref{t4'} is rather long and will be given in the next
subsection. Before that, we give some important consequences of
Theorem \ref{t4'}. The following corollary is obtained from the
construction of a canonical symplectic structure on certain mapping
derived stacks (\cite[Theorem~2.5]{ptvv}).

\begin{cor}\label{poissmoduli}
  Let $Y$ be a derived Artin stack locally of finite presentation and
  endowed with an $n$-shifted symplectic structure. Let $X$ be an
  $\OO$-compact and oriented derived stack of dimension $d$ in the
  sense of \emph{\cite{ptvv}}. We assume that the derived stack
  $\Map(X,Y)$ is a Artin derived stack. Then, $\Map(X,Y)$ carries a
  canonical $(n-d)$-shifted Poisson structure.
\end{cor}

\

\noindent
The main context of application of the above corollary is when $Y=BG$
for $G$ a reductive group endowed with a non-degenerate $G$-invariant
scalar product on its Lie algebra $\frak{g}$. The corollary implies
existence of natural shifted Poisson structures on derived moduli
stacks of $G$-bundles on oriented spaces of various sorts: projective
CY manifolds, compact oriented topological manifolds, de Rham shapes
of smooth and projective varieties, etc. (see \cite{ptvv} for a
discussion of these examples).

Theorem \ref{t4'} together with \cite[Theorem~2.12]{ptvv} yield the
following

\begin{cor}\label{cparf}
  The derived stack $\mathbf{Perf}$ of perfect complexes carries a
  natural $2$-shifted Poisson structure.
\end{cor}

\

\noindent
More generally, via Theorem \ref{t4'}, all the examples of shifted
symplectic derived stacks constructed in \cite{ptvv}, admit
corresponding shifted Poisson structures.

\begin{rmk} More generally we expect suitable generalizations of the
  main results in \cite{ptvv} to hold in the (not necessarily
  non-degenerate) shifted Poisson case.  For example, Theorem
  \ref{poissmoduli} should hold when the target is a general
  $n$-shifted Poisson derived stack, yielding a canonical
  $(n-d)$-shifted Poisson structure on $\Map(X,Y)$.  The same result
  should be true for derived intersections of coisotropic maps (see \S
  \ref{coisosection} for a definition of coisotropic structure on a
  map) into a general shifted Poisson Artin derived stack locally of
  finite presentation over $k$.  Both of these problems are currently
  being investigated by V. Melani.
\end{rmk}
 
\subsection{Proof of Theorem \ref{t4'}} \label{dim}

The proof of this theorem will take us some time and will occupy the rest 
 of this section. Before going into the details of the proof, we present 
its basic steps.

\begin{enumerate}
\item The map $\psi$ induces an isomorphisms on all
homotopy sheaves $\pi_i$ for $i>0$.
\item The derived stacks $\mathsf{Poiss}(\B_X(\s)/\D_{X_{DR}}(\s),n)$ and
$\mathsf{Symp}(\B_X(\s)/\D_{X_{DR}}(\s),n)$ are formal derived stacks in the
sense of Definition~\ref{d10-}. 
\item When $A$ is reduced, the $\pi_0$-sheaves of
  $q_*(\mathsf{Poiss}^{nd}(-,n))$ and $q_*(\mathsf{Symp}(-,n))$,
  restricted to $(\Spec\, A)_{Zar}$, can be described in terms of
  pairing and co-pairing on $L_{\s}$-algebras.
\item $(2)$ and $(3)$ imply that the morphism $\psi$ also induces an
  isomorphism on the sheaves $\pi_0$, by reducing to the case of a
  reduced base.
\end{enumerate}

In the remaining subsections, we will give the proof of
Theorem~\ref{t4'}, following the above outline.

\subsubsection{Derived stacks associated with graded dg-Lie and graded
  mixed complexes} \label{associated}

We will discuss here the general form of the derived stacks
$q_*(\mathsf{Poiss}(-,n))$, and $q_*(\mathsf{Symp}(-,n))$ on
$\dAff_k/X_{DR}$ (where $q: X \to X_{DR}$ is the canonical map). We will see that this will easily lead us to proving
that the morphism $\psi$ of Theorem \ref{t4'} induces isomorphisms on
all higher homotopy sheaves. The case of the sheaves $\pi_0$ will
require more work: it will be a consequence of the results of this
subsection together with a Darboux type statement proved in Lemma
\ref{l12}.

\

\medskip

\noindent \textbf{Derived stacks associated with graded dg-Lie algebras.} 
We work over the $\s$-site $\dAff_{k}/Y$, of derived affine schemes
over some base derived stack $Y$ (it will be $Y=X_{DR}$ later on). We assume 
given a stack of $\OO_{Y}$-linear graded dg-Lie algebras $\mathcal{L}$ on 
$\dAff_{k}/Y$. Here we do not assume $\LL$ to be quasi-coherent, so $\LL$ is
a graded dg-Lie algebra inside the $\s$-category $\mathrm{L}(\OO_Y)$ of 
all (not necessarily quasi-coherent) $\OO_Y$-modules on $\dAff_{k}/Y$. 

We define the \emph{stack associated with $\LL$} to be the $\s$-functor
$$\V(\LL) : (\dAff_{k}/Y)^{op} \longrightarrow \T$$
sending $(\Spec\, A \to Y)$ to the space
$$\V(\LL)(A):=\mathsf{Map}_{\dglie_k^{gr}}(k(2)[-1],\LL(A)).$$
Note that as $\LL$ is a stack of graded dg-Lie algebras, the
definition above makes $\V(\LL)$ into a stack of spaces on
$\dAff_{k}/Y$, because $\mathsf{Map}_{\dglie_k^{gr}}(k(2)[-1],-)$
preserves limits.

We are now going to describe the tangent spaces to the derived stack
$\V(\LL)$. For this, let
$$p : k(2)[-1] \longrightarrow \LL(A)$$
be an $A$-point of $\V(\LL)$ that is given by a strict morphism in the usual
(non $\s$-) category of graded dg-Lie algebras over $k$.  
% By passing to strict models, we may always represent the morphism
% $p$ by an strict morphism in the (usual, non-$\s$) category of
% graded dg-Lie algebras over $k$.
Such a morphism $p$ is thus completely characterized by an element $p
\in \LL(A)(2)^{1}$, of cohomological degree $1$ and weight $2$,
satisfying $[p,p]=dp=0$. We associate to such a $p$ a graded mixed
$A$-dg-module $(\LL(A),p)$ as follows. The underlying graded complex
will be $\LL(A)$ together with its cohomological differential, while
the mixed structure is defined to be $[p,-]$. We will write, as usual,
$$
T^{i}_p(\V(\LL)(A)) := \mathrm{hofib}( \V(\LL)(A\oplus A[i])
\longrightarrow \V(\LL)(A); p).
$$ 
The graded mixed complex $(\LL(A),p)$ is then directly related to the
tangent space of the derived stack $\V(\LL)$ at $p$, as shown by the
following lemma.

\begin{lem}\label{l5}
Assume that for all $i$, the natural
morphism 
$$
\LL(A) \otimes_{A}(A\oplus A[i]) \longrightarrow \LL(A\oplus A[i])
$$
is an equivalence of graded dg-Lie algebras. Then, there is a canonical 
equivalence of spaces 
$$
T^{i}_p(\V(\LL))(A) \simeq
\mathsf{Map}_{\medg_k}(k(2)[-1],(\LL(A),p)[i]).
$$
\end{lem}
\noindent \textbf{Proof.} This is a direct check, using the explicit
way  of \cite{mel} to describe elements in  
$\mathsf{Map}_{\dglie^{gr}_{k}}(k(2)[-1],\LL(A))$.  With
such a description, we see that the space of lifts
$$
k(2)[-1] \longrightarrow \LL(A \oplus A[i]) \simeq \LL(A) \oplus
\LL(A)[i]
$$
of the morphism $p$, consists precisely of the data giving a morphism
of graded mixed complexes $k(2)[-1] \longrightarrow (\LL(A),p)[i]$.
Namely, any such a lift is given by a family of elements $(q_0,\ldots,
q_j, \ldots)$, where $q_j$ is an element of cohomological degree
$(1+i)$ and weight $(2+i)$ in $\LL(A)$, such that the equation
$$
[p,q_j]+d(q_{j+1})=0
$$
holds for all $j\geq 0$. \hfill $\Box$ 

\

\medskip

\noindent \textbf{Derived stacks associated with graded mixed
  complexes.} We work in the same context as before, over the
$\s$-site $\dAff_{k}/Y$, but now we start with a stack of
$\OO_Y$-linear graded mixed dg-modules $\E$ on $\dAff_{k}/Y$. We
define the derived stack associated to $\E$ as
$$\V(\E) : (\dAff_{k}/Y)^{op} \longrightarrow \T$$
sending $\Spec\, A \mapsto Y$ to the space
$$\V(\E)(A):=\Map_{\medg_{k}}(k(2)[-1],\E(A)).$$
Let 
$$\omega : k(2)[-1] \longrightarrow \E(A)$$
be an $A$-point of $\V(\E)$, and 
$$
T^{i}_\omega(\V(\E)(A)) := \mathrm{hofib}( \V(\LL)(A\oplus A[i])
\longrightarrow \V(\E)(A); \omega).
$$ 
Lemma \ref{l5} has the following version in this case, 
with a straightforward proof. 

\begin{lem}\label{l6}
  With the notations above, and assuming that for all $i\geq 0$ the
  natural morphism
$$
\E(A) \otimes_{A}(A\oplus A[i]) \longrightarrow \E(A\oplus A[i]) $$ is
an equivalence of graded mixed $A$-dg-modules. Then, there is a
canonical equivalence of spaces
$$ 
T^{i}_\omega(\V(\E))(A) \simeq
\mathsf{Map}_{\medg_k}(k(2)[-1],\E(A)[i]). $$
\end{lem}

\noindent \textbf{Trivial square zero extensions.} Here is an easy
variation on the two previous lemmas \ref{l5} and \ref{l6}.

\begin{lem}\label{l7}
  Let $\LL$ be a graded dg-Lie algebra over $\Spec\, A$ and $p :
  k(2)[-1] \longrightarrow \LL$ a strict morphism of graded dg-Lie
  algebras over $k$. For all $i\in \mathbb{Z}$, we have a natural
  equivalence of derived stacks over $\Spec\, A$
$$
\V(\LL \oplus \LL[i])\times_{\V(\LL)}\Spec\, A \simeq \V((\LL,p)[i]),
$$
where $(\LL,p)$ is the graded mixed dg-module associated to $\LL$ and
$p$.
\end{lem}

\subsubsection{Higher automorphisms groups}\label{higheraut}

In this subsection we use the descriptions of the tangent spaces given
in \S \ref{associated} in order to conclude that the morphism $\psi$
of Theorem \ref{t4'} induces an isomorphisms on all $\pi_{i}$-sheaves,
for $i>0$.

Let $\Spec\, A \longrightarrow X_{DR}$ and let us fix a non-degenerate
$n$-shifted Poisson structure $p$ on the corresponding base change
$X_A$ of $q: X \to X_{DR}$. We already know that $p$ corresponds to a
non-degenerate $\mathbb{P}_{n+1}$-structure on $\B_X(\s)(A)$ relative
to $\D_{X_{DR}}(\s)(A)=\D(A)(\s)$. We first compute the derived stack of
loops of $\mathsf{Poiss}(X,n)$ based at $p$.

We represent $\B_X(\s)(A)$ by a strict $\mathbb{P}_{n+1}$-algebra $C$,
inside the category of $\D(A)(\s)$-modules (note that everything here
is happening inside the category of Ind-objects in $\medg_{k}$). The
Poisson structure $p$ is then given by a strict morphism of graded
dg-Lie algebras
$$
k(2)[-1] \longrightarrow \Pol(C/\D(A)(\s),n+1)[n+1].
$$ 
Moreover, the derived stack 
$q_*(\mathsf{Poiss}(-,n))$is, by definition 
of n-shifted Poisson structures, given by 
$$
q_*(\mathsf{Poiss}(-,n))_{|\Spec\, A}\simeq 
\V(\underline{\Pol}(\B_X(\s)/\D_{X_{DR}}(\s),n+1)[n+1]),
$$
where $\underline{\Pol}(\B_X(\s)/\D_{X_{DR}}(\s),n+1)$ is the sheafified
version of $\Pol(\B_X(\s)/\D_{X_{DR}}(\s),n+1)$ on $X_{DR}$. i.e.
$$
\underline{\Pol}(\B_X(\s)/\D_{X_{DR}}(\s),n+1) : (\Spec\, A \rightarrow
X_{DR}) \mapsto \Pol(\B_X(\s)(A)/\D_{X_{DR}}(\s)(A),n+1).
$$

We consider the based loop stack
$$\Omega_p  q_*(\mathsf{Poiss}(-,n)),$$
which is a derived stack over $\Spec\, A$.  The strict morphism 
$p$ induces a graded mixed structure
on the complex 
$$
\Pol(C/\D(A)(\s),n+1)[n+1]\simeq
\underline{\Pol}(\B_X(\s)/\D_{X_{DR}}(\s),n+1)(A),
$$ 
and we denote the corresponding graded mixed
complex by $(\LL,p)$. 

\begin{lem}\label{l8}
There is a natural equivalence of derived stacks over $\Spec\, A$
$$\Omega_p  \pi_*(\mathsf{Poiss}(-,n)) \simeq 
\V((\LL,p)[-1]).$$
\end{lem}
\noindent \textbf{Proof.}  This is a general fact. If $\LL$ is a
graded dg-Lie over $\Spec\, A$, then there is a natural equivalence
$$Map(S^1,\V(\LL)) \simeq \V(\LL^{S^1}),$$
where $\LL^{S^1}$ is the $S^1$-exponentiation in the $\s$-category 
of graded dg-Lie algebras. As a graded dg-Lie algebra this
exponentiation is equivalent to $\LL \otimes_{k}C^*(S^1)$, where
$C^*(S^1)$ is the cdga of cochains on $S^1$. As 
$C^*(S^1)$ is naturally equivalent to $k\oplus k[-1]$, we find that 
$$
Map(S^1,\V(\LL)) \simeq \V(\LL \oplus \LL[-1]).
$$
The statement now follows from Lemma \ref{l7}.
\hfill $\Box$ 

\

\begin{cor}\label{cl8}
The morphism $\psi$ of Theorem \ref{t4'}
induces an equivalence on based loop stacks, i.e. for each 
$$p : \Spec\, A \longrightarrow q_*(\mathsf{Poiss}^{nd}(X,n)),$$
the induced morphism
$$\Omega_{p}q_*(\mathsf{Poiss}^{nd}(X,n))  \longrightarrow 
\Omega_{\psi(p)}q_*(\mathsf{Symp}(X,n))$$
is an equivalence of derived stacks over $\Spec\, A$.
\end{cor}

\noindent \textbf{Proof.} Lemma \ref{l8} describes 
$\Omega_{p}q_*(\mathsf{Poiss}^{nd}(X,n))$ as  
$\V(\LL,p)[-1]$, where $(\LL,p)$ is the graded mixed complex
given by $\Pol(C/\D(A)(\s),n+1)[n+1]$ with the mixed
structure being $[p,-]$ (and where as above
$C$ is a strict $\mathbb{P}_{n+1}$-algebra over $\D(A)(\s)$
representing $p$). The strict morphism $p$ induces a morphism 
of graded mixed complexes
$$
\phi_{p} : \DR(C/\D(A)(\s)) \longrightarrow
\Pol(C/\D(A)(\s),n+1)[n+1].
$$
But, $p$ being non-degenerate, this morphism is an equivalence. 
By Lemma \ref{l8}, we get
$$
\Omega_{p}q_*(\mathsf{Poiss}^{nd}(X,n)) \simeq
\V(\DR(C/\D(A)(\s))[-1]).
$$
Now, we have a canonical identification (see Lemma \ref{l8})
$$
\V(\DR(C/\D(A)(\s))[-1]) \simeq \Omega_{\psi(p)}
q_*(\mathsf{Symp}(X,n)).
$$
Thus we find an equivalence of derived stacks over $\Spec\, A$
$$
\Omega_{p}q_*(\mathsf{Poiss}^{nd}(X,n)) \simeq
\Omega_{\psi(p)}q_*(\mathsf{Symp}(X,n)),
$$
which can be easily checked to be exactly the morphism induced by the
map $\psi$ in Theorem \ref{t4'}. \hfill $\Box$

\

\begin{cor}\label{cl8'}
The morphism  
$$
\psi : \mathsf{Poiss}^{nd}(X,n) \longrightarrow \mathsf{Symp}(X,n)
$$ 
of Theorem \ref{t4} 
has discrete homotopy fibers.
\end{cor}

So, we are left to proving that $\psi$ of theorem \ref{t4'} induces an isomorphism also on
$\pi_{0}$-sheaves. In order to do this, we will need some preliminary
reductions.

\subsubsection{Infinitesimal theory of shifted Poisson and symplectic
  structures}

%Recall that we have a derived Artin stack $X$ and its projection
%$$q : X \longrightarrow X_{DR}.$$
%We have two derived stacks on the big $\s$-site of derived affine
% schemes over $X_{DR}$, and a morphism $\psi$ between them $$\psi :
% \pi_*(\mathsf{Poiss}^{nd}(-,n)) \longrightarrow
% \pi_*(\mathsf{Symp}(-,n)).$$
 
In this section we prove a result that enables us to reduce Theorem
\ref{t4'} to a question over reduced base rings. Let
$\dAff_{k}^{red}/X_{DR}$ be the sub $\s$-site of $\dAff_{k}/X_{DR}$ consisting of $\Spec\, A
\longrightarrow X_{DR}$ with $A=A_{red}$.  The $\s$-site
$\dAff_{k}^{red}/X_{DR}$ is equivalent to the big $\s$-site $\dAff_{k}^{red}/X_{red}$ of reduced
affine schemes over $X_{red}$, and it comes equipped with an inclusion
$\s$-functor
$$
j : \dAff_{k}^{red}/X_{red} \hookrightarrow \dAff_{k}/X_{DR}.
$$

\noindent The result we need is then the following

\begin{prop}\label{p5}
The morphism 
$$
\psi : q_*(\mathsf{Poiss}^{nd}(-,n)) \longrightarrow
q_*(\mathsf{Symp}(-,n))
$$
of Theorem \ref{t4'} is an equivalence of stacks if and only if the
induced morphism 
$$
j^*\psi : j^* q_*(\mathsf{Poiss}^{nd}(-,n)) \longrightarrow
j^*q_*(\mathsf{Symp}(-,n))
$$
is an equivalence of stacks over $\dAff_{k}^{red}/X_{red}$.
\end{prop}
\noindent \textbf{Proof.} We will use a deformation theory
argument. We have to prove that if $\Spec\, A \longrightarrow X_{DR}$
is an object in $\dAff_{k}/X_{DR}$, then
$$
\psi_{A} : \pi_*(\mathsf{Poiss}(-,n))(A) \longrightarrow
\pi_*(\mathsf{Symp}(-,n))(A)
$$
is an equivalence as soon as 
$$
\psi_{A_{red}} :  \pi_*(\mathsf{Poiss}(-,n))(A_{red}) \longrightarrow 
\pi_*(\mathsf{Symp}(-,n))(A_{red})
$$
is an equivalence.

\begin{lem}\label{l9}
  The two derived stacks $q_*(\mathsf{Symp}(-,n))$ and
  $q_*(\mathsf{Poiss}^{nd}(-,n))$ are nilcomplete and infinitesimally
  cohesive in the sense of Definition~\ref{d10-}.
\end{lem}
\noindent \textit{Proof of the lemma.}  Remind that nilcomplete and
infinitesimally cohesive for $F$ a derived stack over $X_{DR}$, means
the following two conditions.

\begin{enumerate}
\item For all $\Spec\, B \longrightarrow X_{DR} \in \dAff_{k}/X_{DR}$,
  the canonical map
$$
F(B)\longrightarrow \lim_{k}F(B_{\leq k}),
$$
where $B_{\leq k}$ denotes the $k$-th Postnikov truncation of $B$, is
an equivalence in $\T$.
\item For all 
fibered product of almost finite presented $k$-cdgas in non-positive
degrees 
$$
\xymatrix{
B \ar[r] \ar[d] & B_{1} \ar[d] \\
B_{2} \ar[r] & B_{0},}
$$
such that each $\pi_{0}(B_i) \longrightarrow \pi_{0}(B_0)$ is
surjective with nilpotent kernels, and all morphism $\Spec\, B
\longrightarrow X_{DR}$, the induced square
$$
\xymatrix{
F(B) \ar[r] \ar[d] & F(B_{1}) \ar[d] \\
F(B_{2}) \ar[r] & F(B_{0}),}
$$
is cartesian in $\T$.
\end{enumerate}

To prove the lemma we write the two derived stacks
$q_*(\mathsf{Poiss}(-,n))$ and $q_*(\mathsf{Symp}(-,n))$ in the form
(see \S \ref{associated})
$$
q_{*}(\mathsf{Poiss}(-,n))\simeq \V(\LL) \qquad 
q_* (\mathsf{Symp}(-,n)) \simeq \V(\E).
$$  
Here, 
$$
\LL = \underline{\Pol}(\B_{X}(\s)/\D_{X_{DR}}(\s),n+1)[n+1]
$$
is the stack of ($\OO_{X_{DR}}$-linear) graded dg-algebras
of $(n+1)$-shifted polyvectors on $\B_{X}(\s)$ relative
to $\D_{X_{DR}}(\s)$, and 
$$
\E = \underline{\DR}(\B_{X}(\s)/\D_{X_{DR}}(\s))[n+1].
$$
The fact that $\V(\LL)$ and $\V(\E)$ are both nilcomplete and
infinitesimally cohesive will result from the fact that both $\LL$ and
$\E$, considered as stacks of complexes, are themselves nilcomplete
and infinitesimally cohesive. By looking at weight graded components,
this will follow from the fact that the two stacks of complexes on
$X_{DR}$
$$
q_*(Sym^p(\mathbb{T}_{X}[-n-1]) \quad \text{and} \quad
q_*(Sym^p(\mathbb{L}_{X}[-1])
$$
are themselves nilcomplete and infinitesimally cohesive. Let us prove
that this is the case for $$q_*(Sym^p(\mathbb{T}_{X}[-n-1]),$$ the
other case being established by the same argument (since
$\mathbb{T}_{X}$ is perfect).

The stack $q_*(Sym^p(\mathbb{T}_{X}[-n-1])$
can be described explicitly as follows.
Given a map $\Spec\, A \longrightarrow X_{DR}$, we let, as usual, 
$$
X_{A}:=X\times_{X_{DR}}\Spec\, A.
$$
The derived stack $X_{A}$ is the formal completion of $\Spec\, A_{red}
\longrightarrow X\times \Spec\, A$, and it comes
equipped with a natural morphism $u : X_A \longrightarrow X$.\\ The
value of the derived stack $q_*(Sym^p(\mathbb{T}_{X}[-n-1])$ at $A$ is
then
$$
q_*(Sym^p(\mathbb{T}_{X}[-n-1])(A)=\Gamma(X_A,u^{*}(
Sym^p(\mathbb{T}_{X}[-n-1])).
$$
The lemma then follows from the following elementary fact whose proof
we leave to the reader.

\begin{sublem}\label{sl9}
  Let $f : Y \longrightarrow \Spec\, A$ be any derived stack over
  $\Spec\, A$ and $E \in \mathrm{L}_{\mathsf{Perf}}(F)$ be a perfect
  complex over $Y$. Then, the stack of complexes $f_*(E)$ over
  $\Spec\, A$ is nilcomplete and infinitesimally cohesive.
\end{sublem}

\

\noindent The Sub-Lemma achieves the proof of Lemma \ref{l9}. \hfill $\Box$ 

\

\bigskip

We are now able to finish the proof of Proposition \ref{p5}. By Lemma
\ref{l9} and the standard Postnikov decomposition argument, we will be
done once we prove the following statement. Suppose that $\Spec\, A
\longrightarrow X_{DR}$ is such that the induced morphism
$$
\psi_{A} : q_*(\mathsf{Poiss}(-,n))(A) \longrightarrow
q_*(\mathsf{Symp}(-,n))(A)
$$
is an equivalence. Let $M$ be a module of finite type over $A_{red}$,
$i\geq 0$ and $A\oplus M[i]$ the trivial square zero extension of $A$
by $M[i]$. We have to prove that the induced morphism
$$
\psi_{A\oplus M[i]} : 
q_*(\mathsf{Poiss}(-,n))(A\oplus M[i]) \longrightarrow
q_*(\mathsf{Symp}(-,n))(A\oplus M[i])
$$
is again an equivalence. This morphism fibers over the morphism
$\psi_{A}$, which is an equivalence by assumption and it is then
enough to check that the morphism induced on the fibers is an
equivalence. But this is identical to the computation carried out in
subsection \ref{higheraut}. \hfill $\Box$

\subsubsection{Completion of the proof of Theorem \ref{t4'}}

We are now in a position to conclude the proof of Theorem \ref{t4'}.
We consider the morphism
$$
\psi : q_*(\mathsf{Poiss}^{nd}(-,n)) \longrightarrow 
q_*(\mathsf{Symp}(-,n))
$$
of the theorem. This is a morphism of derived stacks over the big
$\s$-site $\dAff_{k}/X_{DR}$, of derived affine schemes over $X_{DR}$,
and, by Corollary~\ref{cl8}, we know that it induces equivalences on all
based loop stacks, hence on all higher homotopy sheaves. It remains to
prove that the induced morphism
$$
\pi_0(q_*(\mathsf{Poiss}^{nd}(-,n))) \longrightarrow 
\pi_0(q_*(\mathsf{Symp}(-,n)))
$$
is an isomorphism of sheaves of sets on $\dAff_{k}/X_{DR}$. By
Proposition~\ref{p5} it is enough to show that the restriction of this
morphism to \emph{reduced} affine schemes over $X_{DR}$ is an
isomorphism of sheaves of sets.

We thus fix a reduced affine scheme $S=\Spec\, A$ with a morphism $S
\longrightarrow X_{DR}$; by definition of $X_{DR}$, this corresponds
to a morphism $u : S \longrightarrow X$. We consider
$$
X_A:=X\times_{X_{DR}}\Spec\, A,
$$
which is naturally identified with the formal completion of the graph
morphism $S \longrightarrow X \times S$ (Proposition~\ref{fibersDR}). We
have natural projection
$$
q^{A} : X_A \longrightarrow S,
$$
and we consider the induced sheaves of sets on the small Zariski site
$S_{Zar}$
$$
\pi_0(q^{A}_{*}(\mathsf{Poiss}^{nd}(-,n))) \quad
\text{and} \quad 
\pi_0(q^{A}_{*}(\mathsf{Symp}(-,n)))
$$
as well as the morphism induced by $\psi$
$$
\psi_{A} : \pi_0(q^{A}_{*}(\mathsf{Poiss}^{nd}(-,n))) \longrightarrow 
\pi_0(q^{A}_{*}(\mathsf{Symp}(-,n))).
$$
We will prove that $\psi_{A}$ is an isomorphism of sheaves on
$S_{Zar}$.  This will be achieved by using certain \emph{minimal
  models} for graded mixed cdgas over $A$ in order to reconstruct
$\mathbb{P}_{n+1}$-structures out of symplectic structures. We start
by discussing such models.

The perfect formal derived stack $X_A$ has a corresponding graded
mixed cdga $\D(X_A)$. Since $A$ is reduced, we note that $\D(X_A)$
here is an $A$-linear graded mixed cdga which, as a non-mixed graded
cdga, is of the form (see Proposition~\ref{p3})
$$
\D(X_A) \simeq Sym_{A}(u^{*}(\mathbb{L}_{X})),
$$
where $u^{*}(\mathbb{L}_{X})$ is the pull back of the cotangent
complex of $X$ along the morphism  
$u : S \longrightarrow X$ (note that
$\mathbb{L}_{(X_{A})_{red}/A}$ is trivial here, so $u^{*}(\mathbb{L}_{X}) \simeq
\mathbb{L}_{(X_A)_{red}/X_{A}}[-1]$).

We introduce a strict model for $\D(X_A)$ as follows.  We choose a
model $L$ for $u^{*}(\mathbb{L}_{X})$ as a bounded complex of projective $A$-modules of finite
rank, and we consider the graded cdga $B:=Sym_A(L)$. We also fix a
strict model $C$ for $\D(X_A)$, as a cofibrant graded mixed cdga.  As
$B$ is a cofibrant graded cdga (and $C$ is automatically fibrant), we
can chose an equivalence of graded cdgas
$$
v : B \longrightarrow C.
$$
The mixed structure on $C$ can be transported to a \emph{weak} mixed
structure on $B$ as follows.  The equivalence $v$ induces a canonical
isomorphism inside the homotopy category
$\mathrm{Ho}(\dglie_{k}^{gr})$ of graded dg-Lie
$$
v : Der^{gr}(B,B) \simeq Der^{gr}(C,C),
$$
where $Der^{gr}$ denotes the graded dg-Lie algebra of graded
derivations.  The mixed structure on $C$ defines a strict morphism of
graded dg-Lie algebras
$$
k(1)[-1] \longrightarrow Der^{gr}(C,C),
$$
which can be transported by the equivalence $v$ into a morphism in
$\mathrm{Ho}(\dglie_{k}^{gr})$
$$
\ell : k(1)[-1] \longrightarrow Der^{gr}(B,B).
$$
The morphism $\ell$ determines the data of 
an $\mathcal{L}_{\s}$-structure on $L^{\vee}[-1]$, that is a family of 
morphisms of complexes of $A$-modules
$$
[\cdot,\cdot]_{i} : L \longrightarrow Sym^{i}_{A}(L),
$$
for $i\geq 2$ satisfying the standard equations (see
e.g. \cite[4.3]{ko}).

We thus consider $L$ equipped with this $\LL_{\s}$-structure. It
induces a Chevalley differential on the commutative cdga $B$ making it
into a mixed cdga. Note that the mixed structure is not strictly
compatible with the weight grading, so $B$ is not a graded mixed cdga
for the Chevalley differential, it is however a filtered mixed cdga
for the natural filtration on $B$ associated to the weight grading. By
taking the total differential, sum of the cohomological and and the
Chevalley differential, we end up with a well defined commutative
$A$-cdga
$$
|B|:=\prod_{i\geq 0}Sym_A^{i}(L).
$$
Note that $|B|$ is also the completed Chevalley complex
$\widehat{C}^{*}(\LL^{\vee}[-1])$ of the $\LL_{\s}$-algebra
$\LL^{\vee}[-1]$.

We define explicit de Rham and polyvector objects, which are
respectively a graded mixed complex and a graded dg-Lie algebra over
$k$, as follows. We let
$$
\DR^{ex}(B):=\bigoplus_{p} |B|\otimes_{A}Sym_{A}^{p}(L[-1]).
$$
The object $\DR^{ex}(B)$ is first of all a graded dg-module over $k$,
by using the total differential sum of the cohomological and Chevalley
differential. Put differently, each $|B|\otimes_{A}Sym_{A}^{p}(L[-1])$
can be identified with the Chevalley complex with coefficient in the
$\LL_{\s}$-$L^{\vee}[-1]$-module $Sym_{A}^{p}(L[-1])$. Moreover,
$\DR^{ex}(B)$ comes equipped with a de Rham differential
$$
dR : |B|\otimes_{A}Sym_{A}^{p}(L[-1]) \longrightarrow
|B|\otimes_{A}Sym_{A}^{p+1}(L[-1]), 
$$
making it into a graded mixed complex over $A$.

The case of polyvectors is treated similarly. We set 
$$
\Pol^{ex}(B,n):=\bigoplus_{p} |B|\otimes_{A}Sym_{A}^{p}(L^{\vee}[-n]).
$$
We consider $\Pol^{ex}(B,n)$ endowed with the total differential, sum
of the cohomological and the Chevalley differential for the
$\LL_{\s}$-$L$-module $L^{\vee}[-n]$.  Moreover, $\Pol^{ex}(B,n)$ is
also equipped with a natural bracket making it into a a graded
$\mathbb{P}_{n+1}$-algebra. In particular,
$\Pol^{ex}(B,n)[n]$ has a natural structure of graded dg-Lie algebra
over $A$. 

\

The next Lemma shows that $\DR^{ex}(B)$ and $\Pol^{ex}(B)$ provide
strict models. 

\begin{lem}\label{l11}
We have natural equivalences of
\begin{enumerate}
\item $\DR^{ex}(B)\simeq \DR(\D_X(A)/A)$
\item $\Pol^{ex}(B) \simeq \Pol^{t}(\D_X(A)/A).$
\end{enumerate}
\end{lem}
\noindent \textbf{Proof.} We consider $k(1)[-1]$ (i.e. $k$ sitting in pure
weight $1$ and in pure cohomological degree $1$), as a graded dg-Lie
algebra with zero differential, and with bracket of weight $0$. Beware
that this is different from the standard convention used in the rest
of the paper.  Note that the graded Lie dg-modules over $k(1)[-1]$ are
exactly graded mixed complexes.

We now consider the canonical quasi-free resolution of $k(1)[-1]$ as
graded dg-Lie algebras $k[f_*] \simeq k(1)[-1]$ described in
\cite{mel}. Here for $i\geq 0$, $f_0$ is a generator of cohomological
degree $-1$ (set $f_i=0$ for $i<0$), pure of weight $(i+1)$. We
moreover impose equations for all $i\geq -1$
$$
df_{i+1} + \frac{1}{2}\sum_{a+b=i}[f_a,f_b]=0.
$$
The graded dg-Lie $k[f_*]$ is a cofibrant model for $k(1)[-1]$. 
The $\s$-category of graded $k(1)[-1]$-dg-modules is thus 
equivalent to the $\s$-category of graded Lie-$k[f_*]$-dg-modules. We denote
this second $\s$-category by
$$
w-\medg:=k[f_*]-\dg_{k}^{gr}.
$$ 
Objects in this second $\s$-category will be simply called \emph{weak
  graded mixed dg-modules}, where \emph{weak} refers here to the mixed
structure. In concrete terms, an object in $w-\medg$ consists of a
graded complex $E=\oplus_{p}E(p)$, together with family of morphism of
complexes (for $i\geq 0$)
$$
\epsilon_i : E(p) \longrightarrow E(p+i+1)[1],
$$
such that 
$$
d\epsilon_{i+1} + \frac{1}{2}\sum_{a+b=i}[\epsilon_a,\epsilon_b]=0
$$
holds inside $\underline{End}^{gr}(E)$, the graded dg-Lie algebra of
graded endomorphisms of $E$.

We can now do differential calculus inside the $\s$-category $\mathcal M := w-\medg$
as we have done in \S 1, and more precisely inside the model category
of weak graded mixed dg-modules. By construction, our cdga
$B=Sym_A(L)$ in the lemma is endowed with a structure of weak graded
mixed cdga over $A$. As such, its de Rham object is precisely given by
our explicit complex $\DR^{ex}(B)$. In the same way, $\Pol^{ex}(B,n)$
identifies with the polyvector objects of $B$ considered as a weak
graded mixed cdga over $A$. Moreover, $B$ is, as a weak graded
$A$-cgda, equivalent to $\D_X(A)$, so the lemma holds simply because
the natural inclusion from graded mixed complexes to weak graded mixed
complexes induces an equivalence of symmetric monoidal model
categories.  \hfill $\Box$\\

Because of Lemma \ref{l11} we can now work with the explicit de Rham
and polyvector objects $\DR^{ex}(B)$ and $\Pol^{ex}(B,n)$ constructed
above. Now, Corollary~\ref{cpoissnd2} provides a morphism of spaces
$$
\psi : \mathsf{Map}_{\dglie_{k}^{gr}}(k(2)[-1],\Pol^{ex}(B,n+1)[n+1])
\longrightarrow \mathsf{Map}_{\medg}(k(2)[-n-2],\DR^{ex}(B)).
$$
This morphism can be stackified over $S_{Zar}$, where $S=\Spec\, A$, by 
sending an open $\Spec\, A' \subset \Spec\, A$ to the map
$$
\xymatrix{
\mathsf{Map}_{\dglie_{k}^{gr}}(k(2)[-1],\Pol^{ex}(B,n+1)[n+1]\otimes_{A}A')
\ar[d]^-{\psi_{A'}} \\
\mathsf{Map}_{\medg}(k(2)[-n-2],\DR^{ex}(B)\otimes_A
A').
}
$$
We already know that this morphism of stacks induces equivalences on
all higher homotopy sheaves, so it only remains to show that it also
induces an isomorphism on the sheaf $\pi_0$.

In order to prove this, we start by the following strictification
result. Recall that a morphism of graded dg-Lie algebras
$$
p : k(2)[-1] \longrightarrow \Pol^{ex}(B,n+1)[n+1]
$$
is non-degenerate if the morphism induced by using the augmentation
$|B| \rightarrow A$
$$
k \rightarrow |B|\otimes_{k}Sym^{2}(L^{\vee}[-n-1])[n] \longrightarrow 
Sym^{2}(L^{\vee}[-n-1])[n] 
$$
induces an equivalence of complexes of $A$-modules
$L\simeq L^{\vee}[-n-2]$. 

The following lemma is an incarnation of the Darboux lemma for shifted
symplectic and shifted Poisson structures. It is inspired by the
Darboux lemma for $\LL_{\s}$-algebras of Costello-Gwilliam
\cite[Lemma~11.2.0.1]{factorization}.

\begin{lem}\label{l12}
  Assume that the complex $L$ is minimal at a point $p\in \Spec\,
  A$, in the sense that its differential vanishes on
  $L\otimes_{A}k(p)$.
\begin{enumerate}
\item Any morphism in the $\s$-category of graded mixed complexes
$$
\omega : k(2)[-2-n] \longrightarrow \DR^{ex}(B),
$$ 
is homotopic to a strict morphism of graded mixed complexes.
\item For any morphism in the $\s$-category of graded dg-Lie algebras
$$
\pi : k(2)[-1] \longrightarrow \Pol^{ex}(B,n+1)[n+1],
$$
which is non-degenerate at $p$, there is a Zariski open neighborood
$\Spec\, A' \subseteq \Spec\, A$ with $p\in \Spec\, A'$, such that
$$
\pi_A' : k(2)[-1] \longrightarrow \Pol^{ex}(B,n+1)[n+1]\otimes_A A'
$$
is homotopic to a strict morphism of graded dg-Lie algebras. 
\end{enumerate}
\end{lem} 
\noindent \textbf{Proof.} $(1)$ The de Rham cohomology of 
the weak graded mixed cdga $B$ is acyclic, because $B$ is a free
cdga. In other words, the natural augmentation
$$
|\DR^{ex}(B)| \longrightarrow A
$$
is an equivalence (where $|\DR^{ex}|$ denotes the standard realization
of the graded mixed complex $\DR^{ex}$).  By using the Hodge
filtration, we find an equivalence of spaces
$$
\mathsf{Map}_{\medg_k}(k(2)[-2-n],\DR^{ex}(B)) \simeq
\mathsf{Map}_{\dg_k}(k,|\DR^{ex}(B)/A|^{\leq 1}[1+n].
$$
To put things differently, any closed 
$2$-form of degree $n$ on $B$ can be represented by an element $\omega'$ 
of the form $dR(\eta)$ for
$\eta \in (|B|\otimes_{k}L)^{n}$, such that there exists  
 $f \in (|B|/A)^{n-1}$ with $d(f)+dR(\eta)=0$. 
In particular, $\omega'$ is an element of cohomological degree $(n+2)$
in $\DR^{ex}(B)$ which is both $d$ and $dR$-closed. It is thus determined
by a strict morphism of graded dg-modules
$$
k(2)[-2-n] \longrightarrow \DR^{ex}(B).
$$

$(2)$ Let $\pi : k(2)[-1] \longrightarrow \Pol^{ex}(B,n+1)[n+1]$ be
non-degenerate at $p$. We represent $\pi$ by a strict morphism of
graded dg-Lie algebras
$$
p : k[f_*] \longrightarrow \Pol^{ex}(B,n+1)[n+1].
$$
As $L$ is minimal at $p$, there is a Zariski 
open $p \in \Spec\, A' \subset \Spec\, A$ such that 
$\pi_A'$ is strictly non-degenerate, i.e. the induced morphism
$$
L\otimes_A A'\simeq L^{\vee}\otimes_A A'[-n-2]
$$
is an isomorphism. By replacing $A$ by $A'$, we can assume that $\pi$ is
in fact strictly non-degenerate over $A$. 

The morphism $\pi$ consists of a family of elements
$$
\{ p_i \in \Pol^{ex}(B,n+1)^{n+2}\}_{i\geq 0},
$$
of cohomological degree $(n+2)$, with $p_i$ pure of weight $(i+2)$,
satisfying the equation
$$
dp_{i+1} + \frac{1}{2}\sum_{a+b=i}[p_a,p_b]=0.
$$
We consider 
$$
p_0 \in |B|\otimes_{k}Sym^{2}(L^{\vee}[-n])^{n+2},
$$
and we write it as $p_0=q+p_0'$, with respect to the direct sum
decomposition coming from $|B|\simeq A \oplus |B|\geq 1$. The element
$q$ of $|B|\otimes_{k}Sym^{2}(L^{\vee}[-n])^{n+2}$ has now constant
coefficients, and satisfies $d(q)=[q,q]=0$. Therefore, it defines a
strict morphism of graded dg-Lie algebras
$$
q : k(2)[-1] \longrightarrow \Pol^{ex}(B,n+1)[n+1],
$$
which is the leading term of $\pi$. 

The strict morphism $q$ defines a strict $\mathbb{P}_{n+1}$-structure
on the weak graded mixed cdga $B$, which is strictly
non-degenerate. It induces, in particular, an isomorphism of graded
objects
$$\phi_q : \DR^{ex}(B) \simeq \Pol^{ex}(B,n+1).$$
The isomorphism $\phi_q$ is moreover an isomorphism of graded mixed
objects where the mixed structure on the right hand side is given by
$[q,-]$. After Tate realization, we obtain a filtered isomorphism of
filtered complexes
$$
|\phi_q|^{t} : |\DR^{ex}(B)|^{t}[n+1] \longrightarrow
|(\Pol^{ex}(B,n+1),[q,-])|^{t}[n+1].
$$
We will only be interested in the part of weight higher than 2, that
is the induced isomorphism
$$
|\phi_q|^{t} : |\DR^{ex}(B)^{\geq 2}|^{t}[n+1] \longrightarrow 
|(\Pol^{ex}(B,n+1)^{\geq 2},[q,-])|^{t}[n+1].
$$
We are now going to modify the filtrations on $|\DR^{ex}(B)|^{t}$ and
$\Pol^{ex}(B,n+1)$ by also taking into account the natural filtration
on $|B|$ induced by the augmentation ideal $I \subset |B|$. We have
$$
|\DR^{ex}(B)|^{t}=\bigoplus_{p}|B|\otimes_{A}Sym^{p}(L[-1]),
$$
and we set 
$$
F^{i}|\DR^{ex}(B)|^{t}:=\bigoplus_{p\geq 0}I^{i-p} \otimes_{A} Sym^{p}(L[-1])
\subset |\DR^{ex}(B)|^{t}.
$$
This defines a descending filtration on $|\DR^{ex}(B)|^{t}$ which is
complete. In the same way, we have
$$
|\Pol^{ex}(B,n+1)|^{t}=\bigoplus_{p}|B|\otimes_{A}Sym^{p}(L^{\vee}[-n])
$$
and we set 
$$
F^{i}|(\Pol^{ex}(B,n+1),[q,-])|^{t}:= \bigoplus_{p\geq 0}I^{i-p}
\otimes_{A} Sym^{p}(L^{\vee}[-n]) \subset
|(\Pol^{ex}(B,n+1),[q,-])|^{t},
$$
which is a complete filtration of $\mathbb{P}_{n+2}$-algebras.
The isomorphism $|\phi_q|^{t}$ constructed above is compatible with 
these filtrations $F*$, and thus induces a filtered isomorphisms
$$
f_1 : F^{3}|\DR^{ex}(B)^{\geq 2}|^{t}[n+1] \longrightarrow 
F^{3}(\Pol^{ex}(B,n+1)^{\geq 2},[q,-])|^{t}[n+1].
$$
Note that we have
$$
F^{3}|\DR^{ex}(B)^{\geq 2}|^{t} = I\otimes_{A}Sym^{2}(L[-1]) \oplus 
\bigoplus_{p\geq 3}|B| \otimes_{A} Sym^{p}(L[-1]),
$$
and as well for the polyvector sides.

By the results of Fiorenza-Manetti \cite[Corollary~4.6]{fido}, the
morphism $f_1$ is the leading term of a filtered
$\LL_{\s}$-isomorphism
$$
f_* : F^{3}|\DR^{ex}(B)^{\geq 2}|^{t}[n+1] \longrightarrow
F^{3}(\Pol^{ex}(B,n+1)^{\geq 2},[q,-])|^{t}[n+1]
$$
of dg-lie algebras, where the lie bracket on the left hand side is
taken to be zero. This $\LL_{\s}$-isomorphism is moreover obtained by
exponentiating an explicit bilinear operator obtained as the
commutator of the cup product of differential forms and of the
contraction by the Poisson bivector $q$.  In particular, the
$\LL_{\s}$-isomorphism $f_*$ induces an isomorphism on the spaces of
Mauer-Cartan elements (here we use that the filtrations are complete,
see \cite{yal})
$$
MC(F^{3}|\DR^{ex}(B)^{\geq 2}|^{t}[n+1]) \simeq 
MC(F^{3}(\Pol^{ex}(B,n+1)^{\geq 2},[q,-])|^{t}[n+1]).
$$
The MC elements on the left hand side are simply 1-cocycles in
$F^{3}|\DR^{ex}(B)^{\geq 2}|^{t}[n+1]$, and thus are closed $2$-forms
of degree $n$ with no constant terms in $Sym^{2}(L[-1]) \subset
|B|\otimes_{A}Sym^{2}(L[-1])$.  Moreover, by the explicit form of the
$\LL_{\s}$-isomorphism $f_*$ we see that closed $2$-forms of degree
$n$ which are strict (i.e. pure of weight $2$), corresponds in
$MC(F^{3}(\Pol^{ex}(B,n+1)^{\geq 2},[q,-])|^{t}[n+1])$ to MC elements
which are also pure of weight $2$.

We are now back to our Poisson structure $\pi$, given by the family of
elements $p_i$. Recall that $q$ is the constant term of $p_0$, let us
write $p_0=q+p_0'$.  The family of elements $p_0',p_1,\dots,p_n,\dots$
defines an element in $MC(F^{3}(\Pol^{ex}(B,n+1)^{\geq
  2},[q,-])|^{t}[n+1])$, denoted by $\pi'$. In other terms we have
$$
d\pi'+[q,\pi']+\frac{1}{2}[\pi',\pi']=0,
$$
which is another way to write the original MC equation satisfied by
$\pi$.  By the $\LL_{\s}$-isomorphism above this element $\pi'$
provides a closed $2$-form $\omega'$. By the point $(1)$ of the lemma
\ref{l12}, $\omega'$ is equivalent to a strict closed $2$-form
$\omega''$, which by the $\LL_{\s}$-isomorphism gives a new MC element
$\pi''$ in $F^{3}(\Pol^{ex}(B,n+1)^{\geq 2},[q,-])|^{t}[n+1])$. This
MC element is pure of weight $2$, so the equation
$$
d\pi''+[q,\pi'']+\frac{1}{2}[\pi'',\pi'']=0
$$
implies that 
$$
d\pi''=0 \qquad [q,\pi'']+\frac{1}{2}[\pi'',\pi''].
$$
In other words, $q+\pi''$ is a strict $\mathbb{P}_{n+1}$-structure on
$B$, which by construction is equivalent to the original structure
$\pi$.  \hfill $\Box$

\

\medskip

\noindent
We come back to our morphism
$$
\psi_{A} : \pi_0(\pi_*(\mathsf{Poiss}^{nd}(-,n))) \longrightarrow 
\pi_0(\pi_*(\mathsf{Symp}(-,n)))
$$
of sheaves on the small Zariski site of $S=Spec\, A$.  Lemma \ref{l12}
$(1)$ easily implies that this morphism has local sections. Indeed,
locally on $S_{Zar}$ any $n$-shifted symplectic structure can be
represented by a strictly non-degenerate strict symplectic structure,
which can be dualized to a strict
$\mathbb{P}_{n+1}$-structure. Moreover, the point $(2)$ of the lemma
\ref{l12} implies that these local sections are locally
surjective. This implies that $\psi_{A}$ is an isomorphism of sheaves
of sets.

This, finally proves Theorem \ref{t4'}.

\subsection{Coisotropic structures}\label{coisosection}

In this Subsection, we propose a notion of \emph{coisotropic
  structure} in the shifted Poisson setting. Our approach here is
based on the so-called \emph{additivity theorem}, a somewhat folkloric
operadic result which should be considered as a Poisson analogue of
Deligne's conjecture as proved in \cite{lualg}. N. Rozenblyum has
comunicated to us a very nice argument for a proof of this
additivity theorem, based on the duality between chiral and
factorization algebras. For future reference we state the additivity
theorem as Theorem \ref{cadd} below. Since the details of Rozenblyum's
argument are not yet publicly available we also give some conceptual
explanations of why such a statement should be true (see Remark
\ref{reasons}). More recently, another proof of this additivity theorem has been given by P. Safronov in \cite{pavel}.

The dg-operad $\mathbb{P}_{n}$ is a Hopf operad, i.e. it comes
equipped with a comultiplication morphism
$$
\nabla : \mathbb{P}_{n} \longrightarrow \mathbb{P}_{n} \otimes_{k}
\mathbb{P}_{n}, 
$$
making it into a cocommutative coalgebra object inside the category of
dg-operads over $k$. We recall that $\mathbb{P}_{n}$ is the homology
of the $E_{n}$-operad (for $n>1$),
%guess it should be $E_n$ here, and n>1
and the morphism $\nabla$ is simply defined by the diagonal morphism
of $E_{n}$.  For our base model category $M$ (as in \S 1.1), this
implies that the category of $\mathbb{P}_{n}$-algebra objects in $M$
has a natural induced symmetric monoidal structure. The tensor
product of two $\mathbb{P}_{n}$-algebras $A$ and $B$ is defined as
being the tensor product in $M$ together with the
$\mathbb{P}_{n}$-structure induced by the following compositions
$$
\xymatrix{ \mathbb{P}_{n}(p) \otimes (A\otimes B)^{\otimes p}
  \ar[r]^-{\nabla} & (\mathbb{P}_{n}(p)\otimes_{k} \mathbb{P}_{n}(p))
  \otimes (A^{\otimes p}\otimes B^{\otimes p}) \ar[r]^-{a\otimes b} &
  A\otimes B,}
$$
where $a$ and $b$ are the $\mathbb{P}_{n}$-structures of $A$ and $B$
respectively.

This construction defines a natural symmetric monoidal structure on
the $\s$-category $\mathbb{P}_{n}-\cdga_{\C}$ for $\C=L(M)$, the
$\s$-category associated to $M$, such that the forgetful $\s$-functor
$$
\mathbb{P}_{n}-\cdga_{\C} \longrightarrow \C
$$
has a natural structure of symmetric monoidal $\s$-functor. In
particular, it makes sense to consider the $\s$-category
$\mathsf{Alg}(\mathbb{P}_{n}-\cdga_{\C})$ of unital and associative
monoids in $\mathbb{P}_{n}-\cdga_{\C}$ (in the sense of
\cite[4.1]{lualg}).

\ 

\noindent
The additivity property of Poisson operads, proven by N. Rozenblyum, and more recently, and independently by P. Safronov (see \cite{pavel}), 
can then be stated as follows. 

\begin{thm}\label{cadd}
  For any $n \geq 1$ and any $\s$-category $\C=L(M)$ as in Section
  \ref{1.1}, there exists an equivalence of $\s$-categories
$$
\mathrm{Dec}_{n+1} : \mathbb{P}_{n+1}-\cdga_{\C} \longrightarrow
\mathsf{Alg}(\mathbb{P}_{n}-\cdga_{\C})
$$
satisfying the following two properties
\begin{enumerate}
\item The $\s$-functor $\mathrm{Dec}_{n+1}$ is natural, with respect
  to symmetric monoidal $\s$-functors, in the variable $\C$.
\item The $\s$-functor $\mathrm{Dec}_{n+1}$ commutes with the
  forgetful $\s$-functors to $\C$.
\end{enumerate}
\end{thm}

\

\begin{rmk}\label{reasons} Before going further, we make a few
  comments about this theorem.  As the additivity for the operad Lie
  is rather straightforward, Theorem \ref{cadd} can be made even more
  precise by requiring the compatibility of $\mathrm{Dec}_{n+1}$ with
  respect to the forgetful maps induced from the inclusion of the
  (shifted) Lie operad inside $\mathbb{P}_{n}$.  We can, moreover,
  require compatibility with respect to the inclusion of the
  commutative algebras operad $\mathbf{Comm}$ into $\mathbb{P}_{n}$,
  as, again, the additivity property for $\mathbf{Comm}$ is
  straightforward. Indeed, the main difficulty in proving
  Theorem~\ref{cadd} is in constructing the $\s$-functor $Dec$. Once
  it is constructed and it is shown to satisfy these various
  compatibilities, it is rather easy to check that it has to be an
  equivalence.

  As a second comment, we should mention that there is an indirect
  proof to this theorem based on formality. Indeed, as we are in
  characteristic zero, we are entitled to chose equivalences of
  dg-operads
$$
\alpha_{n} : \mathbb{E}_n \simeq \mathbb{P}_{n}
$$
for each $n>1$. These equivalences can be actually chosen as
equivalences of Hopf dg-operads. Now, the solution to the Deligne's
conjecture given in \cite{lualg} implies the existence of a natural
equivalence of $\s$-categories
$$
\mathrm{Dec}_{n+1}^{\mathbb{E}} : \mathbb{E}_{n+1}-\cdga_{\C}
\longrightarrow \mathsf{Alg}(\mathbb{E}_{n}-\cdga_{\C}), 
$$
satisfying all the required properties. Then, we can simply define
$\mathrm{Dec}$ by transporting $\mathrm{Dec}_{n+1}^{\mathbb{E}}$
through the equivalences $\alpha_n$ and $\alpha_{n+1}$. This proof
is however not explicit, depends on the choices of the
$\alpha_n$'s, and thus is not very helpful for us. 
%However, it clearly
%indicates that the conjecture is formally true as stated, but we
%strongly believe that there exists a direct solution, totally
%independent of formality.
\end{rmk}

For our purposes, the importance of Theorem~\ref{cadd} is that it
allows for a notion of $\mathbb{P}_{n+1}$-structure on a
\emph{morphism} between cdgas. Indeed, we can consider the
$\s$-category $\mathbb{P}_{(n+1,n)}-\cdga_{\C}$, whose objects consist
of pairs $(A,B)$ where $A$ is an object in
$\mathsf{Alg}(\mathbb{P}_{n}-\cdga_{\C})$ and $B$ is an $A$-module in
$\mathbb{P}_{n}-\cdga_{\C}$. Theorem \ref{cadd} implies that this
$\s$-category $\mathbb{P}_{(n+1,n)}-\cdga_{\C}$ comes equipped with
two forgetful $\s$-functors
$$
\xymatrix{ \mathbb{P}_{n+1}-\cdga_{\C} & \ar[r] \ar[l]
  \mathbb{P}_{(n+1,n)}-\cdga_{\C} & \mathbb{P}_{n}-\cdga_{\C}.}
$$
Moreover, $\mathbb{P}_{(n+1,n)}-\cdga_{\C}$ has a forgetful
$\s$-functor to the $\s$-category $\mathbb{E}_{(1,0)}(\cdga_{\C})$ of
pairs $(A,B)$, where $A \in \mathsf{Alg}(\cdga_{\C})$ and $B$ is an
$A$-module in $\cdga_{\C}$. It is easy to see that the $\s$-category
$\mathbb{E}_{(1,0)}(\cdga_{\C})$ is equivalent to the $\s$-category
$\mathsf{Mor}(\cdga_{\C})$ of morphisms between cdgas in $\C$.  We are
then able to give the following definition of
$\mathbb{P}_{(n+1,n)}$-structure on a given morphism between cdgas.

\begin{df}\label{coisoaff}
  Let $f : A \longrightarrow B$ be a morphism between cdgas in
  $\C$. The \emph{space of $\mathbb{P}_{(n+1,n)}$-structures on $f$}
  is the fiber at $f$ of the forgetful $\s$-functor constructed above
$$\mathbb{P}_{(n+1,n)}-\cdga_{\C} \longrightarrow \mathsf{Mor}(\cdga_{\C}).$$
It will be denoted by 
$$
\mathbb{P}_{(n+1,n)}\mathsf{-Str}(f):=\mathbb{P}_{(n+1,n)}-\cdga_{\C}
\times_{\mathsf{Mor}(\cdga_{\C})}\{f\}.
$$
\end{df}

\

\medskip

\noindent
Note that, for a morphism $f:A \to B$, the space
$\mathbb{P}_{(n+1,n)}\mathsf{-Str}(f)$ has two natural projections
$$
\xymatrix{
\mathbb{P}_{n+1}(A) & \ar[r] \ar[l] \mathbb{P}_{(n+1,n)}\mathsf{-Str}(f)
& \mathbb{P}_{n}(B),}
$$
where $\mathbb{P}_{n+1}(A)$ (respectively, $\mathbb{P}_{n}(B)$)
denotes the space of $\mathbb{P}_{n+1}$-structures
(resp. $\mathbb{P}_{n}$-structures) on the given cdga $A$
(resp. $B$). Loosely speaking, a $\mathbb{P}_{(n+1,n)}$-structure on a
given $f$ consists of a $\mathbb{P}_{n+1}$-structure on $A$, a
$\mathbb{P}_{n}$-structure on $B$, together with some compatibility
data between these structures. These data not only express the fact
that $B$ is an $A$-module in $\mathbb{P}_{n}$-algebras, through the
$\s$-equivalence $\mathrm{Dec}_{n+1}$ of Theorem~\ref{cadd}, but also
that this module structure induces the given morphism $f$
between the corresponding cdgas.

We are now able to use Definition~\ref{coisoaff} in order to introduce
the important notion of \emph{shifted coisotropic structures}. Let $f
: X \longrightarrow Y$ be a morphism of derived Artin stacks locally
of finite presentation over $k$. Recall (Definition~\ref{d11}) that we have
constructed stacks of graded mixed cdgas $\D_{X_{DR}}$ and
$\D_{Y_{DR}}$, the mixed crystalline structure sheaves of,
respectively, $X$ and $Y$. These are stacks of graded mixed cdgas on
$X_{DR}$ and $Y_{DR}$, respectively.  The morphism $f$ obviously
induces a pull-back morphism (where we simply write $f^*$ for
$f_{DR}^*$)
$$
f^{*}(\D_{Y_{DR}}) \longrightarrow \D_{X_{DR}}
$$
which is an equivalence of stacks of graded mixed cdgas over
$X_{DR}$. 

By Definition~\ref{d11}, we also have the mixed principal parts
$\B_X$ and $\B_Y$, of, respectively, $X$ and $Y$, which are stacks of
graded mixed $\D_{X_{DR}}$ and $\D_{Y_{DR}}$
algebras. respectively. The morphism $f$ induces a pull-back map
$$
f^{*}(\B_{Y}) \longrightarrow \B_{X},
$$
which is a morphism of graded mixed $\D_{X_{DR}}$-cdgas.
Over an 
affine $\Spec\, A \longrightarrow X_{DR}$, corresponding to 
a morphism $\Spec\, A_{red} \longrightarrow X$, the morphism
$$
f^{*}(\B_{Y})(A) \longrightarrow \B_{X}(A)
$$
is the image by the $\s$-functor $\D$ of the morphism of perfect
formal derived stacks over $\Spec\, A$
$$
X_A \longrightarrow Y_A,
$$
where $X_A$ is the formal completion of the morphism $\Spec\, A_{red}
\longrightarrow \Spec\, A \times X$, and, similarly, $Y_A$ is the
formal completion of the morphism $\Spec\, A_{red} \longrightarrow
\Spec\, A \times Y$. By tensoring with $k(\s)$, we obtain a morphism
of stacks of Ind-objects in graded mixed $\D_{X_{DR}}$-cdgas on
$X_{DR}$
$$
f^{*}(\B_Y(\s)) \longrightarrow \B_{X}(\s).
$$
If we suppose that $Y$ is endowed with an $n$-shifted Poisson
structure, then $\B_Y(\s)$ comes equipped with a
$\mathbb{P}_{n+1}$-structure, and is thus a stack of graded mixed
$\D_{Y_{DR}}(\s)-\mathbb{P}_{n+1}$-cdgas on $Y_{DR}$.  The pull-back
$f^{*}(\B_Y(\s))$ is therefore a stack of graded mixed
$\D_{X_{DR}}(\s)-\mathbb{P}_{n+1}$-cdgas on $X_{DR}$

\begin{df}\label{coiso}
  Let $f : X \longrightarrow Y$ be a morphism of derived Artin stacks
  locally of finite presentation over $k$, and assume that $Y$ is
  equipped with an $n$-shifted Poisson structure $p$. We denote by
$$f_{\B}^* : f^{*}(\B_Y(\s)) \longrightarrow \B_X(\s)$$
the induced morphism of $\D_{X_{DR}}(\s)$-algebras.  The \emph{space
  of coisotropic structures on $f$ relative to $p$} is defined as
$$
\mathsf{Coiso}(f,p):= \mathbb{P}_{(n+1,n)}\mathsf{-Str}(f_{\B}^*)
\times_{\mathbb{P}_{n+1}\mathsf{-Str}(f_{\B}^{*}
  (\B_Y(\s))/\D_{X_{DR}}(\s))}\{p\}.
$$
\end{df}

\

\noindent
In the above definition, $f^{*}(\B_Y(\s))$ acquires an induced
$\D_{X_{DR}}$-linear $\mathbb{P}_{n+1}$-structure coming from the
$n$-shifted Poisson structure $p$. The datum of a coisotropic
structure on $f$ consists of the datum of a $\D_{X_{DR}}(\s)$-linear
$\mathbb{P}_{n}$-structure on $\B_{X}(\s)$ together with a suitably
compatible structure of module over $f^{*}(\B_Y(\s))$, inside the
$\s$-category of $\D_{X_{DR}}(\s)$-linear graded mixed
$\mathbb{P}_{n}$-algebras on $X_{DR}$. We note, in particular, that a
coisotropic structure on $f:X \to Y$, with $Y$ $n$-shifted Poisson,
trivially induces an $(n-1)$-shifted Poisson structure on the target
$X$ itself.

We end this subsection by the following statement, which is a relative
version of our comparison Theorem \ref{t4'}. We state it now as a
conjecture as we have not yet carried out all the details.

\begin{conj}
  Let $Y$ be a derived Artin stack with an $n$-shifted symplectic
  structure $\omega$, and $f : X \longrightarrow Y$ be a morphism of
  derived Artin stacks. Let $p$ denote the $n$-shifted Poisson
  structure corresponding to $\omega$ via Theorem \ref{t4'}.  Then,
  there exists a natural equivalence of spaces
$$
\mathsf{Lag}(f,\omega) \simeq \mathsf{Coiso}(f,p)^{nd},
$$
between the space of Lagrangian structures on $f$ with respect to
$\omega$ (in the sense of \cite[2.2]{ptvv}) and an appropriate space
of \emph{non-degenerate} coisotropic structures on $f$ relative to
$p$.
\end{conj}

Note that the above conjecture recovers Theorem \ref{t4'}, by taking
$Y=\Spec\, k$ (and $\omega=0$). \\
{\bf Added in proofs.} This conjecture has recently been proven for an 
a priori different notion of coisotropic structure in \cite[Thm. 6.11]{ms}. 
Our notion of coisotropic structure and the one from \cite{ms} have been proven to 
be equivalent in the even more recent preprint \cite{pavel}. 

\begin{rmk} We expect the Lagrangian intersection theorem
  \cite[Theorem~2.9]{ptvv} to extend to shifted Poisson structures as
  follows. Let $(X, p)$ be a $n$-shifted Poisson Artin stack locally
  of finite presentation over $k$, and $f_i: Y_i \to X$, $i=1,2$ be
  maps of derived Artin stacks, each endowed with a coisotropic
  structure relative to $p$. Then, we expect the existence of a
  $(n-1)$-shifted Poisson structure on the derived pullback $Y_1
  \times_{X} Y_2$, suitably compatible with the given coisotropic
  structures on $f_1$ and $f_2$. A first evidence of this result comes
  from \cite{gb}, which basically treats (on the cohomological level)
  the case $n=0$, for $X, Y_1$ and $Y_2$ smooth schemes. \\
  {\bf Added in proofs.} This result has recently been proven for an 
a priori different notion of coisotropic structure in \cite[Thm. 5.2]{ms}. 
Our notion of coisotropic structure and the one from \cite{ms} have then been proven to 
be equivalent in the even more recent preprint \cite{pavel}. 

\end{rmk}

\subsection{Existence of quantization}

We propose here a notion of \emph{deformation quantization} of $n$-shifted Poisson
structures on derived Artin stacks, and prove that they always exist
as soon as $n\neq 0$.  The special case of $n=0$ would require further
investigations and will not be treated in this paper. Also, the more
general, and more delicate, problem of quantization of coisotropic
structures will not be addressed here.

\subsubsection{Deformation quantization problems}

Let $FM_n$ be Fulton-MacPherson's topological operad: given a finite set $I$, the space of 
operations with entries labelled by $I$ is the compactified configuration space $FM_n(I)$ 
of $I$-indexed configurations of points in $\mathbb{R}^n$. The corresponding chain $k$-dg-operad is a model 
for $\mathbb{E}_n$: 
$$
\mathbb{E}_n=C_{-*}(FM_n,k)\,.
$$
If $n\geq1$ one can construct a filtration on $\mathbb{E}_n$ with associated graded being the $n$-shifted Poisson graded $k$-dg-operad $\mathbb{P}_n$. 

Let us start with the case $n=1$, which is a bit special. First of all observe that $\mathbb{E}_1$ is equivalent to its cohomology: 
$$
\mathbb{E}_1\cong H_0(FM_1,k)=\mathbb{A}s\,,
$$
where $\mathbb{A}s$ is the associative $k$-operad. Note that the $k$-module $\mathbb{A}s(I)$ of operations with entries labelled by $I$ is the 
free associative algebra in $I$ generators: 
$$
\mathbb{A}s(I)=<x_i|i\in I>=U\big(Lie(x_i|i\in I)\big)\,,
$$
where $Lie(-)$ denoted the ``free Lie algebra generated by''. 
Being a universal envelopping algebra, $\mathbb{A}s(I)$ is filtered with associated graded 
$$
Sym\big(Lie(x_i|i\in I)\big)=\mathbb{P}_1(I)\,.
$$ 
Applying the Rees construction we get a $k[\hbar]$-dg-operad $\mathbb{BD}_1$.  

\medskip

If $n\geq2$ then we consider the filtration of $\mathbb{E}_n$ given by a (functorial) Postnikov Tower of $FM_n$. 
Applying the Rees construction one gets a $k[\hbar]$-linear-dg-operad $\mathbb{BD}_n$. 
Note that in this case the associated graded is the homology operad of $FM_n$, which is $\mathbb{P}_n$. 

\medskip

So, for $n\geq1$ we get a $k[\hbar]$-dg-operad operad $\mathbb{BD}_n$ such that $\mathbb{BD}_n\otimes_{k[\hbar]}k\cong\mathbb{P}_n$ 
and $\mathbb{BD}_n\otimes_{k[\hbar]}k[\hbar,\hbar^{-1}]\cong \mathbb{E}_n[\hbar,\hbar^{-1}]$. 

\begin{rmk}
The story in the case $n=0$ is even more special than for $n=1$, and is discussed extensively in \cite{factorization}. 
Namely, one introduces a $k[\hbar]$-dg-operad $\mathbb{BD}_0$ defined as follows: its underlying graded operad is $\mathbb{P}_0$ 
and the differential sends the degree $1$ generating operation $\{-,-\}$ to $0$ and the degree $0$ generating operation $-\cdot-$ to $\hbar\{-,-\}$. 
Again, one has that $\mathbb{BD}_0\otimes_{k[\hbar]}k\cong\mathbb{P}_0$ and $\mathbb{BD}_0\otimes_{k[\hbar]}k[\hbar,\hbar^{-1}]\cong \mathbb{E}_0[\hbar,\hbar^{-1}]$. 
\end{rmk}

Whenever one has a $\mathbb{P}_n$-algebra object $A_0$ in a symmetric monoidal $k$-linear (i.e.~$\mathbb{E}_\infty$-monoidal) $\infty$-category $\mathcal C$, 
the {\it deformation quantization problem} reads as follows: 
\begin{Q}[Deformation quantization problem]\label{defquanpro}
Does there exists a $\mathbb{BD}_n$-algebra object $A$ in $\mathcal C\otimes_kk[\hbar]$ such that $A\otimes_{k[\hbar]}k\cong A_0$  as $\mathbb{P}_n$-algebra objects in $\mathcal C$? 
If this happens we say that $A$ is a \emph{deformation quantization} of $A_0$. 
\end{Q} 
Observe that one can also consider the {\it formal deformation quantization problem}, where one replaces $k[\hbar]$ by $k[\![\hbar]\!]$ everywhere. 

\medskip

Let $(X,p)$ be a derived Artin stack locally of finite presentation
over $k$, endowed with an $n$-shifted Poisson structure $p$, with
$n\geq-1$. By Thm. \ref{tmel}, $p$ corresponds to a
$\D_{X_{DR}}(\s)$-linear $\mathbb{P}_{n+1}$-structure on the stack
(over $X_{DR}$) $\B_{X}(\s)$ of Ind-objects in graded mixed
$k$-cdgas. 

In this context one can 
\begin{itemize}
\item either pose the question of the existence of a (formal) deformation quantization of the $\mathbb{P}_{n+1}$-algebra object $\B_{X}(\s)$, whenever $n\geq-1$ (here $\mathcal C$ is the $\infty$-category of prestacks of Ind-objects in mixed graded $\D_{X_{DR}}(\s)$-dg-modules).
\item or pose the question of the existence of a (formal) deformation quantization of the $\mathbb{P}_n$-monoidal $\infty$-category 
$$
\mathsf{Perf}(X)\cong\Gamma(X_{DR},\B_{X}(\s)-\mathsf{Mod}^{\mathsf{Perf}})\,,
$$
whenever $n\geq0$ (here $\mathcal C$ is the $\infty$-category $\mathbf{dgCat}_{/k}$ of dg-categories over $k$ or, equivalently, the $\infty$-category of stable $k$-linear $\infty$-categories).
\end{itemize}
\begin{rmk}
The second variant makes use of the equivalence $\mathrm{Dec}_{n+1}$ from Theorem \ref{cadd} as well as the version of $\mathbb{P}_n$ for $k$-dg-categories from \cite{toenbranes}. 
\end{rmk}

\subsubsection{Solution to the deformation quantization problem}

Let $(X,p)$ be a derived Artin stack locally of finite presentation
over $k$, endowed with an $n$-shifted Poisson structure $p$. 

\begin{thm}
If $n>0$ then there exists a deformation quantization of the stack $\B_{X}(\s)$ of $\D_{X_{DR}}(\s)$-linear graded mixed $\mathbb{P}_{n+1}$-algebras over $X_{DR}$ from Thm. \ref{tmel}. 
\end{thm}
\textbf{Proof.} Since $n>0$ we can choose a \emph{formality} equivalence of
$k$-dg-operads
$$
\alpha_{n+1} : \mathbb{E}_{n+1} \simeq \mathbb{P}_{n+1}\,.
$$
It induces an equivalence $\mathbb{BD}_{n+1}\simeq\mathbb{P}_{n+1}\otimes_kk[\hbar]$ which is the identity mod $\hbar$. 

Therefore one can consider $\B_{X}(\s)\otimes_kk[\hbar]$ as a stack of $\D_{X_{DR}}(\s)$-linear graded mixed $\mathbb{BD}_{n+1}$-algebras on $X_{DR}$. 
It is a deformation quantization of $\B_{X}(\s)$. \hfill $\Box$

\

Putting $\hbar=1$ in the above Theorem (or directly using the formality equivalence $\alpha_{n+1}$) we can consider 
$\B_{X}(\s)$ as a stack of $\D_{X_{DR}}(\s)$-linear graded mixed $\mathbb{E}_{n+1}$-algebras on $X_{DR}$. 

We denote by $\B_{X}(\s)-\mathsf{Mod}^{\mathsf{Perf}}_{p}$ the stack
of perfect $\B_{X}(\s)$-modules on $X_{DR}$, where $\B_{X}(\s)$ is viewed as a stack of
$\D_{X_{DR}}(\s)$-linear graded mixed $\mathbb{E}_{n+1}$-algebras on
$X_{DR}$. By \cite[5.1.2.2 and 5.1.2.8]{lualg},
$\B_{X}(\s)-\mathsf{Mod}^{\mathsf{Perf}}_{p}$ is endowed with the
structure of a stack of $\mathbb{E}_{n}$-monoidal $\s$-categories on
$X_{DR}$. We denote this stack by
$\B_{X}(\s)-Mod^{\mathsf{Perf}}_{\mathbb{E}_n, \, p}$.

\begin{df}\label{dquant}
  With the notation above, and $n>0$, \emph{the quantization of $X$ with respect to $p$} 
  is the $\mathbb{E}_{n}$-monoidal $\s$-category
  $$
  \mathsf{Perf}(X,p):=\Gamma(X_{DR},
  \B_{X}(\s)-\mathsf{Mod}^{\mathsf{Perf}}_{\mathbb{E}_n, 
  \, p}).
  $$
\end{df}

\begin{rmk} 
  Technically speaking $\mathsf{Perf}(X,p)$
  also depends on the choice of the formality equivalence
  $\alpha_{n+1}$. However, as $\alpha_{n+1}$ can be chosen
  independently of all $X$ and $p$, we simply assume that such a
  choice has been made and will omit to mention it in our notation.
\end{rmk}

Now observe that the underlying $\s$-category of $\mathsf{Perf}(X,p)$
is exactly the category of sections  
$\Gamma(X_{DR},\B_{X}(\s)-\mathsf{Mod}^{\mathsf{Perf}})$
which coincides with the $\s$-category $\mathsf{Perf}(X)$ of perfect
$\OO_X$-modules on $X$. In other words, $\mathsf{Perf}(X,p)$ consists 
of the \emph{datum of a $\mathbb{E}_{n}$-monoidal 
structure on} $\mathsf{Perf}(X)$.

This $\mathbb{E}_{n}$-monoidal structure can also be understood as a
\emph{deformation} of the standard symmetric monoidal
(i.e. $\mathbb{E}_{\s}$-) structure on $\mathsf{Perf}(X)$ by
considering the family, parametrized by the affine line
$\mathbb{A}_{k}^1$, of $n$-shifted
Poisson structure $\hbar \cdot p$, with $\hbar \in k$. 

\begin{conj}
The quantization of $X$ with respect to $p$ is indeed a deformation quantization of the $\mathbb {P}_n$-monoidal structure on $\mathsf{Perf}(X)$. 
\end{conj}
This conjecture is actually a consequence of a result that has been announced by N.~Rozenblyum, stating that $\mathrm{Dec}_{n+1}$ can be lifted to an equivalence of $\infty$-categories 
$$
\mathbb{BD}_{n+1}-\cdga_{\C} \longrightarrow
\mathsf{Alg}(\mathbb{BD}_{n}-\cdga_{\C})
$$
such that one recovers the Dunn-Lurie additivity from \cite[5.1.2.2]{lualg} by evaluating at $\hbar=1$ and the equivalence $\mathrm{Dec}_{n+1}$ by evaluating at $\hbar=0$. 

\

\noindent \textbf{Quantization for negative $n$.} 

Let us now treat the case of a $n$-shifted Poisson structure $p$ on
$X$, with $n<0$. Let $\hbar_{2n}$ be a formal variable of cohomological
degree $2n$, and consider
$$
\B_{X}(\s)[\hbar_{2n}], 
$$ 
which is now a stack, on $X_{DR}$, of Ind-objects in graded
$k(\s)[\hbar_{2n}]$-linear mixed cdgas. It comes equipped with a
natural $k(\s)[\hbar_{2n}]$-linear $\mathbb{P}_{1-n}$-structure,
induced by $\hbar_{2n}\cdot p$. Since $n<0$, we are back to the situation of positively shifted Poisson structures. 

Namely, we may choose a formality equivalence of $k$-dg-operads
$$
\alpha_{1-n} : \mathbb{E}_{1-n} \simeq \mathbb{P}_{1-n},
$$
and thus view $\B_{X}(\s)[\hbar_{2n}]$ as a an
$k(\s)[\hbar_{2n}]$-linear $\mathbb{E}_{1-n}$-algebra. Again by using
\cite[5.1.2.2 and 5.1.2.8]{lualg}, the associated stack
$\B_{X}(\s)[\hbar_{2n}]-\mathsf{Mod}^{\mathsf{Perf}}_{p}$ of perfect
$\B_{X}(\s)[\hbar_{2n}]$-modules comes equipped with a natural
$\mathbb{E}_{-n}$-monoidal structure which will be denoted by
$$
\B_{X}(\s)[\hbar_{2n}]-\mathsf{Mod}^{\mathsf{Perf}}_{\mathbb{E}_{-n},\,p}.
$$

\begin{df}\label{dquant2}
  With the notation above, and $n<0$, the \emph{quantization of $X$
    with respect to $p$} is the $\mathbb{E}_{-n}$-monoidal
  $\s$-category
$$
\mathsf{Perf}(X,p):=\Gamma(X_{DR},
\B_{X}(\s)[\hbar_{2n}]-\mathsf{Mod}^{\mathsf{Perf}}_{\mathbb{E}_{-n},\,p}).
$$
\end{df}

\

\noindent
Now observe that, by construction, the underlying $\s$-category of
$\mathsf{Perf}(X,p)$
is $$\mathsf{Perf}(X)\otimes_{k}k[\hbar_{2n}]=:\mathsf{Perf}(X)[\hbar_{2n}].$$
The quantization of Definition \ref{dquant2} consists then of the
\emph{datum of a $\mathbb{E}_{-n}$-monoidal structure on}
$\mathsf{Perf}(X)[\hbar_{2n}]$. As above, such a quantization can be
considered as a deformation of the standard symmetric monoidal
(i.e. $\mathbb{E}_{\s}$-) structure on
$\mathsf{Perf}(X)[\hbar_{2n}]$. Note that this standard symmetric
monoidal structure on $\mathsf{Perf}(X)[\hbar_{2n}]$ recovers the
standard symmetric monoidal structure on $\mathsf{Perf}(X)$ after base
change along the canonical map $k[\hbar_{2n}] \rightarrow k$.

\subsection{Examples of quantizations}

\subsubsection{Quantization formally at a point}\label{sec:quantpoint}

Let $X$ be an Artin derived stack and $x:*:=\Spec\, k\to X$ a closed
point. We start with an obvious observation.
\begin{lem}
  $\Pol(\widehat{X}_x,n+1)=\Pol^t(\mathcal
  B_X/\mathbb{D}_{X_{DR}},n+1)(x)$
\end{lem}
\textbf{Proof.} Observe that $(\widehat{X}_x)_{DR}=\Spec\, k$. Therefore, 
$$
\Pol(\widehat{X}_x,n+1)=\Pol^t(\mathcal
B_{\widehat{X}_x},n+1)=\Pol^t(\mathbb{D}_{\widehat{X}_x},n+1)=\Pol^t(\mathcal
B_X/\mathbb{D}_{X_{DR}},n+1)(x)\,.
$$
\ \hfill $\Box$

\

In particular, we get a dg-lie algebra morphism 
$$
\Pol(X,n+1)=\Gamma\big(X_{DR},\Pol^t(\mathcal
B_X/\mathbb{D}_{X_{DR}},n+1)\big)\to \Pol^t(\mathcal
B_X/\mathbb{D}_{X_{DR}},n+1)(x)=\Pol(\widehat{X}_x,n+1)\,.
$$
Therefore, any $n$-shifted Poisson structure on $X$ induces an
$n$-shifted Poisson structure on the formal completion $\widehat{X}_x$
at $x$.

Recall from Theorem \ref{t1} that, as a (non-mixed) graded cdga over
$k$, $\mathcal B_{\widehat{X}_x}$ is equivalent to
$$
Sym(\mathbb{L}_{*/\widehat{X}_x}[-1])\cong
Sym(x^*\mathbb{L}_{\widehat{X}_x})\cong Sym(x^*\mathbb{L}_X)\,.
$$
We therefore get a graded mixed $\mathbb P_{n+1}$-algebra structure on
$Sym(x^*\mathbb{L}_X)$, whose underlying graded mixed cdgas is the one
from $\mathcal B_{\widehat{X}_x}$.  After a choice of formality
$\alpha_{n+1}$, we get a graded mixed $\mathbb E_{n+1}$-structure on
$Sym(x^*\mathbb{L}_X)$ whenever $n>0$.

\medskip

We would like to make the above $\mathbb E_{n+1}$-structure on
$Sym(x^*\mathbb{L}_X)$ rather explicit for a large class of examples.

Before doing so, let us recall very briefly Kontsevich's construction
of an equivalence $\alpha_{n+1}$ \cite{KoOpMo}.  Let $FM_{n+1}$ be the
Fulton-MacPherson operad of compactified configuration spaces of
points in $\mathbb{R}^{n+1}$ (which is a topological model for the
operad $\mathbb E_{n+1}$: $\mathbb{E}_{n+1}=C_{-*}(FM_{n+1},k)$ and
$\mathbb{P}_{n+1}=H_{-*}(FM_{n+1},k)$).  The equivalence
$\alpha_{n+1}$ comes from a zig-zag of explicit equivalences, which
can be easily understood on the dual cooperads:
$$
C^*(FM_{n+1},k)\longleftarrow Graphs_{n+1}\longrightarrow
H^*(FM_{n+1},k)\,. 
$$
Here $Graphs_{n+1}$ is a certain cooperad in quasi-free cdgas:
generators of $Graphs_{n+1}(I)$ are certain connected graphs, with
external and internal vertices, having their external vertices labeled
by $I$. The morphism $Graphs_{n+1}(I)\to H^*(FM_{n+1}(I),k)$ sends
\begin{itemize}
\item the connected graph without internal vertex and linking $i$ to
  $j$, to the pull-back $a_{ij}$ of the fundamental class of
  $FM_{n+1}(2)\cong S^n$ along the map $FM_{n+1}(I)\to FM_{n+1}(2)$
  that forgets all points but $i$ and $j$.
\item all other generators, to zero. 
\end{itemize}
The morphism $Graphs_{n+1}(I)\to C^*(FM_{n+1}(I),k)$ is transcendental
in nature: it sends a graph $\Gamma$ to the form
$$
\int_{\textrm{internal vertices}}\bigwedge_{\textrm{edges
  }(i,j)}\omega_{ij}\,, 
$$
where $\omega_{ij}$ is the pull-back of the $SO(n+1)$-invariant volume
form on $FM_{n+1}(2)\cong S^n$ along the map $FM_{n+1}(I)\to
FM_{n+1}(2)$ that forgets all points but $i$ and $j$.

Let us now chose a minimal model $L$ for $x^*\mathbb{L}_X$.  As we
already observed, we get a weak mixed structure on the graded cdga
$B:=Sym(L)$, that is equivalent to $\mathcal B_{\widehat{X}_x}$.  This
weak mixed structure induces (and is actually equivalent to) the data
of an $\mathcal L_\infty$-structure on $L^{\vee}[-1]$.

If we further assume that the $n$-shifted Poisson structure on $X$ we
started with is non-degenerate at $x$, then Lemma \ref{l12} tells us
that the induced Poisson structure on $\widehat{X}_x$ is homotopic to
a strict morphism of graded dg-lie algebras
$$
k(2)[-1]\longrightarrow \Pol^{ex}(B,n+1)\,.
$$
Let us assume for simplicity that the strict degree $-n$ Poisson
bracket $q$ we get that way on $B$ is constant (meaning, as in the
proof of Lemma \ref{l12}, that $q$ is a degree $n+2$ element in
$Sym^2(L^\vee[-n-1])\subset |B|\otimes Sym^2(L^\vee[-n-1])$).  In this
case the corresponding strict $\mathbb{P}_{n+1}$-structure on $B$ has
the following remarkable description: structure maps
$$
B^{\otimes I}\longrightarrow B\otimes H^*(FM_{n+1}(I),k)
$$
are given by 
$$
B^{\otimes I}\overset{exp(a)}{\longrightarrow}B^{\otimes I}\otimes
H^*(FM_{n+1}(I),k)\overset{m\otimes id}{\longrightarrow}B\otimes
H^*(FM_{n+1}(I),k)\,,
$$
where $m$ is the multiplication on $B$ and 
$$
a:=\sum_{i\neq j}\partial_p^{i,j}\otimes a_{ij}\,.
$$
It can be checked that this formula lifts to graphs without
modification whenever $p$ is constant, and thus the induced
$\mathbb{E}_{n+1}$-structure on $B$ can be described by structure maps
$$
B^{\otimes I}\overset{exp(A)}{\longrightarrow}B^{\otimes
  I}\widehat\otimes C^*(FM_{n+1}(I),k)\overset{m\otimes
  id}{\longrightarrow}B\widehat\otimes C^*(FM_{n+1}(I),k)\,,
$$
where 
\begin{equation}\label{eq:n-weyl}
A:=\sum_{i\neq j}\partial_p^{i,j}\otimes\omega_{ij}\,.
\end{equation}
Of course $A$ is a formal sum, but when evaluated on chains it becomes
finite and makes perfect sense.

We recover that way the Weyl $n$-algebras that were recently defined
by Markarian (see \cite{markarian}).

\subsubsection{Quantization of $BG$}

Let now $X=BG$, where $G$ is an affine group scheme, and observe that
$X_{DR}=B(G_{DR})$. Let $x:*\to BG$ be the classifying map of the unit
$e:*\to G$. We have a fiber sequence of groups
$$
\widehat{G}_e\longrightarrow G\longrightarrow G_{DR}\,,
$$
so that $\widehat{BG}_x\simeq B(\widehat{G}_e)$. 

We have already seen in the previous \S~that the pull-back of
$\mathcal B_X$ along $x_{DR}:*\to BG_{DR}$ is $\mathcal
B_{\widehat{X}_x}$. Therefore we get that the symmetric monoidal
$\infty$-category
$$
\mathsf{Perf}(BG)\simeq \mathcal
B_X-Mod_{\epsilon-\mathbf{dg}^{gr}}^{Perf}
$$
is equivalent to the symmetric monoidal $\infty$-category of
$G_{DR}$-equivariant objects in
$$
\mathcal
B_{\widehat{X}_x}-Mod_{\epsilon-\mathbf{dg}^{gr}}^{Perf}
\simeq\mathsf{Perf}(B\widehat{G}_e)\,. 
$$

Therefore, given an $n$-shifted Poisson structure $p$ on $BG$, the
quantization we get is completely determined by the
$G_{DR}$-equivariant graded mixed $\mathbb{E}_{n+1}$-algebra structure
on $\mathcal B_{\widehat{X}_x}$ obtained from the equivalence
$\alpha_{n+1}:\mathbb{P}_{n+1}\simeq\mathbb{E}_{n+1}$.  This shall
have a fairly explicit description as $\mathcal
B_{\widehat{X}_x}\simeq\mathbb{D}(B\widehat{G}_e)$ is equivalent to
$Sym(x^*\mathbb{L}_{BG})\simeq Sym(\mathfrak g^\vee[-1])$ as a graded
(non-mixed) cdga, where $\mathfrak g:=e^*T_G$.

\medskip

Before going further, let us prove that $\mathbb{D}(B\widehat{G}_e)$
is actually equivalent to the Chevalley-Eilenberg graded mixed cdga of
the Lie algebra $\mathfrak g$.  The proof mainly goes in two steps:
\begin{itemize}
\item we first prove that equivalences classes graded mixed cdga
  structures on $Sym(V^\vee[-1])$, for $V$ a discrete projective
  $k$-module of finite type, are in bijection with isomorphisms
  classes of strict Lie algebra structures on $V$.
\item we then show that the Lie algebra structure on $\mathfrak g$
  coming from the above mixed structure on $Sym(\mathfrak g^\vee[-1])$
  is isomorphic to the standard Lie algebra structure on $\mathfrak
  g=e^*T_G$.
\end{itemize}

For $C \in \cdga^{gr}_{k}$, we will denote by $\mecdga_{k} (C)$ the
fiber product
$$
\xymatrix{ \mecdga_{k} (C) \ar[r] \ar[d] & \mecdga_{k} \ar[d]^-{U} \\
  \{* \} \ar[r]_-{C} & \cdga^{gr}_{k} }
$$
where $U$ denotes the forgetful functor, and $C$ the given graded cdga
structure. We then define $\epsilon-\mathrm{cdga}_{k}^{gr}(C) := \pi_0
(\mecdga_{k}{C})$.  For $V$ a $k$-module, we write
$\mathrm{LieAlg}^{\mathrm{str}}(V)$ for the set of isomorphism classes
of Lie algebra structures on $V$.

\begin{prop}\label{Lieb}
Let $V$ be a discrete projective $k$-module of finite type. 
\begin{enumerate}
\item for $B \in \epsilon-\mathrm{cdga}_{k}^{gr}(\mathrm{Sym}
  (V^{\vee} [-1]))$, let $H(B)$ be the graded mixed cdga defined by
$$
H(B)(p):= H^p (B(p))[-p]\, , \, p\in \mathbb{Z}
$$ 
with mixed differential induced by $H^*(\epsilon_B)$. Then there is a
canonical equivalence $B \simeq H(B)$ in $\mecdga_{k}$ (i.e. $B$ is
\emph{formal} as a graded \emph{mixed} cdga).
\item there is a bijection 
$$
\mathrm{Lie}: \epsilon-\mathrm{cdga}_{k}^{gr}(\mathrm{Sym} (V^{\vee}
[-1])) \longrightarrow \mathrm{LieAlg}^{\mathrm{str}}(V)
$$
whose inverse
$$
\mathrm{Mix}: \mathrm{LieAlg}^{\mathrm{str}}(V) \longrightarrow
\epsilon-\mathrm{cdga}_{k}^{gr}(\mathrm{Sym} (V^{\vee} [-1]))
$$
is given by the (strict) Chevalley-Eilenberg construction.
\end{enumerate}
\end{prop}

\noindent \textbf{Proof.} (1) Let $B \in
\epsilon-\mathrm{cdga}_{k}^{gr}(\mathrm{Sym} (V^{\vee} [-1]))$, and
$u: B \simeq \mathrm{Sym} (V^{\vee} [-1])$ an equivalence in
$\cdga^{gr}_{k}$.  Since the differential in $\mathrm{Sym} (V^{\vee}
[-1])$ is zero, $\mathrm{Sym} (V^{\vee} [-1])$ is a formal graded
cdga, and we have
$$
H^*(B(p))=0\, , \,\, \textrm{for any}\, p<0 ,
$$
$$
H^i (B(p))=0\, , \,\,\, \textrm{for any} \, p\geq 0\, , i \neq p ,
$$
and $u$ induces $k$-module isomorphisms 
$$
H^p (B(p)) \simeq \wedge^{p} V^{\vee} \, , \,\,\, \textrm{for any} \,
p\geq 0\,.
$$
We may also consider $\tau_{\leq}(B)$ as 
$$
\tau_{\leq}(B)(p):= \tau_{\leq p}(B(p))\, , \, p\in \mathbb{Z}\,,
$$ 
where $\tau_{\leq p}(E)$ denotes the good truncation of a dg-module
$E$.  One can check that the graded mixed cdga structure on $B$
induces a graded mixed cdga structure on $\tau_{\leq}(B)$, and that
the obvious dg-modules maps define a strict diagram of graded mixed
cdgas
$$
\xymatrix{B & \tau_{\leq}(B) \ar[r]^-{g} \ar[l]_-{h} & H(B). }
$$ 
By our computation of $H(B)$ above, we deduce that both $g$ and $h$
are graded quasi-isomorphisms, hence that $B$ is equivalent to $H(B)$
in $\mecdga_k,$ i.e. any $B\in
\epsilon-\mathrm{cdga}_{k}^{gr}(\mathrm{Sym} (V^{\vee} [-1]))$ is
formal as a graded \emph{mixed} cdga.

(2) For $B$ as above, we now consider the mixed differential
$\epsilon_1 : B(1) \to B(2)[1]$, for $p\geq 0$.  It induces on $H^1$ a
map
$$
V^{\vee} \simeq H^1 (B(1)) \to H^{2}(B(2)) \simeq \wedge^2 V^{\vee}
$$ 
whose dual 
$$
\langle\,,\,\rangle_{u}: \wedge^{2}V \to V
$$ 
can easily be checked to define a Lie bracket on $V$.  If $B' \in
\epsilon-\mathrm{cdga}_{k}^{gr}(\mathrm{Sym} (V^{\vee} [-1]))$, $u':
B' \simeq \mathrm{Sym} (V^{\vee} [-1])$ an equivalence in
$\cdga^{gr}_{k}$, and $B\simeq B'$ in
$\epsilon-\mathbf{cdga}_{k}^{gr}$, then $\langle\, , \, \rangle_{u}$
and $\langle\, , \, \rangle_{u'}$ defines the same element in
$\mathrm{LieAlg}^{\mathrm{str}}(V)$.  Thus, we have a well defined map
$$
\mathrm{Lie} : \epsilon-\mathrm{cdga}_{k}^{gr}(\mathrm{Sym} (V^{\vee}
[-1])) \longrightarrow \mathrm{LieAlg}^{\mathrm{str}}(V)\,.
$$ 

Let us show that $\mathrm{Lie}$ is injective.  Let us recall
(e.g. \cite[Lemma 2.2]{xu}) that the map $\mathrm{Lie}^{str}$ sending
a strict graded mixed cdga structure $\{\epsilon_p : \wedge^p
V^{\vee}[-p] \to \wedge^{p+1}V^{\vee}[-p] \}$ to $(\epsilon_1)^{\vee}:
\wedge^2 V \to V$ defines a bijection between strict isomorphism
classes of (strict) graded mixed cdga structures on $\mathrm{Sym}
(V^{\vee} [-1])$ and $\mathrm{LieAlg}^{\mathrm{str}}(V)$.  We denote
its inverse by $\mathrm{str Mix}$.  Let $B$ and $B'$ be such that
$\mathrm{Lie}(B)= \mathrm{Lie}(B')$. By definition of $\mathrm{Lie}$,
and the bijection just mentioned, we have strict isomorphisms of
graded mixed cdgas
$$
H(B) \simeq (\mathrm{Sym} (V^{\vee} [-1]),
\mathrm{strMix}(\mathrm{Lie}(B)))$$ $$H(B') \simeq (\mathrm{Sym}
(V^{\vee} [-1]), \mathrm{strMix}(\mathrm{Lie}(B')))\,.
$$ 
But $\mathrm{Lie}(B)= \mathrm{Lie}(B')$, so we get a strict
isomorphism of graded mixed cdgas $H(B) \simeq H(B')$. Since we have
proved that $B$ and $H(B)$ (respectively, $B'$ and $H(B')$) are
equivalent as graded mixed cdgas, we conclude that $\mathrm{Lie}$ is
injective.

\

\noindent
Now, the (strict) Chevalley-Eilenberg construction yields a map 
$$
\mathrm{Mix}: \mathrm{LieAlg}^{\mathrm{str}}(V) \longrightarrow \epsilon-\mathrm{cdga}_{k}^{gr}(\mathrm{Sym} (V^{\vee} [-1]))
$$ 
which is easily checked to be a left inverse to $\mathrm{Lie}$; therefore $\mathrm{Lie}$ is surjective, hence bijective with inverse $\mathrm{Mix}$.
\hfill $\Box$\\

Recall that $\mathfrak{g}$ is the Lie algebra of $G$, and denote by $[\,,\,]$ its Lie bracket. As we have already seen in, we have a canonical equivalence 
$$
u: \D(B\widehat{G}_e) \simeq \mathrm{Sym} (\mathfrak{g}^{\vee} [-1])
$$
in $\cdga^{gr}_{k}$. Since $\D(B\widehat{G})$ has a canonical structure of graded mixed cdga, let $\langle\,,\,\rangle_{u}$ the Lie bracket induced on $\mathfrak{g}$ according to Proposition \ref{Lieb}.

\begin{prop}
With the above notation, and assume that $k$ is a field, we have
\begin{enumerate}
\item $(\mathfrak{g},[\,,\,])$ and  $(\mathfrak{g}, \langle\,,\,\rangle_u)$ are isomorphic Lie algebras.
\item There is an equivalence 
$$
\D(B\widehat{G}) \simeq \big(\mathrm{Sym} (\mathfrak{g}^{\vee} [-1]), \epsilon:=d_{\mathrm{CE}, [\, , \,]}\big)=: \mathsf{CE}\big(\mathfrak{g},[\,,\,]\big)
$$ 
in $\mecdga_{k}$.
\end{enumerate}
\end{prop}

\noindent \textbf{Proof.} (1) Recall the equivalence of symmetric monoidal $\infty$-categories 
$$
\mathbb{D}(B\widehat{G}_e)-Mod_{\epsilon-\mathbf{dg}^{gr}}^{Perf}\simeq \mathsf{Perf}(B\widehat{G}_e
$$
Let $\mathbb{D}(B\widehat{G}_e)-Mod_{\epsilon-\mathbf{dg}^{gr}}^{qffd}$ be the full sub-$\infty$-category of $\mathbb{D}(B\widehat{G}_e)-Mod_{\epsilon-\mathbf{dg}^{gr}}^{Perf}$ consisting of \emph{quasi-free finite dimensional modules}; i.e. those $\mathbb{D}(B\widehat{G}_e)$-modules which are equivalent as graded modules to $\mathbb{D}(B\widehat{G}_e)\otimes V$, where $V$ is a discrete finite dimensional $k$-vector space that is concentrated in pure weight $0$. 
The above equivalence then restricts to an equivalence of tensor $k$-linear (discrete) categories 
$$
\mathbb{D}(B\widehat{G}_e)-Mod_{\epsilon-\mathbf{dg}^{gr}}^{qffd}\simeq \mathsf{Rep}^{fd}(\widehat{G}_e)\,,
$$
where $\mathsf{Rep}^{fd}(\widehat{G}_e)$ is the tensor $k$-linear category of finite dimensional representations of $\widehat{G}_e$. 
Observe that this equivalence commutes with the obvious fiber functors to $\mathsf{Vect}(k)$ (whose geometric origin is simply the pull-back $x^*$ along the point $x:*\to B\widehat{G}_e$), where $\mathsf{Vect}(k)$ is the category of vector spaces. In particular, the above equivalence is an equivalence of neutral Tannakian categories, and we therefore have the following chain of equivalences between neutral Tannakian categories: 
$$
\mathsf{Rep}^{fd}\big(\mathfrak g, \langle\,,\,\rangle_u\big)\simeq \mathsf{CE}\big(\mathfrak{g},\langle\,,\,\rangle_u\big)-Mod_{\epsilon-\mathbf{dg}^{gr}}^{qffd}
\simeq  \mathbb{D}(B\widehat{G}_e)-Mod_{\epsilon-\mathbf{dg}^{gr}}^{qffd}\simeq \mathsf{Rep}^{fd}(\widehat{G}_e)\simeq \mathsf{Rep}^{fd}\big(\mathfrak g, [\,,\,]\big)\,.
$$
We refer to \cite{de-mi} for general facts about the Tannakian formalism, which tells us that we therefore have the following sequence of Lie algebra morphisms: 
\begin{equation}\label{eq:tannakian}
(\mathfrak g, \langle\,,\,\rangle_u\big)\longrightarrow \mathbf{End}(f_{\langle\,,\,\rangle_u})\cong\mathbf{End}(f_{[\,,\,]})\longleftarrow (\mathfrak g, [\,,\,]\big)\,,
\end{equation}
where $\mathbf{End}(f)$ is the Lie $k$-algebra of natural transformations of a given fiber functor $f$ (endowed with the commutator as Lie bracket), and $f_{\langle\,,\,\rangle_u}$ and $f_{[\,,\,]}$ are the fiber functors of $\mathsf{Rep}^{fd}\big(\mathfrak g, \langle\,,\,\rangle_u\big)$ and $\mathsf{Rep}^{fd}\big(\mathfrak g, [\,,\,]\big)$, respectively. 
It is a general fact that the leftmost and rightmost morphisms in \eqref{eq:tannakian} are injective. Moreover, $(\mathfrak g,[\,,\,])$ being algebraic, the leftmost morphism is actually an isomorphism. 
Therefore we get an injective Lie algebra morphism $(\mathfrak g, \langle\,,\,\rangle_u\big)\to (\mathfrak g, [\,,\,]\big)$, which must be an isomorphism for obvious dimensional reasons. \\
(2) To ease notations, we will write $B:= \D(B\widehat{G})$ as a graded mixed cdga, and $\epsilon_B$ its mixed differential. 
Since $B \in \epsilon-\mathrm{cdga}_{k}^{gr}(\mathrm{Sym} (V^{\vee} [-1]))$, by Proposition \ref{Lieb} we have  
$$\langle\, , \, \rangle_u= \mathrm{Lie}(H(B))= \mathrm{Lie}(B)\,.
$$
By $(1)$, and,  again, Proposition \ref{Lieb}, we get 
$$
\mathsf{CE}\big(\mathfrak{g},[\,,\,]\big) =\mathrm{Mix}([\, , \,])= \mathrm{Mix}(\langle\, , \, \rangle_{u}) = H(B) = B\,,
$$
where the equalities are in $\epsilon-\mathrm{cdga}_{k}^{gr}(C) := \pi_0 (\mecdga_{k}{C})$. 
In particular, $B$ and $\mathsf{CE}\big(\mathfrak{g},[\,,\,]\big)$ are equivalent in $\mecdga_{k}$.
\hfill $\Box$\\
\begin{rmk}
Let us give an alternative,  less elementary but direct proof of $(2)$. As observed in \S  \ref{t4'}, an equivalence of graded cdgas $v: B \simeq \mathrm{Sym} (\mathfrak{g}^{\vee} [-1])$ induces a \emph{weak} mixed structure (see proof of Lemma \ref{l11}) on $C:=\mathrm{Sym} (\mathfrak{g}^{\vee} [-1])$, i.e. a family of strict maps 
$$
\epsilon_i : C(p) \longrightarrow C(p+i+1)[1]\, , \, i\geq 0
$$
satisfying a Maurer-Cartan-like equation. In our case 
$$
\epsilon_i : (\wedge^p \mathfrak{g}^{\vee})[-p] \longrightarrow (\wedge^{p+i+1} \mathfrak{g}^{\vee})[-p-i]
$$
hence $\epsilon_i=0$ for $i>0$, because $\mathfrak{g}$ sits in cohomological degree $0$. 
The only non-trivial remaining map is $\epsilon_0$, and the Maurer-Cartan equation tells us exactly that  it defines a strict graded mixed cdga structure on $\mathrm{Sym} (\mathfrak{g}^{\vee} [-1])$, and that, with such structure, the equivalence  $v: B \simeq \mathrm{Sym} (\mathfrak{g}^{\vee} [-1])$ is indeed an equivalence of graded \emph{mixed} cdgas.
\end{rmk}

\paragraph{The case $n=1$ for a reductive $G$.}

We have seen in \S \ref{sec:3.1} that equivalences classes of $1$-shifted Poisson structures on $BG$, for a reductive group $G$, are in bijection with elements $Z\in\wedge^3(\mathfrak g)^G$. 
The induced $1$-shifted Poisson structure on the graded mixed cdga $\mathsf{CE}(\mathfrak g)$ is then very explicit in terms of a so-called semi-strict $\mathbb{P}_{n+1}$-structure (see \cite{mel}): 
all structure $2$-shifted polyvectors are trivial except for the $3$-ary one which is constant and given by $Z$. 

Our deformation quantization in particular leads to a deformation of $\mathsf{Rep}^{fd}(\mathfrak{g})$ as a monoidal category. 

\begin{ex}\label{ex:Z}
Given a non-degenerate invariant pairing $<\,,\,>$ on $\mathfrak g$, such an element can be obtained from the $G$-invariant linear form 
$$
\wedge^3\mathfrak g\longrightarrow k\quad,\quad(x,y,z)\longmapsto <x,[y,z]>\,.
$$
Alternatively, any invariant symmetric $2$-tensor $t\in Sym^2(\mathfrak g)^G$ leads to such an element $Z=[t^{1,2},t^{2,3}]\in\wedge^3(\mathfrak g)^G$. 
In this case the deformation of $\mathsf{Rep}^{fd}(\mathfrak{g})$ as a monoidal category can be obtained by means of a deformation of the associativity constraint only (see \cite{Dr1}), which then looks like   
$$
\Phi=1^{\otimes3}+\hbar^2Z+o(\hbar^2)\in U(\mathfrak{g})^{\otimes3}[[\hbar]]\,.
$$
\end{ex}

\begin{rmk}
Note that even in the case when $G$ is not reductive, every element $Z\in\wedge^3(\mathfrak g)^G$ lead to a $1$-shifted Poisson structure on $BG$ as well (but we have a map $\wedge^3(\mathfrak g)^G\to\pi_0\mathsf{Pois}(BG,1)$ rather than a bijection). The above reasoning works as well for these $1$-shifted Poisson structures. 
\end{rmk}

\paragraph{The case $n=2$ for a reductive $G$.}

We have seen in \S \ref{sec:3.1} that equivalences classes of $2$-shifted Poisson structures on $BG$, for a reductive group $G$, are in bijection with elements $t\in Sym^2(\mathfrak g)^G$. 
The induced $2$-shifted Poisson structure on the graded mixed cdga $\mathsf{CE}(\mathfrak g)$ is strict and constant. 
The graded mixed $\mathbb{E}_3$-structure on $\mathsf{CE}(\mathfrak g)$ given by our deformation quantization then takes the form of a Weyl $3$-algebra, as described in \S \ref{sec:quantpoint} (one simply has to replace $p$ by $t$ in \eqref{eq:n-weyl}). 

Note that, as we already mentioned, this graded mixed $\mathbb{E}_3$-structure is $G_{DR}$-equivariant by construction, so that it leads to an $\mathbb{E}_2$-monoidal deformation of $\mathsf{Perf}(BG)$. 
This in particular leads to a braided monoidal deformation of $\mathsf{Rep}^{fd}(\mathfrak{g})$. 

\begin{rmk}
Note that even in the case when $G$ is not reductive, elements $t\in Sym^2(\mathfrak g)^G$ are exactly $2$-shifted Poisson structure on $BG$ (i.e. we have a map $Sym^2(\mathfrak g)^G\cong\pi_0\mathsf{Pois}(BG,2)$). 
The above reasoning works as well for these $2$-shifted Poisson structures.  
\end{rmk}

Such deformation quantizations of $BG$ have already been constructed: 
\begin{itemize}
\item when $\mathfrak g$ is reductive and $t$ is non-degenerate, by
  means of purely algebraic methods: the quantum group
  $U_\hbar(\mathfrak g)$ is an explicit deformation of the enveloping
  algebra $U(\mathfrak g)$ as a quasi-triangular Hopf algebra.
\item without any assumption, by Drinfeld \cite{Dr2}, using
  transcendental methods similar to the ones that are crucial in the
  proof of the formality of $\mathbb{E}_2$.
\end{itemize} 
It is known that Drinfeld's quantization is equivalent to the quantum
group one in the semi-simple case (see e.g. \cite{kas2} and references
therein).

\begin{rmk}
  It is remarkable that our quantization relies on the formality of
  $\mathbb{E}_3$ rather than on the formality of $\mathbb{E}_2$.  It
  deserves to be compared with Drinfeld's one, but this task is beyond
  the scope of the present paper.
\end{rmk}

%%%%%%%%%%%%%%%%%%%

\

\bigskip

\bigskip

\appendix

\Appendix{\ } \label{appendix:A}

\setcounter{subsection}{1}
\setcounter{thm}{0}
\setcounter{equation}{0}

\

This Appendix contains a few technical results needed in Sect. \ref{RDCsection}.

\begin{prop}\label{Mstable}
Any $C(k)$-model category is a stable model category. 
\end{prop}
\textbf{Proof.} Let $N$ be a $C(k)$-model category, and let
$\underline{\mathrm{Hom}}_{k} (-,-)$ be its enriched
hom-complex. There is a unique map $0 \to \underline{\mathrm{Hom}}_{k}
(*,\emptyset)$ in $C(k)$, where $*$ (respectively, $\emptyset$) is the
final (respectively, initial) object in $N$.  By Composing with the
map $k \to 0$ in $C(k)$, we get a map in $N$ from its final to its
initial object: hence $N$ is pointed. Let us denote by $\Sigma:
\mathrm{Ho}(N) \to \mathrm{Ho}(N)$ the corresponding suspension
functor.  For $X \in N$ cofibrant we have that $X\otimes_{k} k[1]
\simeq \Sigma (X)$ (since $X\otimes_{k} (-)$ preserves homotopy
pushouts and $k[1]$ is the suspension of $k$ in $C(k)$). Therefore,
the suspension functor $\Sigma$ is an equivalence, its quasi inverse
being given by $(-)\otimes^{\mathbb{L}}_{k} k[-1]$.  \hfill $\Box$ 

\

\begin{prop}\label{1-5formod}Let $M$ be a symmetric monoidal
  combinatorial model category  
satisfying the standing assumptions $(1)-(5)$ of Section \ref{1.1},
and let $A \in Comm(M)$.  Then the symmetric monoidal combinatorial
model category $A-Mod_{M}$ also satisfies the standing assumptions
$(1)-(5)$.
\end{prop}
\textbf{Proof.} Left to the reader. \hfill $\Box$ 

\

\begin{prop}\label{modinvariance} Let $M$ be a symmetric monoidal combinatorial 
model category satisfying the standing assumptions $(1)-(5)$ of
Section \ref{1.1}.  If $w: A \to B$ is a weak equivalence in
$Comm(M)$, then the Quillen adjunction
$$
w^*= -\otimes_{A}B \,\, : A-Mod_{M} \longleftrightarrow B-Mod_{M} :
w_* 
$$ 
is a 
Quillen equivalence.
\end{prop}
\textbf{Proof.} Since $w_*$ reflects weak equivalences, $w^*$ is a
Quillen equivalence iff for any cofibrant $A$-module $N$, the natural
map $i: \mathrm{id}_{N}\otimes w: N \simeq N\otimes_{A}A \to
N\otimes_{A} B$ is a weak equivalence. Since $N$ is cofibrant, we may
write it as $\textrm{colim}_{\beta \leq \alpha} N_{\beta}$ (colimit in
$A-Mod_{M}$) where $\alpha$ is an ordinal, $N_{0}= 0$ and each map
$N_{\beta} \to N_{\beta +1} $ is obtained as a pushout in $A-Mod_{M}$
$$
\xymatrix{ A\otimes X \ar[r]^-{\textrm{id}\otimes u} \ar[d] & A
  \otimes Y \ar[d] \\  
N_{\beta} \ar[r] & N_{\beta +1}}
$$ 
where $u: X \to Y$ belongs to the set $I$ of generating
cofibrations of $M$ (all assumed with $M$-cofibrant domain, by
standing assumption (3)). In order to prove that $i: N \simeq
N\otimes_{A}A \to N\otimes_{A} B$ is a weak equivalence, we will
prove, by transfinite induction, that each $i_{\beta}: N_{\beta}
\simeq N_{\beta}\otimes_{A}A \to N_{\beta} \otimes_{A} B$ is a weak
equivalence. 

\

\medskip

\noindent
Since $N_{0}=0$, the induction can start. Let us suppose that
$i_{\beta}$ is a weak equivalence, and consider the pushout diagram
$\mathsf{P}$ defining $N_{\beta} \to N_{\beta +1}$
$$
\xymatrix{ A\otimes X \ar[r]^-{\textrm{id}\otimes u} \ar[d] & A
  \otimes Y \ar[d] \\ 
N_{\beta} \ar[r] & N_{\beta +1}.}
$$ 
Now, let us apply the functor $w^*$ to this pushout. We obtain the diagram  
$\mathsf{P'}$ 
$$
\xymatrix{ B\otimes X \ar[r]^-{\textrm{id}\otimes u} \ar[d] & B
  \otimes Y \ar[d] \\  
N_{\beta} \otimes_{A} B \ar[r] & N_{\beta +1} \otimes_{A} B}
$$ 
which is again a pushout in $B-Mod_{M}$ (since $w^*$ is left
adjoint).  There is an obvious map of diagrams from $\mathsf{P}$ to
$\mathsf{P'}$ induced by the maps $w\otimes \textrm{id}_{X}: A\otimes
X \to B\otimes X$, $i_{\beta}: N_{\beta} \to N_{\beta} \otimes_{A} B$,
and $w \otimes \textrm{id}_{Y}: A\otimes Y \to B\otimes Y$.  All these
three maps are weak equivalences ($i_{\beta}$ by induction hypothesis,
and the other two by standing assumption $(3)$, since $X$ is
cofibrant, and so is $Y$, $u$ being a cofibration).  Since the
forgetful functor $A-Mod_{M} \to M$ has right adjoint the internal
hom-functor $\underline{Hom}_{M}(A, -)$, both $\mathsf{P}$ and
$\mathsf{P'}$ are pushouts in $M$, too. Thus
(\cite[Proposition~13.5.10]{hir}) also the induced map $i_{\beta +1}:
N_{\beta +1} \to N_{\beta +1} \otimes_{A} B$ is a weak equivalence (in
$M$) as the two diagrams $\mathsf{P}$ and $\mathsf{P'}$ are also
homotopy pushouts, by standing assumption $(2)$ on $M$. We are done
with the successor ordinal case and left to prove the limit ordinal
case. The family of maps $\{i_{\beta} \}$ are all weak equivalences
and define a map of sequences $\{ N_{\beta} \} \to \{N_{\beta}
\otimes_{A} B \} $, where each map $N_{\beta} \to N_{\beta +1}$ is a
cofibration (as pushout of a cofibration), and the same is true for
each map $ N_{\beta} \otimes_{A} B \to N_{\beta + 1} \otimes_{A} B$
(since $w^*$ is left Quillen). Moreover, each $N_{\beta}$ is cofibrant
(since $N_0 =0$ is and each $N_{\beta} \to N_{\beta +1}$ is a
cofibration), and the same is true for each $N_{\beta} \otimes_{A} B$
(since $w^*$ is left Quillen). Therefore the map induced on the
(homotopy) colimit is a weak equivalence too.

\hfill $\Box$ 

\

\begin{prop}\label{cofibforgettocofib} Let $M$ be a symmetric monoidal combinatorial model 
category satisfying the standing assumptions $(1)-(5)$ of Section
\ref{1.1}.  Then the forgetful functor $Comm(M) \to M$ preserves
fibrant-cofibrant objects.
\end{prop}
\textbf{Proof.} The forgetful functor is right Quillen, so it
obviously preserves fibrant objects.  The $C(k)$-enrichment, together
with $\mathrm{char}(k)=0$, implies that $M$ is freely powered in the
sense of \cite[Definition~4.5.4.2]{lualg}.  By
\cite[Lemma~4.5.4.11]{lualg}, $M$ satisfies the strong commutative
monoidal axiom of \cite[Definition~3.4]{white}. Then, the statement
follows from our standing assumption $(1)$ and from
\cite[Corollary~3.6 ]{white}.

\hfill $\Box$

\

\bigskip

\bigskip

\Appendix{\ } \label{appendix:B}

\setcounter{subsection}{1}
\setcounter{thm}{0}
\setcounter{equation}{0}

\

\

\noindent
We prove here several technical statement about differential forms and
formal completions in the derived setting, needed in Sect. \ref{FLsection}.

\begin{lem}
Let $X \longrightarrow U \longrightarrow Y$ be morphisms of derived
algebraic $n$-stacks. Let $U_*$ be the nerve of the morphism
$U\longrightarrow Y$. Then, for all $p$ there is a natural equivalence
$$\Gamma(X,\wedge^{p}\mathbb{L}_{X/Y}) \simeq
\lim_{n \in \Delta}\Gamma(X,\wedge^{p}\mathbb{L}_{X/U_{n}}).$$
\end{lem}
\noindent \textbf{Proof.}  
For $F\in\dSt_k$ we consider the shifted tangent derived stack
$$T^1(F):=\Map(\Spec\, k[\epsilon_{-1}],F),$$
the internal Hom object, where $k[\epsilon_{-1}]=k\oplus k[1]$ is
the free cdga over one generator in degree $-1$. The natural 
augmentation $k[\epsilon_{-1}] \rightarrow k$ induces
a projection $T^1(F) \longrightarrow F$. Moreover, if $F$ is
an algebraic derived n-stack then $T^1(F)$ is an algebraic
derived $(n+1)$-stack.

For a morphism
$F \longrightarrow G$, we let
$$T^1(F/G):=T^1(F)\times_{T^1(G)}G,$$ as a derived stack over $F$. The
multiplicative group $\mathbb{G}_{m}$ acts on $T^1(F/G)$, and thus we
can consider $\Gamma(T^1(F/G),\mathcal{O})$ as a graded complex. As
such, its part of weight $p$ is
$$\Gamma(F,\wedge^{p}\mathbb{L}_{F/G})[-p].$$
In order to conclude, we observe that the induced morphism, which is
naturally $\mathbb{G}_{m}$-equivariant
$$T^1(X/U) \longrightarrow T^1(X/F)$$ is an epimorphism of derived
stacks. The nerve of this epimorphism is the simplicial object $n
\mapsto T^1(X/U_n)$. By descent for functions of weight $p$ we see
that the natural morphism
$$\Gamma(X,\wedge^p\mathbb{L}_{X/F}) \longrightarrow 
\lim_n\Gamma(X,\wedge^p\mathbb{L}_{X/U_n})$$
is an equivalence. 
\hfill $\Box$

\

\noindent
For the next lemma, we will use Koszul commutative dg-algebras. For a
commutative k-algebra $B$, and $f_1,\dots,f_p$ a family of elements in
$B$, we let $K(B,f_{1},\dots,f_p)$ be the commutative dg-algebra
freely generated over $B$ by variables $X_1,\dots,X_p$ with
$deg(X_i)=-1$, and with $dX_i=f_i$. When $f_1,\dots,f_p$ form a
regular sequence in $B$, then $K(B,f_1,\dots,f_p)$ is a cofibrant
model for $B/(f_1,\dots,f_p)$ considered as a $B$-algebra. In general,
$\pi_i(K(B,f_1,\dots,f_p))\simeq Tor_{i}^{B}(B/(f_1),\dots,B/(f_p))$
are possibly non zero only when $i\in [0,p]$.

\begin{lem}
Let $B$ be a commutative (non-dg) $k$-algebra of finite type and
$I \subset B$ an ideal generated by $(f_1,\dots,f_p)$. 
Let $f : X=\Spec\, B/I \longrightarrow Y=\Spec\, B$ 
be the induced morphism of affine schemes, 
and $X_n:=\Spec\, K(B,f_1^n,\dots,f_p^n)$. Then, the natural morphism
$$\mathrm{colim}_{n}X_n \longrightarrow \widehat{Y}_{f}$$ is an
equivalence of derived prestacks: for all $\Spec\, A \in \dAff_{k}$ we
have an equivalence
$$\mathrm{colim}_{n}(X_n(A)) \simeq \widehat{Y}_{f}(A).$$
\end{lem}
\noindent \textbf{Proof.} We let $F$ be the colimit prestack
$\mathrm{colim}_n X_n$. There is a natural morphism of derived
prestacks
$$\phi : F \longrightarrow \widehat{Y}_{f}.$$
For any $k$-algebra $A$ of finite type, the induced morphism of sets
$$F(A) \longrightarrow \widehat{Y}_{f}(A)$$ is bijective. Indeed, the
left hand side is equivalent to the colimit of sets $\mathrm{colim}_n
Hom_{k-Alg}(B/I(n),A)$, where $I(n)$ is the ideal generated by the
$n$-th powers of the $f_i$'s, whereas the right hand side consists of
the subset of $Hom_{k-Alg}(B,A)$ of maps $f : B \longrightarrow A$
sending $I$ to the nilpotent radical of $A$.  In order to prove that
the morphism $\phi$ induces an equivalences for all $\Spec\, A \in
\dAff_{k}$ we use a Postnikov decomposition of $A$
$$\xymatrix{A \ar[r] & \dots \ar[r] & A_{\leq k} \ar[r] & A_{\leq k-1}
\ar[r] & \dots \ar[r] & A_{\leq 0}=\pi_0(A).}$$
As prestacks, i.e. as $\s$-functors on $\dAff_{k}^{op}$, both
$F$ and $\widehat{Y}_{f}$ satisfy the following two properties.

\begin{itemize}

\item For all $\Spec\, A \in \dAff_{k}$, we have equivalences 
$$F(A) \simeq \lim_{k}F(A_{\leq k}) \qquad \widehat{Y}_{f}(A) \simeq \lim_{k} 
\widehat{Y}_{f}(A_{\leq k})$$

\item For all 
fibered product of almost finite presented $k$-cdgas in non-positive degrees
$$\xymatrix{ B \ar[r] \ar[d] & B_{1} \ar[d] \\ B_{2} \ar[r] &
  B_{0},}$$ such that $\pi_{0}(B_i) \longrightarrow \pi_{0}(B_0)$ are
surjective with nilpotent kernels, the induced square
$$\xymatrix{
F(B) \ar[r] \ar[d] & F(B_{1}) \ar[d] \\
F(B_{2}) \ar[r] & F(B_{0}),}$$
is cartesian in $\T$.
\end{itemize}
 
The above two properties are clear for $\widehat{Y}_{f}$, because
$\widehat{Y}_{f}$ is a formal stack. The second property is also clear
for $F$ because filtered colimits preserve fiber products. Finally,
the first property is satisfied for $F$ because for each fixed $n$,
and each fixed $i\geq 0$ the projective system of homotopy groups
$$\xymatrix{\pi_i(X_{n}(A)) \ar[r] & \dots \ar[r] & 
\pi_i(X_{n}(A_{\leq k})) \ar[r] & 
\pi_i(X(A_{\leq k-1}))
\ar[r] & \dots \ar[r] & \pi_i(X(A_{\leq 0}))}$$
stabilizes (this is because $K(B,f_1^n,\dots,f_p^n)$ are cell $B$-cdga with 
finitely many cells and thus with a perfect cotangent complex).

By these above two properties, and by Postnikov decomposition,  
we are reduced to prove that for any
non-dg $k$-algebra $A$ of finite type, any $A$-module $M$ of 
finite type, and any $k\geq 1$ the induced morphism
$$F(A\oplus M[k]) \longrightarrow \widehat{Y}_{f}(A\oplus M[k])$$ is
an equivalence. We can fiber this morphism over $F(A) \simeq
\widehat{Y}_{f}(A)$ and thus are reduced to compare cotangent
complexes of $F$ and $\widehat{Y}_{f}$.

By replacing $X$ by one of the $X_n$, we can assume that $\Spec\, A =
X$ and thus that $A=B/I$.  We thus consider the morphism induced on
cotangent complexes for the morphism $X \longrightarrow F
\longrightarrow \widehat{Y}_{f}$
$$\mathbb{L}_{X/F} \longrightarrow \mathbb{L}_{X/\widehat{Y}_{f}}.$$
Here, $\mathbb{L}_{X/F}$ is not quite an $A$-dg-module but 
is a pro-object in $L_{coh}^{\leq 0}(A)$ which represents
the adequate $\s$-functor. This pro-object is explicitly given by
$$\mathbb{L}_{X/F} \simeq "\lim_{n}"\mathbb{L}_{X/X_n}.$$ We have to
prove that the morphism of pro-objects
$$"\lim_{n}"\mathbb{L}_{X/X_n} \longrightarrow
\mathbb{L}_{X/\widehat{Y}_{f}},$$ where the right hand side is a
constant pro-object, is an equivalence. Equivalently, using various
exact triangles expressing cotangent complexes we must prove that the
natural morphism
$$"\lim_{n}"u_n^*(\mathbb{L}_{X_n/Y}) \longrightarrow
u^{*}(\mathbb{L}_{\widehat{Y}_{f}/Y})$$ is an equivalence of
pro-objects, where $u_n : X \longrightarrow X_n$ and $u : X
\longrightarrow Y$ are the natural maps. The right hand side vanishes
because $\widehat{Y}_{f} \longrightarrow Y$ is formally \'etale.
Finally, the left hand side is explicitly given by the projective
systems of $A=B/I$-dg-modules $"\lim_n(A^p[1])$ (because
$K(B,f_1^n,\dots,f_p^n)\otimes_{B}A$ is freely generated over $A$ by
$p$ cells of dimension $1$).  Here the transition morphisms are
obtained by multiplying the $i$-th coordinate of $A^p$ by $f_i$ and
thus are the zero morphisms. This pro-object is therefore equivalent
to the zero pro-object, and this finishes the proof of the lemma.
\hfill $\Box$

\begin{lem}\label{lappend}
Let $X$ be an affine formal derived stack. We assume that, as a derived prestack $X$ is
of the form $X\simeq \mathrm{colim}_{n\geq 0}X_n$,
with $X_n \in \dAff_{k}$ for all $n$.
Then, for all $p$, the natural morphism
$$\wedge^p \mathbb{L}_{X_{red}/X} \simeq \lim_{n}\wedge^p \mathbb{L}_{X_{red}/X_n}$$
is an equivalence in $\mathsf{L}_{\mathsf{Qcoh}}(X_{red})$.
\end{lem}

\noindent \textbf{Proof:} We
consider the $\s$-functor co-represented by 
$\mathbb{L}_{X_{red}/X}$
$$\mathsf{Map}(\mathbb{L}_{X_{red}/X},-) : L_{coh}^{\leq 0}(X_{red}) \longrightarrow \T.$$
Note that because $X$ is a colimit of derived schemes its cotangent complex
$\mathbb{L}_{X_{red}/X}$ sits itself in $L_{coh}^{\leq 0}(X_{red})$. Moreover, 
as $X$ is the colimit of the $X_n$ as derived prestacks, the $\s$-functor 
$\mathsf{Map}(\mathbb{L}_{X_{red}/X},-)$
is also pro-representable by the pro-object in $L_{coh}^{\leq 0}(X_{red})$
$$"\lim_{n}"\{\mathbb{L}_{X_{red}/X_{n}}\}.$$
Therefore, this pro-object is equivalent, in the $\s$-category of 
pro-objects in $L_{coh}^{\leq 0}(X_{red})$, to the
constant pro-object $\mathbb{L}_{X_{red}/X}$. Passing to 
wedge powers, we see that for all $p$ the pro-object  
$"\lim_{n}"\{\wedge^p\mathbb{L}_{X_{red}/X_{n}}\}$ is also 
equivalent to the constant pro-object $\wedge^p\mathbb{L}_{X_{red}/X}$, and
the lemma follows.
\hfill $\Box$

\

\

\smallskip

\noindent
Damien Calaque, {\sc Universit\'e de
  Montpellier \& Institut Universitaire de France}, \\ damien.calaque@univ-montp2.fr

\smallskip

\noindent
Tony Pantev, {\sc University of Pennsylvania}, tpantev@math.upenn.edu

\smallskip

\noindent
Bertrand To\"{e}n, {\sc  Universit\'e Paul Sabati\'er \& CNRS},
Bertrand.Toen@math.univ-toulouse.fr 

\smallskip

\noindent
Michel Vaqui\'{e},  {\sc  Universit\'e Paul Sabati\'er \& CNRS},
michel.vaquie@math.univ-toulouse.fr

\smallskip

\noindent
Gabriele Vezzosi, {\sc Universit\`a di Firenze},
gabriele.vezzosi@unifi.it

\end{document}